\documentclass[11pt,a4paper,titlepage,pdf]{book}
\usepackage[latin1]{inputenc}
\usepackage[T1]{fontenc}
\usepackage{numprint}
\usepackage{makeidx}
\usepackage[francais]{babel}
\usepackage{amssymb,amsmath,stmaryrd,amscd,bbm,marvosym}
\usepackage{theorem} 

\usepackage[francais]{babel} 

\renewcommand{\cleardoublepage}{%
  \clearpage%
  \if@twoside%
    \ifodd\c@page%
    \else%
      \thispagestyle{empty}\hbox{}\newpage%
      \if@twocolumn \hbox{}\newpage \fi%
    \fi%
  \fi%
}

\newcounter{claim}

\newenvironment{preuve}{{\setcounter{claim}{0}\noindent\sc preuve
    ---}}{\hfill$\Box$\vspace{2ex}} 

\newenvironment{claim}[1][]%
{\refstepcounter{claim}\vspace{1ex}\noindent{(\it\arabic{claim}){#1}{}}\it}{\vspace{1ex}}

\newenvironment{preuveclaim}[1][]%
	{\noindent {}{#1}{}}{ Ceci prouve~(\arabic{claim}).\vspace{1ex}}

 {\noindent {\setcounter{claim}{0}\sc proof ---
    }{#1}{}}{\hfill$\Box$\vspace{2ex}}

\newcommand{\bp}{}
\newcommand{\tp}{\!-\!}
\newcommand{\ep}{}

\newtheorem{defeng}{Définition}[chapter]
\newtheorem{theoreme}[defeng]{Théorème}

\newtheorem{lemme}[defeng]{Lemme}

\newtheorem{conjecture}[defeng]{Conjecture}
\newtheorem{corollaire}[defeng]{Corollaire}
\usepackage{graphicx} 

\newcommand{\n}{\numprint}

\newcounter{outcomes}
\newenvironment{outcomes}
	       {
		 \begin{list}
		   {(\alph{outcomes})}
		   {
		     \usecounter{outcomes}
		     \setlength{\leftmargin}{1.5em}
		     \setlength{\labelwidth}{.5em}
		   }
	       }
	       {\end{list}}

{\theoremstyle{break}\theorembodyfont{\rmfamily} \newtheorem{algorithme}[defeng]{Algorithme}}
{\theoremstyle{break}\theorembodyfont{\rmfamily} \newtheorem{probleme}[defeng]{Probleme}}
{\theorembodyfont{\rmfamily} \newtheorem{question}[defeng]{Question}}

\begin{document}

\thispagestyle{empty}
{\textwidth18cm
\oddsidemargin0cm
\topmargin0cm
\footskip0cm
\textheight27cm
\huge
{\center 

\noindent \rule{1.5cm}{0cm}Thèse\\}

\vspace{4ex}

\Large
\noindent Présentée par {\sc Nicolas Trotignon}
}\\

\parbox{13cm}{

\noindent pour obtenir le grade de docteur \\de l'université Grenoble I,
Joseph Fourier\\

\noindent Spécialité : Mathématiques Informatique

\noindent Formation doctorale : Recherche Opérationnelle,\\ Combinatoire
et Optimisation
}

{\center
\vspace{20ex}\mbox{}\\
\Huge
\noindent \rule{1.5cm}{0cm}Graphes parfaits~:\\
\noindent \rule{1.5cm}{0cm}Structure et algorithmes\\
\huge
\vspace{5ex}\mbox{}\\
}

{
\Large

\vfill

\vspace{11ex}
\noindent Soutenue le 28 septembre 2004\\

\noindent {\sc Michel Burlet} (co-directeur de thèse)\\
\noindent {\sc Gérard Cornuéjols} (rapporteur)\\
\noindent {\sc Jean Fonlupt} (examinateur)\\
\noindent {\sc Jean-Luc Fouquet} (examinateur)\\
\noindent {\sc Frédéric Maffray} (directeur de thèse)\\
\noindent {\sc Bruce Reed} (rapporteur) \\
}

\newpage
\thispagestyle{empty}
\mbox{}
\newpage

\setcounter{page}{1}
\tableofcontents
\nocite{berge.chvatal:topics}
\nocite{perec:cantatrix}
   
\bibliographystyle{plain}

\chapter*{Remerciements}
\addcontentsline{toc}{chapter}{Remerciements}

Je  tiens tout  d'abord à  remercier Gérard  Cornuéjols et  Bruce Reed
d'avoir accepté  de rapporter cette thèse.  Je  remercie également mes
deux examinateurs, Jean Fonlupt et Jean-Luc Fouquet.

Je tiens à remercier mon employeur, l'Université Pierre Mendès-France,
et  plus particulièrement  son Conseil  Scientifique qui  m'a déchargé
d'une partie de mes enseignements pour me permettre de mener à bien ce
travail.

\vspace{4ex}

Je  remercie les  institutions suivantes  qui m'ont  permis,  par leur
soutien financier, d'assister à des conférences internationales~:

\vspace{4ex}

\begin{itemize}
  \item 
  Mathematisches   Forschungsinstitut   à   Oberwolfach~:   Conférence
  ``Geometric  Convex   Combinatorics''  organisée  par   B.  Gerards,
  A. Seb\H o et R. Weismantel en juin 2002.

  \item
  American Institute  of Mathematics,  à Palo Alto~:  ``Perfect Graphs
  Workshop'' organisé par P. Seymour et R. Thomas en novembre 2002.

  \item 
  Caesarea  Edmond  Benjamin de  Rothschild  Foundation Institute  for
  Interdisciplinary Applications of Computer Science, Haifa~: ``Expert
  workshop on graph classes and algorithms'' organisé par M. Golumbic,
  A. Berry et F. Maffray en avril 2004.
\end{itemize}

\vspace{4ex}

Je remercie  l'équipe du laboratoire Leibniz, et  tout d'abord Sylvain
Gravier  qui  m'a  fait  faire  mes premiers  pas  au  laboratoire  en
encadrant  mon stage de  DEA.  J'ai  pu bénéficier  du soutien  et des
conseils d'autres membres de son équipe, Michel Mollard, Charles Payan
et  Xuong  Nguyen,  ainsi  que  des membres  de  l'équipe  ``Recherche
Opérationnelle'', Nadia Brauner, Gerd Finke et Maurice Queyranne.

L'équipe  ``graphes  et  optimisation combinatoire''  m'a  accueilli~:
Andr\'as Seb\H o et son enthousiasme constant et communicatif~; Myriam
Preissmann, avec qui j'ai passé de  longues heures à lire la preuve du
théorème fort  des graphes parfaits~; Michel Burlet,  qui a co-encadré
cette thèse, et m'a orienté vers les graphes parfaits.

La  bonne  ambiance   du  laboratoire  ne  le  serait   pas  sans  mes
condisciples thésards,  présents et passés~:  Prakash Countcham, David
Defossez, Éric Duchène, Haris  Gravanovic, Vincent Jost, Yann Kieffer,
Mohamed  Kobeissi,  Marie Lalire,  Pierre  Lemaire, Benjamin  Lévêque,
Medhi  Mallah, Frédéric  Meunier, Julien  Moncel, Simon  Perdrix, Éric
Tannier.   J'ai eu  des  contacts enrichissants  avec des  visiteurs~:
Nicola Apollonio, Hein Van  der Holst, Celina de Figueiredo, Cl\'audia
Linhares Sales, Giacomo Zambelli notamment.

\vspace{4ex}

Je remercie aussi mon ami  Jérôme Renault, toujours prêt à attaquer de
nouvelles énigmes, et souvent à les résoudre. Le travail effectué avec
lui pendant la durée de cette thèse m'a beaucoup aidé.

\vspace{4ex}

Je remercie tous mes relecteurs~: en premier lieu Frédéric Maffray qui
a participé à tous les  résultats présentés ici.  Bruce Reed a proposé
des  améliorations   de  fond  qui   ont  permis  de   raccourcir  des
démonstrations  fastidieuses,  et  d'accélérer  certains  algorithmes.
Michel Burlet a lui aussi proposé de nombreuses améliorations.  Myriam
Preissmann, Benjamin  Lévêque, David Defossez,  Yann Kieffer, Isabelle
Rave  et  Jérôme  Renault  ont   corrigé  pas  mal  de  défauts.   Mes
discussions  avec  Pierre  Lemaire  sur  \LaTeX,  Bibtex  et  Metapost
(excellent langage gratuit pour  dessiner des graphiques), ont apporté
beaucoup à la présentation.

\vspace{4ex}

Enfin, mon épouse Christelle  Petit, non mathématicienne mais lettrée,
a eu le courage de relire l'ensemble du document pour l'expurger d'une
quantité de  fautes d'orthographe et de  maladresses d'expression.  Je
veux lui dire  que je n'aurais jamais pu mener à  bien ce travail sans
son soutien, ni sans nos enfants Émile, Alice et Coline.

\vspace{4ex}

Je tenais à  dédier ce travail aux nombreux  excellents professeurs de
mathématiques que j'ai rencontrés tout  au long de mes longues études,
au nombre desquels~: Mme Ma et Pelissou, MM Moscovici, Danset, Cori et
Burlet.

Dernier d'entre eux en date, Frédéric Maffray, mon directeur de thèse,
ne saura jamais  assez ma gratitude pour la manière  dont il a encadré
mon travail,  sa patience et sa  disponibilité.  J'ai eu  la chance de
bénéficier de sa connaissance profonde des graphes parfaits durant une
période critique. Entre mille bienfaits,  je tiens à le remercier dans
le  détail  d'un épisode  de  notre  travail  commun, pour  sa  valeur
exemplaire  et morale  (ne  lui  en déplaise)~:  son  insistance à  ce
qu'avec Myriam Preissmann  nous lisions les 148 pages  de la preuve du
théorème fort des graphes parfaits~: travail fastidieux mais hautement
bénéfique.

\markboth{}{}

\chapter*{Introduction}
\addcontentsline{toc}{chapter}{Introduction}
\markboth{}{}

Considérons un  pays dont la station nationale  de radiodiffusion veut
affecter  une  fréquence  à  chaque  ville.  Quand  deux  villes  sont
éloignées de  moins de  50 km, il  faut leur attribuer  des fréquences
différentes  afin d'éviter  les  interférences. On  peut modéliser  le
problème  par un  \emph{graphe}~: on  représente chaque  ville  par un
point (en  théorie des  graphes, on dit  plutôt \emph{sommet}),  et on
relie deux sommets par une  ligne (on dit plutôt \emph{arête}) dès que
les deux villes correspondantes sont éloignées de moins de 50 km. Cela
pourrait donner l'un des graphes représentés ci-dessous~:

\begin{figure}[ht]
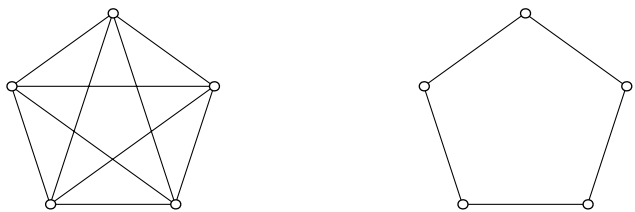

  \center
  \begin{tabular}{cc}
    \parbox{5cm}{\center\includegraphics{fig.base.27}} & 
    \parbox{5cm}{\center\includegraphics{fig.base.28}}
  \end{tabular}
  \caption{Deux graphes}
\end{figure}

Il  faut bien  noter que  sur la  figure ci-dessus,  l'endroit  où est
représenté un  sommet et  la forme effective  des arêtes  n'ont aucune
importance. Seul  compte le fait que  deux sommets sont  reliés ou non
reliés  par une  arête.  Notre  problème d'affectation  des fréquences
radio revient alors  à associer à chaque sommet  d'un graphe un nombre
entier  (correspondant à  une  fréquence) de  sorte  que deux  sommets
reliés  par une  arête  ne reçoivent  jamais  le même  nombre. Il  est
toujours  possible d'affecter  à  chaque sommet  du  graphe un  nombre
différent, mais on souhaite généralement minimiser le nombre d'entiers
utilisés --- les fréquences radio sont précieuses.

En  théorie des  graphes, plutôt  que  d'associer un  nombre à  chaque
sommet,       il        est       coutumier       d'associer       une
\emph{couleur}\index{couleur}~:  rouge, bleu~\dots\ Cette  habitude un
peu curieuse provient sans  doute de la première apparition historique
du  problème~: en  1852 un  mathématicien amateur,  Francis Guthrie,
avait demandé  combien de  couleurs étaient nécessaires  pour colorier
les  pays  sur  un  atlas   afin  que  les  cartes  fussent  lisibles,
c'est-à-dire que deux  pays ayant une frontière commune  ne soient pas
de  la même  couleur.  Cette  question difficile  n'a reçu  sa réponse
définitive qu'en 1977,  quand Appel et Haken ont  démontré que dans ce
cas  particulier  des  graphes  planaires, quatre  couleurs  suffisent
toujours.   Sur la  figure ci-dessous,  on montre  une  coloration des
graphes de la figure~1 avec un nombre minimal de couleurs~:

\begin{figure}[ht]
  \center
  \begin{tabular}{cc}
    \includegraphics{fig.base.27} & \includegraphics{fig.base.28} \\
    \parbox{5cm}{\center Ce graphe nécessite 5 couleurs}&
    \parbox{5cm}{\center Ce graphe nécessite 3 couleurs}
  \end{tabular}
  \caption{Exemples de graphes coloriés optimalement}
\end{figure}

Le  problème général  de la  coloration\index{coloration}  des graphes
consiste donc en l'affectation à chaque sommet d'une couleur, de sorte
que  deux sommets reliés  par une  arête ne  reçoivent jamais  la même
couleur, et en utilisant un nombre minimal de couleurs.  Ce problème a
de nombreuses applications  pratiques~: conception d'emplois du temps,
affectation   de  fréquences   radio~\dots\  Depuis   une  quarantaine
d'années, on soupçonne qu'il  n'existe aucune méthode pour le résoudre
efficacement et en toute  généralité, même à l'aide d'ordinateurs très
rapides.   Pour contourner cette  difficulté, plusieurs  approches ont
été  essayées~:   résolution  du  problème  de   manière  rapide  mais
approximative, accéléreration des méthodes lentes, ou enfin résolution
efficace en se restreignant  à des graphes particuliers. Notre travail
se place dans le cadre de cette dernière approche.

\section*{Les graphes parfaits}

Au  début  des  années 1960,  Claude  Berge  a  défini la  classe  des
\emph{graphes parfaits} qui est apparue  par la suite comme une classe
assez générale pour laquelle le problème de la coloration pouvait être
résolu  efficacement.   De  nombreux  travaux ont  été  consacrés  aux
graphes parfaits~: on rencense  plus de 500 articles scientifiques sur
le sujet.  En 1984, Martin Grötschel, L\'aszl\'o Lov\'asz et Alexander
Schrijver  ont démontré  qu'il est  possible de  colorier  les graphes
parfaits en  temps polynomial.  Ils  utilisent à cette fin  la méthode
dite des ellipsoïdes,  difficile à implémenter en pratique  à cause de
problèmes  d'instabilité numérique. 

En même  temps que Claude  Berge définissait les graphes  parfaits, il
proposait  une  conjecture   les  caractérisant  par  interdiction  de
sous-graphes  induits~:  la  conjecture  forte des  graphes  parfaits.
Celle-ci s'est avérée récalcitrante à  la preuve~: pendant 40 ans, les
approches les plus diverses ont été tentées sans succès.  L'une d'elle
fut promue  par Va\v sek  Chv\'atal~: pour démontrer la  perfection de
certains types de graphes, il proposait de les ``casser'' en plusieurs
sous-graphes plus simples par  des strucures ne pouvant pas apparaître
dans des  graphes imparfaits minimaux, voire quand  c'est possible, de
``recoller  les morceaux''  par des  moyens préservant  la perfection.
Cette approche  a permis à  de nombreux auteurs, (Michel  Burlet, Jean
Fonlupt, et d'autres ont suivi)  de prouver la perfection de certaines
classes de graphes,  mais ce n'est qu'à la fin  des années 1990, qu'un
groupe  de trois  chercheurs (Michele  Conforti, Gérard  Cornuéjols et
Kristina  Vu\v skovi\'c)  a proposé  un plan  précis et  réaliste pour
appliquer cette méthode en  toute généralité, espérant ainsi démontrer
la  conjecture forte des  graphes parfaits.   Ces trois  chercheurs ne
sont toutefois  pas parvenus  à faire aboutir  leurs idées,  malgré de
nombreux  résultats.  Lorsque j'ai  commencé ma  thèse, les  choses en
étaient là~: la  conjecture forte des graphes parfaits  était l'un des
grands problèmes ouverts de la théorie des graphes.

En  mai  2002,  un   groupe  de  chercheurs  (Maria  Chudnovsky,  Neil
Robertson,  Paul Seymour  et Robin  Thomas) a  démontré  la conjecture
forte  des  graphes  parfaits  en  suivant  la  méthode  proposée  par
Cornuéjols {\it et al.}. Cette percée théorique n'a pas eu de retombée
directe  sur le  problème de  la coloration  des graphes  parfaits. En
novembre 2002,  Chudnovsky, Cornuéjols, Liu, Seymour  et Vu\v skovi\'c
ont résolu une  autre question ouverte~: ils ont  décrit un algorithme
en temps polynomial qui  décide si un graphe est parfait.

C'est une chance exceptionnelle d'avoir  pu réaliser une thèse dans un
climat de recherche aussi actif et riche en idées nouvelles, au moment
précis  où  des questions  difficiles  et  anciennes trouvaient  leurs
réponses.  Les nouvelles idées  dégagées par l'équipe de Paul Seymour,
associées  aux  tentatives  précédentes  de mon  directeur  de  thèse,
Frédéric Maffray, nous ont  permis de résoudre une conjecture proposée
au  début des  années 1990  par Hazel  Everett et  Bruce  Reed.  Cette
conjecture  affirme  qu'une  sous-classe  des  graphes  parfaits,  les
graphes  d'Artémis,  peut être  coloriée  rapidement  en utilisant  la
notion de  paire d'amis d'un  graphe.  Nous avons également  trouvé un
algorithme de reconnaissance pour les graphes d'Artémis, et aussi pour
d'autres classes  similaires. On trouvera ci-après un  résumé de cette
thèse.

\vspace{1ex}

L'ensemble de ce travail a  été réalisé en collaboration avec Frédéric
Maffray.

\vspace{1ex}

\begin{description}
\item[Chapitre  1  :  Notions  de  base]\mbox{}

  Nous rappelons toutes les  définitions de théorie des graphes utiles
  pour  la suite.  Nous  donnons un  aperçu rapide  et informel  de la
  théorie de la complexité  algorithmique. Le lecteur connaissant bien
  ces sujets  pourra consulter la  dernière section où  sont rappelées
  les quelques définitions utiles à notre travail ne figurant pas dans
  les ouvrages classiques.  L'index en fin de document permet aussi de
  retrouver les définitions de tous les termes employés.

\item[Chapitre 2  : Graphes parfaits]\mbox{}

  Nous donnons  un historique des principaux  résultats concernant les
  graphes parfaits,  des premiers travaux  de Claude Berge  jusqu'à la
  preuve de  la conjecture forte des graphes  parfaits par Chudnovsky,
  Robertson, Seymour et Thomas.

\item[Chapitre  3 : Paires d'amis]\mbox{}

  Nous  rappelons les définitions  et résultats  principaux concernant
  une notion importante pour les algorithmes efficaces de coloration~:
  la notion de  paire d'amis d'un graphe.  Nous  donnons un algorithme
  original de  détection des paires d'amis dans  les line-graphes.  La
  dernière   section  présente   des  résultats   originaux   sur  des
  généralisations possibles de la notion de paire d'amis.

\item[Chapitre 4 : Le lemme de Roussel et Rubio]\mbox{}

  Nous  rappelons  un résultat  important~:  le  lemme  de Roussel  et
  Rubio.  Nous en  donnons  une preuve  originale,  puis nous  donnons
  quelques corollaires et variantes qui nous serviront par la suite.

\item[Chapitre  5 : Graphes d'Artemis]\mbox{}

  Nous prouvons  une conjecture de  Reed et Everett,  affirmant qu'une
  certaine  classe de  graphes parfaits  (les graphes  d'Artémis) peut
  être coloriée  efficacement en utilisant la notion  de paires d'amis
  (complexité~: $O(mn^2)$).  Nous donnons un théorème de décomposition
  des  graphes  d'Artémis  (nous  n'avons pas  démontré  nous-même  ce
  théorème~: il est énoncé implicitement et démontré dans la preuve de
  la conjecture forte des  graphes parfaits par Chudnovsky, Robertson,
  Seymour et Thomas).   Les résultats de ce chapitre  sont en révision
  favorable   au  {\it   Journal  of   Combinatorial   Theory,  Series
  B}~\cite{nicolas:artemis}.

\item[Chapitre 6 : Problèmes de reconnaissance]\mbox{}

  Grâce à des  techniques utilisées par Chudnovsky et  Seymour pour la
  reconnaissance des  graphes de Berge, nous donnons  un algorithme en
  temps  polynomial  de reconnaissance  des  graphes d'Artémis.   Nous
  donnons également des algorithmes pour les autres classes de graphes
  intéressantes du  point de  vue de la  notion de paire  d'amis. Nous
  montrons que nos algorithmes de reconnaissance peuvent être utilisés
  à  des fins  de coloration.   Nos  algorithmes sont  fondés sur  des
  méthodes  de détection de  sous-graphes induits  par calcul  de plus
  courts chemins dans les graphes sans trou impair.  Nous montrons que
  ces mêmes  problèmes de détection sont NP-complets  dans les graphes
  généraux.  Les  résultats de  ce chapitre sont  soumis au  {\it SIAM
  Journal of Discrete Mathematics}~\cite{nicolas:reco}.

\end{description}


\chapter{Notions de base}
Nous présentons ici les notions  de bases utilisées dans cette thèse~:
les graphes et  les algorithmes.  Le lecteur habitué  au sujet peut se
dispenser  de lire  ce chapitre,  mais devrait  consulter  sa dernière
section  qui rappelle  les  rares points  où  nos conventions  peuvent
paraître inhabituelles. Ce chapitre n'a aucune prétention pédagogique,
et bien que toutes les définitions  de théorie des graphes utiles à la
compréhension de cette thèse soient données rigoureusement, le lecteur
novice     préfèrera    sans     doute     consulter    un     ouvrage
classique~\cite{berge:85,  diestel:graph}.    Pour  ce  qui   est  des
algorithmes, nous  ne donnons pas  d'exposé rigoureux, mais  un simple
aperçu des définitions et des résultats les plus classiques.

\section{Les graphes}

Si $V$  est un ensemble  quelconque et si  $k$ est un  entier positif,
alors  on  note $V  \choose  k$\index{0@\rule{0cm}{3ex}$V \choose  k$}
l'ensemble des  parties de  $V$ qui ont  exactement $k$  éléments.  Un
\emph{graphe}\index{graphe} est un couple d'ensembles $(V, E)$ tel que
$E\subset  {V \choose  2}$.   Les  éléments de  $V$  sont appelés  les
\emph{sommets}\index{sommet}  de  $G$, et  les  éléments  de $E$  sont
appelés  les \emph{arêtes}\index{arête}  de $G$.   Afin  d'alléger les
notations, on  notera $uv$ la paire  contenant $u$ et $v$  (au lieu de
$\{u, v\}$).   Si $G$ est un  graphe, alors $V(G)$\index{V@$V(\dots)$}
désignera  l'ensemble  de  ses sommets  et  $E(G)$\index{E@$E(\dots)$}
l'ensemble de ses arêtes.

Soit $G$ un graphe.  Soient $u$ et $v$ deux sommets de $G$. Si $uv \in
E(G)$     alors     on     dit     que     $u$     et     $v$     sont
\emph{adjacents}\index{adjacent},   que   $u$    et   $v$   sont   les
\emph{extrémités}\index{extrémités!d'une  arête} de l'arête  $uv$, que
$u$ est un  \emph{voisin}\index{voisin} de $v$~; on dit  aussi que $u$
\emph{voit}\index{voir} $v$. Si $uv \notin E(G)$, alors on dit que $u$
et $v$ sont \emph{non adjacents},  que $u$ est un \emph{non-voisin} de
$v$~; on  dit aussi que $u$ \emph{manque}\index{manquer}  $v$.  Si $u$
est   un   sommet   et  $e$   une   arête,   on   dit  que   $e$   est
\emph{incidente}\index{incident}  à  $u$ si  et  seulement  si $u  \in
e$. Soient $e$ et $f$ deux arêtes  de $G$.  On dit que $e$ et $f$ sont
\emph{incidentes}  si et  seulement si  $e$ et  $f$ ont  une extrémité
commune, autrement dit si $e\cap f \neq \emptyset$.

Notons   que,   par    définition,   nos   graphes   sont   \emph{sans
boucle}\index{boucle},  c'est-à-dire qu'un sommet  $v$ ne  peut jamais
être son propre voisin (car $\{v,v\}$ ne possédant qu'un seul élément,
il n'appartient pas à ${V  \choose 2}$). Notons que par définition nos
graphes ne  sont pas \emph{orientés}\index{orienté},  c'est-à-dire que
pour tous sommets $u$ et $v$, $uv = vu$.  Notons enfin que nos graphes
sont \emph{simples}\index{simple},  c'est-à-dire qu'entre deux sommets
$u$  et $v$,  il existe  0 ou  1 arête,  et non  un  nombre quelconque
d'arêtes.  Notons  enfin que tous  les graphes que  nous considérerons
seront  \emph{finis}\index{fini},  c'est-à-dire  que nous  supposerons
toujours  que $V$  possède un  nombre fini  d'éléments.  On  notera en
général  $n$ le  nombre  de sommets  d'un  graphe, et  $m$ son  nombre
d'arêtes.   On  appelle  \emph{taille}\index{taille!d'un~graphe}  d'un
graphe son nombre de sommets.

Soit   $G$   un   graphe.   Si   $v   \in   V(G)$,   alors   on   note
$N(v)$\index{N@$N(\dots)$} l'ensemble des  voisins de $v$.  Notons que
$v \notin  N(v)$.  On note $d(v)$ le  \emph{degré}\index{degré} de $v$
défini par $d(v) = |N(v)|$.  Si $A \subset V(G)$, alors $N(A)$ désigne
l'ensemble  des sommets  de $G$  ayant au  moins un  voisin  dans $A$.
Notons que $A$ et $N(A)$ ne sont pas nécessairement disjoints.

Soit         $G$         un         graphe.         On         appelle
\emph{complémentaire}\index{complémentaire}  de  $G$  le  graphe  noté
$\overline{G}$\index{0@$\overline{\dots\rule{0cm}{1ex}}$~dans~$\overline{G}$}
défini par $V(\overline{G}) = V(G)$ et $E(\overline{G}) = \{uv ; u \in
V(G) \text{ et } v\in V(G)  \text{ et } uv \notin E(G)\}$.  Soient $F$
et         $G$        deux         graphes.          On        appelle
\emph{isomorphisme}\index{isomorphisme}   entre  $F$   et   $G$  toute
fonction  $\varphi$ bijective,  associant à  chaque sommet  de  $F$ un
sommet   de  $G$   et   telle  que   $uv   \in  E(F)   \Leftrightarrow
\varphi(u)\varphi(v) \in E(G)$.  S'il existe un isomorphisme entre $F$
et  $G$, alors  on  dit que  $F$  et $G$  sont \emph{isomorphes}.   On
appelle   graphe   \emph{autocomplémentaire}\index{autocomplémentaire}
tout graphe isomorphe à son complémentaire.

Soient $G$  un graphe et $W$ un  ensemble de sommets de  $G$.  On note
$G[W]$\index{0@$[\dots]$  dans $G[\dots]$}  le graphe  dont l'ensemble
des sommets est  $W$ et dont l'ensemble des arêtes  est ${W \choose 2}
\cap E(G)$.  On dit que  $G[W]$ est le \emph{sous-graphe de $G$ induit
par $W$}\index{induit}\index{sous-graphe}.  Pour ne pas trop alourdir
les notations, on s'autorise certains abus~: si $F$ est un sous-graphe
induit de  $G$, on note $G\setminus  F$ le graphe  qu'on devrait noter
$G[V(G) \setminus V(F)]$\index{0@$\setminus$, dans $G\setminus F$}. Si
$F$  est un  ensemble de  sommets de  $G$, on  note $G\setminus  F$ le
graphe qu'on devrait noter $G[V(G) \setminus F]$.

On     appelle     \emph{clique}\index{clique}     ou     \emph{graphe
complet}\index{complet!graphe ---} tout graphe vérifiant $E(G) = {V(G)
\choose 2}$.  On  appelle \emph{triangle}\index{triangle} toute clique
ayant trois  sommets.  Il est coutumier  de noter $K_n$\index{K@$K_n$}
une clique  avec $n$ sommets.   On appelle \emph{stable}\index{stable}
tout  graphe  dont  le  complémentaire  est une  clique.   On  appelle
\emph{clique de $G$} tout graphe $K$  qui est une clique et qui est un
sous-graphe  induit  de $G$.   Plus  généralement,  chaque fois  qu'on
définira un type de graphe ``bidule'' (clique, stable, trou~\dots), on
appellera \emph{``bidule''  de $G$} tout  graphe qui est  un ``bidule''
(ou  qui est  isomorphe à  un ``bidule'')  et qui  est  un sous-graphe
induit  de $G$.   On dira  que $G$  \emph{contient}\index{contient} un
``bidule'' si  et seulement s'il  existe un sous-graphe induit  de $G$
qui  est  un  ``bidule'' de  $G$.   On  dira  que $G$  est  \emph{sans
``bidule''}\index{sans} si et seulement s'il n'existe aucun ``bidule''
qui soit un sous-graphe induit  de $G$.
 
Soit   $v$    un   sommet    de   $G$.    On    dit   que    $v$   est
\emph{simplicial}\index{simplicial} si  et seulement si  $N(v)$ induit
une    clique   de    $G$.    On    dit   qu'un    graphe    $G$   est
\emph{biparti}\index{biparti}  si et  seulement s'il  est  possible de
partitionner $V(G)$  en deux ensembles $A$  et $B$ tels  que $G[A]$ et
$G[B]$ sont  des stables. Si  de plus il  existe une arête  entre tout
sommet de  $A$ et  tout sommet  de $B$, alors  on dit  que $G$  est un
graphe  \emph{biparti  complet}\index{complet!biparti~---}.   Dans  ce
cas, la bipartition  est unique. Si $G$ est  un graphe biparti complet
avec  $|A|  =  a$ et  $|B|  =  b$,  alors  on le  note  habituellement
$K_{a,b}$\index{K@$K_{a,b}$}.

Soit $k$ un entier. On appelle $k$-\emph{coloration}\index{coloration}
de $G$  toute partition de $V(G)$  en $k$ ensembles  $A_1, \dots, A_k$
tels que pour  tout $1\leq i \leq k$, $G[A_i]$  est un stable.  Notons
que  tout graphe  avec $n$  sommets possède  une  $n$-coloration.  Les
ensembles $A_i$  sont appelés des  \emph{couleurs}\index{couleur}.  Si
$v\in A_i$, on dit qu'on a \emph{donné à $v$ la couleur $i$}. Parfois,
on désigne les  numéros des couleurs par des  adjectifs bizarres comme
``rouge'', ``bleu''  \dots\ ou même ``rose'' ou  ``vert''.  On appelle
\emph{nombre chromatique}\index{nombre chromatique}\index{chromatique}
de  $G$   le  plus   petit  entier  $k$   tel  que  $G$   possède  une
$k$-coloration.   On  note  $\chi(G)$\index{0chi@$\chi(G)$} le  nombre
chromatique de $G$.  Une coloration  de $G$ avec $\chi(G)$ sommets est
dite    \emph{optimale}\index{optimal!coloration   ---}.     On   note
$\omega(G)$\index{0omega@$\omega(G)$} la taille d'une clique de $G$ de
taille  maximale.   On  note $\alpha(G)$\index{0alpha@$\alpha(G)$}  la
taille d'un  stable de  $G$ de taille  maximale. Le lemme  suivant est
trivial~:

\begin{lemme} \mbox{}
\label{base.l.trivbip}
Soit $G$ un graphe.
\begin{itemize}
\item 
  $G$ est un stable si et seulement si $\chi(G) = 1$.
\item 
  $G$ est un stable si et seulement si $\omega(G) = 1$.
\item 
  $G$ est une clique si et seulement si $\chi(G) = n$. 
\item 
  $G$ est une clique si et seulement si $\omega(G) = n$. 
\item 
  $G$ est biparti si et seulement si $\chi(G) \leq 2$.
\end{itemize}
\end{lemme}


\vspace{2ex}

Soit $G$ un graphe.  On appelle \emph{line-graphe} \index{line-graphe}
\index{L@$L(\dots)$} de $G$ le  graphe noté $L(G)$ dont l'ensemble des
sommets  est  $E(G)$,  et  dont  deux sommets  sont  adjacents  si  et
seulement  si les  deux  arêtes  de $G$  qui  leur correspondent  sont
incidentes.  Formellement,  $V(L(G)) = E(G)$ et  $E(L(G)) = \{  ef ; e
\in E(G) \text{ et } f \in E(G) \text{ et } e\cap f \neq \emptyset\}$.

Si $G$ est un  graphe et si $v$ est un sommet  de $G$ ayant exactement
deux   voisins  non   adjacents   $v'$  et   $v''$,   alors  on   note
$G/v$\index{0@$/$, dans  $G/v$} le graphe dont  l'ensemble des sommets
est $V(G)  \setminus \{v\}$, et  dont l'ensemble des arêtes  est $(E(G)
\setminus \{vv',  vv''\}) \cup \{v' v''\}$.   Intuitivement, le graphe
$G/v$ est le graphe obtenu en ``lissant'' le sommet $v$. Si $G$ et $F$
sont    des   graphes,    alors    on   dit    que    $G$   est    une
\emph{subdivision}\index{subdivision}  de  $F$  si et  seulement  s'il
existe une suite de graphes $G =  G_0$, $G_1$, \dots, $G_k = F$ et une
suite $v_0, \dots v_{k-1}$ tel que pour tout $0\leq i \leq k-1$, $v_i$ est
un sommet de $G_i$ de degré 2 et $G_{i+1} = G_i/v_i$.

\index{représentation  d'un graphe}Il  est coutumier  de  dessiner les
graphes  en  associant à  chaque  sommet un  point  du  plan que  l'on
symbolise par un petit cercle (vide ou plein), puis en reliant par une
ligne  toutes  les  paires  de   points  qui  sont  des  arêtes.  Pour
représenter un graphe biparti, on représente parfois les sommets de la
première  partie  de  la  partition  en ``plein''  et  les  autres  en
``creux''.

\subsection{Chaînes,  chemins, cycles, trous et antitrous}
\label{base.ss.ccctt}

Dans tout ce travail, il  sera beaucoup question de ``chemins''. Or il
existe en  théorie des graphes  plusieurs notions de  chemins, variant
selon  les auteur. Il  nous semble  qu'aucune notation  n'est vraiment
standard. Nous présentons ici nos propres notations, que nous espérons
claires, intuitives et rigoureuses.

Soit  $G$ un  graphe.  On  appelle \emph{chaîne}\index{chaîne}  de $G$
toute  suite $C  = (v_1,  v_2,  \dots, v_k)$  de sommets  deux à  deux
distincts de $G$ vérifiant pour tout $1\leq i \leq k-1$~: $v_i v_{i+1}
\in E(G)$. Les arêtes de $G$ de  la forme $v_i v_{i+1}$ avec $1 \leq i
\leq   k-1$,    sont   alors   appelées   les    \emph{arêtes   de   la
chaîne}\index{arête!d'une chaîne}. Les arêtes  de $G$ de la forme $v_i
v_j$  avec   $|i-j|  >  1$   sont  appelées  les  \emph{cordes   de  la
chaîne}\index{corde!d'une  chaîne}.  Les sommets  $v_1$ et  $v_k$ sont
appelés  les   \emph{extrémités  de  la  chaîne}\index{extrémités!d'une
chaîne}.   La longueur\index{longueur!d'une chaîne}  de la  chaîne est
égale à son nombre d'arêtes.
  
On  appelle \emph{cycle}\index{cycle} de  $G$ toute  suite $C  = (v_1,
v_2, \dots,  v_k)$ de sommets deux  à deux distincts  de $G$ vérifiant
pour tout  $1 \leq i \leq  k$~: $v_i v_{i+1} \in  E(G)$, où l'addition
des indices  est entendue modulo $n$.   Les arêtes de $G$  de la forme
$v_i  v_{i+1}$  avec  $1  \leq  i  \leq k$,  sont  alors  appelées  les
\emph{arêtes du  cycle}\index{arête!d'un cycle}. Les arêtes  de $G$ de
la forme $v_i  v_j$ avec $|i-j| > 1$ sont  appelées les \emph{cordes du
cycle}\index{corde!d'un cycle}.

On appelle \emph{chemin}\index{chemin} tout graphe $P$ dont l'ensemble
des  sommets peut être  ordonné de  manière à  former une  chaîne sans
corde de  $P$, et on  appelle \emph{trou}\index{trou} tout  graphe $H$
avec au moins quatre sommets dont les sommets peuvent être ordonnés de
manière à former  un cycle sans corde de $H$.   Notons qu'il existe au
plus deux  manières d'ordonner  les sommets d'un  chemin de  manière à
former une chaîne.  Pour ces ordres, les extrémités sont les mêmes, et
on les  appelle les \emph{extrémités  du chemin}\index{extrémités!d'un
chemin}~;    les    autres     sommets    du    chemin    sont    dits
\emph{intérieurs}\index{intérieur}.               On              note
$P^*$\index{0@$*$,~dans~$P^*$}  le sous-chemin de  $P$ induit  par les
sommets intérieurs de $P$.  On note $\lg(P)$\index{L0@$\lg(\dots)$} la
\emph{longueur}\index{longueur!d'un  chemin} de $P$,  c'est-à-dire son
nombre d'arêtes.   On appelle \emph{antichemin}\index{antichemin} tout
graphe    qui    est     le    complémentaire    d'un    chemin,    et
\emph{antitrou}\index{antitrou} tout graphe  qui est le complémentaire
d'un trou.  On dit qu'un  trou (ou un antitrou) est \emph{pair} (resp.
\emph{impair}) s'il contient un  nombre pair d'arêtes (resp.  impair).
On dit  qu'un trou ou  un antitrou est  \emph{long}\index{long!trou ou
antitrou ---} s'il contient au moins 5 sommets.
 
Si $n\geq 1$ est un entier et  si $v_1, v_2 \dots v_n$ sont des objets
mathématiques quelconques, on note $\bp v_1 \tp v_2 \tp \cdots \tp v_n
\ep$ le  graphe dont  les sommets sont  $v_1$, $v_2$, \dots,  $v_n$ et
dont  l'ensemble des  arêtes est  $\{ v_i  v_{i+1}; 1\leq  i  \leq n-1
\}$. Le graphe  $\bp v_1 \tp v_2  \tp \cdots \tp v_n \ep$  est donc un
chemin d'extrémités $v_1$, $v_n$ et  de longueur $n-1$. Si $P$ est un
chemin et si $u$, $v$ sont des sommets de $P$, alors, on note $\bp u
\tp  P \tp  v  \ep$\index{0@$-$, dans  $\bp u  \tp  P \tp  v \ep$}  le
\emph{sous-chemin}\index{sous-chemin} de $P$  d'extrémités $u$ et $v$,
c'est-à-dire  l'unique sous-graphe  induit de  $P$ qui  est  un chemin
d'extrémités   $u$    et   $v$.     Il   est   coutumier    de   noter
$P_n$\index{P@$P_n$} un  chemin avec $n$  sommets, $\overline{P_n}$ un
antichemin  avec $n$  sommets, $C_n$\index{C@$C_n$}  un trou  avec $n$
sommets et  $\overline{C_n}$ un antitrou avec $n$  sommets. On appelle
\emph{carré}\index{carré} le trou avec quatre sommets, noté $C_4$.

On définit une relation entre les sommets d'un graphe $G$~: on note $u
\sim v$ si et seulement s'il  existe un chemin de $G$ d'extrémités $u$
et  $v$. On  vérifie que  $\sim$ est  une relation  d'équivalence.  On
appelle  \emph{composantes   connexes}\index{composante}  de  $G$  les
graphes  induits par  les  classes d'équivalence  de  $\sim$.  Si  $G$
possède  une  seule  composante  connexe  alors on  dit  que  $G$  est
\emph{connexe}\index{connexe}. Si $\overline{G}$ est connexe, alors on
dit  que $G$  est  \emph{anticonnexe}\index{anticonnexe}.  On  appelle
\emph{composantes  anticonnexes}   de  $G$  les   complémentaires  des
composantes connexes  de $\overline{G}$.  Notons que pour  tout graphe
$G$, l'un  au moins  de $G$, $\overline{G}$  est connexe.   On appelle
\emph{ensemble   d'articulation}\index{articulation}   de   $G$   tout
ensemble de sommets $F$ tel  que $G\setminus F$ n'est pas connexe.  Si
pour tout ensemble $F$ de  $k-1$ sommets de $G$, le graphe $G\setminus
F$     est     connexe,    alors     on     dit     que    $G$     est
\emph{$k$-connexe}\index{connexe!$k$-connexe}.

Soit  $G$   un  graphe,  et  $u$,   $v$  deux  sommets   de  $G$.   La
\emph{distance}\index{distance}  de $u$  à  $v$ est  la longueur  d'un
chemin de  $G$ de longueur  minimale (on dit souvent  \emph{plus court
chemin})\index{chemin!plus~court~---}, d'extrémités  $u$ et  $v$.  Si
$P$  est un chemin  de $G$  ayant une  extrémité $u$,  si on  parle du
sommet $u'$  de $P$  \emph{le plus proche}  de $u$ ayant  une certaine
propriété,  alors il  s'agira du  sommet $u'$  minimisant  la distance
entre  $u$ et $u'$  calculée dans  le chemin  $P$ et  non pas  dans le
graphe $G$.

Nous terminons par un théorème bien connu~:

\begin{theoreme}[Cf \cite{diestel:graph}]
  Un graphe est biparti si et seulement s'il n'a pas de cycle impair.
\end{theoreme}

\subsection{Pyramides, prismes et quasi-prismes}
\label{base.ss.pps}

\begin{figure}[ht]
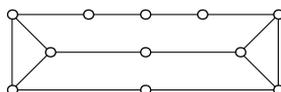
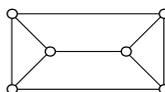

  \center
  \begin{tabular}{lc}
  Un prisme impair~: &\parbox[c]{7cm}{\center
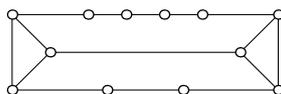
\includegraphics{fig.base.3}}\\
  \rule{0cm}{1.2cm}&\\
  Un prisme pair~: & \parbox[c]{7cm}{\center\includegraphics{fig.base.4}} \\
  \rule{0cm}{1.2cm}&\\
  Le plus petit prisme~: & \parbox[c]{7cm}{\center\includegraphics{fig.base.7}} \\
  \rule{0cm}{.5cm}&\\
  \end{tabular}
  \caption{Des prismes}
\end{figure}

On  appelle  \emph{prisme}\index{prisme}  tout  graphe  $F$  dont  les
sommets  se partitionnent  en  trois ensembles  induisant des  chemins
$P_1$, $P_2$, $P_3$  tels que pour $i=1, 2, 3$,  $P_i$ est de longueur
au moins 1, d'extrémités $a_i$ et $b_i$, tels que $\{ a_1, a_2, a_3\}$
et $\{b_1, b_2, b_3\}$ induisent  des triangles, et tels qu'il n'y ait
aucune arête  entre des sommets de $P_i$  et $P_j$ ($1\leq i  < j \leq
3$) autre  que celles  des triangles.  Si  $F$ est  sous-graphe induit
d'un graphe $G$,  alors on dit que les trois  chemins $P_1$, $P_2$, et
$P_3$ \emph{forment}\index{former!un prisme} un prisme de $G$.  On dit
que     les    sommets    $a_i,     b_i$,    $i=1,2,3$     sont    les
\emph{coins}\index{coin!d'un prisme}  du prisme.   On dit que  $F$ est
\emph{pair}\index{pair!prisme                ---}               (resp.
\emph{impair}\index{impair!prisme  ---}) si  les  trois chemins  $P_1,
P_2,  P_3$ sont  de longueur  paire  (resp.  impaire).   On dit  qu'un
prisme  est \emph{long}\index{long!prisme ---}  s'il possède  au moins
sept sommets,  c'est-à-dire si l'un  au moins  de ses chemins  est de
longueur supérieure ou égale à 2.

\begin{figure}[ht]
  \center
  \includegraphics{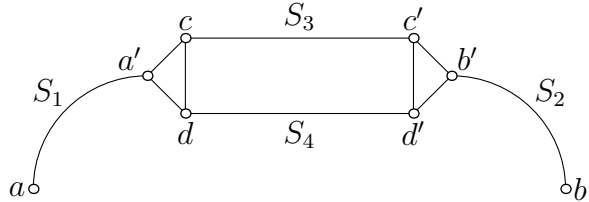}
  \caption{Un quasi-prisme\label{base.fig.quasiP}}
\end{figure}

On  appelle \emph{quasi-prisme}\index{quasi-prisme}\index{prisme!quasi
  - ---} tout graphe  $S$ dont les sommets se  partitionnent en quatre
ensembles induisant  des chemins  $S_1= \bp a  \tp \cdots  \tp a'\ep$,
$S_2= \bp b \tp \cdots \tp b'  \ep$, $S_3=\bp c \tp \cdots \tp c' \ep$
et $S_4= \bp d \tp \cdots \tp  d' \ep$, où $S_1$ et $S_2$ peuvent être
de longueur~$0$  mais où $S_3, S_4$ sont de longueur  au moins $1$, et
tels que  $E(S)= E(S_1)\cup E(S_2) \cup E(S_3)  \cup E(S_4) \cup\{a'c,
a'd, cd, b'c', b'd', c'd'\}$. Notons que $\{a', c , d\}$ et $\{b', c',
d'\}$  sont  des  triangles  de  $S$.   On  appelle  $a$  et  $b$  les
\emph{extrémités}\index{extrémités!d'un quasi-prisme} du quasi-prisme.
Le  quasi-prisme est dit  \emph{strict}\index{strict!quasi-prisme ---}
si l'un au moins de $S_1$ et $S_2$ est de longueur au moins $1$.

\begin{figure}[ht]
  \center
  \begin{tabular}{lc}
  Une pyramide~: & \parbox[c]{7cm}{\center\includegraphics{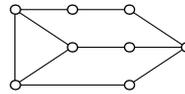}} \\
  \rule{0cm}{1.2cm}&\\
  La plus petite pyramide~: & \parbox[c]{7cm}{\center\includegraphics{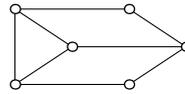}} \\
  \rule{0cm}{.5cm}&\\
  \end{tabular}
  \caption{Des pyramides}
\end{figure}

On  appelle  \emph{pyramide}\index{pyramide}   tout  graphe  $F$  dont
l'ensemble des  sommets est l'union  de trois ensembles  induisant des
chemins $P_1$,  $P_2$, $P_3$ tels que  pour $i=1, 2, 3$,  $P_i$ est de
longueur au  moins 1, d'extrémités $a$  et $b_i$, et  tels que $\{b_1,
b_2,  b_3\}$ induit  un triangle,  tels  que $a$  est l'unique  sommet
commun aux  trois chemins, et tels  qu'il n'y ait  aucune arête entre
des sommets de $P_i$ et $P_j$ ($1\leq  i < j \leq 3$) autre que celles
du triangle et celles d'extrémité  $a$.  De plus, un seul chemin parmi
$P_1, P_2,  P_3$ est autorisé  à être de  longueur 1, les  deux autres
doivent être de longueur au  moins $2$.  Si $F$ est sous-graphe induit
d'un graphe $G$,  alors on dit que les trois  chemins $P_1$, $P_2$, et
$P_3$   \emph{forment}\index{former!une  pyramide}  une   pyramide  de
$G$. On  dit que le  sommet $a$ est  le \emph{coin}\index{coin!d'une
pyramide} de la pyramide.

Le lemme suivant  est souvent très pratique pour  prouver qu'un graphe
contient un trou impair~:

\index{pyramide!contient un trou impair}
\begin{lemme}
  \label{base.l.pyramide}
  Soit $F$ une pyramide.  Alors $F$ contient un trou impair.
\end{lemme}

\begin{preuve}
  Il existe deux chemins de $F$ qui sont de même parité, sans perte de
  généralité, $P_1$  et $P_2$. Parmi eux, il  y en a au  maximum un de
  longueur 1 et $V(P_1) \cup V(P_2)$ induit un trou impair.
\end{preuve}

Le lemme suivant  montre que dans les graphes  sans trou impair, les
prismes sont tous pairs ou impairs~:

\begin{lemme}
  \label{base.l.prismepairimpair}
  Soit $F$  un prisme ni pair  ni impair.  Alors $F$  contient un trou
  impair.
\end{lemme}

\begin{preuve}
  Soit $F$ un prisme. Si l'un des chemins $P_1$  de $F$ est de
  longueur paire et  un autre chemin $P_2$ de  longueur impaire, alors
  $V(P_1) \cup V(P_2)$ induit un trou impair.
\end{preuve}

\subsection{Quelques graphes particuliers}
\label{base.ss.gparti}

Nous présentons ici quelques  graphes particuliers qui vont tous jouer
un rôle dans la suite de ce travail. De manière assez curieuse (et non
préméditée),  beaucoup de ces  graphes sont  autocomplémentaires. Pour
une étude  des graphes autocomplémentaires, nous renvoyons  à la thèse
de Master très complète de Alastair Farrugia~\cite{farrugia:these}.

\begin{figure}[ht]
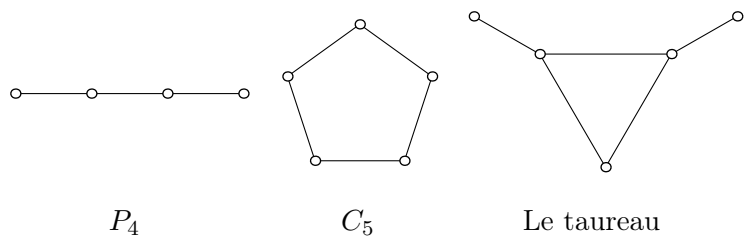

  \center
  \begin{tabular}{ccc}
  \parbox[c]{0cm}{\rule{0cm}{3cm}}
  \includegraphics{fig.base.11}&
  \parbox[c]{2cm}{\includegraphics{fig.base.12}}&
  \parbox[c]{3.3cm}{\includegraphics{fig.base.13}}\\
   $P_4$&  $C_5$& Le taureau
  \end{tabular}
\caption{Les   trois   graphes  autocomplémentaires   avec   4  ou   5
sommets\label{base.fig.45}}
\end{figure}

\begin{figure}[ht]
  \center
  \includegraphics{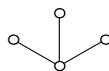}
\caption{La griffe\label{base.fig.griffe}}
\end{figure}

Le graphe $\bp a \tp b \tp c \tp d \ep$ est souvent appelé $P_4$ (voir
figure~\ref{base.fig.45}).     C'est     le    plus    petit    graphe
autocomplémentaire (à part le cas  trivial du graphe vide et du graphe
réduit à un sommet). C'est aussi, avec le carré et son complémentaire,
le plus grand  graphe qui soit à la fois  biparti et complémentaire de
biparti. Notons qu'il est aussi le line-graphe d'un graphe biparti, et
le   complémentaire   du   line-graphe   d'un  graphe   biparti.    La
\emph{griffe}\index{griffe},  parfois  notée  $K_{1,3}$, est  le  plus
petit    graphe   qui    ne   soit    pas   un    line-graphe~:   voir
figure~\ref{base.fig.griffe}.

Le  $C_5$ est  représenté figure~\ref{base.fig.45}.   C'est  un graphe
autocomplémentaire.    C'est   le   plus   petit  trou   impair.    Le
\emph{taureau}\index{taureau}  est  le  graphe contenant  un  triangle
représenté   figure~\ref{base.fig.45}.    C'est   encore   un   graphe
autocomplémentaire.   Le taureau et  le $C_5$  sont les  seuls graphes
autocomplémentaires avec cinq sommets (voir \cite{farrugia:these}).

\begin{figure}[ht]
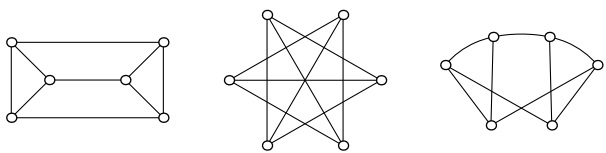

  \center
  \begin{tabular}{ccccc}
   \parbox[c]{0cm}{\rule{0cm}{2cm}}
   \parbox[c]{2cm}{\includegraphics{fig.base.7}}&\parbox{1cm}{}&
   \parbox[c]{2cm}{\includegraphics{fig.base.14}}&\parbox{1cm}{}& 
   \parbox[c]{2cm}{\includegraphics{fig.base.15}}
  \end{tabular}
\caption{Trois représentations de $\overline{C_6}$\label{base.fig.c6b}}
\end{figure}

Le graphe $\overline {C_6}$ est représenté de diverses manières sur la
figure~\ref{base.fig.c6b}.    La  première   représentation   le  fait
apparaître comme le plus petit prisme. Sur la deuxième représentation,
on le voit comme le plus petit antitrou pair connexe. Sur la troisième
représentation, on  le voit apparaître  comme un chemin de  longueur 3
(en haut)  dont les extrémités  voient deux sommets non  adjacents (en
bas), et cela sans que le chemin de longueur 3 ne possède d'arête dont
les extrémités voient les deux sommets en bas.

Le \emph{diamant}\index{diamant}  est le graphe à  quatre sommets dont
le    complémentaire   ne   contient    qu'une   seule    arête.    Le
\emph{double-diamant}\index{double-diamant}  est le  graphe  avec huit
sommets   représenté   figure~\ref{base.fig.dd}.    C'est  un   graphe
autocomplémentaire.

\begin{figure}[ht]
  \center
  \includegraphics{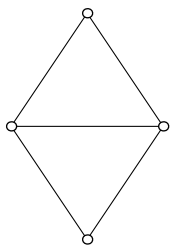}
  \caption{Le diamant\label{base.fig.diam}}
\end{figure}

\begin{figure}[ht]
  \center
  \includegraphics{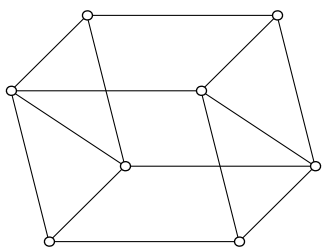}
  \caption{Le double-diamant\label{base.fig.dd}}
\end{figure}

\begin{figure}[h]
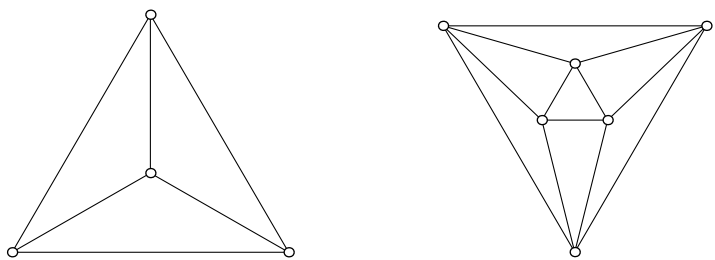

  \center
  \begin{tabular}{ccc}
  \includegraphics{fig.base.18}&\rule{1cm}{0cm}&
  \includegraphics{fig.base.19}
  \end{tabular}
  \caption{$K_4$ et $L(K_4)$\label{base.fig.k4}}
  \index{L0@$L(K_4)$}
\end{figure}

Figure~\ref{base.fig.k4},  on a représenté  $K_4$ et  son line-graphe,
noté $L(K_4)$.  Figure~\ref{base.fig.k33}, on a représenté $K_{3, 3}$.

\vspace{3cm}

\begin{figure}[h]
  \center
  \mbox{}\vspace{2cm}\\
  \includegraphics{fig.base.20}
  \caption{$K_{3,3}$\label{base.fig.k33}}
\end{figure}

\begin{figure}[h]
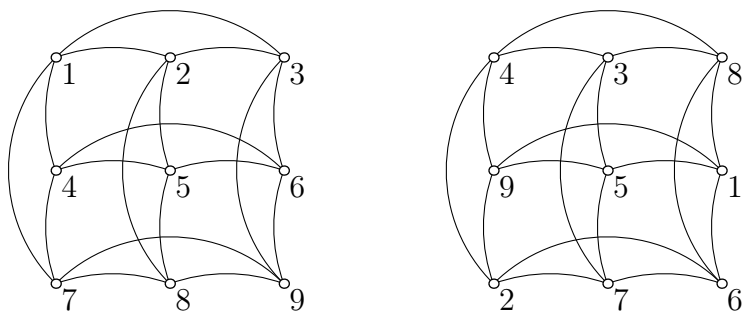

  \center
  \begin{tabular}{ccc}
  \includegraphics{fig.base.8}&\rule{1cm}{0cm}&
  \includegraphics{fig.base.9}
  \end{tabular}
  \caption{$L(K_{3,3})$ et son complémentaire\label{base.fig.lk33}}
  \index{L0@$L(K_{3,3})$}
\end{figure}

\begin{figure}[h]
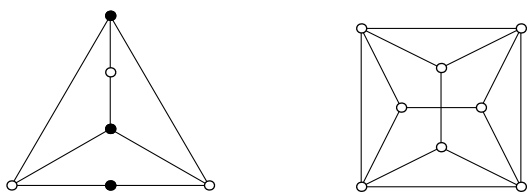

  \center
  \begin{tabular}{ccc}
  \includegraphics{fig.base.21}&\rule{1cm}{0cm}&
  \includegraphics{fig.base.22}
  \end{tabular}
  \caption{$K_{3,3}\setminus e$ et $L(K_{3,3}\setminus e)$\label{base.fig.lk33-e}}
  \index{L0@$L(K_{3,3} \setminus  e)$}
\end{figure}

Figure~\ref{base.fig.lk33},  on   a  représenté  le   line  graphe  de
$K_{3,3}$ avec une  numérotation naturelle des sommets.  À  côté, on a
représenté  son  complémentaire  de  manière  à  faire  apparaître  un
isomorphisme avec le graphe de départ, montrant ainsi que $L(K_{3,3})$
est  autocomplémentaire.  Il  faut  bien noter  que  les deux  graphes
représentés côte  à côte sont  \emph{différents}, que chaque  arête de
l'un est une non-arête de l'autre.  On remarquera (fait anecdotique~!)
que  les  numéros  des  sommets  du complémentaire  forment  alors  un
\emph{carré  magique}\index{carré!---  magique},  c'est-à-dire que  si
l'on  additionne les  trois chiffres  d'une colonne,  d'une  ligne, ou
d'une   diagonale  du   carré,  on   obtient  toujours   la  constante
``magique''~: 15.   Pour plus d'informations sur  les carrés magiques,
voir~\cite{descombes:carres}.     Il    est    habituel    de    noter
$K_{3,3}\setminus     e$      le     graphe     biparti     représenté
figure~\ref{base.fig.lk33-e}.  C'est le  plus petit graphe biparti qui
soit une subdivision de $K_4$.  Son line-graphe représenté sur la même
figure  est  un  graphe   autocomplémentaire,  que  l'on  désigne  par
$L(K_{3,3}\setminus e)$.  Beineke~\cite{beineke:comp} a prouvé que les
seuls  line-graphes autocomplémentaires sont  $K_1$, $P_4$,  $C_5$, le
taureau, $L(K_{3,3}\setminus e)$ et $L(K_{3,3})$.

\section{Les algorithmes}

\index{algorithme}
La notion d'algorithme est  malaisée à définir formellement, mais cela
est  tout  à   fait  possible  (cf.   \cite{papadimitriou:comp})~:  un
algorithme  est un objet  mathématique tout  aussi clair  qu'un nombre
entier ou un  graphe.  La notion de rapidité  d'un algorithme est très
importante  si l'on  envisage  des applications  pratiques.  De  plus,
depuis des temps très reculés,  on sait que la recherche d'algorithmes
rapides enrichit globalement les  mathématiques pures --- nous voulons
dire  ici  pures  d'algorithmes.   Citons comme  exemple  l'algorithme
d'Euclide pour  le calcul du PGCD de  deux entiers qui en  plus de son
intérêt pratique  permet de prouver des  théorèmes importants (Bezout,
théorème   de   Gauss   \dots).    Toutefois,   jusqu'au   milieu   du
XX\ieme~siècle, la  notion de rapidité d'un algorithme  était vague et
ne faisait  pas l'objet de  recherches pour elle-même.   Finalement le
développement de  la théorie de la complexité  algorithmique depuis le
milieu des années  1960 a placé la notion  de rapidité d'un algorithme
au c\oe ur de questions très difficiles et profondes.

\subsection{Problèmes, algorithmes, complexité}

\index{algorithme}
Un \emph{problème}\index{problème} est défini par la donnée d'un objet
mathématique   représenté    par   un   nombre    fini   de   symboles
(\emph{l'instance}\index{instance}  ou  \emph{l'entrée}\index{entrée})
et  d'une  \emph{question}\index{question}   dont  la  réponse  dépend
uniquement de l'instance. Si la réponse à la question ne peut être que
``oui''    ou    ``non'',    on    parle    de    \emph{problème    de
décision}\index{problème de  décision}\index{décision}.  Naïvement, un
algorithme   est  une   description   de  longueur   finie,  dans   un
\emph{langage  de  programmation},   d'une  suite  finie  d'opérations
élémentaires  (\emph{le  calcul}\index{calcul})  ne dépendant  que  de
l'instance  et  fournissant  une   réponse  à  la  question  (\emph{la
sortie}\index{sortie}).  Tous nos algorithmes auront pour instance des
graphes, ou des graphes avec quelques objets supplémentaires (sommets,
chemins du graphe \dots).

On               appelle              \emph{taille}\index{taille!d'une
instance}\index{instance!taille d'une ---}  de l'instance le nombre de
symboles nécessaires à sa  représentation.  Pour chaque algorithme, et
pour   chaque   entier   $n$,    nous   chercherons   à   évaluer   la
\emph{complexité}\index{complexité}  de l'algorithme,  c'est-à-dire le
nombre   d'opérations  élémentaires   nécessaires  à   l'exécution  de
l'algorithme pour  une instance  de taille $n$  dans le pire  des cas.
Nos estimations de complexité  seront toujours données à une constante
multiplicative     près    à    l'aide     de    la     notation    de
Landau\index{Landau~(notation~de~---)}~:
$O(\dots)$\index{O@$O(\dots)$}.    Dire   qu'un   algorithme  est   de
complexité  $O(n^{23})$  par  exemple,  veut  dire  qu'il  existe  une
constante  $c$   telle  que  le  nombre   d'opérations  pour  exécuter
l'algorithme sur une  instance de taille $n$ est  inférieur à~$c\times
n^{23}$.
 
Nous considérons  qu'un algorithme est \emph{efficace}\index{efficace}
s'il                 est                 en                \emph{temps
polynomial}\index{temps!polynomial}\index{polynomial}, c'est-à-dire si
sa  complexité  est  majorée  par   un  polynôme  dans  la  taille  de
l'instance.  Un  algorithme dont la complexité est  $O(2^n)$ sera donc
par définition non-efficace.  Le premier à avoir remarqué que beaucoup
d'algorithmes    classiques   s'exécutent    en    temps   polynomial,
indépendamment  de  tout   choix  raisonnable  de  représentation  des
données,  semble être  Cobham~\cite{cobham:comput},  mais l'idée  d'en
faire    un   critère    de    classification   est    due   à    Jack
Edmonds~\cite{edmonds:ptf}.   Ce dernier  a aussi  le  premier formulé
l'idée  souvent vérifiée que  la découverte  d'un algorithme  en temps
polynomial pour un problème va généralement de pair avec une meilleure
compréhension de celui-ci.

Pour         certains         problèmes         (les         problèmes
\emph{indécidables}\index{indécidable}), il n'existe aucun algorithme.
Mais pour  tous les  problèmes que nous  traiterons dans  cette thèse,
l'existence d'un  algorithme non-efficace sera  triviale.  Notre souci
sera donc de rechercher des algorithmes efficaces.

On  considère généralement que  la taille\index{taille!d'un  graphe en
tant qu'instance}  d'un graphe,  en tant qu'instance  d'un algorithme,
est la  somme de  son nombre  de sommets et  de son  nombre d'arêtes~:
$n+m$.  Si  l'on admet  qu'un algorithme dont  l'entrée est  un graphe
doit au  moins prendre  connaissance de l'instance,  on voit  que pour
pour tout  problème non-trivial, il est vain  d'espérer une complexité
meilleure que $O(n+m)$.  Un  algorithme ayant une telle complexité est
dit   \emph{en   temps   linéaire}\index{temps!linéaire}  ou   parfois
\emph{optimal}\index{optimal!algorithme  ---}.   Nous nous  soucierons
assez peu  de la représentation  de nos données.   Les représentations
classiques des  graphes par matrice d'incidence  ou liste d'adjacences
sont   essentiellement  équivalentes   pour  la   complexité   de  nos
algorithmes, dans la  mesure où nous ne donnerons  aucun algorithme en
temps linéaire.

Nous tenons  à préciser que notre  point de vue sur  la complexité des
algorithmes,  qui est  d'ailleurs  le point  de  vue conventionnel  et
majoritaire, est  par certains aspects critiquable.   Tout d'abord, un
algorithme mauvais  dans le  pire des cas  peut être  intéressant pour
beaucoup  d'instances, voire  en  moyenne.  De  plus,  même si  chacun
convient qu'un algorithme en temps exponentiel sera toujours plus lent
qu'un  algorithme en  temps polynomial  \emph{à partir  d'une certaine
taille $t$ de l'instance}, il  se peut que $t$ soit suffisamment grand
pour que seules  les instances de taille très  inférieure à $t$ soient
traitables  en  pratique  par  l'une  ou  l'autre  des  méthodes.   En
pratique,  un  algorithme  en  temps  $O({{2}^n})$  peut  donc  être
préférable à un algorithme en $O(n^{23})$ par exemple.  En outre, pour
les  applications  pratiques,  la  constante multiplicative  est  loin
d'être négligeable.  Finalement, la complexité algorithmique telle que
nous la donnons n'est qu'une indication de la performance effective de
l'algorithme~:  l'implémentation finale de  l'algorithme et  des tests
empiriques  sont  à  mon  avis  indispensables si  l'on  se  préoccupe
d'applications   concrètes.    En   dépit   de   ces   critiques,   la
classification  des  algorithmes présentée  ci-dessus  reste un  outil
simple  et  puissant,  qui  a  l'avantage  d'être  peu  dépendant  des
représentations de données et des langages de programmation.

Nous utiliserons  souvent l'algorithme fondamental  ci-dessous. Notons
que cet algorithme parfois  appelé ``parcours en largeur d'abord'' est
souvent attribué à  Dijkstra, parce qu'il a été  le premier à résoudre
un problème plus général (le chemin de poids minimum).  Nous renvoyons
au livre de  Alexander Schrijver (\cite{schrijver:opticombA}, page 87)
pour des éclaircissements historiques.

\index{chemin!plus~court~---}
\begin{algorithme}[Berge~\cite{berge:graphesapp}, Moore~\cite{moore:shortest}]
\label{base.a.pluscourt}
\begin {itemize}
\item[\sc Entrée~:] Un graphe $G$ et deux sommets $s$ et $t$ de $G$.

\item[\sc Sortie~:]  Un plus court  chemin de $G$ d'extrémités  $s$ et
$t$.

\item[\sc Complexité :] $O(n+m)$.
\end{itemize}

\end{algorithme}

\subsection{Les problèmes P, NP, CoNP et NP-complets}

\label{base.ss.npcomplet}

On dit qu'un problème appartient à la \emph{classe $P$} s'il existe un
algorithme en  temps polynomial pour le  résoudre. Un \emph{certificat
du oui}\index{certificat}  (resp.  du \emph{non}) pour  un problème de
décision est un  objet représentable par une suite  finie de symboles,
dépendant seulement de l'instance et  qui existe si et seulement si la
réponse  au  problème  est  ``oui'' (resp.   ``non'').   Un  \emph{bon
certificat}\index{bon~certificat} est un certificat dont la taille est
bornée  par  un  polynôme en  la  taille  de  l'instance, et  qui  est
\emph{vérifiable en temps polynomial}. C'est-à-dire qu'il doit exister
un algorithme en temps polynomial,  prenant en entrée l'instance et le
certificat, et qui  répond ``oui'' (resp. ``non'') si  et seulement si
l'entrée  est un  certificat du  oui  (resp.  du  ``non'') pour  cette
instance.   On dit  qu'un problème  appartient à  la  \emph{classe NP}
(resp.    \emph{CoNP})\index{NP}\index{Co-NP}  s'il  existe   pour  ce
problème un  bon certificat du ``oui'' (resp.   ``non'').  Les notions
de  certificats  pour  un   problème  ont  été  introduites  par  Jack
Edmonds~\cite{edmonds:matpar}.   Il  est clair  qu'un  problème P  est
toujours NP et CoNP.  Une conjecture célèbre affirme que la réciproque
est également vraie~:

\begin{conjecture}[P $=$ NP$\cap$coNP]
  Si un problème  admet un bon certificat du oui  et un bon certificat
  du non,  alors il existe un  algorithme en temps  polynomial pour le
  résoudre.
\end{conjecture}

Rappelons qu'une fonction booléenne\index{booléen} à $n$ variables est
une fonction $f$ de $\{0,  1\}^n$ dans $\{0, 1\}$.  Un vecteur booléen
$\xi\in\{0,    1\}^n$   \emph{satisfait}\index{satisfaire}    $f$   si
$f(\xi)=1$.   Pour toute  variable booléenne  $x$ sur  $\{0,  1\}$, on
écrit  $\overline{x}:=1-x$, et  on appelle  $x$ et  $\overline{x}$ des
\emph{littéraux}.   Une instance  de  {\sc $3$-sat}  est une  fonction
booléenne  $f$ donnée  comme  un produit  booléen  de clauses,  chaque
clause étant  la somme booléenne  de trois littéraux.  La  question de
{\sc  $3$-sat}\index{sat@{\sc 3-sat}}\index{3sat@{\sc  3-sat}}  est de
décider  si  on  peut  satisfaire  $f$.   Steve  Cook  a  démontré  un
extraordinaire théorème~:

\begin{theoreme}[Cook, \cite{cook:np}]
  Si l'on  dispose d'un algorithme  en temps polynomial  pour résoudre
  {\sc $3$-sat}, alors, pour n'importe   quel problème NP, on peut
  donner un algorithme en temps polynomial.
\end{theoreme}

On     dit      que     le     problème      {\sc     $3$-sat}     est
\emph{NP-complet}\index{complet!NP-    ---}\index{NP-complet}.    Plus
généralement,   tout  problème  $\Pi$   de  la   classe  NP   est  dit
\emph{NP-complet} si pour chaque problème  $\Pi'$ de NP, il existe une
\emph{réduction de  Turing} au  problème $\Pi$. C'est-à-dire  que pour
chaque problème $\Pi'$ de NP, on doit pouvoir fournir un algorithme en
temps polynomial  invoquant une  sous-routine de résolution  de $\Pi$.
Si un problème n'appartenant pas à la classe NP peut être réduit de la
sorte, on  dit qu'il est  \emph{NP-difficile}\index{NP-difficile}.  En
pratique, pour montrer qu'un  problème est NP-complet ou NP-difficile,
il  suffit  de  supposer   qu'on  dispose  d'un  algorithme  en  temps
polynomial pour le  résoudre, puis de montrer qu'on  peut utiliser cet
algorithme  comme sous-routine  pour résoudre  en temps  polynomial un
problème  qu'on  sait déjà  être  NP-complet,  comme  {\sc 3-sat}  par
exemple.  À la  suite du théorème de Cook, fondateur  de la théorie de
la complexité algorithmique, on a  prouvé que de nombreux problèmes de
décision sont  NP-complets (voir~\cite{garey.johnson:np}).  Le premier
exemple en théorie des graphes a été donné par Cook lui-même~:

\index{clique!NP-complétude}
\begin{probleme}[Clique]
\begin{itemize}
  \item[\sc Entrée~:] Un graphe $G$ et un entier $k$.
  \item[\sc Question~:] Existe-t-il une clique de $G$ de taille $k$ ?
  \item[\sc Complexité~:] NP-complet (Cook \cite{cook:np}).
\end{itemize}
\end{probleme}

Évidemment,  l'existence  d'un  algorithme  en temps  polynomial  pour
résoudre  {\sc  3-sat}  est  une question  fondamentale,  mais  encore
ouverte  à ce  jour. L'opinion  majoritaire est  qu'un  tel algorithme
n'existe sans doute pas~:

\begin{conjecture}[P $\neq$ NP]
  Il n'existe pas d'algorithme  en temps polynomial pour résoudre {\sc
  3-sat}.
\end{conjecture}

\section{Résumé des définitions les plus importantes}

\begin{description}
\item[Contient :]  On dit qu'un graphe  $G$ contient un  graphe $H$ si
  $H$ est un sous-graphe induit de $G$.
  
\item[Sans~:]  On dit  qu'un  graphe $G$  est  sans ``bidule''  (trou,
  diamant ...) s'il ne contient aucun ``bidule''.

\item[Chaîne~:] Une chaîne est une  suite de sommets d'un graphe telle
  que chaque sommet (sauf le dernier) voit le suivant. Dans le graphe,
  il  peut  y  avoir  d'autres  arêtes  que  celles  nécessaires  pour
  satisfaire la définition~: les \emph{cordes}.

\item[Chemin~:] Un chemin est un  graphe dont les sommets peuvent être
  ordonnés de  manière à  former une chaîne  sans corde.   La notation
  $\bp v_1  \tp \cdots \tp v_n  \ep$ désigne le chemin  sur lequel les
  sommets apparaissent  dans l'ordre  $v_1$, \dots, $v_n$.   Ce graphe
  est un  chemin par définition, mais  même si les  sommets $v_i$ sont
  dans un graphe $G$, ce n'est  pas forcément un chemin de $G$. Si $P$
  est un chemin, si  $u$ et $v$ sont des sommets de  $P$, alors $\bp u
  \tp P \tp  v \ep$ désigne le sous-chemin de  $P$ d'extrémités $u$ et
  $v$.

\item[Cycle~:] Un cycle est une suite de sommets d'un graphe telle que
  chaque  sommet voit  le  suivant  (le suivant  du  dernier étant  le
  premier).  Dans  le graphe,  il peut y  avoir d'autres  arêtes~: les
  \emph{cordes}.

\item[Trou~:] Un trou  est un graphe avec au moins  4 sommets dont les
  sommets  peuvent être  ordonnés de  manière à  former un  cycle sans
  corde . On appelle trou long tout trou ayant au moins 5 sommets.

\item[Anti~:]  Le  préfixe  ``anti''  désigne  une  propriété  ou  une
structure  du complémentaire  (à une  exception près~:  les partitions
antisymétriques).

\item[Prisme~:] Parfois appelé 3PC$(\Delta,\Delta)$, ou stretcher. Les
  prismes  sont  des  graphes  obtenus  à  partir  de  deux  triangles
  disjoints $\{a_1, a_2, a_3\}$ et $\{b_1, b_2, b_3\}$ en reliant pour
  tout $i$ les sommets $a_i$ et $b_i$ par un chemin $P_i$.  Les seules
  arêtes sont  celles des chemins  et des triangles.  Les  sommets des
  triangles sont  appelés les \emph{coins}  du prisme.  Un  prisme est
  dit  \emph{pair} (resp.   \emph{impair}) si  les trois  chemins sont
  pairs (resp.   impairs). Un prisme  est dit \emph{long} si  au moins
  l'un des  trois chemins est de  longueur au moins 2  (le seul prisme
  qui n'est pas long est $\overline{C_6}$).

\item[Quasi-prisme~:]  Les quasi-prismes  sont des  graphes  obtenus à
  partir  de deux  triangles disjoints  $\{a', c,  d\}$ et  $\{c', d',
  b'\}$ en  reliant les  sommets $c$  et $c'$ par  un chemin  $S_3$ de
  longueur au moins  1 et les sommets $d$ et $d'$  par un chemin $S_4$
  de longueur au moins 1. Il y  a en plus un chemin $S_1$ reliant $a'$
  à un sommet $a$ et un chemin  $S_2$ reliant $b'$ à un sommet $b$. Il
  n'y a  aucune autre  arête que celles  des chemins et  des triangles
  ci-dessus.  Les sommets $a$ et $b$ sont appelés \emph{extrémités} du
  quasi-prisme.  Les  chemins $S_1$ et $S_2$ peuvent  être de longueur
  0.  Si l'un de  $S_1$, $S_2$ est de longueur au moins  1, on dit que
  le quasi-prisme est \emph{strict}.

\item[Pyramide~:] Parfois appelée 3PC$(\Delta, P)$. Les pyramides sont
  des  graphes obtenus  à  partir  d'un sommet  $a$  et d'un  triangle
  $\{b_1,  b_2, b_3\}$ en  reliant pour  tout $i$  les sommets  $a$ et
  $b_i$  par un  chemin  $P_i$.   Les seules  arêtes  sont celles  des
  chemins  et du  triangle. Par  définition, un  seul des  chemins est
  autorisé à être  de longueur 1, les autres  doivent être de longueur
  au moins 2.  Le sommet $a$ est appelé le \emph{coin} de la pyramide.
  L'intérêt des  pyramides est  qu'elles contiennent toujours  un trou
  impair.
\end{description}

\chapter{Graphes parfaits}
Les graphes parfaits ont été  introduits par Claude Berge au début des
années 1960 à  la suite d'un cheminement assez  complexe qu'il raconte
lui-même  dans un  article  agréable~\cite{berge.r:origin}.  Pour  une
synthèse  complète  des  recherches  sur  les  graphes  parfaits  nous
renvoyons  à   quatre  références~:  un  article  un   peu  ancien  de
L.~Lov\'asz~\cite{lovasz:pgsurvey},  un ouvrage collectif  publié sous
la direction  de C.~Berge et V.~Chv\'atal~\cite{berge.chvatal:topics},
un  ouvrage  collectif  plus  récent,  publié  sous  la  direction  de
L.~Ram\'irez  Alfons\'in et  B.~Reed~\cite{livre:perfectgraphs}  et un
article   récent  de   M.~Chudnovsky,   N.~Robertson,  P.~Seymour   et
R.~Thomas~\cite{chudvovsky.r.s.t:progress}.    De   manière   un   peu
anachronique,  on peut ajourd'hui  motiver l'introduction  des graphes
parfaits par la théorie de la complexité.  Le problème qui nous occupe
est celui de la coloration,  reformulé ci-dessous comme un problème de
décision~:

\index{coloration!NP-complétude}
\begin{probleme}[Coloration] 
\begin{itemize}
  \item[\sc Entrée~:] Un graphe $G$ et un entier $k$.
  \item[\sc Question~:] Existe-t-il une $k$-coloration de $G$ ? 
  \item[\sc Complexité~:] NP-complet (Karp, \cite{karp:np}).
\end{itemize}
\end{probleme}

Il est clair  que si un graphe $G$ contient une  clique de taille $k$,
il faut au moins $k$ couleurs pour le colorier. Donc, pour tout graphe
on a~: $\chi(G) \geq \omega(G)$.  Une plus grande clique apparaît donc
comme un bon candidat pour être un certificat du non du problème de la
coloration.   Il est donc  naturel de  se demander  s'il y  a toujours
égalité entre  le nombre  chromatique et la  taille d'une  plus grande
clique.   Ce n'est  pas le  cas et  le plus  petit  contre-exemple est
$C_5$~:  $\chi(C_5) =  3$ et  $\omega(C_5) =  2$.   Plus généralement,
n'importe quel trou  impair $H$ vérifie $\chi(H) =  3$ et $\omega(H) =
2$.  Alain Ghouila-Houri~\cite{gouilathouri:note} a également remarqué
que n'importe quel antitrou  impair $\overline{H}$ avec $2n+1$ sommets
vérifie   $\chi(\overline{H})  =   n+1$  et   $\omega(\overline{H})  =
n$. Donc,  l'égalité $\chi  = \omega$ est  fausse en général,  mais on
peut  se demander  pour quels  graphes  elle demeure  exacte. Dès  les
origines, Claude  Berge s'est intéressé aux graphes  vérifiant $\chi =
\omega$ (ou  ayant comme il  disait ``la belle propriété'')  pour tous
leurs  sous-graphes   induits.  En  effet,   sans  cette  restriction,
n'importe  quel graphe  peut être  muni artificiellement  de  la belle
propriété par  l'ajout d'une grande  clique par exemple.   L'étude des
graphes ayant la belle propriété risque donc de se ramener à celle des
graphes ayant quelque part une grande clique.

\index{parfait~(graphe~---)}
\begin{defeng}[Berge, \cite{berge:61}]
  Un graphe $G$ est dit parfait si et seulement si pour chacun de ses
  sous-graphes induits $G'$ on a $\chi(G') = \omega(G')$.
\end{defeng}

Pour la classe des graphes parfaits, on dispose d'un bon certificat du
non et  aussi du oui pour le  problème de la coloration.   Il est donc
naturel de se  demander s'il existe un algorithme  en temps polynomial
pour colorier les graphes parfaits. C'est bien le cas~:

\index{coloration!graphes parfaits}
\begin{algorithme}[Gr\"ostchel, Lov\'asz, Schrijver, \cite{grostchel.l.s:color}]
  \label{graphpar.a.coloration}
  \begin{itemize}
  \item[\sc Entrée~:] Un graphe parfait $G$.
  \item[\sc Sortie~:] Une coloration optimale de $G$.
  \item[\sc Complexité~:] Polynomiale.  
  \end{itemize}
\end{algorithme}

Il  faut  noter  que  les  recherches  ayant  conduit  à  l'algorithme
ci-dessus ont joué  un rôle central dans le  développement de branches
importantes  de l'optimisation,  comme  la programmation  semi-définie
(Cf.~\cite{reed:gentle}). Il utilise la méthodes des ellipsoïdes qui a
la  réputation  d'être  très  difficile  à implémenter  en  raison  de
problèmes  d'instabilité  numérique.   La  recherche  d'un  algorithme
performant  ou  conceptuellement  simple  pour  colorier  les  graphes
parfaits reste donc d'actualité.

Nous avons vu que les trous  impairs et les antitrous impairs sont des
graphes imparfaits. Donc, dès qu'un graphe contient l'un d'eux, il est
lui aussi  imparfait. Les  graphes sans trou  impair et  sans antitrou
impair sont donc la plus grande classe de graphes dont on peut espérer
la perfection. Cela justifie la définition suivante~:

\index{Berge~(graphe de ---)}
\begin{defeng}
  Un graphe $G$  est dit de Berge si et seulement  s'il ne contient ni
  trou impair ni antitrou impair.
\end{defeng}

Au  début des  années 1960,  Claude  Berge a  énoncé deux  conjectures
célèbres~\cite{berge:indian}.   La   première,  dite  \emph{conjecture
faible    des     graphes    parfaits},    a     été    étudiée    par
Fulkerson~\cite{fulkerson:antiblocking},  ce  qui  l'a  conduit  à  la
théorie des  polyèdres antibloquants (que nous  n'évoquerons pas ici).
Finalement, L\'aszl\'o Lov\'asz l' a démontrée en 1971~:

\begin{theoreme}[Théorème  des graphes parfaits, Lov\'asz, \cite{lovasz:pg,lovasz:nh}]\mbox{}
  \label{graphespar.t.faible}

  Un  graphe $G$  est parfait  si et  seulement si  $\overline{G}$ est
  parfait.
\end{theoreme}

La deuxième conjecture de  Claude Berge, \emph{la conjecture forte des
graphes parfaits},  a fait l'objet de nombreuses  recherches (Va\v sek
Chv\'atal   recense    plus   de   500   articles    sur   le   sujet,
voir~\cite{chvatal:bib})   avant   d'être   finalement  démontrée   en
mai~2002~:

\begin{theoreme} 
  {\rm\bf(Théorème   fort  des   graphes   parfaits,  \\   \rule{3cm}{0cm}Chudnovsky,
    Robertson, Seymour, Thomas \cite{chudvovsky.r.s.t:spgt})}

  Un graphe est parfait si et seulement s'il est de Berge.
\end{theoreme}

À  la fin  de l'année  2002, la  principale question  ouverte  sur les
graphes  parfaits  restait  celle  de  leur  reconnaissance  en  temps
polynomial. Ce problème a  été finalement résolu par Maria Chudnovsky,
Gérard  Cornuéjols,  Xinming  Liu,  Paul Seymour  et  Kristina  Vu\v
skovi\'c~\cite{chudnovsky.c.l.s.v:cleaning,chudnovsky.seymour:reco,cornuejols.liu.vuskovic:reco}.
Au chapitre~\ref{reco.chap}, nous  donnerons quelques détails sur
leurs méthodes.

Finalement,  les   \emph{trois  grands  problèmes}   sur  les  graphes
parfaits, à savoir  la conjecture forte, le problème  de la coloration
et celui  de la reconnaissance des graphes  parfaits, sont aujourd'hui
résolus.   De  nombreuses   questions  demeurent  cependant  ouvertes,
l'existence  d'algorithmes   de  coloration  plus   rapide  que  celui
utilisant  la  méthode  de  l'éllipsoïde notamment.   On  souhaiterait
également des algorithmes qui tirent parti de propriétés combinatoires
des  graphes  de  Berge.  On  considère donc  que  la  recherche  d'un
algorithme ``combinatoire'' de coloration des graphes parfaits est une
question ouverte.  Notons que la  notion d'algorithme ``combinatoire''
n'a pas de définition formelle.

Nous présentons  dans ce chapitre  un certain nombre de  résultats sur
les  graphes  parfaits,  parfois   assez  anciens.   Notre  choix  est
subjectif et  ne prétend pas à l'exhaustivité.   Outre notre ignorance
et nos  besoins spécifiques  pour la suite  de ce travail,  le critère
retenu est celui de la pertinence des résultats au vu de la résolution
des trois  problèmes, et notamment de  la preuve du  théorème fort des
graphes parfaits par Chudnovsky {\it et al.}.

\section{Les graphes basiques}

\index{line-graphe~de~biparti} Nous mentionnons ici les quatre classes
dites \emph{basiques}\index{basique} de graphes parfaits~: les graphes
bipartis, les line-graphes de  bipartis, et leurs complémentaires.  La
perfection des graphes bipartis est triviale (elle résulte par exemple
du lemme~\ref{base.l.trivbip}).  La  perfection des complémentaires de
bipartis,  des  line-graphes de  bipartis  et  des complémentaires  de
line-graphes   de   bipartis   résultent   de   trois   théorèmes   de
K\"onig~(voir~\cite{lovasz:pgsurvey,gallai:triangule,konig:16,konig:31}).

Dans notre  travail, nous aurons parfois besoin  d'informations sur la
structure  des   line-graphes  de  bipartis.    Signalons  que  Lowell
Beineke~\cite{beineke:linegraphs}  a prouvé un  théorème caractérisant
les  line-graphes en général  (il affirme  que sa  caractérisation par
sous-graphes  induits interdits  était  en fait  déjà  connue de  Neil
Robertson).   Indépendamment  l'un   de  l'autre,  Philippe  Lehot  et
Nicholas Roussopoulos ont ensuite donné l'algorithme suivant~:

\index{line-graphe!calcul de la racine}
\begin{algorithme}[Lehot \cite{lehot:root}, Roussopoulos \cite{roussopoulos:linegraphe}]
  \label{graphespar.a.lehot}
  \begin{itemize}
    \item[\sc Instance~:] Un graphe $G$.
    \item[\sc Sortie~:]  Un graphe  $R$ tel que  $G= L(R)$. Si  un tel
    graphe n'existe pas, l'algorithme  retourne ``$G$ n'est pas un
    line-graphe''.
    \item[Complexité~:] $O(n+m)$. 
  \end{itemize}
\end{algorithme}

Nul   n'a  semble-t-il   songé  à   spécialiser  cet   algorithme  aux
line-graphes  de bipartis.   Il  y a  de  bonnes raisons  à cela.   Un
théorème ancien de Hassler Whitney~\cite{whitney:graphs} affirme qu'en
dehors  du cas  trivial  de la  griffe  et du  triangle, deux  graphes
différents ont toujours des line-graphes différents.  Donc si on donne
à  l'algorithme  de  Lehot  un  line-graphe de  biparti  différent  du
triangle,  on est  certain que  le graphe  $R$ retourné  sera biparti.
Autrement  dit, l'algorithme de  Lehot se  spécialise de  lui-même aux
line-graphes de bipartis,  et comme il est déjà  en temps linéaire, il
est  impossible d'améliorer  sa  complexité.  Il  existe cependant  un
théorème sur la structure des line-graphes de graphes bipartis~:

\index{line-graphe~de~biparti!caractérisation}
\begin{theoreme}[Harary et Holzmann, \cite{harary.holzmann:lgbip}]
  \label{graphespar.t.lgb}
  Soit  $G$  un  graphe. $G$  est  un  line-graphe  de biparti  si  et
  seulement si $G$ est sans griffe, sans diamant et sans trou impair.
\end{theoreme}

\section{Les travaux de Lov\'asz}

Nous  avons  déjà rencontré  les  principaux  théorèmes de  L\'aszl\'o
Lov\'asz,  donnés au  début des  années  1970, et  fondamentaux pour  la
compréhension des  graphes parfaits.  Nous avons besoin  pour la suite
de notre travail de détails supplémentaires sur ses travaux~: un outil
essentiel  dans la  preuve du  théorème  des graphes  parfaits est  le
\emph{lemme de réplication}\index{réplication~(lemme~de~---)}~:

\begin{lemme}[Lov\'asz, \cite{lovasz:nh}]
  \label{graphespar.l.rep}
  Soit $G$  un graphe parfait  et $v \in  V(G)$.  Soit $G'$  le graphe
  obtenu en ajoutant un nouveau sommet  $v'$ et en le reliant à $v$ et
  à tous les voisins de $v$. Alors, $G'$ est parfait.
\end{lemme}

Pour démontrer  le théorème des  graphes parfaits, Lov\'asz  montre en
fait  le  théorème suivant  qui  est  plus  fort. Mentionnons  que  
G.~Gasparian~\cite{gasparian:minimp}  a donné une  preuve simplifiée
de ce théorème~:

\begin{theoreme}[Lov\'asz, \cite{lovasz:pg}]
  \label{graphespar.t.lovasz}
  Un graphe $G$  est parfait si et seulement  si pour tout sous-graphe
  $G'$ de $G$ on a $\alpha(G') \omega(G') \geq |V(G')|$.
\end{theoreme}

Lov\'asz a également introduit  une notion importante. Soient $\alpha,
\omega  \geq 1$ des  entiers et  $G$ un  graphe.  On  dit que  $G$ est
\emph{$(\alpha,
\omega)$-partitionnable}\index{partitionnable~(graphe~---)}    si   et
seulement si  pour tout sommet $v$  de $G$, le graphe  $G \setminus v$
peut être  partitionné en  $\alpha$ cliques de  taille $\omega$  et en
$\omega$  stables   de  taille  $\alpha$.    On  appelle  \emph{graphe
minimalement      imparfait}\index{minimalement~imparfait~(graphe~---)}
\index{imparfait|see{minimalement~imparfait}}  tout  graphe qui  n'est
pas parfait et tel que chacun de ses sous-graphes induits est parfait.
Notons que la conjecture forte des graphes parfaits dit simplement que
les  trous impairs  et les  antitrous impairs  sont les  seuls graphes
minimalement  imparfaits.   Le théorème  suivant  est une  conséquence
directe du théorème~\ref{graphespar.t.lovasz}~:

\begin{theoreme}[Lov\'asz, \cite{lovasz:pg}]
  \label{graphespar.t.minimparfait}
  Soit $G$ un graphe imparfait minimal. Alors $G$ est partitionnable.
\end{theoreme}

Les  graphes   partitionnables  ont  de   très  nombreuses  propriétés
intéressantes  (voir~\cite{preissmann.sebo:minimal}).
Notons  que pour  prouver le  théorème fort  des graphes  parfaits, il
suffit  de montrer  que tout  graphe partitionnable  contient  un trou
impair  ou  un  antitrou   impair.  Le  théorème  suivant  montre  les
principales propriétés des graphes partitionnables. Notons que Padberg
\cite{padberg:74}  avait déjà  démontré  ce théorème  dans  le cas  plus
restreint des graphes minimalement imparfaits.

\begin{theoreme}[Bland, Huang, Trotter \cite{bland.h.t:79}]
  \label{graphespar.t.partitionnable}
  Soit $G$ un  graphe $(\alpha, \omega)$-partitionnable avec $n=\alpha
  \omega +1$ sommets. Alors~:
  \begin{outcomes}
    \item 
      $G$ possède exactement $n$ cliques de taille $\omega$.
    \item
      $G$ possède exactement $n$ stables de taille $\alpha$.
    \item 
      Chaque sommet de $G$ appartient à exactement $\omega$ cliques de
      taille $\omega$.
    \item 
      Chaque sommet de $G$ appartient à exactement $\alpha$ stables de
      taille $\alpha$.
    \item 
      \label{graphespar.o.cs}
      Chaque clique de $G$ de taille $\omega$  est disjointe
      d'exactement un stable de $G$ de taille $\alpha$.
     \item 
      \label{graphespar.o.sc}
      Chaque stable de $G$ de taille $\alpha$  est disjoint
      d'exactement une clique de $G$ de taille $\omega$.
    \item 
      \label{graphespar.o.unicolor}
      Pour chaque sommet  $v$ de $G$, il existe  une unique coloration
      de $G \setminus v$ avec $\omega$ couleurs.
  \end{outcomes}
\end{theoreme}

\section{Les théorèmes de décomposition}
\label{graphespar.ss.structure}

\index{décomposition}  Par théorème  de décomposition,  nous entendons
ici tout  résultat affirmant que  pour tout graphe $G$  d'une certaine
classe $\cal C$, ou bien $G$ appartient à une sous-classe $\cal C'$ ou
bien $G$  peut être  cassé (décomposé) d'une  manière {\it ad  hoc} en
plusieurs ``morceaux''. Un tel résultat sera d'autant plus intéressant
que  la  sous-classe  $\cal  C'$ sera  simple  et  que l'opération  de
décomposition  préservera  de  bonnes  propriétés, comme  par  exemple
l'appartenance à la classe~$\cal C$.

On peut débattre  longtemps pour savoir si tel ou  tel théorème est ou
non un véritable théorème de décomposition.  Certains résultats ont de
bonnes      propriétés       algorithmiques      (2-joint,      clique
d'articulation~\dots),  d'autres  permettent  de montrer  beaucoup  de
théorèmes~(étoile~d'articulation),  d'autres enfin  ont  l'avantage de
s'appliquer  à   des  classes   de  graphes  très   larges  (partition
antisymétrique, étoile double~\dots). Nous avons choisi la solution de
facilité, consistant à  ne pas entrer dans ce débat  et à inclure dans
cette  section tous  les  résultats nous  semblant  nécessaires. À  la
section~\ref{graphespar.s.poq}  consacrée aux problèmes  ouverts, nous
expliquerons pourquoi  la partition  antisymétrique n'a pas  de bonnes
propriétés algorithmiques.

L'utilisation de théorèmes de décomposition pour l'étude de classes de
graphes  parfaits  a  commencé   avant  même  que  Berge  ait  formulé
clairement    sa   célèbre    conjecture   (voir~\cite{berge.r:origin}
et~\cite{reed:chap2} pour un  historique détaillé des origines).  Mais
ce  n'est qu'à  la fin  des  années 1970  qu'est apparu  le projet  de
prouver la conjecture forte des graphes parfaits par des techniques de
décomposition. Lors de la conférence  en l'honneur de Claude Berge qui
s'est tenue à Paris en juillet  2004, Va\v sek Chv\'atal a raconté que
la première personne ayant à sa connaissance formulé ce projet est Sue
Whitesides, à l'automne 1977. Le premier résultat dans cette direction
fut  l'étude de la  décomposition par  \emph{amalgame} des  graphes de
Meyniel par Burlet et Fonlupt~\cite{burlet.fonlupt:meyniel}. Une autre
étape importante fut la conférence  de Princeton en 1993, où Chv\'atal
et Reed  ont compris que  les théorèmes de décomposition  des matrices
équilibrées de Conforti, Cornuéjols, Kapoor et Vu\v skovi\'c pouvaient
avoir  des  analogies avec  les  graphes  parfaits.   Puis sont  venus
d'autres résultats que nous allons examiner ici.

L'exemple  le   plus  ancien  et  le  plus   simple  d'utilisation  de
décompositions   provient  d'un   théorème  de   Dirac.    On  appelle
\emph{graphe triangulé}\index{triangulé}  tout graphe qui  ne contient
aucun  trou.    Si  $G$  est   un  graphe,  on   appelle  \emph{clique
d'articulation}\index{clique!d'articulation}  de   $G$  tout  ensemble
d'articulation $K$  de $G$, qui induit  une clique de  $G$. On appelle
alors \emph{pièces}\index{pièce}  de $G$ les graphes  $G[G_1 \cup K]$,
\dots, $G[G_k  \cup K]$  où $G_1$, \dots,  $G_k$ sont  les composantes
connexes de $G \setminus K$.

\index{décomposition!graphes triangulés}
\index{décomposition!par~clique~d'articulation}
\begin{theoreme}[Dirac, \cite{dirac:chordal}]
  \label{graphpar.t.dirac}
  Soit $G$ un  graphe triangulé. Alors ou bien $G$  est une clique, ou
  bien $G$ possède une clique d'articulation.
\end{theoreme}

\begin{theoreme}[Gallai, \cite{gallai:triangule}]
  \label{graphpar.t.gallai}
  Soit  $G$ un  graphe  possédant une  clique  d'articulation. Si  les
  pièces de $G$ sont parfaites, alors $G$ est parfait.
\end{theoreme}

Ce  dernier théorème  a  pour conséquence  immédiate  que les  graphes
minimalement  imparfaits n'ont pas  de clique  d'articulation.  Claude
Berge~\cite{berge:61}  a  donc  remarqué  la  perfection  des  graphes
triangulés, que  l'on peut prouver  par une récurrence  évidente~: les
cliques sont des  graphes parfaits.  Si un graphe  triangulé n'est pas
une    clique,   il    a    une   clique    d'articulation   par    le
théorème~\ref{graphpar.t.dirac},   ses  pièces   sont   parfaites  par
hypothèse    de   récurrence,    et    il   est    parfait   par    le
théoreme~\ref{graphpar.t.gallai}.  Depuis les années 1960, de nombreux
théorèmes sur  les graphes  parfaits ont été  démontrés sur  ce modèle
(voir \cite{rusu:cutsets}), l'idée étant de trouver des décompositions
interdites  dans  les   graphes  minimalement  imparfaits,  ou  mieux,
préservant         la          perfection.          On         appelle
\emph{étoile}\index{etoile@étoile} tout  graphe qui possède  un sommet
voyant tous les autres  sommets du graphe. L'étoile d'articulation est
l'une des décompositions les plus connues.  Son intérêt provient entre
autres du théorème suivant~:

\index{minimalement~imparfait~(graphe~---)!étoile~d'articulation}
\begin{theoreme}[Chv\'atal, \cite{chvatal:starcutset}]
  \label{graphespar.t.scs}
  Dans un  graphe partitionnable (et donc dans  un graphe minimalement
  imparfait) il n'y a pas d'étoile d'articulation.
\end{theoreme}

\index{décomposition!par~partition~antisymétrique}    L'étoile
d'articulation a permis de prouver la perfection de nombreuses classes
de graphes.  Chv\'atal~\cite{chvatal:starcutset} a aussi introduit une
autre  décomposition qui  s'est avérée  importante par  la  suite~: on
appelle                                                 \emph{partition
antisymétrique}\index{partition~antisymétrique}\index{antisymétrique|see{partition~antisymétrique}}
d'un graphe $G$  toute partition des sommets de  $G$ en deux ensembles
$A$  et $B$  tels que  $G[A]$ n'est  pas connexe  et $G[B]$  n'est pas
anticonnexe.    Notons  qu'une  étoile   d'articulation  est   un  cas
particulier de  partition antisymétrique.  Chv\'atal  a conjecturé que
les  graphes  minimalement  imparfaits   n'ont  pas  de  de  partition
antisymétrique~:

\index{minimalement~imparfait~(graphe~---)!partition~antisymétrique}
\index{décomposition!par~partition~antisymétrique}
\begin{conjecture}[Chv\'atal, \cite{chvatal:starcutset}]
  \label{graphespar.conj.sp}
  Soit $G$ un graphe minimalement imparfait. Alors $G$ ne possède pas
  de partitition antisymétrique. 
\end{conjecture}

Cette  conjecture   s'est  révélée  difficile,  et   de  nombreux  cas
particuliers      ont       été      étudiés.       Cornuéjols      et
Reed~\cite{cornuejols.reed:sp}  ont montré  qu'un  graphe minimalement
imparfait ne peut pas  contenir d'ensemble d'articulation induisant un
graphe multiparti complet avec  au moins deux composantes anticonnexes
(un   graphe  \emph{multiparti   complet}  est   un  graphe   dont  le
complémentaire a pour  composantes connexes des cliques).  Mentionnons
les  travaux de  Ch\'inh Hoàng~\cite{hoang:minimp}  ainsi que  ceux de
Florian  Roussel  et  Philippe  Rubio  qui ont  eu  des  prolongements
inattendus (voir chapitre~\ref{rr.chap})~:

\begin{theoreme}[Roussel et Rubio, \cite{roussel.rubio:01}]
  Soit $G$ un graphe minimalement imparfait possédant une partition
  antisymétrique $(X,Y)$. Alors aucune composante anticonnexe de
  $G[Y]$ n'est un stable.
\end{theoreme}

\index{décomposition!par~2-joint}  Une autre  décomposition importante
est le  \emph{2-joint}\index{2-joint}\index{joint@2-joint} définit par
G.~Cornuéjols  et W.~H.~Cunningham~\cite{cornuejols.cunningham:2join}.
Notons  que   la  définition  précise  du  2-joint   varie  selon  les
auteurs.  Nous donnons  une définition  provisoire, qui  sera précisée
ensuite.  Soit $G$  un graphe. On appelle \emph{2-joint}  de $G$ toute
partition des  sommets de  $G$ en deux  ensembles $X_1$ et  $X_2$ tels
qu'il existe $A_i, B_i \subseteq  X_i$ ($i= 1, 2$) vérifiant certaines
conditions parmi les suivantes~:

\vspace{1ex}

\label{graphespar.p.2joint}
\begin{enumerate}
  \item
    \label{2joint1}
    Pour $i=1,2$, $A_i \cap B_i = \emptyset$~;
  \item
    \label{2joint2}
    Chaque sommet de $A_1$ voit chaque sommet de $A_2$~; 
  \item
    \label{2joint3}
    Chaque sommet de $B_1$ voit chaque sommet de $B_2$~;
  \item 
    \label{2joint4}
    Il  n'y a  aucune  arête entre  $X_1$  et $X_2$  autres que  celles
    entre~$A_i$ et $B_i$ ($i=1, 2$);
  \item 
    \label{2joint5}
    Pour $i=  1, 2$,  si $|A_i|  = |B_i| =  1$ et  si $G[X_i]$  est un
    chemin  joignant l'unique  sommet de  $A_i$ à  l'unique  sommet de
    $B_i$, alors ce chemin est de longueur au moins 3~;
  \item
    \label{2joint6}
    Pour $i= 1, 2$, chaque composante connexe de $G[X_i]$ comporte des
    sommets  de $A_i$ et de $B_i$~;
  \item 
    \label{2joint7}
    Pour $i= 1,  2$, si $|A_i| = |B_i| = 1$,  alors $G[X_i]$ n'est pas
    un chemin joignant  l'unique sommet de $A_i$ à  l'unique sommet de
    $B_i$~;
  \item
    \label{2joint8}
    Pour  $i =  1,  2$, il  existe  un chemin  de  $G[X_i]$ ayant  une
    extrémité dans $A_i$ et une extrémité dans $B_i$~;
  \item
    \label{2joint9}
    Les ensembles $X_1$ et $X_2$ sont de cardinal au moins~3.
\end{enumerate}

\vspace{2ex}

Dans les articles de Chudnvovsky {\it et
al.}    dans~\cite{chudvovsky.r.s.t:spgt}, les 2-joints doivent
seulement vérifier les conditions~\ref{2joint1} à~\ref{2joint6}. En
l'absence d'indications supplémentaire, c'est cette définition que
nous utilisons.   Dans   les   articles  de
Cornuéjols  {\it  et al},  des définitions  légerement différentes  sont
utilisées, que nous précisons pour chaque théorème.

Les  2-joints sont  en  un  certain sens  interdits  dans les  graphes
minimalement  imparfaits.  Plus précisément,  ils sont  interdits dans
les  contre-exemple  minimaux  à   la  conjecture  forte  des  graphes
parfaits.  Notons  que  ce   théorème  est  valable  pour  toutes  les
définitions du 2-joints.

\index{minimalement~imparfait~(graphe~---)!2-joint}
\begin{theoreme}[Cornuéjols  et Cunningham~\cite{cornuejols.cunningham:2join}]
  \label{graphespar.t.2joint}
  Soit  $G$  un  graphe  minimalement  imparfait. Si  $G$  possède  un
  2-joint, alors $G$ est un trou impair.
\end{theoreme}

En  2001, Gérard  Cornuéjols, Michelangelo  Conforti et  Kristina Vu\v
skovi\'c ont formulé une conjecture qui, si l'on admet que les graphes
minimalement imparfaits n'ont pas de partition antisymétrique, implique
la conjecture forte des graphes parfaits~:

\index{Berge~(graphe de ---)!décomposition}
\index{décomposition!graphes de Berge}
\begin{conjecture}[Conforti, Cornuéjols, Vu\v skovi\'c~\cite{conforti.c.v:square}]
  \label{graphespas.conj.ccv}
  Soit $G$ un graphe de Berge.  Alors ou bien $G$ est basique, ou bien
  l'un de $G$, $\overline{G}$ possède  un 2-joint, ou bien $G$ possède
  une partition antisymétrique.
\end{conjecture}

Cette conjecture apparaît aujourd'hui comme l'idée décisive conduisant
à la preuve du théorème  fort des graphes parfaits.  Ses trois auteurs
ont pu la démontrer dans  des cas particuliers, ou en relaxant quelque
peu  la conclusion.  Dans  le théorème  suivant, les  2-joints doivent
seulement   vérifier  les   conditions~\ref{2joint1}  à~\ref{2joint4},
\ref{2joint7} et~\ref{2joint8}.

\index{décomposition!graphes de Berge sans carré} 
\index{carré!graphes~de~Berge~sans~---,~décomposition}
\begin{theoreme}[Conforti, Cornuéjols, Vu\v skovi\'c~\cite{conforti.c.v:square}]
  \label{graphespar.t.sansc4}
  Soit $G$ de Berge sans $C_4$. Alors ou bien $G$ est basique, ou bien
  $G$   possède  un   2-joint   ou  bien   $G$   possède  une   étoile
  d'articulation. En conséquence, $G$ est parfait.
\end{theoreme}

On                         appelle                        \emph{étoile
double}\index{double|see{étoile~double}}\index{etoile~double@étoile~double}
tout graphe  $G$ réduit à un  sommet ou possédant une  arête $e$ telle
que chaque sommet de $G$ est  adjacent à au moins l'une des extrémités
de  $e$. Dans  le  théorème suivant,  les  2-joints doivent  seulement
vérifier  les conditions~\ref{2joint1}  à~\ref{2joint4}, \ref{2joint7}
et~\ref{2joint9}.

\index{décomposition!graphes~sans~trou~impair}
\index{décomposition!par~étoile~double}     \index{Berge~(graphe    de
---)!décomposition}
\begin{theoreme}[Conforti, Cornuéjols, Vu\v skovi\'c~\cite{conforti.c.v:dstrarcut}]
  \label{graphespar.t.sansti}
  Soit $G$  sans trou impair. Alors  ou bien $G$ est  basique, ou bien
  $G$  possède  un 2-joint  ou  bien  $G$  possède une   étoile double
  d'articulation.
\end{theoreme}

Les  antitrous impairs  d'au  moins 7  sommets  possèdent des  étoiles
doubles  d'articulation. Il  est donc  impossible de  montrer  que les
graphes   minimalement   imparfaits    n'ont   pas   d'étoile   double
d'articulation. On pourrait en déduire  un peu vite que le théorème de
décomposition des graphes sans trou impair est totalement inutile pour
prouver le théorème  fort des graphes parfaits. Cela  n'est pas tout à
fait  exact~: on  sait aujourd'hui  que la  conjecture  ci-dessous est
vraie  ---  c'est  une   conséquence  du  théorème  fort  des  graphes
parfaits. Une preuve  directe (dont nous n'avons pas  la moindre idée)
fournirait une alternative à l'article de Chudnovsky et al.~:

\index{minimalement~imparfait~(graphe~---)!double~étoile}
\begin{conjecture}
  \label{graphespar.conj.dsc}
  Soit $G$ un graphe minimalement imparfait. Alors l'un de $G$ et
  $\overline{G}$ ne possède pas d'étoile double d'articulation.
\end{conjecture}

Voici  comment on  montre le  théorème  fort des  graphes parfaits  en
admettant la conjecture ci-dessus~: soit $G$ un contre-exemple minimal
au théorème fort  des graphes parfaits.  Notons que  $G$ est de Berge.
Donc  $G$  et  $\overline{G}$  sont  sans  trou  impair.   D'après  le
théorème~\ref{graphespar.t.faible},   $\overline{G}$   est  aussi   un
contre-exemple            minimal.             Appliquons           le
théorème~\ref{graphespar.t.sansti} à  $G$ et $\overline{G}$.   Si l'un
de  $G$, $\overline{G}$ est  basique, on  contredit la  perfection des
graphes basiques~; si l'un  de $G$, $\overline{G}$ possède un 2-joint,
on  contredit  le  théorème~\ref{graphespar.t.2joint}~;  donc  $G$  et
$\overline{G}$  possèdent  une étoile  double  d'articulation, ce  qui
contredit la conjecture~\ref{graphespar.conj.dsc}.

\section{Le théorème fort des graphes parfaits}
\label{sec.graphpar.preuve}

Maria Chudnovsky,  Neil Robertson, Paul  Seymour et Robin  Thomas sont
parvenus  à  prouver  des   versions  légèrement  plus  faibles,  mais
suffisantes,          des         conjectures~\ref{graphespar.conj.sp}
et~\ref{graphespas.conj.ccv},  prouvant  ainsi  le théorème  fort  des
graphes parfaits  et répondant à l'une des  grandes questions ouvertes
de  la  théorie des  graphes.   Leur  preuve  est longue  (le  premier
manuscrit diffusé,  daté du  28 octobre 2002,  comporte 148  pages) et
technique.  Nous en présentons très brièvement les étapes principales.
Nous mettons l'accent sur les résultats intermédiaires susceptibles de
nous  servir  par  la  suite,  qui  ne sont  pas  forcément  les  plus
importants selon d'autres critères.

Nous avons  besoin de quelques  notions supplémentaires.  Soit  $G$ un
graphe.  Si  $X, Y  \subset V(G)$,  on dit que  le couple  $(X,Y)$ est
\emph{complet}\index{complet!ensemble~---~à} \index{complet!couple~---} (ou que $X$ est complet à
$Y$) si  chaque sommet de  $X$ voit chaque  sommet de $Y$. On  dit que
$(X,Y)$ est \emph{anticomplet}\index{anticomplet}  si chaque sommet de
$X$  manque chaque  sommet de  $Y$. On  appelle  \emph{paire homogène}
\index{paire~homogène}\index{homogène|see{paire~homogène}}\index{décomposition!par~paire~homogène}
toute partition des sommets de $G$ en 6 ensembles non vides $(A, B, C,
D, E, F)$ tels que~:

\vspace{1ex}

\begin{itemize}
\item
  Chaque sommet de $A$ a un voisin et un non-voisin dans $B$. Chaque
  sommet de $B$ a un voisin et un non-voisin dans $A$. 
\item
  Les couples $(C, A)$, $(A, F)$, $(F, B)$, $(B, D)$ sont complets.
\item
  Les  couples   $(D,  A)$,  $(A,   E)$,  $(E,  B)$,  $(B,   C)$  sont
  anticomplets.
\end{itemize}

\vspace{2ex}

Notons que notre définition de la paire homogène est celle donnée sous
le   nom   de  $M$-joint\index{Mj@$M$-joint|see{paire~homogène}}   par
Chudnovsky  {\it  et  al}. Mais  que  la  notion  a été  inventée  par
Chv\'atal  et Sbihi~\cite{chvatal.sbihi:bullfree},  qui ont  montré le
théorème suivant~:

\begin{theoreme}[Chv\'atal et Sbihi~\cite{chvatal.sbihi:bullfree}]  
  \label{grapphespar.t.paireh}
  Soit $G$ un graphe minimalement imparfait. Alors $G$ ne contient pas
  de paire homogène.
\end{theoreme}

\label{graphespar.ss.SP}
\index{pair!partition~antisymétrique~---}
\index{partition~antisymétrique~paire} Nous  allons maintenant énoncer
un certain  nombre de  théorèmes démontrés par  Chudnovsky, Robertson,
Seymour et Thomas  dans~\cite{chudvovsky.r.s.t:spgt}.  Nous donnons en
référence  le numéro  du théorème  et la  page à  laquelle on  peut le
trouver  dans le manuscrit  daté du  28 octobre  2002.  On  dit qu'une
partition antisymétrique $(A,B)$ est  \emph{paire} si tous les chemins
ayant leurs  extrémités dans  $B$ et leur  intérieur dans $A$  sont de
longueur paire et si tous  les antichemins ayant leurs extrémités dans
$A$ et  leur intérieur  dans $B$ sont  de longueur paire.  Le théorème
suivant  est  une version  affaiblie  de  la  conjecture de  Chv\'atal
(conjecture~\ref{graphespar.conj.sp})~:

\index{minimalement~imparfait~(graphe~---)!partition~antisymétrique}
\begin{theoreme}[\cite{chudvovsky.r.s.t:spgt}, \no 4.9, page 16]
  \label{graphespar.t.bsp}
  Soit $G$ un graphe de Berge,  non parfait, et de taille minimum avec
  ces propriétés. Alors $G$ ne possède pas de partition antisymétrique
  paire.
\end{theoreme}

Nous  devons maintenant  définir une  nouvelle classe  de  graphes. On
appelle  \emph{bicographe}\index{bicographe} tout  graphe $G$  tel que
$V(G)$ peut se partitionner en quatre ensembles de taille au moins 2~:
$\{a_1, \dots,  a_m\}$, $\{b_1, \dots, b_m\}$,  $\{c_1, \dots, c_n\}$,
$\{d_1, \dots, d_n\}$, tels que~:

\vspace{1ex}

\begin{itemize}
\item 
  Pour tout $1\leq i \leq m$,  $a_i$ voit $b_i$ et pour tout $1 \leq j
  \leq n$, $c_j$ manque $d_j$.
\item
  Pour tout $1 \leq i < i' \leq m$, il n'y a aucune arête entre $\{a_i,
  b_i\}$ et $\{a_{i'}, b_{i'}\}$. Pour tout $1 \leq j < j' \leq m$, il y
  a  toutes les  arêtes possibles  entre $\{c_j,  d_j\}$  et $\{c_{j'},
  d_{j'}\}$.
\item 
  Pour  tout $1\leq  i \leq  m$ et  tout $1  \leq j  \leq n$,  il  y a
  exactement deux  arêtes entre $\{a_i,  b_i\}$ et $\{c_j,  d_j\}$, et
  ces deux arêtes sont non incidentes.
\end{itemize}

\vspace{2ex}

Quels sont  les plus petits bicographes~? Les  plus petits bicographes
ont 8 sommets~:  $a_1, a_2, b_1, b_2, c_1,  c_2, d_1, d_2$. L'ensemble
$\{c_1, c_2,  d_1, d_2\}$ induit un  $C_4$.  Le sommet  $a_1$ voit une
arête de ce $C_4$, et $b_1$  voit l'arête opposée. De même, $a_2$ voit
une arête du $C_4$ et  $b_2$ l'arête opposée.  À isomorphisme près, il
n'y a donc que  2 cas~: $a_1$, $a_2$ voient la même  arête du $C_4$ ou
$a_1$, $a_2$  voient des arêtes  incidentes du $C_4$. Dans  le premier
cas, on  obtient un  double-diamant, et  dans le  deuxième, $L(K_{3,3}
\setminus e)$. Le théorème suivant est facile~:
 
\begin{theoreme}[\cite{chudvovsky.r.s.t:spgt}, page 1]
  \label{graphespar.t.bicographe}
  Soit $G$ un bicographe. Alors $G$ est parfait.
\end{theoreme}

La preuve  du théorème de structure  des graphes de  Berge utilise une
technique classique,  d'ailleurs déjà mise  en \oe uvre  par Conforti,
Cornuéjols et  Vu\v skovi\'c pour prouver  le cas sans  carré. On part
d'un graphe  de Berge $G$, et  on suppose que $G$  contient un certain
sous-graphe $G'$ d'une  certaine classe $\cal C$. On  montre alors que
tout  le  graphe  s'organise  autour  de $G'$,  et  qu'il  possède  la
décomposition  souhaitée.  On  recommence  alors avec  les graphes  ne
contenant aucun sous-graphe de la  classe $\cal C'$\dots\ Voyons cela
plus en détails avec la suite de théorèmes ci-dessous~:

\index{line-graphe~de~subdivision~bipartie~de~$K_4$!décomposition}
\begin{theoreme}[\cite{chudvovsky.r.s.t:spgt}, \no 9.7, page 58]
  \label{graphespar.t.spgtlsbk4}
  Soit $G$ un  graphe de Berge.  Si $G$  contient le line-graphe d'une
  subdivision  bipartie  de  $K_4$,   alors  ou  bien  l'un  de  $G$,
  $\overline{G}$ est  un line-graphe  de biparti, ou  bien $G$  est un
  bicographe, ou bien l'un  de $G$, $\overline{G}$ possède un 2-joint,
  ou bien $G$ possède une partition antisymétrique paire.
\end{theoreme}

\begin{theoreme}[\cite{chudvovsky.r.s.t:spgt}, \no 10.6, page 62]
  \label{graphespar.t.spgtprismepair}
  Soit  $G$ un  graphe  de  Berge ne  contenant  aucun line-graphe  de
  subdivision  bipartie de  $K_4$. Si  $G$ contient  un  prisme pair,
  alors ou  bien $G$ est  un prisme pair  avec 9 sommets, ou  bien $G$
  possède un 2-joint, ou bien $G$ possède une partition antisymétrique
  paire.
\end{theoreme}

\begin{theoreme}[\cite{chudvovsky.r.s.t:spgt}, \no 13.4, page 86]
  \label{graphespar.t.prismeimp}
  Soit $G$  un graphe  de Berge  tel que ni  $G$ ni  $\overline{G}$ ne
  contient  de line-graphe de  subdivision bipartie  de $K_4$.  Si $G$
  contient  un  prisme  impair  long,  alors  ou  bien  l'un  de  $G$,
  $\overline{G}$ possède un 2-joint, ou bien $G$ possède une partition
  antisymétrique paire, ou bien $G$ possède une paire homogène.
\end{theoreme}

\index{double-diamant!décomposition}
\begin{theoreme}[\cite{chudvovsky.r.s.t:spgt}, \no 14.3, page 92]
  Soit $G$  un graphe  de Berge  tel que ni  $G$ ni  $\overline{G}$ ne
  contient  de line-graphe de  subdivision bipartie  de $K_4$  ou de
  prisme long. Si  $G$ contient un double-diamant,  alors ou bien l'un
  de $G$, $\overline{G}$  possède un 2-joint, ou bien  $G$ possède une
  partition antisymétrique paire.
\end{theoreme}

\begin{theoreme}[\cite{chudvovsky.r.s.t:spgt}, \no 1.3.7 à 1.3.11, page~4]
  \label{graphespar.t.sbipar}
  Soit $G$  un graphe  de Berge  tel que ni  $G$ ni  $\overline{G}$ ne
  contient  de  line-graphe de  subdivision  bipartie  de $K_4$,  de
  prisme long  ou de  double-diamant.  Alors  ou bien $G$  possède une
  partition  antisymétrique  paire,  ou  bien $G$  ne  contient  aucun
  antitrou  long, ou  bien $\overline{G}$  ne contient  aucun antitrou
  long.
\end{theoreme}

\begin{theoreme}[\cite{chudvovsky.r.s.t:spgt}, \no 1.3.7 à 1.3.12, page 4]
  \label{graphespar.t.sartemis}
  Soit $G$ un graphe de Berge sans  antitrou long et tel que ni $G$ ni
  $\overline{G}$ ne contient  de line-graphe de subdivision bipartie
  de $K_4$,  de prisme long ou  de double-diamant.  Alors  ou bien $G$
  possède  une partition  antisymétrique paire,  ou bien  $G$  est une
  clique ou bien $G$ est biparti.
\end{theoreme}

De  cette liste  de théorèmes,  on  déduit facilement  le théorème  de
décomposition des graphes de Berge~:

\index{décomposition!graphes de Berge}
\index{Berge~(graphe de ---)!décomposition}
\begin{theoreme} 
  \label{graphespar.t.structure}
  Soit $G$ un graphe de Berge.  Alors ou bien $G$ est basique, ou bien
  $G$ est un  bicographe, ou bien l'un de  $G$, $\overline{G}$ possède
  un 2-joint, ou bien  $G$ possède une partition antisymétrique paire,
  ou bien $G$ possède une paire homogène.
\end{theoreme}

 En  admettant tous les  théorèmes de  ce chapitre,  il est  facile de
 démontrer  le  théorème  fort  des  graphes parfaits~:  soit  $G$  un
 contre-exemple  de  taille  minimum  au  théorème  fort  des  graphes
 parfaits.    Notons    que   $G$   est   de    Berge.    D'après   le
 théorème~\ref{graphespar.t.faible},   $\overline{G}$  est   aussi  un
 contre-exemple     de      taille     minimum.      Appliquons     le
 théorème~\ref{graphespar.t.structure} à $G$.   Si $G$ est basique, on
 contredit  la  perfection  des   graphes  basiques,  si  $G$  est  un
 bicographe,  on contredit  le théorème~\ref{graphespar.t.bicographe},
 si  $G$  ou  $\overline{G}$  possède  un  2-joint,  on  contredit  le
 théorème~\ref{graphespar.t.2joint},   si   $G$   possède  une   paire
 homogène, on  contredit le théorème~\ref{grapphespar.t.paireh}  et si
 $G$  possède  une partition  antisymétrique  paire,  on contredit  le
 théorème~\ref{graphespar.t.bsp}.   Dans  tous   les  cas,  on  a  une
 contradiction.

\subsection*{Trigraphes}
\label{graphespar.ss.trigraphes}

\index{trigraphe} La thèse de Maria Chudnovsky~\cite{chudnovsky:these}
est consacrée  à la preuve  de la conjecture~\ref{graphespas.conj.ccv}
proprement    dite,   qui,    notons   le,    n'est    pas   démontrée
dans~\cite{chudvovsky.r.s.t:spgt}. Le premier  problème à résoudre est
trivial~:  tous les  bicographes\index{bicographe}  ont une  partition
antisymétrique (pas  forcément paire), et on peut  donc facilement les
éliminer         du          théorème         de         décomposition
(théorème~\ref{graphespar.t.structure}).    Reste   l'élimination  des
paires homogènes qui semble être un détail, puisqu'elles sont très peu
utilisées  dans la  preuve  du théorème~\ref{graphespar.t.structure}~:
elles           n'apparaissent          que           dans          le
théorème~\ref{graphespar.t.prismeimp}.    Toutefois  leur  élimination
semble être un problème difficile.

Les explications données ici  proviennent de notes personnelles prises
lors d'un exposé de Maria Chudnovsky à Palo Alto en novembre 2002.  La
technique  qu'elle utilise n'est  rien moins  qu'une relaxation  de la
notion de  graphe~: un trigraphe $G$  est un triplet  $(V, E_1, E_2)$.
Les  ensembles $E_1$  et $E_2$  sont des  sous-ensembles  disjoints de
l'ensemble $V \choose 2$.   L'ensemble $E_1$ représente l'ensemble des
arêtes \emph{obligatoires} de $G$, tandis que $E_2$ est l'ensemble des
arêtes    \emph{optionnelles}.    Plus   formellement,    on   appelle
\emph{réalisation}  de  $G$ tout  graphe  $G'=(V,  E)$ vérifiant  $E_1
\subseteq E  \subseteq E_1 \cup  E_2$.  Un trigraphe est  dit \emph{de
Berge}\index{Berge~(trigraphe~de~---)} si toutes  ses réalisations sont
des graphes de Berge. Notons qu'un  graphe de Berge peut être vu comme
un  trigraphe de  Berge  ne comportant  aucune  arête optionnelle,  et
n'ayant donc  pour seule  réalisation que lui-même.   Maria Chudnovsky
démontre alors un théorème de  structure pour les trigraphes de Berge,
qui,  appliqué  au  cas   particulier  d'un  graphe  de  Berge,  donne
exactement la  conjecture~\ref{graphespas.conj.ccv}.  Formellement, le
théorème de structure des trigraphes  est plus difficile à prouver que
le théorème de  structure des graphes, mais, comme  il arrive parfois,
démontrer  quelque   chose  de  plus  difficile   permet  de  disposer
d'hypothèses  d'induction  plus  fortes,  et s'avère  finalement  plus
facile.  Au total,  la preuve en passant par  les trigraphes est aussi
longue  et complexe  que celle  par les  graphes, mais  elle  donne un
résultat un peu plus fort~:

\index{décomposition!graphes de Berge}
\index{Berge~(graphe de ---)!décomposition}
\begin{theoreme}[Chudnovsky \cite{chudnovsky:these,chudnovsky:trigraphs}] 
  \label{graphespar.t.structuresM}
  Soit $G$ un graphe de Berge.  Alors ou bien $G$ est basique, ou bien
  l'un de $G$, $\overline{G}$ possède  un 2-joint, ou bien $G$ possède
  une partition antisymétrique.
\end{theoreme}

\section{Problèmes ouverts et questions}
\label{graphespar.s.poq}

Comme on l'a déjà dit,  de nombreuses questions concernant les graphes
parfaits sont  encore ouvertes, la plus importante  étant peut-être la
recherche  d'un   algorithme  combinatoire  de   coloration.  D'autres
questions apparaîtront  au fil de ce  travail. Mais le  problème de la
décomposition des graphes de Berge est en un sens encore ouvert.

En effet,  on attend en  général d'un théorème de  décomposition qu'il
permette de  démontrer des théorèmes  par induction.  En ce  sens, les
théorèmes~\ref{graphespar.t.structure}
et~\ref{graphespar.t.structuresM}  remplissent leur mission.   Mais on
attend  aussi que  les décompositions  fournissent des  algorithmes en
temps  polynomial fonctionnant  sur le  modèle suivant~:  on  donne un
graphe  $G$, on  vérifie  s'il est  basique.   S'il ne  l'est pas,  on
cherche une  décomposition qui permette de ``casser''  le graphe, puis
on réapplique  l'algorithme récursivement sur  les ``morceaux''.  Pour
autant  qu'on  sache,  le  théorème~\ref{graphespar.t.structuresM}  ne
permet  pas   de  mettre  en   \oe  uvre  cette  stratégie   en  temps
polynomial.  Pour  s'en  convaincre,  on  va essayer  de  résoudre  le
problème de la reconnaissance des graphes de Berge~:

On  nous donne  un graphe  $G$. Il  est facile  de tester  si  $G$ est
basique ---  en utilisant  entre autre l'algorithme  de reconnaissance
des line-graphes  (algorithme~\ref{graphespar.a.lehot}). S'il ne l'est
pas,  on recherche un  2-joint dans  $G$. Cela  est possible  en temps
polynomial  grâce   à  un  algorithme  de   Cornuéjols  et  Cunningham
\cite{cornuejols.cunningham:2join}\index{2-joint!détection}\index{détection!2-joints}.
S'il y  a un 2-joint  $(X,Y)$, on parvient  à ``casser'' le  graphe en
deux    morceaux   et    les    choses   se    passent   bien    (voir
\cite{cornuejols.liu.vuskovic:reco}  par  exemple  pour des  détails).
S'il n'y a pas de  2-joint, on recherche une partition antisymétrique.
Il existe  pour cela  un algorithme  en temps polynomial  dû à  C.  de
Figueiredo,      S.      Klein,      Y.      Kohayakawa      et     B.
Reed~\cite{figuereido.k.k.r:sp}. 
\index{partition~antisymétrique!détection}
\index{détection!partition~antisymétrique}

\index{décomposition~et~reconnaissance}
\index{décomposition!par~partition~antisymétrique}
L'ensemble $V(G)$ est  donc partitionné en deux ensembles  $X$ et $Y$,
où  $G[X]$ n'est  pas connexe  et  $G[Y]$ n'est  pas anticonnexe.   On
choisit $A$, composante connexe de $G[X]$  et on pose $B = X \setminus
A$.  On choisit $C$, composante anticonnexe  de $G[Y]$, et on pose $D =
Y \setminus  C$. Il n'y a aucune  arête entre $A$ et  $B$, mais toutes
les arêtes possibles  entre $C$ et $D$.  Il se peut  très bien que $G$
contienne  un trou  impair alors  que $G[A]$,  $G[B]$,  $G[C]$ $G[D]$,
$G[X]$ et $G[Y]$ sont de  Berge. Donc, pour reconnaître les graphes de
Berge, il  est inutile de relancer récursivement  l'algorithme sur $A,
B, C, D, X,  Y$. Par contre, il est facile de  vérifier que $G$ est de
Berge si  et seulement si  $G[A\cup B \cup  C]$, $G[A\cup B  \cup D]$,
$G[A\cup C \cup D]$ et $G[B\cup C \cup D]$ sont de Berge. Donc on peut
relancer l'algorithme récursivement sur ces quatre graphes. Finalement
dans le cas extrême où trois des ensembles $A, B, C, D$ sont de taille
1, on voit que si on  donne à l'algorithme un graphe avec $n$ sommets,
il  risque de  lancer 3  appels récursifs  sur des  graphes  de taille
$n-1$. Au pire, il rique d'y  avoir $3^{n-3}$ appels, ce qui donne une
complexité exponentielle.

Comme on le voit, la  partition antisymétrique n'est pas une véritable
décomposition,  au  sens  algorithmique  et  pour le  problème  de  la
reconnaissance des  graphes de  Berge tout du  moins. Pour  autant que
nous  le sachions,  nul  n'a la  moindre  idée aujourd'hui  de ce  que
pourrait être un théorème de  décomposition des graphes de Berge ayant
de bonnes propriétés algorithmiques.

\chapter{Paires d'amis} 
Nous présentons ici  la notion de paire d'amis  d'un graphe introduite
par Meyniel~\cite{meyniel:87}~:

\index{amis|see{paire~d'amis}}
\index{paire~d'amis}
\begin{defeng}
  On appelle \emph{paire d'amis} d'un graphe $G$ toute paire $\{x,y\}$
  de sommets de $G$ telle que  tous les chemins reliant $x$ à $y$ sont
  de longueur paire.
\end{defeng}

Nous verrons que les paires d'amis sont un outil efficace pour prouver
la  perfection  de  certaines   classes  de  graphes  de  Berge,  bien
qu'historiquement nous n'ayons  connaissance d'aucune classe de graphe
dont la perfection ait d'abord  été prouvée par les paires d'amis. Les
graphes de  Berge eux-même  ne font pas  exception à la  règle puisque
Chudnovsky \emph{et al.}  n'utilisent pas les paires d'amis, peut-être
par insuffisance de  l'outil, mais peut-être aussi parce  que leur but
était  de  prouver   la  conjecture~\ref{graphespas.conj.ccv}  sur  la
structure des graphes de Berge.  Dans cette optique, il ne fallait pas
recourir  aux  paires  d'amis.   Notons  toutefois que  la  notion  de
partition antisymétrique paire fait  appel à des contraintes de parité
de  chemins, que  l'on peut  formuler en  termes de  paires  d'amis de
sous-graphes.   Dans  un  sens  Chudnovsky \emph{et  al.}   font  donc
peut-être  un usage  implicite de  la  notion de  paire d'amis.   Mais
jusqu'à présent,  aucun de leurs  résultats n'a permis de  trouver une
paire d'amis dans un graphe de  Berge.  Il reste que les paires d'amis
pourraient  être un  outil  pour simplifier  certaines  parties de  la
preuve du théorème fort des graphes parfaits.

Dans  ce  chapitre, nous  nous  contenterons  de mentionner  certaines
classes de graphes dont la perfection peut se prouver grâce aux paires
d'amis.   Nous donnerons  des résultats  sur la  détection  des paires
d'amis  et  présenterons à  cette  occasion  un  nouvel algorithme  de
détection de paires d'amis  dans les line-graphes.  Nous mentionnerons
une conjecture de R. Thomas et F. Maffray sur les graphes bipartisans,
dont une démonstration pourrait remplacer les 50 dernières pages de la
preuve  de  la  conjecture  forte  des graphes  parfaits.   Puis  nous
montrerons  l'intérêt   des  paires  d'amis  pour   la  coloration  en
présentant   les  graphes   parfaitement   contractiles.   Nous   nous
concentrerons sur les  résultats les plus éclairants pour  la suite de
ce  mémoire et nous  renvoyons au  chapitre sur  les paires  d'amis du
livre \emph{Perfect graphs} publié sous  la direction de Bruce Reed et
Jorge   Ram\'irez   Alfons\'in~\cite{everett.f.l.m.p.r:ep}  pour   une
synthèse plus générale (mais  un peu moins récente).  Nous terminerons
par des tentatives de généralisation de la notion de paire d'amis.

\section{Graphes de quasi-parité}

Henri Meyniel a prouvé le théorème suivant~:
\index{minimalement~imparfait~(graphe~---)!paire~d'amis}
\index{paire~d'amis!graphes~minimalement~imparfaits}
\begin{theoreme}[Meyniel, \cite{meyniel:87}]
  Un graphe minimalement imparfait n'a pas de paire d'amis.
\end{theoreme}\index{triviale}
\index{subdivision!---~triviale}

\index{quasi-parité~(graphe~de~---)}
\index{quasi-parité~stricte~(graphe~de~---)}
\index{strict!quasi-parité~---|see{quasi-parité~stricte}}
Les trous impairs  n'ont pas de paire d'amis,  les antitrous longs non
plus (entre  deux sommets non-adjacents, il existe  toujours un chemin
de  longueur~3).  Le  théorème  de Meyniel  est  donc aujourd'hui  une
conséquence  triviale  du théorème  fort  des  graphes parfaits  mais,
évidemment,  Meyniel  l'avait  prouvé  directement.   Le  théorème  de
Meyniel   et   le   théorème   des  graphes   parfaits   de   Lov\'asz
(théorème~\ref{graphespar.t.faible})  ont pour conséquence  directe la
perfection des classes de graphes  suivantes~: on dit qu'un graphe $G$
est de \emph{quasi-parité}  si tout sous-graphe induit $H$  de $G$ non
réduit  à un  sommet contient  une paire  d'amis de  $H$ ou  une paire
d'amis de  $\overline{H}$.  Par exemple,  les graphes bipartis  et les
complémentaires de  bipartis sont des graphes de  quasi-parité. On dit
qu'un graphe  est de  \emph{quasi-parité stricte} si  tout sous-graphe
induit $H$ de $G$ différent  d'une clique contient une paire d'amis de
$H$.   Par   exemple,  les  graphes  bipartis  sont   des  graphes  de
quasi-parité  stricte.  Tout  graphe  de quasi-parité  stricte est  un
graphe de quasi-parité.

On peut se demander s'il est  facile de trouver des paires d'amis dans
un  graphe. Dan  Bienstock a  montré que  ce problème  est  en général
difficile~:

\index{paire~d'amis!détection~en~général}
\index{détection!paire~d'amis}
\begin{probleme}
  \label{pair.prob.bienstock}
  \begin{itemize}
  \item[\sc Instance  :] Un graphe $G$  et deux sommets $a$  et $b$ de
    $G$.
  \item[\sc Question~:] La paire  $\{a, b\}$ est-elle une paire d'amis
    de $G$~?
  \item[\sc        Complexité~:]        CoNP-complet       (Bienstock,
  \cite{bienstock:evenpair}).
  \end{itemize}
\end{probleme}

Mais  la  recherche de  paires  d'amis  peut  être réalisée  en  temps
polynomial  si  l'on  se  restreint  aux graphes  de  Berge,  grâce  à
l'algorithme      de      reconnaissance      des      graphes      de
Berge~(algorithme~\ref{reco.a.berge},    page~\pageref{reco.a.berge})~:
si $G$ est un graphe de Berge  et si $\{a,b\}$ est une paire de sommet
de $G$, alors on construit le  graphe $G'$ obtenu en ajoutant à $G$ un
sommet qui  ne voit que $a$ et  $b$. On constate alors  que $\{a, b\}$
est une paire d'amis de $G$ si et seulement si $G'$ est de Berge.

Nous  allons maintenant présenter  une classe  de graphes  parfaits de
quasi-parité,  et d'autres  classes que  l'on soupçonne  telles.  Nous
présenterons plus tard des  classes de graphes de quasi-parité stricte
car il s'avère que toutes les classes que nous mentionnerons seront en
fait  contenues  dans une  classe  aux  propriétés  plus riches~:  les
graphes parfaitement contractiles.

\subsection{Les graphes sans taureau}
\index{taureau!graphe sans ---}
On rappelle  que le  \emph{taureau} est un  graphe particulier  à cinq
sommets représenté ci-dessous~:

\begin{figure}[ht]
  \center
  \includegraphics{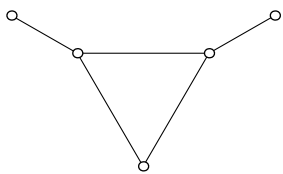}
\caption{Le taureau\label{pair.fig.taureau}}
\end{figure}

\index{quasi-parité~(graphe~de~---)!graphes~sans~taureau}
Va\v sek Chv\'atal  et Najiba Sbihi~\cite{chvatal.sbihi:bullfree} ont
démontré que les graphes de Berge sans taureau sont parfaits.  Par la
suite,  Celina de  Figueiredo,  Frédéric Maffray  et  Oscar Porto  ont
montré  que les graphes  de Berge  sans taureau  sont des  graphes de
quasi-parité~\cite{figuereido.m.p:bullfree}.  C'est  la classe la plus
générale de graphes de quasi-parité connue à ce jour. 

Bruce  Reed et  Najiba Sbihi~\cite{reed.sbihi:bullfree}  ont  donné un
algorithme en  temps polynomial pour reconnaître les  graphes de Berge
sans     taureau.      Celina     de    Figueiredo     et     Frédéric
Maffray~\cite{figuereido.m:bullfreeopt}  ont donné des  algorithmes en
temps polynomial pour résoudre les problèmes classiques d'optimisation
sur   ces  graphes,  la   coloration  notamment.    Leurs  algorithmes
n'utilisent pas les paires d'amis.

\subsection{Les conjectures de Hougardy}

Stephan  Hougardy~\cite{hougardy:95} a étudié  les paires  d'amis dans
les  line-graphes de  bipartis.  Il  a  donné un  algorithme en  temps
polynomial pour décider si  le line-graphe d'un graphe biparti possède
une                                                               paire
d'amis\index{paire~d'amis!détection~dans~line-graphes~de~biparti}\index{détection!paire~d'amis}.
Le théorème suivant montre que de nombreux line-graphes de bipartis ne
possèdent pas de paires d'amis.

\begin{theoreme}[Hougardy, \cite{hougardy:95}]
  \label{pair.t.hougardy}
  Soit   $R$   un   graphe   biparti  3-connexe.   Alors   $L(R)$   et
  $\overline{L(R)}$ n'ont pas de paire d'amis.
\end{theoreme}

Hougardy  a également  proposé des  conjectures pour  caractériser les
graphes de quasi-parité et les graphes de quasi-parité stricte~:

\index{quasi-parité~stricte~(graphe~de~---)!caractérisation}
\begin{conjecture}[Hougardy, \cite{hougardy:these}]
  Tout  graphe qui  n'est  pas  de quasi-parité  stricte,  et qui  est
  minimal pour l'inclusion  avec cette propriété, est ou  bien un trou
  impair, ou bien un antitrou long, ou bien le line-graphe d'un graphe
  biparti.
\end{conjecture}

La  conjecture ne  précise  pas quels  line-graphes  de bipartis  sont
minimalement non  de quasi-parité  stricte.  Elle est  encore ouverte,
mais un certain  nombre de cas particuliers ont  été résolus~: dans le
cas  des graphes  sans  taureau  par C.   de  Figueiredo, Maffray  et
Porto~\cite{figuereido.m.p:bullfree},  dans le  cas  des graphes  sans
diamant  par K\'ezdy  et Scobee~\cite{kezdy.scobee:diamond},  dans le
cas     des     graphes    sans     griffe     par    Linhares     et
Maffray~\cite{linhares.maffray:98}, dans le cas des graphes planaires
par  Linhares,  Maffray  et Reed~\cite{linhares.m.r:99}.   Hougardy  a
également proposé la conjecture suivante~:

\index{quasi-parité~(graphe~de~---)!caractérisation}
\begin{conjecture}[Hougardy, voir~\cite{everett.f.l.m.p.r:ep}]
  Soit $G$ un graphe qui n'est  pas de quasi-parité et qui est minimal
  pour  l'inclusion   avec  cette  propriété.   Alors   l'un  de  $G$,
  $\overline{G}$  est un  trou impair  ou le  line-graphe  d'un graphe
  biparti.
\end{conjecture}

\subsection{Les graphes bipartisans}

\index{bipartisan~(graphe~---)}
Les graphes  bipartisans sont définis implicitement dans  la preuve du
théorème  fort des  graphes  parfaits~\cite{chudvovsky.r.s.t:spgt}, et
explicitement  dans  l'article  de  synthèse  de  Chudnovsky  {\it  et
al.}~\cite{chudvovsky.r.s.t:progress}.       On     appelle     graphe
\emph{bipartisan}  tout  graphe  de   Berge  sans  prisme  long,  sans
complémentaire de prisme long,  sans double-diamant et sans $L(K_{3,3}
\setminus e)$.   Une conjecture, proposée  indépendamment par Frédéric
Maffray et Robin Thomas (à la conférence de Palo-Alto, novembre 2002),
stipule que les graphes  bipartisans sont des graphes de quasi-parité.
À la dernière section du chapitre~\ref{artemis.chap}, nous verrons que
les théorèmes~\ref{graphespar.t.sbipar} et~\ref{graphespar.t.sartemis}
page~\pageref{graphespar.t.sbipar} ne sont rien d'autre qu'un théorème
de  décomposition  des  graphes   bipartisans.   Une  preuve  de  leur
quasi-parité, et  donc de leur  perfection, permettrait de  prouver la
conjecture  forte  des  graphes  parfaits  sans  recours  à  ces  deux
théorèmes,   dont   la   preuve   occupe  les   50   dernières   pages
de~\cite{chudvovsky.r.s.t:spgt}.   Au  chapitre~\ref{reco.chap},  nous
donnerons  un  algorithme de  reconnaissance  des graphes  bipartisans
(algorithme~\ref{reco.a.bipartisan}   page~\pageref{reco.a.bipartisan},
complexité $O(n^9)$).

\index{bipartisan~(graphe~---)!paire~d'amis}
\index{paire~d'amis!graphes~bipartisans}
\index{bipartisan~(graphe~---)!conjecture}
\begin{conjecture}[Maffray et Thomas]
  Soit $G$ un graphe bipartisan. Alors $G$ est de quasi-parité.
\end{conjecture}

Notons  que la  réciproque de  la  conjecture est  fausse~: le  double
diamant n'est  pas bipartisan, et  pourtant, il est  de quasi-parité~;
les prisme longs sont des  graphes de quasi-parité stricte, et donc de
quasi-parité, sans être bipartisans.

\subsection{Détection des paires d'amis dans les line-graphes}

\index{paire~d'amis!détection~dans~les~line-graphes}
\index{détection!paire~d'amis}

Nous nous proposons ici d'étendre quelque peu l'algorithme de Hougardy
qui détecte les paires d'amis  dans les line-graphes de bipartis. Nous
donnons un  algorithme qui décide si  un line-graphe \emph{quelconque}
contient une  paire d'amis.   Il suffit pour  cela de  rassembler deux
résultats assez  éloignés~: l'algorithme~\ref{graphespar.a.lehot}, qui
permet de reconstruire la  racine d'un line-graphe, et une application
peu  connue  des  algorithmes  de  couplage due  à  Jack  Edmonds,  et
qu'Andr\'as Seb\H o m'a signalée~:

\index{chaîne!algorithme~de~plus~courte~---~de~longueur~paire}
\begin{algorithme}[Edmonds, Cf.  \cite{schrijver:opticombA},
    pages~515 et 458]
  \label{pair.a.edmonds}
  \begin{itemize}
  \item[\sc Instance  :] Un graphe $G$  et deux sommets $a$  et $b$ de
    $G$.
  \item[\sc Sortie  :] Une chaîne  de $G$, de longueur  paire, reliant
  $a$ et  $b$ et de longueur  minimale avec ces propriétés  --- si une
  telle chaîne existe.  Sinon, ``pas de chaîne de longueur paire''.
  \item[\sc Complexité :] $n^3$.
  \end{itemize}
\end{algorithme}

Notre algorithme repose aussi sur  le lemme suivant qui traduit l'idée
bien connue que les chaînes d'un graphe, même si elles ont des cordes,
deviennent des chemins quand on passe au line-graphe~:

\begin{lemme}
  \label{pair.l.chainechemin}
  Soit $R$ un  graphe et $G = L(R)$. Soit $F$  un ensemble d'arêtes de
  $R$. Alors $F$  est l'ensemble des arêtes d'une chaîne  de $R$ si et
  seulement si $F$ induit un chemin de $G$.
\end{lemme}

\begin{preuve}
  Si $F$ est l'ensemble des arêtes  d'une chaîne de $R$, alors on note
  $e_1,  e_2,  \dots,  e_{k}$   ses  arêtes,  dans  l'ordre  où  elles
  apparaissent sur la chaîne.  Ces arêtes de $R$ constituent aussi une
  suite de sommets de $G$ qui est clairement une chaîne de $G$.  Cette
  chaîne est sans corde car  une éventuelle corde entre des sommets de
  $G$, $e_i$ et $e_j$ avec $|i-j| > 1$, impliquerait que les arêtes de
  $R$,  $e_i$ et $e_j$,  soient incidentes  bien que  non consécutives
  dans la chaîne,  ce qui est impossible.  On a  montré que $F$ induit
  bien un chemin de $G$.
    
  Réciproquement,  si $F$  (vu comme  un ensemble  de sommets  de $G$)
  induit  un chemin  de  $G$, alors  on  note $v_1,  \dots v_{k}$  les
  sommets  de $F$,  dans l'ordre  où  ils apparaissent  sur $F$.   Les
  sommets $v_i$, $i= 1, \dots, k$, sont aussi des arêtes de $R$.  Pour
  tout $i \in \{1, \dots, k-1\}$,  il est clair que les arêtes de~$R$,
  $v_i$ et  $v_{i+1}$, sont incidentes.   On note $u_i \in  V(R)$ leur
  extrémité commune.  On note $u_0$ l'extrémité de $v_1$ qui n'est pas
  $u_1$  et  on  note  $u_{k}$  l'extrémité de  $v_k$  qui  n'est  pas
  $u_{k-1}$.   S'il existe $0  \leq i  < j  \leq k$  avec $u_i  = u_j$
  alors, ou bien il y a au  moins trois arêtes de $R$ incidentes en un
  même sommet qui donnent un triangle dans $P$ (absurde), ou bien $u_0
  = u_k$,  et il y a un  cycle dans $P$ (également  absurde).  On voit
  donc que  les sommets  $u_i$, $0 \leq  i \leq  k$, sont deux  à deux
  distincts.   Autrement dit, $F$  peut être  vu comme  l'ensemble des
  arêtes d'une chaîne de~$R$.
\end{preuve}

Voici notre algorithme. Son principe est très simple, l'idée étant que
si on  sait faire quelque chose  avec des chaînes dans  les graphes en
général, alors on  sait faire la même chose avec  des chemins dans les
line-graphes.

\index{chemin!plus~court~---~de~longueur~impaire}
\index{paire~d'amis!détection~dans~les~line-graphes}

\begin{algorithme}
  \label{pair.a.eplg}
  \begin{itemize}
  \item[\sc Instance :] Un line-graphe  $G$ et deux sommets $a$ et $b$
    de $G$.
  \item[\sc Sortie  :] Un  chemin de longueur  impaire reliant  $a$ et
    $b$, s'il existe un tel  chemin.  Sinon, ``$\{a,b\}$ est une paire
    d'amis de $G$''.
    \item[\sc Calcul~:] Grâce à l'algorithme~\ref{graphespar.a.lehot},
    calculer un graphe $R$ tel que  $G=L(R)$. On peut donc voir $a$ et
    $b$ comme des  arêtes de $R$.  On note  $u_{a,1}$ et $u_{a,2}$ les
    deux extrémités de l'arête $a$, on note $u_{b,1}$ et $u_{b,2}$ les
    deux     extrémités      de     l'arête     $b$.       Grâce     à
    l'algorithme~\ref{pair.a.edmonds}, pour  tout $i \in  \{1, 2\}$ et
    pour tout $j\in\{1, 2\}$, on  calcule une  chaîne de $R
    \setminus  \{  u_{a, i+1},  u_{a,  j+1}  \}$  notée $C_{i,j}$,  de
    longueur paire  et reliant  $u_{a,i}$ à $u_{b,j}$  (l'addition des
    indices s'entend  modulo 2).  Si  aucune de ces  chaînes n'existe,
    stopper et  répondre ``$\{a, b\}$  est une paire d'amis  de $G$''.
    Sinon,  considérer  l'ensemble $P$  des  arêtes  de  l'une de  ces
    chaînes, qui est aussi un  ensemble de sommets de $G$.  Stopper et
    retourner le sous-graphe de $G$ induit par $Q = P\cup \{a, b\}$.

  \item[\sc Complexité :] $n^3$. 
  \end{itemize}
\end{algorithme}

\begin{preuve}
  Si l'algorithme  d'Edmonds trouve une chaîne  $C_{i,j}$ de longueur
  paire, alors on note $e_1, e_2, \dots, e_{2\alpha}$ ses arêtes, dans
  l'ordre où  elles apparaissent  sur $C_{i,j}$.  On  note $e_0  = a$,
  $e_{2\alpha+1} = b$  et $Q = \{ e_0,  e_1, \dots, e_{2\alpha+1} \}$.
  Comme les  sommets de  $C_{i,j}$ sont choisis  dans $R  \setminus \{
  u_{a, i+1}, u_{a, j+1} \}$, l'ensemble $Q$ est l'ensemble des arêtes
  d'une chaîne de  $R$ qui comporte $2\alpha +  2$ arêtes.  D'après le
  lemme~\ref{pair.l.chainechemin},  $Q$  est   aussi  un  ensemble  de
  sommets  de $G$  qui induit  un chemin  de $G$  de  longueur impaire
  reliant $a$ à $b$~: l'algorithme donne la bonne réponse.
    
  Réciproquement, si  l'algorithme d'Edmonds  ne trouve pas  de chaîne
  $C_{i,j}$, alors on va montrer que $\{a,b\}$ est une paire d'amis de
  $G$.  Pour  cela, on suppose en vue  d'une contradiction l'existence
  dans $G$ d'un chemin $P'$ de  longueur impaire entre $a$ et $b$.  On
  note $v_0  = a$, $v_1$, \dots,  $v_{2\alpha}$, $v_{2\alpha+1}=b$ les
  sommets de  $P'$, dans  l'ordre où ils  apparaissent sur  $P'$.  Les
  sommets $v_k$,  $k= 0, \dots,  2\alpha+1$, sont aussi des  arêtes de
  $R$  qui, d'après  le lemme~\ref{pair.l.chainechemin},  peuvent être
  vues comme les  arêtes d'une chaîne de longueur  paire de $R$.  Mais
  alors $\{v_1, \dots, v_{2\alpha} \}$ est l'ensemble des arêtes d'une
  chaîne  $C_{i,j}$  de  longueur  paire que  l'algorithme  aurait  dû
  détecter~: une contradiction.
    
  Donc l'algorithme  donne la bonne  réponse. Sa complexité  est bien
  dominée par celle de l'algorithme~\ref{pair.a.edmonds}.
\end{preuve}

\section{Les graphes parfaitement contractiles}
\label{pair.s.pc}

Étant  donnés   deux  sommets  $x,y$  d'un  graphe   $G$,  on  définit
l'opération                                                          de
\emph{contraction}\index{contraction}\index{paire~d'amis!contraction}
de $x$ et  $y$~: on enlève $x$ et~$y$ de~$V$, et  on ajoute un nouveau
sommet noté $xy$ et relié à chaque sommet de $V \setminus \{x,y\}$ qui
voit l'un  au moins de  $x$ et~$y$.  On note  $G/xy$\index{0@$/$, dans
$G/xy$}     le     graphe     obtenu.      On     dit     que     $xy$
\emph{représente}\index{représenter~un~sommet} $x$ dans $G/xy$ (notons
que $xy$ représente aussi $y$).  Il est clair que pour tout graphe $G$
et toute paire de sommets non adjacents $\{x, y\}$, on a $\chi(G) \leq
\chi(G/xy)$.   En effet, étant  donnée une  $\chi(G/xy)$-coloration de
$G/xy$,  on  peut facilement  obtenir  une $\chi(G/xy)$-coloration  de
$G$~: on garde  la couleur des sommets différents de  $x$ et $y$, puis
on donne à  $x$ et $y$ la couleur du sommet  contracté.  On obtient de
même  l'inégalité  $\omega(G)   \leq  \omega(G/xy)$.   Ces  inégalités
peuvent être strictes~:  si $G$ est un chemin  avec quatre sommets, si
$x$  et $y$ sont  les extrémités  de ce  chemin, alors  $\chi(G)=2$ et
$\omega(G)=2$  tandis que  $\chi(G/xy)=3$  et $\omega(G/xy)=3$.   Jean
Fonlupt  et  Jean-Pierre  Uhry~\cite{fonlupt.uhry:82}  ont  prouvé  le
théorème suivant dont nous rappelons  la preuve classique, car il nous
semble que l'argument d'échange  bichromatique utilisé pour prouver la
conclusion~(\ref{pair.o.fonlupt0}) a des  chances d'être exploité dans
un contexte  plus général  (nous nous expliquerons  sur ce point  à la
section~\ref{pair.ss.generalisation.pp}   consacrée   aux   partitions
paires).

\index{paire~d'amis!coloration}
\begin{theoreme}[Fonlupt et Uhry~\cite{fonlupt.uhry:82}]
  \label{pair.t.fonlupt}

  Soit $G$ un graphe et $\{x, y\}$ une paire d'amis de $G$. Alors~:

  \begin{outcomes}
  \item 
    \label{pair.o.fonlupt0}  
    Il existe une coloration optimale de $G$ qui donne la même couleur
    à $x$ et $y$.
  \item 
    \label{pair.o.fonlupt1}
    $G$ et $G/xy$ ont le même nombre chromatique.
  \item
    \label{pair.o.fonlupt2}
    $G$ et $G/xy$ ont des cliques maximum de même taille.
  \item 
    \label{pair.o.fonlupt3}
    Si $G$ est parfait, alors $G/xy$ est parfait.
  \end{outcomes}

\end{theoreme}

\begin{preuve}
  Soit  $A$ une  coloration optimale  de $G$  ne donnant  pas  la même
  couleur à  $x$ et $y$  (disons que $x$  est rouge et $y$  bleu).  On
  appelle  $X$ l'ensemble  des sommets  rouges et  $Y$  l'ensemble des
  sommets bleus.  L'ensemble $X\cup  Y$ induit dans $G$ un sous-graphe
  biparti, et les  sommets $x$ et $y$  ne sont pas du même  côté de la
  bipartition $(X,  Y)$.  Cela signifie  que tous les chemins  de $G[X
  \cup Y]$ reliant $x$ à $y$ sont de longueur impaire.  Or, $\{x, y\}$
  étant  une paire  d'amis de  $G$,  on conclut  qu'il n'existe  aucun
  chemin reliant $x$  à $y$ dans $G[X \cup Y]$.  Donc  $x$ et $y$ sont
  dans deux composantes connexes distinctes  de $G[X \cup Y]$, et dans
  celle  de $x$,  on peut  échanger les  couleurs rouge  et  bleue.  On
  obtient alors  une $\chi(G)$-coloration de $G$ dans  laquelle $x$ et
  $y$ sont bleus.  Ceci prouve~(\ref{pair.o.fonlupt0}).

  Pour prouver~(\ref{pair.o.fonlupt1}), il suffit de vérifier $\chi(G)
  \geq   \chi(G/xy)$.    Pour  cela,   il   suffit   de  trouver   une
  $\chi(G)$-coloration de $G/xy$. Soit $A$ une $\chi(G)$-coloration de
  $G$  donnant la  même couleur  à $x$  et $y$,  dont  l'existence est
  garantie  par~(\ref{pair.o.fonlupt0}). On  obtient  alors facilement
  une $\chi(G)$ coloration de $G/xy$ en conservant la couleur de tous
  les sommets et en donnant  au sommet contracté la couleur commune de
  $x$ et $y$.

  Pour   prouver~(\ref{pair.o.fonlupt2}),  il  suffit   de  vérifier~:
  $\omega(G) \geq  \omega(G/xy)$. Pour cela, il suffit  de trouver une
  $\omega(G/xy)$-clique   dans    $G$.    Si   $G/xy$    possède   une
  $\omega(G/xy)$-clique ne  contenant pas le sommet  contracté, il est
  clair que  c'est aussi  une $\omega(G/xy)$-clique dans  $G$.  Sinon,
  soit  $K$ une  $\omega(G/xy)$-clique de  $G/xy$ contenant  le sommet
  contracté $xy$.   Le graphe $K  \setminus \{xy\}$ est une  clique de
  $G$. Si $x$ possède un non-voisin $a$ dans $K\setminus \{xy\}$ et si
  $y$ possède  un non-voisin  $b$ dans $K\setminus \{xy\}$,  alors $x$
  voit $b$ et $y$ voit $a$ car $K$ est une clique de $G/xy$. Donc $\bp
  x \tp  b \tp a  \tp y \ep$  est chemin de  longueur 3 reliant  $x$ à
  $y$~:  une  contradiction.   L'un  au  moins  de  $x$  et  $y$  voit
  entièrement $K\setminus  \{xy\}$ (dans $G$),  ce qui montre  que $G$
  possède bien une $\omega(G/xy)$-clique.

  Pour  prouver~(\ref{pair.o.fonlupt3}), il  suffit  de remarquer  que
  pour tout sous-graphe induit $H$ de $G/xy$, il existe un sous-graphe
  induit  $H'$  de  $G$  vérifiant~:  $H  =  H'/xy$  (si  $x$  ou  $y$
  n'appartient pas à $V(H')$, on adopte la convention naturelle $H'/xy
  =  H'$).    Comme  $H'$  est  parfait  par   hypothèse,  il  résulte
  par~(\ref{pair.o.fonlupt1})     et~(\ref{pair.o.fonlupt2})     que~:
  $\chi(H)  = \chi(H')  = \omega(H')  =  \omega(H)$ ce  qui montre  la
  perfection de $G/xy$.
\end{preuve}

\begin{figure}
  \center
  \begin{tabular}{cc}
    \includegraphics{fig.pair.1} & \includegraphics{fig.pair.2} \\
    \parbox{5cm}{\center Un graphe $G$ imparfait avec une paire d'amis
    $\{x,y\}$}&
    \parbox{5cm}{\center Pourtant le graphe $G/xy$ est parfait.}
  \end{tabular}
  \caption{Exemple de contraction d'une paire d'amis\label{pair.fig.cassenpar}}
\end{figure}

Il faut noter  que la contraction d'une paire  d'amis peut transformer
un     graphe    imparfait    en     un    graphe     parfait    (voir
figure~\ref{pair.fig.cassenpar}).  Ceci montre que les contractions de
paires d'amis  sont apparemment un mauvais outil  pour reconnaître les
graphes parfaits.\index{paire~d'amis!outil~de~reconnaissance}

Le  théorème  ci-dessus permet  de  définir  un  algorithme simple  de
coloration pour  les graphes  parfaits. Si $G$  est parfait et  si $G$
possède une paire d'amis $\{x,y\}$,  ce qui précède montre qu'à partir
d'une  coloration  optimale de  $G/xy$,  on  déduit immédiatement  une
coloration  optimale de  $G$ en  donnant à  $x$ et  $y$ la  couleur du
sommet  contracté  $xy$.  L'algorithme  fonctionne  comme suit~:  tant
qu'on trouve dans le graphe $G$  une paire d'amis, on la contracte. On
peut  donc suivre la  trace de  tout sommet  $v$ de  $G$ au  cours des
contractions successives et chaque sommet de $G$ est représenté par un
unique  sommet  dans  le  graphe  résultant.   On  colorie  le  graphe
résultant  (qui ne possède  pas de  paire d'amis),  et l'on  en déduit
immédiatement une coloration  du graphe de départ en  donnant à chaque
sommet la  couleur de  son représentant dans  le graphe  résultant. On
peut aussi  définir sur le même  modèle un algorithme  de recherche de
clique maximum.   L'idée, qui apparaît  dans la preuve du  théorème de
Fonlupt et Uhry, est que si  $\{a,b\}$ est une paire d'amis de $G$, et
que si  on connait  une clique  maximum $K$ de  $G/ab$, alors  on peut
calculer une clique maximum de $G$~:  si $K$ ne contient pas le sommet
contracté $ab$, alors $K$ est  aussi une clique maximum de $G$.  Sinon
l'un des  deux ensembles $(K\setminus ab) \cup  \{a\}$ ou $(K\setminus
ab) \cup \{b\}$ induit une clique maximum de $G$.

Ces algorithmes posent deux  problèmes~: tout d'abord, comment trouver
une paire  d'amis dans un  graphe~? On a  déjà vu que la  recherche de
paire       d'amis       était       un       problème       difficile
(problème~\ref{pair.prob.bienstock}), mais que leur détection dans les
graphes  de Berge  est réalisable  en temps  polynomial.   Le deuxième
problème posé  par l'algorithme  est encore ouvert~:  comment colorier
efficacement  un graphe  de  Berge sans  paire  d'amis~?  Résoudre  ce
problème ouvert  permettrait de colorier tout graphe  parfait en temps
polynomial.

Évidemment, à la  suite de diverses contractions de  paires d'amis, on
peut espérer tomber par chance  sur un graphe trivial à colorier. Pour
toute classe $\cal  C$ de graphes faciles à  colorier, on peut définir
la classe  des graphes $G$  tels qu'en partant  de $G$ et que  par une
suite de contraction  de paires d'amis on puisse  parvenir à un graphe
de $\cal  C$.  À  notre connaissance une  approche aussi  générale n'a
jamais donné de résultat intéressant, et  seul le cas où $\cal C$ est
la classe des cliques a été sérieusement étudié dans la littérature.

Ainsi Marc Bertschi  \cite{bertschi:90} a-t-il proposé les définitions
suivantes.          Un        graphe        $G$         est        dit
\emph{contractile}\index{contractile~(graphe~---)}  si   $G$  est  une
clique ou  s'il existe une suite  $G_0, \ldots, G_k$  de graphes telle
que $G=G_0$, et  que pour $i=0, \ldots, k-1$,  le graphe $G_i$ possède
une paire d'amis  $\{x_i, y_i\}$ telle que $G_{i+1}  = G_i/x_iy_i$, et
que $G_k$  est une clique.   Le graphe $G$ est  dit \emph{parfaitement
contractile}\index{parfaitement~contractile~(graphe~---)}    si   tout
sous-graphe induit de $G$ est contractile.

Quels sont les  graphes non parfaitement contractiles minimaux  ? On a
déjà vu  que les  trous impairs  et les antitrous  longs n'ont  pas de
paires d'amis.   Les prismes impairs peuvent avoir  des paires d'amis,
mais pour tout  choix de séquence de contraction,  on arrive au graphe
$\overline{C_6}$ qui est un antitrou  long (ce fait n'est pas trivial,
voir~\cite{linhares.m.r:97}  pour des  détails).  Les  prismes impairs
sont donc  des graphes non  parfaitement contractiles, dont  on montre
facilement qu'ils  sont minimaux.  En  résumé, les trous  impairs, les
antitrous  longs et  les prismes  impairs sont  les seuls  graphes non
parfaitement contractiles minimaux connus à ce jour.

Bruce     Reed     a     proposé    d'appeler     \emph{graphes     de
Grenoble}\index{Grenoble~(graphe~de~---)}   les   graphes  sans   trou
impair,  sans  antitrou long  et  sans  prisme  impair.  Craignant  la
concurrence de la noix de  Grenoble, nous proposons pour notre part de
les             appeler             \emph{graphes            d'Artémis
pairs}\index{Artémis~pair~(graphe~d'---)}.   Nous suivons  en revanche
cet         auteur          pour         appeler         \emph{graphes
d'Artemis}\index{Artémis~(graphe~d'---)} les graphes sans trou impair,
sans antitrou long  et sans prisme.  Notons que  tout graphe d'Artémis
est  un graphe d'Artémis  pair.  La  conjecture suivante  peut sembler
triviale alors qu'elle est en fait encore ouverte~:

\begin{conjecture}[Everett et Reed \cite{reed:93}]
  \label{pair.conj.EverettReed2}
  Soit  $G$   un  graphe  parfaitement   contractile  différent  d'une
  clique.  Alors $G$  possède  une paire  d'amis  dont la  contraction
  redonne un graphe parfaitement contractile.
\end{conjecture}

Hazel  Everett et  Bruce Reed  ont proposé  une autre  conjecture pour
caractériser les graphes parfaitement contractiles~:

\index{Artémis~pair~(graphe~d'---)!paire~d'amis}
\index{paire~d'amis!graphes~d'Artémis~pairs}
\index{Artémis~pair~(graphe~d'---)!conjecture}
\begin{conjecture}[Everett et Reed \cite{reed:93}]
  \label{pair.conj.EverettReed}
  Un graphe est  parfaitement contractile si et seulement  s'il est un
  graphe d'Artémis pair.
\end{conjecture}
\index{parfaitement~contractile~(graphe~---)}

Une conjecture plus faible a également été proposée~:

\index{Artémis~(graphe~d'---)!paire~d'amis}
\index{paire~d'amis!graphes~d'Artémis}
\index{Artémis~(graphe~d'---)!conjecture}
\begin{conjecture}[Everett et Reed \cite{reed:93}]
  \label{pair.conj.Artemis}
  Si un graphe est d'Artémis, alors il est parfaitement contractile.
\end{conjecture}

Des cas  particuliers de ces  conjectures ont été démontrés~:  dans le
cas    des    graphes    planaires    par   Linhares,    Maffray    et
Reed~\cite{linhares.m.r:99}, dans  le cas des graphes  sans griffe par
Linhares  et  Maffray~\cite{linhares.maffray:98} et  dans  le cas  des
graphes    sans    taureau    par    de   Figueiredo,    Maffray    et
Porto~\cite{figuereido.m.p:bullfree}.   Irena  Rusu  à  donné  des
conditions suffisantes  pour qu'un graphe d'Artémis  possède une paire
d'amis~\cite{rusu:artemis}.  Mais le pas,  selon nous décisif, vers la
preuve  de  la conjecture~\ref{pair.conj.Artemis}  a  été franchi  par
Linhares et Maffray qui ont prouvé le théorème suivant~:

\index{carré!graphe~d'Artémis~sans~---}
\begin{theoreme}[Linhares et Maffray, \cite{linhares.maffray:evenpairsansc4}]
  Soit  $G$ un  graphe d'Artémis  sans carré.  Alors $G$  possède une
  paire d'amis.
\end{theoreme}

Nous    démontrerons    la    conjecture~\ref{pair.conj.Artemis}    au
chapitre~\ref{artemis.chap}. Notons que nous prouverons par là-même la
perfection des  graphes d'Artémis indépendamment du  théorème fort des
graphes parfaits.  Notre preuve ressemble à celle du cas sans carré, à
l'exception notable de l'usage d'un lemme fameux dû à Roussel et Rubio
--- voir  chapitre~\ref{rr.chap}.   Elle   donnera  un  algorithme  de
coloration   des    graphes   d'Artémis   de    complexité   $O(mn^2)$
(algorithme~\ref{artemis.a.color} page~\pageref{artemis.a.color}).  La
conjecture   la  plus   forte  reste   quant  à   elle   ouverte.   Au
chapitre~\ref{reco.chap},  consacré aux  problèmes  de reconnaissance,
nous donnerons  un algorithme de reconnaissance  des graphes d'Artémis
(algorithme~\ref{reco.a.artemis}         page~\pageref{reco.a.artemis},
complexité $O(n^9)$),  et un algorithme de  reconnaissance des graphes
d'Artémis            pairs           (algorithme~\ref{reco.a.grenoble}
page~\pageref{reco.a.grenoble},  complexité   $O(n^{20})$).   Nous  en
déduirons  un algorithme  en temps  polynomial donnant  une coloration
\emph{non  nécessairement   optimale}  des  graphes   d'Artémis  pairs
(algorithme~\ref{reco.a.colorpc}         page~\pageref{reco.a.colorpc},
complexité  $O(n^{23})$).   Si  les  conjectures de  Reed  et  Everett
(conjectures~\ref{pair.conj.EverettReed}
et~\ref{pair.conj.EverettReed2})  sont vraies,  alors il  sera immédiat
que cet  algorithme donne  une coloration \emph{optimale}  des graphes
d'Artémis pair.

%
%

\begin{figure}[p]
  \center
  \label{pair.fig.inter}
  \includegraphics{fig.pair.15}  \\
  \vspace{3ex}

  \begin{parbox}{5cm}
    { \begin{tabular}{lcl}
 F.    T.  & : & Faiblement  triangulé.\\   P.O. &  : & Parfaitement
      ordonnable.\\  Q.P. & : & Quasi-parité.\\  Q.P.S. & : & Quasi-parité
      stricte.\\  P. C.  & : & Parfaitement  contractile. 
      \end{tabular}
      \vspace{2ex}}
  \end{parbox}
  \caption{Relations  d'intersections et d'inclusions pour
  huit classes de graphes parfaits}
\end{figure}

Nous  allons maintenant passer  en revue  quelques classes  de graphes
parfaitement  contractiles.    Tous  ces  graphes   sont  des  graphes
d'Artémis et notre preuve de la conjecture~\ref{pair.conj.Artemis} est
donc  une alternative  pour prouver  leur parfaite  contractilité.  La
figure~\ref{pair.fig.inter}   page~\pageref{pair.fig.inter}  donne  un
exemple de  graphe pour chaque  possiblilité concernant l'appartenance
aux  classes suivantes  : Berge,  quasi-parité  stricte, quasi-parité,
parfaitement contractile,  Artémis, Meyniel, faiblement  triangulés et
parfaitement ordonnables.

\subsection{Les graphes  de Meyniel}
\label{pair.ss.meyniel}
\index{Meyniel~(graphe~de~---)}

\index{décomposition!graphes de Meyniel} \index{reconnaissance!graphes
de        Meyniel}       \index{Meyniel~(graphe~de~---)!décomposition}
\index{Meyniel~(graphe~de~---)!reconnaissance} Henri  Meyniel a défini
la  classe des graphes  dont tous  les cycles  impairs de  longueur au
moins  5  possèdent  au   moins  deux  cordes.   Ces  graphes  portent
aujourd'hui   son  nom.    Il  a   démontré  que   ces   graphes  sont
parfaits~\cite{meyniel:76}. De manière indépendante, S.E. Markosian et
I.A.  Karapetian~\cite{markosian.karapetian:76}  sont parvenus au même
résultat.   Meyniel a  montré  que  ses graphes  sont  des graphes  de
quasi-parité stricte~\cite{meyniel:87} par  une preuve assez simple et
astucieuse qui repose sur un  petit lemme (nous le verrons au chapitre
suivant~:  lemme~\ref{rr.l.rrmeyniel})  qu'on  peut  voir  aujourd'hui
comme  un  cas  particulier  du  fameux  lemme  de  Roussel  et  Rubio
(lemme~\ref{rr.l.w})         Notre         preuve        de         la
conjecture~\ref{pair.conj.Artemis}  (chapitre~\ref{artemis.chap})  met
en  \oe uvre  des idées  déjà  présentes dans  le cas  des graphes  de
Meyniel --- la preuve de Meyniel sera alors rappelée pour illustrer le
propos  (section~\ref{artemis.ss.meyniel}).   Michel  Burlet  et  Jean
Fonlupt~\cite{burlet.fonlupt:meyniel}   ont  donné   un   théorème  de
décomposition des graphes de Meyniel  qui permet de les reconnaître en
temps  polynomial. Ch\'inh  Ho\`ang~\cite{hoang:87}  a démontré  qu'un
graphe est  de Meyniel  si et seulement  si il est  fortement parfait,
c'est-à-dire si  et seulement si  pour tout sous-graphe induit  $H$ de
$G$, chaque sommet de $H$  appartient à un stable qui rencontre toutes
les cliques maximales de $H$.
 
\begin{figure}
  \center
  \includegraphics{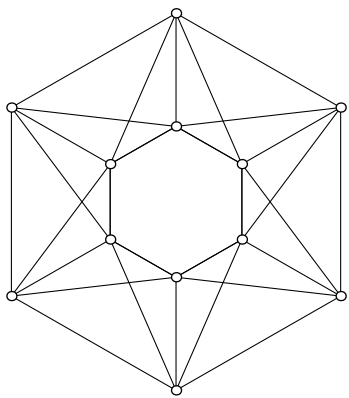}
  \caption{Un graphe de Meyniel ne contenant aucune paire d'amis dont
  la contraction redonne un graphe de  Meyniel.}
  \label{pair.fig.meyniel}
\end{figure}

\index{coloration!graphes                  de                 Meyniel}
\index{Meyniel~(graphe~de~---)!coloration}
\index{Meyniel~(graphe~de~---)!paire~d'amis}
\index{paire~d'amis!graphes~de~Meyniel}  Hertz~\cite{hertz:90}   a
démontré que les graphes  de Meyniel sont parfaitement contractiles et
en a  déduit un algorithme  de coloration de complexité  $O(nm)$.  Une
difficulté mérite d'être mentionnée~: Bertschi~\cite{bertschi:these} a
remarqué que  dans un graphe de  Meyniel, il n'existe  pas toujours de
paire  d'amis dont  la contraction  donne un  graphe de  Meyniel (voir
figure~\ref{pair.fig.meyniel}  un exemple  dû  à Sarkossian  mentionné
dans~\cite{everett.f.l.m.p.r:ep}).   Hertz a  donc  défini une  classe
plus   large   que   les   graphes   de  Meyniel,   les   graphes   de
quasi-Meyniel\index{quasi-Meyniel~(graphe~de~---)}                (voir
~\cite{everett.f.l.m.p.r:ep}),  et   a  montré  que   tout  graphe  de
quasi-Meyniel contient  une paire d'amis dont la  contraction donne un
graphe  de  quasi-Meyniel.   Celina  de Figueiredo  et  Kristina  Vu\v
skovi\'c~\cite{figuereido.vuskovic:meyniel}  ont  donné un  algorithme
qui reconnait les graphes  de quasi-Meyniel.  Florian Roussel et Irena
Rusu~\cite{roussel.rusu:meyniel} ont donné un algorithme de complexité
$O(n^2)$  pour  colorier  les  graphes de  Meyniel.   Tout  récemment,
Benjamin  Lévêque  et Frédéric  Maffray~\cite{leveque.maffray:meyniel}
ont  simplifié  et perfectionné  cet  algorithme  pour  parvenir à  un
algorithme en temps linéaire.

\subsection{Les graphes faiblement triangulés}

Va\v  sek Chv\'atal  et Ryan  Hayward~\cite{hayward:wt} ont  défini la
notion de graphe faiblement  triangulé. Un graphe est \emph{faiblement
triangulé} si et seulement s'il ne contient ni trou long (c'est-à-dire
de longueur au moins 5) ni  antitrou long. Hayward a donné un théorème
de décomposition  pour les graphes faiblement triangulés  qui, à cause
du théorème~\ref{graphespar.t.scs}, entraîne leur perfection~:

\index{faiblement~triangulé~(graphe~---)!décomposition}
\index{décomposition!graphes faiblement triangulés}
\begin{theoreme}[Hayward, \cite{hayward:wt}]
  Un graphe $G$ est faiblement  triangulé si et seulement si pour tout
  sous-graphe induit $H$ de $G$ on a ou bien :
  \begin{itemize}
    \item 
      $H$ est une clique.
    \item 
      $H$ est le complémentaire d'un couplage parfait.
    \item 
      $H$ possède une étoile d'articulation.
  \end{itemize}
\end{theoreme}

Ch\'inh  Ho\`ang  et  Frédéric Maffray~\cite{hoang.maffray:wtsqp}  ont
montré  que les  graphes  faiblement triangulés  sont  des graphes  de
quasi-parité   stricte.    Leur   preuve   utilise  le   théorème   de
décomposition ci-dessus.   Par la  suite, Hayward, Ho\`ang  et Maffray
ont   prouvé  un   théorème  plus   fort,  nécessitant   une  nouvelle
définition~:      on      dit      que     $\{x,y\}$      est      une
\emph{2-paire}\index{2-paire}\index{paire@2-paire}  de   $G$  si  tout
chemin  reliant  $x$  à  $y$  est de  longueur~2.   Les  2-paires  ont
l'avantage par  rapport aux paires d'amis d'être  faciles à détecter~:
$\{x, y\}$ est  une 2-paire d'un graphe $G$ si et  seulement si $x$ et
$y$ sont  dans deux  composantes connexes distinctes  de $G[V\setminus
N(x)      \cap     N(y)]$.       Arikati      et     Pandu      Rangan
\cite{arikati.pandurandan:tp}  ont donné  un algorithme  de complexité
$O(nm)$, qui prend en entrée un  graphe et qui en retourne une 2-paire
s'il                     y                     en                    a
une\index{2-paire!détection}\index{paire@2-paire!détection}\index{détection!2-paires}.
Voici le théorème de Hayward {\it et. al.}~:

\index{faiblement~triangulé~(graphe~---)!paire~d'amis}
\index{paire~d'amis!graphes~faiblement~triangulés}
\index{faiblement~triangulé~(graphe~---)!existence~d'une~2-paire}
\index{2-paire!graphes~faiblement~triangulés}
\index{paire@2-paire!graphes~faiblement~triangulés}
\begin{theoreme}[Hayward, Ho\`ang, Maffray, \cite{hayward.hoang.m:90}]
  \label{pair.t.deuxp}
  $G$  est faiblement triangulé  si et  seulement si  tout sous-graphe
  induit de $G$ différent d'une clique possède une 2-paire.
\end{theoreme}

\index{reconnaissance!graphes~faiblement~triangulés}
\index{faiblement~triangulé~(graphe~---)!reconnaissance}

Ce théorème a pour conséquence immédiate la parfaite contractilité des
graphes faiblement triangulés car la contraction d'une 2-paire dans un
graphe   faiblement    triangulé   redonne   un    graphe   faiblement
triangulé. C'est aussi  à notre connaissance le seul  théorème de type
``paire  d'amis'' qui permette  de reconnaître  une classe  de graphes
parfaits.  En effet, Spinrad et Sritharan~\cite{spinrad.sritharan:awt}
ont remarqué que si un graphe $G$ possède une 2-paire, alors le graphe
obtenu  en  reliant  les  deux  sommets de  la  paire  est  faiblement
triangulé  si et  seulement  si $G$  est  faiblement triangulé.   Pour
reconnaître  un  graphe  faiblement  triangulé,  il  suffit  donc  d'y
rechercher  une  2-paire, et  le  cas  échéant  d'en relier  les  deux
sommets, puis de recommencer sur  le graphe obtenu jusqu'à trouver une
clique.  Si à une  étape il n'y a pas de 2-paire,  c'est que le graphe
de  départ n'était  pas  faiblement triangulé.   En perfectionnant  ce
principe,         Hayward,         Spinrad        et         Sritharan
\cite{hayward.spinrad.sritharan:optwt}  ont proposé  un  algorithme de
reconnaissance de complexité $O(n+m^2)$  qui repose sur l'existence de
2-paires.  Il  faut noter qu'il existe un  algorithme pour reconnaître
les  graphes faiblement  triangulés, moins  performant  (de complexité
$O(n^5)$), mais reposant sur  des idées naïves~: nous rappellerons cet
algorithme  à titre d'illustration  au chapitre  sur les  problèmes de
reconnaissance (sous-section~\ref{reco.s.wt}).  Le meilleur algorithme
à ce jour pour la  coloration des graphes faiblement triangulés a pour
complexité        $O(n+m^2)$        (Hayward,        Spinrad        et
Sritharan~\cite{hayward.spinrad.sritharan:optwt}).

En nous restreignant aux  graphes faiblement triangulés, nous pourrons
donner une version plus  forte (lemme~\ref{rr.l.rrwt}) du fameux lemme
de  Roussel et  Rubio  (lemme~\ref{rr.l.w}).  Ceci  nous permettra  de
spécialiser notre  preuve de la  conjecture~\ref{pair.conj.Artemis} au
cas des  graphes faiblement triangulés, ce qui  donnera une nouvelle
preuve      du       théorème      d'existence      d'une      2-paire
(théorème~\ref{pair.t.deuxp}),        qui        sera        présentée
section~\ref{artemis.ss.ft}.

\subsection{Les  graphes parfaitement ordonnables}
\label{pair.ss.po}

\index{parfaitement~ordonnable~(graphe~---)}
\index{ordonnable~(graphe~parfaitement~---)}
Va\v sek Chv\'atal~\cite{chvatal:84} a défini les graphes parfaitement
ordonnables et  a démontré qu'ils  sont parfaits (voir~\cite{hoang:po}
pour une synthèse sur le sujet).  Étant donné un graphe quelconque $G$
et un  ensemble de  couleurs $\{1,  2, \dots \}$,  il est  possible de
colorier le graphe  de manière gloutonne, c'est à  dire en définissant
un  ordre sur  les sommets,  puis en  coloriant les  sommets  dans cet
ordre, en affectant à chaque  sommet la plus petite couleur disponible
en tenant  compte de ses  voisins déjà coloriés.  Pour  n'importe quel
graphe, il existe d'ailleurs un  ordre sur les sommets qui donne ainsi
une coloration optimale, et  calculer cet ordre est donc NP-difficile.
Mais  existe-t-il un  ordre sur  les sommet  de $G$  qui,  restreint à
n'importe quel  sous-graphe induit $H$, donne  une coloration optimale
de $H$~?  Si tel est le  cas, on dit que l'ordre est \emph{parfait} et
que $G$ est \emph{parfaitement ordonnable}.

\index{parfaitement~ordonnable~(graphe~---)!reconnaissance}
\index{reconnaissance!graphes~parfaitement~ordonnables}   Il  n'existe
pas  de  moyen simple  pour  décider  si  un graphe  est  parfaitement
ordonnable~:   Middendorf   et   Pfeiffer~\cite{middendorf.p:90}   ont
démontré  que la reconnaissance  des graphes  parfaitement ordonnables
est  un  problème NP-difficile.   En  fait, Ho\`ang~\cite{hoang:96}  a
montré que ce problème reste NP-difficile, même si on le restreint aux
graphes faiblement triangulés.  Cependant, Chv\'atal~\cite{chvatal:84}
a  trouvé  un  moyen   de  certifier  qu'un  graphe  est  parfaitement
ordonnable~: un graphe $G$ est parfaitement ordonnable si et seulement
s'il est possible d'orienter les arêtes de $G$ de manière acyclique et
sans obstruction, une  \emph{obstruction} étant un $P_4$ $\bp  a \tp b
\tp  c \tp  d  \ep$, sous-graphe  induit de  $G$,  avec l'arête  $ab$
orientée de $a$ vers $b$ et l'arête $cd$ orientée de $d$ vers $c$.  Ce
critère ne donne évidemment  pas d'algorithme en temps polynomial pour
décider  qu'un  graphe  est  parfaitement  ordonnable,  mais  dans  de
nombreux cas particuliers, il permet de s'en convaincre rapidement par
une suite de forçages. C'est  le cas pour les graphes non-parfaitement
ordonnables  de  la  figure~\ref{pair.fig.inter}.   Si un  graphe  $G$
possède  une orientation  sans  obstruction (mais  avec peut-être  des
cycles   orientés)   alors  on   dit   qu'il  est   \emph{parfaitement
orientable}\index{parfaitement~orientable~(graphe~---)}\index{orientable~(graphe~parfaitement~---)}.
Arikati et Peled~\cite{arikati.pandurandan:po} ont donné un algorithme
de détection des paires d'amis dans les graphes parfaitement orientés.


\index{parfaitement~ordonnable~(graphe~---)!paire~d'amis}
\index{paire~d'amis!graphes~parfaitement~ordonnables}
Meyniel~\cite{meyniel:87}  a  démontré  que les  graphes  parfaitement
ordonnables  sont des graphes  de quasi-parité  stricte.  Hertz  et de
Werra~\cite{hertz.dewerra:88} on démontré que les graphes parfaitement
ordonnables  sont  parfaitement   contractiles.   Il  faut  noter  que
l'existence de la paire d'amis  est démontrée en présupposant connu un
ordre parfait  sur les sommets de  $G$.  Or calculer un  tel ordre est
sans  doute NP-difficile,  puisque, comme  on  l'a vu,  décider si  un
graphe est  parfaitement ordonnable est  NP-difficile. Il y a  donc de
grandes  chances que  la méthode  de Hertz  et de  Werra ne  donne pas
d'algorithme de coloration en  temps polynomial.  Cela n'est toutefois
pas prouvé rigoureusement.

\index{parfaitement~ordonnable~(graphe~---)!coloration}
\index{coloration!graphes~parfaitement~ordonnables}
Notre    algorithme    de    coloration    des    graphes    d'Artémis
(algorithme~\ref{artemis.a.color})  fonctionne en particulier  sur les
graphes parfaitement  ordonnables.  En dehors  de l'algorithme général
de   coloration    des   graphes   parfaits   par    la   méthode   de
l'ellipsoïde~(algorithme~\ref{graphpar.a.coloration}),  c'est  à notre
connaissance le  premier algorithme en temps  polynomial pour colorier
les graphes parfaitement ordonnables.

\subsection{Les graphes avec moins de 10 sommets}

Pour finir,  mentionnons que Stephan  Hougardy \cite{hougardy:class} a
recensé 117 classes de graphes  parfaits, pour chacune desquelles il a
pu calculer par  ordinateur le nombre des graphes à  $k$ sommets, $k =
1, \dots,  10$.  Nous avons  extrait de son article  les dénombrements
concernant   les    classes   étudiées   dans    ce   chapitre   (voir
tableau~\ref{pair.tab.hougardy}). Ceci permet de  se faire une idée de
l'importance  quantitative  des   classes,  idée  peut-être  trompeuse
d'ailleurs pour  qui s'intéresse aux  graphes de plus de  dix sommets.
Selon Hougardy, il y a \n{3065118} graphes parfaitement contractiles à
10 sommets  parmi les  \n{3269264} graphes de  Berge de  cette taille.
Seules  cinq classes,  parmi  lesquelles les  ``quasi-parité'' et  les
``quasi-parité stricte'', sont plus nombreuses.  La classe des graphes
parfaitement  contractiles  réalise donc  un  assez  bon  score à  cet
étrange concours.

\begin{table}
  \center
  \begin{tabular}{|l|rrrrrr|}
    \hline 
    \parbox{10.5em}{Classe  \hfill $n=$}    &  $5$ & $6$ & $7$ & $8$ & $9$ & $10$ \\ \hline 
    Graphes quelconques   & \n{34} & \n{156} & \n{1044} & \n{12346} & \n{274668} & \n{12005168} \\
    Berge   &  \n{33} & \n{148} & \n{906} & \n{8887} & \n{136756} & \n{3269264} \\
    Quasi-parité    & \n{33} & \n{148} & \n{906} & \n{8886} & \n{136735} & \n{3268600} \\
    Quasi-parité stricte   & \n{33} & \n{147} & \n{896} & \n{8684} & \n{131363} & \n{3066504} \\
    Parfaitement contractile   & \n{33} & \n{147} & \n{896} & \n{8683} & \n{131333} & \n{3065118} \\ 
    Parfaitement ordonnable  & \n{33} & \n{147} & \n{896} & \n{8682} & \n{131299} & \n{3062755} \\
    Faiblement triangulé   & \n{33} & \n{146} & \n{886} & \n{8483} & \n{126029} & \n{2866876} \\ 
    Meyniel   & \n{32} & \n{130} & \n{622} & \n{3839} & \n{28614} & \n{258660} \\   \hline 
 \end{tabular}
  \caption{Dénombrements  pour  certaines  classes de  graphes.   Dans
  toutes  les classes  dont il  est question,  il y  a 1  graphe  à un
  sommet, 2  graphes à deux sommets,  4 graphes à trois  sommets et 11
  graphes à quatre sommets.\label{pair.tab.hougardy}}
\end{table}

Selon Hougardy, il y a \n{8887} graphes de Berge avec 8 sommets, parmi
lesquels \n{8886} sont de quasi-parité,  et 8 est le plus petit nombre
de sommets pour  lequel il y ait une différence.  On  en déduit que le
line-graphe      de     $K_{3,3}      \setminus      e$     représenté
figure~\ref{pair.fig.inter} est  le plus petit graphe de  Berge qui ne
soit pas de quasi-parité.

Par  un   raisonnement  similaire  et  en   utilisant  uniquement  les
dénombrements  de Hougardy,  on montre  que  le plus  petit graphe  de
quasi-parité et  non parfaitement contractile  est le $\overline{C_6}$
représenté    figure~\ref{pair.fig.inter},    et    aussi    que    la
conjecture~\ref{pair.conj.EverettReed} est vraie  pour les graphes à 6
sommets~!     Le    prisme    impair    à   8    sommets    représenté
figure~\ref{pair.fig.inter} est  le plus petit  graphe de quasi-parité
stricte qui  n'est pas  parfaitement contractile. Le  prisme pair  à 9
sommets  représenté  figure~\ref{pair.fig.inter}  est  le  plus  petit
graphe parfaitement contractile qui  n'est pas d'Artémis.  Le graphe à
8     sommets    de     la    patate     ``F.      T.''     représenté
figure~\ref{pair.fig.inter} est le  plus petit graphe non parfaitement
ordonnable  et  parfaitement contractile.   Le  plus  petit graphe  de
quasi-parité stricte  qui n'est pas faiblement triangulé  est le $C_6$
représenté figure~\ref{pair.fig.inter},  et on voit que  ce graphe est
le  seul  graphe  à 6  sommets  avec  ces  propriétés.  Donc  il  est
impossible de remplacer le  graphe parfaitement ordonnable à 7 sommets
de la  figure~\ref{pair.fig.inter} par un graphe plus  petit.  Le plus
petit graphe  de Berge  à 5 sommets  qui n'est  pas de Meyniel  est la
maison représentée figure~\ref{pair.fig.inter}.

Finalement,    il    ne   reste    que    trois    graphes   sur    la
figure~\ref{pair.fig.inter}  dont   nous  ne  pouvons   affirmer  avec
certitude qu'ils  soient les plus  petits possibles à leur  place (les
deux graphes avec 18 sommets, et le graphe d'Artémis qui n'appartient
à aucune sous-classe d'Artémis).  Nous croyons qu'ils le sont, mais le
prouver nécessiterait des vérifications selon nous assez pénibles.

\addtocontents{toc}{\mbox{}\newpage\mbox{}\vspace{1cm}\\}
\section{Généralisations de la notion de paire d'amis}
\label{pair.s.generalisation}

Pourquoi consacrer  une section  à la généralisation  de la  notion de
paire  d'amis  tant  il  vrai  que  toute  notion  est  suceptible  de
généralisation   et   qu'aucune   section   n'est   consacrée   à   la
géneralisation  de la  notion  de  2-joint, ni  à  celle de  partition
antisymétrique  par exemple  ?  Les  raisons à  cela  sont multiples.
Tout d'abord  cette question est  souvent posée, par exemple  par Paul
Seymour lors de la session  de problèmes ouverts sur les paires d'amis
à  la conférence de  Palo Alto  sur les  graphes parfaits  en novembre
2002, ou par Yann Kieffer, à  la suite d'un exposé de Frédéric Maffray
à la  conférence ROADEF à Avignon  en mars 2003.   Le lecteur souhaite
peut-être une réponse plus convaincante.

La paire d'amis est une  notion qui permet simultanément de prouver la
perfection et d'envisager des  algorithmes de coloration.  Dans le cas
très particulier des graphes  faiblement triangulés, elle fournit même
un  algorithme  de reconnaissance.   Pour  qui  rêve  d'une preuve  du
théorème fort des graphes parfaits donnant simultanément un algorithme
de coloration, la paire d'amis  est donc un outil prometteur.  Mais on
sait que l'outil fonctionne très mal pour de nombreux graphes parfaits
(certains  lines-graphes de  bipartis par  exemple) qui  n'ont  pas de
paires d'amis.  Il est donc naturel  de se demander si une notion plus
faible permettrait de garder  l'essentiel des propriétés pour tous les
graphes  parfaits.   Nous  donnons  ici  deux  pistes  possibles  (les
partitions  paires  et les  cliques  ennemies),  avec quelques  lemmes
encourageants.

\subsection{Partitions paires}
\label{pair.ss.generalisation.pp}
\index{partition~paire}
\index{pair!partition~---}
La   notion   de  partition   antisymétrique   paire,  déjà   définie
page~\pageref{graphespar.ss.SP}, n'est pas sans rapport avec la notion
de paire d'amis.   Nous proposons ici une notion  encore plus générale
définie   dans    la   preuve    du   théorème   fort    des   graphes
parfaits~\cite{chudvovsky.r.s.t:spgt}~: soit $G$ un graphe et $(A, B)$
une partition des sommets de $V$  en deux ensembles non vides.  On dit
que la partition est \emph{paire} si et seulement si~:

\vspace{1ex}

\begin{itemize}
\item
  Tous les chemins  ayant leurs extrémités dans $B$  et leur intérieur
  dans $A$ sont de longueur paire.
\item
  Tous  les  antichemins ayant  leurs  extrémités  dans  $A$ et  leur
  intérieur dans $B$ sont de longueur paire.
\end{itemize}

\vspace{2ex}

Soit $G$ un graphe quelconque avec au moins 3 sommets et $x$, $y$ deux
sommets non adjacents  de $G$.  Alors $\{x, y\}$  est une paire d'amis
de  $G$ si  et seulement  si $(V\setminus  \{x,y\}, \{x,y\})$  est une
partition paire  de $G$. En ce  sens, la partition paire  est bien une
généralisation de la notion de paire d'amis.

Soit  $G$ est  un graphe  de Berge.   On vérifie  immédiatement qu'une
partition $(A,  B)$ de ses sommets  est paire si  et seulement l'ajout
d'un  sommet $v$ complet  à $B$  et manquant  entièrement $A$  donne à
nouveau un graphe  de Berge.  En fait, par le  lemme de réplication de
Lov\'asz (lemme~\ref{graphespar.l.rep}),  on peut dupliquer  ce sommet
$v$ et même  ajouter à $G$ une clique $K$ en  reliant tous les sommets
de $K$  à tous les  sommets de $B$,  pour encore obtenir un  graphe de
Berge.  Le  lemme suivant utilise cette  construction pour généraliser
la conclusion~(\ref{pair.o.fonlupt0})  du théorème de  Fonlupt et Uhry
(théorème~\ref{pair.t.fonlupt}). Rappelons que ce théorème affirme que
si $\{x,y\}$  est une  paire d'amis d'un  graphe, alors il  existe une
coloration optimale  de $G$ qui  donne la même  couleur à $x$  et $y$,
c'est-à-dire une coloration  de $G$ dont la trace  sur $\{x,y\}$ donne
une coloration optimale de $G[\{x,y\}]$.

\begin{lemme}
  Soit  $G$  un graphe  parfait  et  $(A,B)$  une partition  paire  de
  $G$. Alors  il existe une coloration  optimale de $G$  dont la trace
  sur $B$ donne une coloration optimale de $G[B]$.
\end{lemme}

\begin{preuve}
  On construit  le graphe $G'$ en  ajoutant à $G$ une  clique $K$ avec
  ${\chi(G) -  \chi(G[B])}$ sommets, disjointe de $G$.   On relie tous
  les sommets  de la  clique $K$  aux sommets de  $B$.  On  obtient un
  graphe  de Berge  vérifiant  $\omega(G') =  \omega(G)$.  D'après  le
  théorème fort des graphes parfaits, on sait que $G'$ est parfait, ce
  qui  implique  $\chi(G)  =  \chi(G')$.  Soit  alors  une  coloration
  optimale de $G'$. Sa trace  sur $G$ donne une coloration optimale de
  $G$.  On utilise $\chi(G)  - \chi(G[B])$ couleurs pour colorier $K$,
  et donc $\chi(G[B])$ couleurs pour  colorier $G[B]$.  On a donc bien
  trouvé une  coloration optimale de $G$  dont la trace  sur $B$ donne
  une coloration optimale de $G[B]$.
\end{preuve}

La preuve du lemme précédent peut  très bien remplacer la preuve de la
conclusion~(\ref{pair.o.fonlupt0})  du  théorème  de Fonlupt  et  Uhry
(théorème~\ref{pair.t.fonlupt}).   Évidemment,  cette nouvelle  preuve
est nettement moins intéressante  puisqu'elle utilise le théorème fort
des graphes parfaits.  Son seul mérite est donc de prouver un résultat
plus général.   Nous avons tenté  sans succès de garder  la conclusion
générale tout  en nous  débarrassant du recours  au théorème  fort des
graphes parfaits.

Les  partitions antisymétriques  paires sont  des cas  particuliers de
partitions paires.  Donc, par le théorème de décomposition des graphes
de Berge (théorème~\ref{graphespar.t.structure}), on sait que certains
des graphes de Berge les plus difficiles à colorier, à savoir ceux qui
ne sont pas basiques et qui n'ont ni 2-joint ni complément de 2-joint,
possèdent  de telles  partitions.  Une  preuve algorithmique  du lemme
ci-dessus  pourrait  donc peut-être  ouvrir  une  piste pour  colorier
certains  graphes  parfaits,   reposant  peut-être  sur  des  échanges
bichromatiques,   à   la   manière   de   la   preuve   classique   du
théorème~\ref{pair.t.fonlupt}.   Là   encore,  nous  n'avons   pas  de
résultat.

\subsection{Cliques ennemies}
\label{pair.ss.cliquesennemies}

On  sait que  les line-graphes  de  bipartis sont  souvent sans  paire
d'amis.     Ce     phénomène    est    bien     compris    grâce    au
théorème~\ref{pair.t.hougardy}, dû  à Hougardy.  Michel  Burlet a donc
cherché une nouvelle notion ayant vocation à remplacer la paire d'amis
dans  le cas  des line-graphes  de bipartis.   On sait  que,  sauf cas
trivial d'un  graphe de moins  de trois sommets, les  graphes bipartis
ont toujours des paires d'amis.  Pourquoi alors ne pas regarder ce que
devient cette paire après  passage au line-graphe~?  Michel Burlet est
ainsi  parvenu à  la  notion de  cliques  ennemies~: on  dit que  deux
cliques sont  \emph{ennemies} si  tous les chemins  allant de  l'une à
l'autre  et minimaux  pour l'inclusion  avec cette  propriété  sont de
longueur  impaire.  Plus  précisément,  si $K_1$  et  $K_2$ sont  deux
cliques d'un graphe $G$, on appelle \emph{chemin sortant reliant $K_1$
et
$K_2$}\index{chemin!---~sortant~entre~2~cliques}\index{sortant!chemin~---~entre~2~cliques}
tout chemin $P$ ayant une extrémité dans $K_1$, l'autre extrémité dans
$K_2$  et  tel  que  $V(P^*)  \subseteq V(G)  \setminus  (V(K_1)  \cup
V(K_2))$.   En fait,  deux cliques  sont \emph{ennemies}  si  tous les
chemins sortants allant  de l'une à l'autre sont  de longueur impaire.
On notera  que deux cliques  ennemies sont par  définition disjointes,
car  un éventuel  sommet commun  constituerait  à lui  seul un  chemin
sortant de longueur~0.  On définit  de même la notion de \emph{cliques
amies}~: deux  cliques sont \emph{amies} si tous  les chemins sortants
allant de l'une à l'autre sont de longueur paire.

\subsubsection*{Justification de la notion de cliques ennemies} 
\label{pair.sss.justifennemie}

Nous allons voir que les cliques ennemies jouent dans les line-graphes
de  bipartis un  rôle similaire  à celui  des paires  d'amis  dans les
graphes  bipartis.  Pour  bien faire,  il nous  faut  des informations
précises  sur la  structure  des line-graphes  de  bipartis. Le  lemme
suivant   est   dans   l'esprit   des   travaux   plus   généraux   de
Beineke~\cite{beineke:linegraphs}.  Il appartient au folklore, n'étant
à  notre connaissance  publié nulle  part. Notons  qu'il ne  donne pas
d'algorithme rapide de reconnaissance des line-graphes de graphes sans
triangle.  En  effet, pour  reconnaître de tels  graphes, il  est plus
rapide de  calculer la racine  en temps linéraire par  l'algorithme de
Lehot  ou   Roussopoulos  (algorithme~\ref{graphespar.a.lehot}),  pour
ensuite  y  tester  l'existence   d'un  triangle,  que  de  rechercher
directement des griffes et des diamants.  Notons aussi que notre lemme
est  une étape  naturelle  pour  démontrer le  théorème  de Harary  et
Holzmann  (théorème~\ref{graphespar.t.lgb})  sur  les line-graphes  de
bipartis. Il est donc probable  que Harary et Holzmann ont démontré ce
lemme,  mais malheureusement,  leur article  est très  difficile  à se
procurer,  peut-être parce  que publié  dans le  premier  numéro d'une
revue  chilienne en 1974,  année qui  suit le  coup d'État  du général
Pinochet.  Nous prouvons finalement ce lemme pour ne rien laisser dans
l'ombre.

\index{line-graphe~de~graphe~sans~triangle!caractérisation}
\index{triangle!line-graphe d'un graphe sans ---}
\begin{lemme}
  \label{pair.l.LB}
  Soit $G$ un graphe.  Il existe  un graphe $R$ sans triangle tel que
  $G=L(R)$ si et seulement si $G$ est sans griffe et sans diamant.
\end{lemme}

\begin{preuve}
  Il est clair que le line-graphe d'un graphe sans triangle ne possède
  pas de griffe ni de  diamant. Réciproquement, on considère un graphe
  $G$ sans griffe et sans diamant. On va construire un graphe $R$ sans
  triangle tel que $G$ soit isomorphe à $L(R)$.
  
  \begin{claim}
    \label{pair.c.LBNv}
    Soit $v\in V(G)$. Alors ou  bien $G[N(v)]$ est une clique, ou bien
    $G[N(v)]$   possède  deux  composantes   connexes,  et   ces  deux
    composantes connexes sont des cliques.
  \end{claim}
  
  \begin{preuveclaim}
    Supposons  que $G[N(v)]$ n'est  pas une  clique.  Alors  il existe
    $v_1 \in N(v)$ et $v_2 \in  N(v)$ tels que $v_1$ manque $v_2$.  Si
    $u$ est un sommet de $N(v) \setminus \{v_1, v_2\}$, alors $u$ voit
    exactement  l'un des  sommets $v_1$  et $v_2$.   En effet,  si $u$
    manque  $v_1$ et  $v_2$, alors  $\{v,  u, v_1,  v_2\}$ induit  une
    griffe et si  $u$ voit $v_1$ et $v_2$, alors  $\{v, u, v_1, v_2\}$
    induit un  diamant.  Donc les  voisins de $v$ se  partitionnent en
    deux ensembles  $K_1$ et  $K_2$ avec $K_i  = \{v_i\}  \cup (N(v_i)
    \cap N(v))$,  $i= 1, 2$.  On  voit que $K_1$ induit  une clique de
    $G$ car si  $K_1$ contient deux sommets non  adjacents $u$ et $w$,
    alors l'ensemble $\{v, u, w,  v_{2}\}$ induit une griffe.  De même
    $K_2$ induit une clique. De plus,  il n'y a aucune arête entre les
    sommets de $K_1$  et les sommets de $K_2$. En  effet, si une telle
    arête existe  entre un sommet $u_1\in  K_1$ et un  sommet $u_2 \in
    K_2$, alors on ne peut pas avoir simultanément $u_1 = v_1$ et $u_2
    = v_2$.   On suppose sans perte  de généralité $u_2  \neq v_2$. On
    constate  alors que  $u_1$ voit  $v_2$ car  sinon $\{u_1,  v, v_2,
    u_2\}$ induit  un diamant.   Comme $u_1$ voit  $v_2$, on  sait que
    $u_1 \neq  v_1$. Donc $v_1$ voit  $u_2$ car sinon  $\{u_2, v, v_1,
    u_1\}$ induit  un diamant.  Mais  alors $\{v_1, u_1 ,  v_2, u_2\}$
    induit  un  diamant~:  une  contradiction.  Donc  $K_1$  et  $K_2$
    induisent les deux composantes connexes de $G[N(v)]$.
  \end{preuveclaim}

  De~(\ref{pair.c.LBNv}), on déduit facilement deux petits corollaires
  pratiques~:

  \begin{claim}
    \label{pair.c.LBf}
    Deux cliques  distinctes et maximales de $G$  s'intersectent en au
    plus un sommet.
  \end{claim}

  \begin{claim}
    \label{pair.c.LBg}
    Soit  $v$  un  sommet de  $G$.  Alors  ou  bien $v$  appartient  à
    exactement une  clique maximale  de $G$ ou  bien $v$  appartient à
    exactement deux cliques maximales de $G$.
  \end{claim}

  À partir de $G$, on construit le graphe $R$. Les sommets de $R$ sont
  de  deux types~:  les sommets  de type  \emph{clique}, qui  sont les
  cliques maximales de  $G$ et les sommets de  type \emph{pendant} qui
  sont les sommets de $G$ n'appartenant qu'à une seule clique maximale
  de $G$. On ajoute une arête  entre deux sommets de type clique quand
  les deux cliques qui  leur correspondent dans $G$ s'intersectent. On
  ajoute  une arête  entre  un sommet  de  type pendant  et le  sommet
  représentant l'unique clique maximale  à laquelle il appartient. Les
  sommets de type pendant sont par construction de degré~1 et ont pour
  seul voisin un sommet de type de clique.
  
  \begin{claim}
    Le  graphe  $R$ est  sans  triangle. 
  \end{claim}
  
  \begin{preuveclaim}
    En  effet, si $R$  possède un  triangle $\{a,  b, c\}$,  alors ses
    sommets sont de  type clique et il leur  correspond dans $G$ trois
    cliques maximales distinctes $K_a, K_b, K_c$ s'intersectant deux à
    deux. S'il existe  $u \in K_a \cap K_b  \cap K_c$, alors $G[N(u)]$
    contient    au   moins   trois    cliques   maximales,    ce   qui
    contredit~(\ref{pair.c.LBNv}).   Donc $K_a  \cap K_b  \cap  K_c =
    \emptyset$. Soit alors $u \in K_a \cap K_b$, $v\in K_a\cap K_c$ et
    $w\in  K_b  \cap  K_c$.   D'après~(\ref{pair.c.LBNv}),  l'ensemble
    $N(u)$  se partitionne en  deux cliques  qui sont  les composantes
    connexes de $G[N(u)]$. Ces deux composantes connexes ne peuvent
    être que  $K_b \setminus u$ et  $K_c \setminus u$  et l'arête $vw$
    montre   que    $G[N(u)]$   est   connexe~:    une   contradiction
    avec~(\ref{pair.c.LBNv}).
  \end{preuveclaim}

  On définit  une fonction $f$  de $V(G)$ vers  $E(R)$. Si $v$  est un
  sommet appartenant  à deux  cliques maximales de  $G$, alors  on lui
  associe l'unique  arête de $R$ qui  relie ces deux  cliques.  Si $v$
  n'appartient qu'à  une seule  clique maximale de  $G$, alors  on lui
  associe      l'unique     arête      de     $R$      incidente     à
  $v$.  D'après~(\ref{pair.c.LBg}),  $f$ est  bien  définie pour  tout
  sommet de $G$.

  \begin{claim}
    La fonction $f$ est bijective.  De plus, pour toute  paire de sommets
    $\{u, v\}$ de $G$, $f(u)$ et $f(v)$ sont incidentes dans $R$ si et
    seulement si $u$ et $v$ sont adjacents dans $G$.
  \end{claim}
  
  \begin{preuveclaim}
    La fonction $f$ est injective.  En effet, si $u$ appartient à deux
    cliques  maximales  de  $G$ alors  d'après~(\ref{pair.c.LBf}),  le
    sommet $u$ est  l'unique sommet commun à ces  deux cliques.  Donc
    aucun  autre  sommet de  $G$  n'aura  pour  image $f(u)$.  Si  $u$
    n'appartient qu'à une seule clique, alors $f(u) = u$ et là encore,
    aucun autre sommet de $G$ n'aura pour image $f(u)$.

    La fonction $f$  est surjective. En effet, soit  $ab$ une arête de
    $R$. Si $a$ et $b$ sont de type clique, alors les deux cliques $a$
    et  $b$ de  $G$ ont  un  sommet commun  dont l'image  par $f$  est
    justement $ab$. Si $a$ est de  type clique et $b$ de type pendant,
    alors $ab$ est l'image de $b$ par $f$.

    Supposons $u$ et $v$ adjacents  dans $G$. Si $u$ appartient à deux
    cliques   maximales  de   $G$   notées  $K_u$   et  $K'_u$   alors
    d'après~(\ref{pair.c.LBNv})   on  sait   que   $G[N(u)]$  a   pour
    composantes connexes $K_u \setminus  u$ et $K'_u \setminus u$.  On
    suppose  sans  perte  de   généralité  $v  \in  V(K_u)$.   Si  $v$
    n'appartient qu'à  une seule clique  maximale de $G$,  alors cette
    clique est $K_u$. L'arête de  $R$ égale à $f(u)$ est alors l'arête
    reliant $K_u$ et  $K'_u$ tandis que l'arête de  $R$ égale à $f(v)$
    est  l'arête reliant  $K_u$  à  $v$.  Ces  deux  arêtes sont  bien
    incidentes.  Si  $v$ appartient à  deux cliques maximales  de $G$,
    alors d'après~(\ref{pair.c.LBNv}),  l'une de ces  deux cliques est
    $K_u$ et on note l'autre $K_v$.  L'arête de $R$ égale à $f(u)$ est
    alors l'arête  reliant $K_u$ et  $K'_u$ tandis que l'arête  de $R$
    égale à $f(v)$ est l'arête reliant $K_u$ à $K_v$.  Ces deux arêtes
    sont  bien  incidentes.  Si  $u$  appartient  à  une seule  clique
    maximale  de $G$  notée  $K_u$,  alors on  peut  supposer que  $v$
    n'appartient  aussi qu'à une  seule clique  maximale de  $G$ notée
    $K_v$ (sinon, en échangeant les rôles  de $u$ et $v$, on se ramène
    à un cas déjà  traité).  En fait, d'après~(\ref{pair.c.LBf}), $K_u
    = K_v$  et $f(u)$ est l'arête  de $R$ reliant $u$  à $K_u$, tandis
    que  $f(v)$ est  l'arête de  $R$ reliant  $v$ à  $K_u$~:  ces deux
    arêtes sont bien incidentes.

    Réciproquement, si  les arêtes  $f(u)$ et $f(v)$  sont incidentes,
    alors leur sommet commun est  nécessairement un sommet de degré au
    moins  deux, donc  un sommet  de type  clique correspondant  à une
    clique de $G$  contenant $u$ et $v$, qui  sont donc bien adjacents
    dans $G$.
  \end{preuveclaim}
  
  Les  propriétés  de  la  fonction  $f$ montrent  bien  que  $G$  est
  isomorphe au line-graphe de $R$.
\end{preuve}

Le lemme suivant est à notre connaissance original, mais il ne devrait
pas surprendre un lecteur habitué  au sujet. Il montre que les notions
de cliques amies et ennemies sont pertinentes dans les line-graphes de
bipartis~:

\begin{lemme}
  \label{pair.l.caracennemie}
  Soit $G$  un graphe sans  griffe et sans  diamant. Alors $G$  est un
  line-graphe de biparti si et seulement si ses cliques maximales pour
  l'inclusion peuvent être partitionnées  en deux ensembles $A$ et $B$
  tels que pour toutes cliques maximales $K_1$ et $K_2$ de $G$~:
  \begin{itemize}
  \item
    Si $K_1\in A$ et $K_2\in A$, alors $K_1$ et $K_2$ sont ennemies.
  \item
    Si $K_1\in B$ et $K_2\in B$, alors $K_1$ et $K_2$ sont ennemies.
  \item
    Si $K_1\in A$ et $K_2\in B$, alors $K_1$ et $K_2$ sont amies.
  \end{itemize} 
\end{lemme}

\begin{preuve}
  On  adopte  les notations  de  la  preuve du  lemme~\ref{pair.l.LB}.
  D'après ce lemme, $G$ est  isomorphe au line-graphe d'un graphe sans
  triangle  $R$. On  peut donc  supposer qu'il  existe un  graphe sans
  triangle  $R$ tel  que $G=  L(R)$.  Il  apparaît dans  la  preuve du
  lemme~\ref{pair.l.LB}  que les  cliques  maximales de  $G$ sont  des
  sommets de $R$ appelés sommets de type clique.

  Si  $R$ est  biparti, alors  les sommets  de type  clique de  $R$ se
  partitionnent en  deux stables $A$  et $B$. Cette partition  est aussi
  une partition des cliques maximales de $G$. Soit $K_1 \in A$ et $K_2
  \in A$  deux cliques  maximales de $G$.   Notons que $K_1$  et $K_2$
  sont en même  temps des sommets de $R$ non  adjacents.  On sait donc
  que $K_1$ et $K_2$ sont  des cliques disjointes de $G$.  S'il existe
  un chemin  de $G$,  de longueur paire,  sortant et reliant  $K_1$ et
  $K_2$, alors  d'après le lemme~\ref{pair.l.chainechemin}, les
  sommets intérieurs  de ce  chemin (qui est  de longueur  au moins~2)
  sont  les arêtes  de l'intérieur  d'une  chaîne de  $R$ de  longueur
  impaire reliant le sommet $K_1$  au sommet $K_2$.  Ceci contredit la
  bipartition de $R$.  Donc tous les chemins sortants reliant $K_1$ et
  $K_2$ sont  de longueur impaire~:  $K_1$ et $K_2$ sont  ennemies. On
  démontre de  même que si  $K_1\in B$ et  $K_2\in B$, alors  $K_1$ et
  $K_2$ sont  ennemies. Soient maintenant $K_1  \in A$ et  $K_2 \in B$
  deux  cliques maximales de  $G$. S'il  existe un  chemin de  $G$, de
  longueur impaire,  sortant et reliant  $K_1$ et $K_2$, alors  par un
  raisonnement  similaire, on  trouve une  chaîne de  $R$  de longueur
  paire  reliant  le sommet  $K_1$  au  sommet  $K_2$. Là  encore,  on
  contredit la bipartition de $R$~: $K_1$ et $K_2$ sont amies.

  Si $R$  n'est pas biparti, alors  $R$ possède un trou  impair $H$ de
  longueur au moins  5 (car $R$ est sans  triangle). Notons $v_1, v_2,
  \dots  v_{2k+1}$ les sommets  d'un tel  trou dans  l'ordre où  ils y
  apparaissent. Chaque sommet  de $H$ est de degré  au moins~2, et est
  donc de type clique. Donc chaque sommet $v_i$ est en fait une clique
  maximale de $G$.  Si on arrive à partitionner l'ensemble des cliques
  maximales en deux ensembles $A$  et $B$ comme indiqué dans le lemme,
  il advient que deux cliques  consécutives $v_i$ et $v_{i+1}$ ne sont
  pas disjointes,  ne peuvent donc  être ennemies, et ne  peuvent être
  toutes les deux  dans $A$ ou toutes les deux  dans $B$.  Les cliques
  $v_i$ doivent donc être alternativement  dans $A$ et $B$, ce qui est
  impossible étant donné leur nombre impair.
\end{preuve}

\subsubsection*{Conjectures et résultats partiels}

Si l'on  admet que  les cliques ennemies  sont l'extension  des paires
d'amis aux  line-graphes de bipartis,  alors il est assez  naturel de
proposer la conjecture suivante~:

\begin{conjecture}[Burlet, 2001]
  \label{pair.conj.bubu}
  Soit $G$ un graphe de Berge non réduit à un sommet. Alors ou bien~:
  \begin{outcomes}
  \item
    $G$ ou $\overline{G}$ possède une paire d'amis.
  \item
    $G$ ou $\overline{G}$ possède une paire de cliques ennemies.
  \end{outcomes}
\end{conjecture}

Nous n'avons pas d'idée bien précise pour prouver cette conjecture. Ce
qui précède montre qu'elle est  vraie pour les graphes basiques.  Elle
est vraie également  pour tous les graphes de  quasi-parité.  Pour les
graphes  de  Berge plus  généraux,  le  théorème  de décomposition  de
Chudnovsky,         Seymour,        Robertson         et        Thomas
(théorème~\ref{graphespar.t.structure})  ne donne  rien  de facilement
exploitable. 

Il  est clair  qu'entre deux  cliques  maximales d'un  trou impair,  il
existe toujours  un chemin pair  sortant. Entre deux  cliques maximale
d'un  antitrou  impair,  il   existe  toujous  un  chemin  sortant  de
longueur~2.  Si  l'on admet que les seuls  graphes imparfaits minimaux
sont les trous impairs et  les antitrous impairs, on constate donc que
la conjecture suivante est vraie~:

\index{minimalement~imparfait~(graphe~---)!cliques~ennemies}
\begin{conjecture}[Burlet, 2001]
  \label{pair.conj.cliqueenn}
  Dans  un graphe  imparfait minimal  $G$, il  n'y a  pas de  paire de
  cliques ennemies $\{K_1, K_2\}$, telle  que $K_i$, $i=1, 2$, est une
  clique maximale de $G$
\end{conjecture}

On pourrait souhaiter une preuve  n'utilisant pas le théorème fort des
graphes  parfaits.  Nous présentons  ici  quelques résultats  partiels
allant dans ce sens.

On sait  que la contraction  d'une paire d'amis $\{x,y\}$  d'un graphe
$G$   est  une  opération   intéressante.   Quel   serait  l'opération
équivalente pour  les cliques ennemies~?  Les arêtes  incidentes à $x$
forment une clique $K_x$ de $L(G)$, et celles incidentes à $y$ forment
une clique  $K_y$. Le sommet  contracté $xy$ reçoit toutes  les arêtes
incidentes à $x$ ou $y$, et  donne donc dans $L(G)$ une clique obtenue
en reliant les  sommets de $K_x$ à ceux de $K_y$.   Dans un graphe $G$
avec deux cliques ennemies $K_1$ et $K_2$, on définit donc l'opération
consistant à relier  chaque sommet de $K_1$ à  chaque sommet de $K_2$.
On obtient  un nouveau graphe  noté $G_{K_1\equiv K_2}$.   Nous allons
voir  que cette opération  préserve la  perfection.  Au  préalable, on
vérifie  que l'opération  consistant  à relier  deux cliques  ennemies
$K_1$  et $K_2$  d'un graphe  $G$ ne  peut créer  aucune  autre grande
clique que celle induite par $V(K_1) \cup V(K_2)$. Plus précisément~:

\begin{lemme} \label{pair.l.grosseclique}
  Soit $K_1$ et $K_2$ deux  cliques ennemies d'un graphe $G$. Soit $K$
  une clique de $G_{K_1 \equiv K_2}$. Deux cas seulement se présentent :
  \begin{itemize}
  \item
    $K$ est une clique de $G$.
  \item
    $V(K) \subseteq V(K_1) \cup V(K_2)$.
  \end{itemize}
\end{lemme}

\begin{preuve}
  Si $K$ n'est  pas une clique de $G$, alors $K$  contient au moins un
  sommet $v_1$  de $K_1$ et un sommet  $v_2$ de $K_2$ qui  ne sont pas
  reliés dans $G$.   Si en outre $V(K)$ n'est  pas inclus dans $V(K_1)
  \cup V(K_2)$,  alors $K$  contient un sommet  $v$ qui n'est  ni dans
  $K_1$,  ni dans  $K_2$, et  qui est  relié à  $v_1$ et  $v_2$.  Mais
  alors, $\bp v_1 \tp v \tp v_2  \ep$ est un chemin sortant de $G$, de
  longueur paire, reliant $K_1$ à $K_2$, ce qui est absurde.
\end{preuve}

Voici notre résulat principal sur les cliques ennemies~:

\begin{theoreme}
  \label{pair.l.preserve}
  Soit  $G$ un  graphe parfait  avec  deux cliques  ennemies $K_1$  et
  $K_2$. Alors $G_{K_1 \equiv K_2}$ est lui aussi parfait.
\end{theoreme}

\begin{preuve}
  Soit $H'$ un sous-graphe de  $G_{K_1 \equiv K_2}$.  On considére $H$
  le sous-graphe  de $G$  qui a  les mêmes sommets  que $H'$.   Il est
  clair que  $V(K_1) \cap  V(H)$ et $V(K_2)  \cap V(H)$  induisent des
  cliques ennemies de $H$ et  que $H'=H_{(K_1 \cap H) \equiv (K_2 \cap
  H)}$.  Donc,  pour parvenir à nos  fins, il nous  suffit de vérifier
  que  $\chi(G_{K_1 \equiv  K_2}) =  \omega(G_{K_1 \equiv  K_2})$ sans
  nous  soucier des sous-graphes.   Pour cela,  on montre  qu'à partir
  d'une  coloration de $G$  avec $\omega(G)$  couleurs, on  trouve une
  coloration de $G_{K_1 \equiv K_2}$ avec $\omega(G_{K_1 \equiv K_2})$
  couleurs~:

  On commence  par colorier les sommets  qui ne sont ni  dans $K_1$ ni
  dans $K_2$.   Si $\omega(G_{K_1 \equiv K_2}) >  \omega(G)$, alors on
  sait  d'après le  lemme~\ref{pair.l.grosseclique}  que $V(K_1)  \cup
  V(K_2)$ induit l'unique clique  maximum de $G_{K_1 \equiv K_2}$.  On
  a donc le droit à $\omega(G_{K_1 \equiv K_2}) - \omega(G)$ nouvelles
  couleurs que l'on utilise  pour colorier $\omega(G_{K_1 \equiv K_2})
  - \omega(G)$ sommets quelconques de  $V(K_1) \cup V(K_2)$.  Il reste
  au  maximum $\omega(G)$ sommets  dans $V(K_1)  \cup V(K_2)$  dont on
  garde provisoirement  l'ancienne couleur dans $G$.  On peut supposer
  qu'il existe  un sommet $v_1$ de  $K_1$ et un sommet  $v_2$ de $K_2$
  ayant  la même  couleur (disons  rouge) car  sinon on  a  obtenu une
  $\omega(G_{K_1 \equiv  K_2})$-coloration de $G_{K_1  \equiv K_2}$ et
  la conclusion  du lemme est satisfaite.   Il y a  alors une ancienne
  couleur  (disons bleue)  qui n'est  utilisée ni  dans $K_1$  ni dans
  $K_2$.   On considère  l'ensemble $C$  des sommets  de $G$  qui sont
  rouges ou  bleus.  L'ensemble $C$  induit un sous-graphe  biparti de
  $G$  et on  appelle $C_1$  la composante  connexe de  $v_1$  dans ce
  sous-graphe.  Si $v_2 \in C_1$,  alors un plus court chemin de $C_1$
  reliant  $v_1$ à  $v_2$  est un  chemin  de $G$,  de longueur  paire
  reliant $K_1$  et $K_2$.  Ce chemin  est sortant car il  n'y a aucun
  sommet de  couleur bleue dans $V(K_1) \cup  V(K_2)$.  Ceci contredit
  la  définition des cliques  ennemies, et  on a  montré que  $v_1$ et
  $v_2$ ne  sont pas dans la  même composante connexe  $C_1$.  Donc il
  est possible d'échanger  les couleurs rouge et bleue  dans $C_1$ pour
  affecter à $v_1$ la couleur  bleue, sans changer la couleur du sommet
  $v_2$.  On  recommence cette opération tant qu'il  reste des sommets
  de  même couleur  dans $K_1  \cup K_2$  pour finalement  obtenir une
  $\omega(G_{K_1 \equiv K_2})$-coloration de $G_{K_1 \equiv K_2}$.
\end{preuve}

Si $K_1$ et $K_2$ sont  deux sous-cliques disjointes d'une clique $K$,
alors  elles  sont  ennemies.   Dans  ce  cas,  on  dit  qu'elle  sont
\emph{trivialement} ennemies.  Le lemme suivant  va dans le sens de le
conjecture~\ref{pair.conj.cliqueenn}~:

\index{minimalement~imparfait~(graphe~---)!cliques~ennemies|(}
\begin{lemme} 
  \label{pair.l.omega}
  Soit $G$  un graphe imparfait  minimal.  Soient $K_1$ et  $K_2$ deux
  cliques ennemies  mais non  trivialement ennemies dans  $G$.  Alors
  $|K_1|+|K_2| \neq \omega(G)$.
\end{lemme}

\begin{preuve}
  Supposons  $|K_1|  + |K_2|  =  \omega(G)$.   Remarquons d'abord  que
  d'après  le  lemme~\ref{pair.l.grosseclique}, $\omega(G_{K_1  \equiv
  K_2}) =  \omega(G)$.  Il est  clair d'autre part  que $\alpha(G_{K_1
  \equiv    K_2})     \leq    \alpha(G)$.     Enfin,     d'après    le
  théorème~\ref{graphespar.t.lovasz},  on  a  $\alpha(G)  \omega(G)  <
  |V(G)|$.

  Par définition,  tous les sous-graphes  de $G$ sont  parfaits. Donc,
  d'après le théorème~\ref{pair.l.preserve},  tous les sous-graphes de
  $G_{K_1  \equiv  K_2}$  sont  parfaits.   Pourtant,  comme  $G$  est
  partitionnable, en raison de  la $\omega$-clique $K_1 \cup K_2$ (qui
  n'est pas une clique de $G$ car $K_1$ et $K_2$ sont non trivialement
  ennemies), $G_{K_1  \equiv K_2}$  n'est pas partitionnable.   Par le
  théorème~\ref{graphespar.t.minimparfait},  on  en  déduit qu'il  est
  parfait.  Mais on a $\alpha(G_{K_1 \equiv K_2}) \omega(G_{K_1 \equiv
  K_2}) \leq  \alpha(G) \omega(G) < |V(G)| =  |V(G_{K_1 \equiv K_2})|$
  ce qui contredit le théorème~\ref{graphespar.t.lovasz}.
\end{preuve}

D'après  le lemme  précédent,  si  $K_1$ et  $K_2$  sont deux  cliques
ennemies  d'un  graphe  minimal  imparfait  $G$,  alors  deux  cas  se
présentent :

\begin{enumerate}
\item $|K_1|+|K_2|<\omega(G)$  

  Dans ce  cas, on peut  remarquer que les  arêtes que l'on  ajoute en
  reliant $K_1$  à $K_2$ ne  créent aucune $\omega$-clique  d'après le
  lemme~\ref{pair.l.grosseclique}.  De plus,  ces arêtes ne détruisent
  aucun $\alpha$-stable.  Ce fait mérite une preuve~:

\begin{preuve}
  Supposons   en   vue    d'une   contradiction   qu'on   détruit   un
  $\alpha$-stable de $G$  en reliant $K_1$ et $K_2$.   Il existe alors
  deux sommets $v_1\in  K_1$ et $v_2\in K_2$ qui  sont dans un certain
  $\alpha$-stable $S$.   D'après la conclusion~(\ref{graphespar.o.sc})
  du  théorème~\ref{graphespar.t.partitionnable}, il  existe  donc une
  $\omega$-clique $K$  disjointe de $S$. Soit $v\in  V(K)$. D'après la
  définition des  graphes partitionnables,  $G \setminus v$  peut être
  partitionné en $\omega$ stables de taille $\alpha$. Au moins l'un de
  ces stables (notons le $S'$) est disjoint de $K$, car $K\setminus v$
  contient      $\omega       -      1$~sommets.       D'après      la
  conclusion~(\ref{graphespar.o.cs})                                 du
  théorème~\ref{graphespar.t.partitionnable}, on sait que $S'=S$. Donc
  on  a  trouvé dans  $G$  un sommet  $v$  tel  qu'on puisse  colorier
  optimalement $G  \setminus v$  en donnant à  $v_1$ et $v_2$  la même
  couleur, disons rouge.  Comme $|K_1|+|K_2|<\omega(G)$, il existe une
  couleur  (disons  bleu)  qui  n'est  pas  utilisée  dans  $K_1  \cup
  K_2$. Par un échange bichromatique identique à celui de la preuve du
  théorème~\ref{pair.l.preserve}, on peut trouver une coloration de $G
  \setminus v$ qui donne la couleur  rouge à $v_1$ et la couleur bleue
  à $v_2$ (si un tel  échange échoue, rappelons qu'on trouve un chemin
  sortant  de  longueur  paire  reliant  $K_1$ et  $K_2$  ce  qui  est
  contradictoire).    Finalement,  on   a   trouvé  deux   colorations
  différentes    de   $G\setminus    v$,   ce    qui    contredit   la
  conclusion~(\ref{graphespar.o.unicolor})                           du
  théorème~\ref{graphespar.t.partitionnable}.
\end{preuve}

  Donc, quand on relie $K_1$ à $K_2$, on garde un graphe partitionnable,
  et on ne voit pas apparaître de contradiction.

\item $|K_1|+|K_2|>\omega(G)$ 

  Dans ce cas, en reliant $K_1$  à $K_2$ on obtient un graphe avec une
  clique  maximum  unique  :  $K_1\cup  K_2$. Ce  graphe  n'étant  pas
  partitionnable,  et tous  ses  sous-graphes étant  parfaits, il  est
  parfait. Dans ce cas on parvient donc à un graphe assez extravagant,
  sans qu'aucune contradiction ne semble vouloir apparaître~\dots
\end{enumerate}
\index{minimalement~imparfait~(graphe~---)!cliques~ennemies|)}

Le lemme suivant et son corollaire montrent que si $K_1$ et $K_2$ sont
ennemies dans un graphe partitionnable, alors une sous-clique de $K_1$
a  peu de  chance d'être  ennemie de  $K_2$. Nous  n'avons  pas trouvé
d'application à ce petit fait, peut-être bon à connaître malgré tout~:

\begin{lemme}
  Soient $K_1$  et $K_2$ deux cliques ennemies  maximales et distinctes
  d'un graphe $G$. Soit $K'_1$ une sous-clique (stricte) de $K_1$.  Si
  $K'_1$ et  $K_2$ sont  ennemies, alors  il y a  dans $G$  une étoile
  d'articulation.
\end{lemme}

\begin{preuve}
  Soit $a \in K'_1$ et $b \in K_2$ non reliés (ils existent, car $K_2$
  est maximale).   Soit $c$ un sommet quelconque  de $V(K_1) \setminus
  V(K'_1)$.  On va  montrer que $\{a\} \cup N(a)  \setminus \{c\}$ est
  un ensemble  d'articulation $G$ déconnectant $c$ et  $b$. Pour cela,
  on vérifie  que tout chemin  $P$ reliant $c$  à $b$ sans  passer par
  $K_1$ contient un voisin de $a$ différent de $c$.  En effet~:
  
  Si  $P^*$  ne contient  aucun  sommet de  $K_2$,  alors  $P$ est  de
  longueur impaire car $K_1$ et  $K_2$ sont ennemies.  Comme $K'_1$ et
  $K_2$ sont  ennemies, il y  a une corde  dans la chaîne  de longueur
  paire  $(a, c,  \dots, b)$, et cette corde relie forcément
  $a$ et un sommet de $P^*$, ce qu'on voulait montrer.

  Si  $P^*$  contient  un  sommet   de  $K_2$,  alors  ce  sommet  est
  nécessairement le voisin de $b$ dans  $P$, que l'on note $d$. On voit
  alors que $\bp c \tp P \tp  d \ep$ est de longueur impaire car $K_1$
  et $K_2$  sont ennemies.  Donc la  chaîne $(a, c, \dots,  d)$ est de
  longueur paire, et encore une fois,  on trouve une corde de $a$ vers
  un sommet de $P^*$ (cette corde peut d'ailleurs être $ad$).
\end{preuve}

\begin{corollaire}
  Soit $K_1$  et $K_2$ deux cliques ennemies  maximales, et distinctes
  d'un graphe $G$ partitionnable.  Soit $K'_1$ une sous-clique stricte
  et non vide de $K_1$. Alors $K'_1$ et $K_2$ ne sont pas ennemies.
\end{corollaire}

\begin{preuve}
  Clair d'après le théorème~\ref{graphespar.t.scs}
\end{preuve}

On peut remarquer  que le lemme ci-dessus n'est  pas très encourageant
en vu  de la  conjecture~\ref{pair.conj.bubu}. En effet,  pour prouver
cette dernière, on  pourrait penser suivre la même  démarche que celle
qui a  fonctionné pour décomposer  les graphes de Berge~:  on commence
par  supposer  que $G$  contient  $F$,  line-graphe d'une  subdivision
bipartie de $K_4$.  D'après le lemme~\ref{pair.l.caracennemie}, il est
clair que  $F$ contient des  cliques ennemies.  On peut  alors espérer
faire  ``grossir'' ces cliques  jusqu'à trouver  une paire  de cliques
ennemies de $G$.  On peut continuer ainsi avec  les prismes, puis, une
fois  que les  prismes  sont exclus,  on  est quasiment  raméné à  des
graphes  bipartisans  (aux  double-diamants  près~\dots), et  on  peut
espérer trouver une paire d'amis.  Le lemme ci-dessus montre que cette
idée risque de ne pas fonctionner, car il y a des chances qu'on puisse
pas faire  ``grossir'' les cliques  ennemies, ou alors  seulement dans
des graphes trop simples (ceux qui ont une étoile d'articulation).

\chapter{Le lemme de Roussel et Rubio}
\label{rr.chap}
\index{lemme~de~Roussel~et~Rubio}
\index{Roussel~et~Rubio|see{lemme~de~Roussel~et~Rubio}}
Dans       un        article       consacré       aux       partitions
antisymétriques~\cite{roussel.rubio:01},  Florian Roussel  et Philippe
Rubio ont présenté un lemme d'apparence assez anodine. Ce lemme montre
que d'un certain point de  vue, les ensembles anticonnexes des graphes
de Berge se comportent comme  de simples sommets.  Peu de temps après,
ce lemme a été redécouvert indépendamment par l'équipe de Paul Seymour
(en  collaboration  avec  Carsten  Thomassen), qui  l'a  rebaptisé  le
\emph{wonderful  lemma\footnote{En français~: le  lemme merveilleux.}}
en raison des  applications multiples qu'ils lui ont  trouvées dans la
preuve du  théorème fort  des graphes parfaits  (ils n'ont  pas publié
leur preuve  du lemme).   Nous nous étendrons  un peu  sur différentes
manières de  prouver ce lemme important.  Nous  énoncerons ensuite des
versions plus fortes du  lemme pour des classes particulières (graphes
d'Artémis, graphes  de Meyniel, graphes  faiblement triangulés).  Nous
terminerons  par des corollaires  qui nous  seront utiles  au chapitre
suivant.

\section{Différentes preuves}

Si $T$ est un ensemble de sommets d'un graphe $G$, on dit qu'un sommet
$x$ est  \emph{$T$-complet} si  $x$ voit tous  les sommets de  $T$. On
appelle  \emph{$T$-arête}\index{arête!$T$-arête}  de  $G$ toute  arête
dont les deux extrémités sont $T$-complètes.  Il est trivial que, dans
un graphe de  Berge, si un sommet $t$ voit  les extrémités d'un chemin
$P$ de  longueur impaire (au  moins 3), alors  ce sommet doit  voir au
moins un  sommet de~$P^*$ --- rappelons que  $P^*$ désigne l'intérieur
de  $P$.   En  fait,  si  on  appelle  \emph{intervalle}  de  $P$  les
sous-chemins de $P$ dont les extrémités voient $t$ et dont les sommets
intérieurs manquent $t$, on voit que $P$ se partitionne en intervalles
(au sens des arêtes), que les intervalles sont de longueur paire ou de
longueur   1,  et   donc  que   $P$  possède   un  nombre   impair  de
$\{t\}$-arêtes.   Le  fameux  lemme   de  Roussel  et  Rubio  est  une
généralisation de ce fait~: le sommet $t$ est remplacé par un ensemble
anticonnexe $T$. Nous  donnons une preuve simplifiée par  rapport à la
preuve  originale  qui séparait  déjà  les  cas  ``$G[T]$ stable''  et
``$G[T]$  non  stable''.  La  première  simplification provient  d'une
reformulation  de la  conclusion~(\ref{rr.o.rrodd}) du  lemme~:  là où
Roussel  et  Rubio  écrivaient  ``Le  chemin  $P$  possède  un  sommet
$T$-complet'' nous  écrivons l'assertion plus forte~:  ``Le chemin $P$
comporte  un  nombre impair  de  $T$-arêtes''.   Ce renforcement,  qui
n'avait  évidemment  pas  échappé  au  groupe de  Seymour,  permet  de
disposer d'une hypothèse de récurrence formellement plus forte, ce qui
débarrasse la preuve  de détails fastidieux, notamment pour  le cas où
$G[T]$ n'est pas stable.

Nous  proposons aussi  deux  nouvelles preuves  pour  le cas  ``$G[T]$
stable'',  plus simples  à notre  avis  que la  preuve originale.   La
première  est due  à Ajai  Kapoor, Kristina  Vu\v skovi\'c  et Giacomo
Zambelli~\cite{zambelli:these},  ce  dernier  me  l'ayant  aimablement
communiquée.  Il  a découvert l'argument astucieux de  comptage par la
formule du  crible dans des brouillons  non publiés de  Kapoor et Vu\v
skovi\'c concernant des recherches sur le nettoyage des trous impairs,
et  s'est   rendu  compte  que  cet  argument   pouvait  presque  sans
changement  s'appliquer au  cas  ``stable'' du  lemme  de Roussel  et
Rubio.  Nous présenterons aussi notre propre preuve de ce cas ``$G[T]$
stable'', moins  spectaculaire, mais qui a  peut-être quelques mérites
dont le lecteur jugera.

\begin{figure}[p]
  \center
  \label{rr.fig.oleap}
  \includegraphics{fig.rr.1}  \\  \rule{0cm}{.5cm}\\
  \includegraphics{fig.rr.2} \\  \rule{0cm}{.51cm}\\
  \includegraphics{fig.rr.3} \\
  \rule{0cm}{1cm}\\
  \parbox{9cm}{ Nous avons représenté trois graphes de Berge où le lemme de Roussel et Rubio s'applique~:
    \begin{itemize}
    \item Dans le premier, on a la conclusion~(\ref{rr.o.rrleap}) et
      le graphe entier  est un prisme impair avec  8 sommets.  
    \item Dans le deuxième, on a la conclusion~(\ref{rr.o.rrhop}) et
      le graphe entier  est un antitrou avec 8  sommets.  
    \item   Dans  le   troisième   graphe,  on   a  simultanément   la
      conclusion~(\ref{rr.o.rrleap})               et               la
      conclusion~(\ref{rr.o.rrhop}).   Le   graphe   entier   est   un
      $\overline{C_6}$.
  \end{itemize}} 
  \\ \rule{0cm}{1cm}
  \caption{Lemme de  Roussel et Rubio~: différents cas.}
\end{figure}

\begin{lemme}[Roussel et Rubio, \cite{roussel.rubio:01}]
  \label{rr.l.w}
  Soient $G$ un graphe de Berge,  $T$ un ensemble de sommets de $G$ et
  $P$  un chemin  de  $G$ tels  que  $G[T]$ est  anticonnexe, $T$  est
  disjoint  de $V(P)$  et les  extrémités de  $P$  sont $T$-complètes.
  Alors ou bien~:
  \begin{outcomes}
  \item 
    \label{rr.o.rreven}
    $P$ est de longueur paire et possède un nombre pair de $T$-arêtes.
  \item 
    \label{rr.o.rrodd}
    $P$  est  de longueur  impaire  et  possède  un nombre  impair  de
    $T$-arêtes.
  \item 
    \label{rr.o.rrleap}
    $P$ est de longueur impaire, supérieure  ou égale à 3, et si on le
    note $P= \bp x \tp x' \tp \cdots \tp y' \tp y \ep$, alors il existe
    dans $T$ deux sommets non adjacents $u$ et $v$ tels que $N(u) \cap
    V(P) = \{x, x', y\}$ et $N(v) \cap V(P) = \{x, y', y\}$.
  \item \label{rr.o.rrhop} $P$ est de longueur $3$ et ses deux sommets
    intérieurs  sont  les  extrémités  d'un  antichemin  de  longueur
    impaire dont l'intérieur est dans $T$.
  \end{outcomes}
\end{lemme}

\begin{preuve}
  On démontre le lemme par récurrence sur $|V(P) \cup T|$.  Si $P$ est
  de   longueur    $0$,   $1$   ou   $2$,    l'une   des   conclusions
  (\ref{rr.o.rreven})  ou~(\ref{rr.o.rrodd})  est  satisfaite, ce  qui
  montre que  le lemme est  vrai si $|V(P)  \cup T| \leq 4$.   On peut
  donc supposer que $P$ est de  longueur au moins $3$.  On note $P=\bp
  x  \tp x'\tp  \cdots \tp  y' \tp  y \ep$.   Supposons  qu'aucune des
  conclusions~(\ref{rr.o.rrleap})    et~(\ref{rr.o.rrhop})   ne   soit
  satisfaite.  Ces conclusions ne seront pas non plus satisfaites pour
  un sous-ensemble de $T$.  Donc, par hypothèse de récurrence, on sait
  que pour tout $U$, sous-ensemble  strict et anticonnexe de $T$, il y
  a un nombre de $U$-arêtes dans $P$ qui a même parité que la longueur
  de $P$. On distingue deux cas~:

  \vspace{2ex}

  {\noindent  \it Cas  1: Il  n'y a  pas de  sommet  $T$-complet dans
  $P^*$.}    Si  $P$   est   de   longueur  paire,   on   a  bien   la
  conclusion~(\ref{rr.o.rreven}).  Donc on  peut supposer que $P$ est
  de longueur impaire.  On distingue deux sous-cas qui conduisent tout
  deux à des contradictions~:

  {\vspace{2ex}}

  {\noindent  \it Cas  1.1: $T$  induit un  stable (preuve  de Kapoor,
  Vu\v{s}kovi\'{c} et  Zambelli).}  On marque  les sommets de  $P$ qui
  ont au moins un voisin  dans $T$.  On appelle \emph{intervalle} tout
  sous-chemin de $P$, de longueur au moins 1, dont les extrémités sont
  marquées et  dont les sommets  intérieurs sont non  marqués. Puisque
  $x$  et  $y$  sont  marqués,  les intervalles  de  $P$  forment  une
  partition de $P$ (au sens des arêtes).

  Nous prétendons que  chaque intervalle de $P$ est  de longueur paire
  ou  de longueur  1.   Car supposons  qu'il  y ait  un intervalle  de
  longueur impaire (au moins 3) noté  $P' = \bp x'' \tp \cdots \tp y''
  \ep$ avec $x''$ du côté de $x$  et $y''$ du côté de $y$. Soit $u$ un
  voisin de $x''$ dans $T$ et $v$ un voisin de $y''$ dans $T$.  Si $u=
  v$ alors  $V(P') \cup \{ u  \}$ induit un trou  impair.  Donc $u\neq
  v$. Si $x''\neq  x'$ alors $V(P') \cup \{ u, x,  v\}$ induit un trou
  impair.  Donc $x''=x'$.  De même,  $y''=y'$.  Mais alors $u$ et $v$
  permettent  de  satisfaire  la conclusion~(\ref{rr.o.rrleap})~:  une
  contradiction.

  Puisque  $P$  est de  longueur  impaire, il  y  a  un nombre  impair
  d'intervalles  de  longueur  $1$.   De  plus,  pour  chacun  de  ces
  intervalles  $\bp x''  \tp P  \tp y''  \ep$, $x''$  et $y''$  ont un
  voisin commun dans  $T$. Car sinon, soit $u \in  T$, voisin de $x''$
  manquant $y''$,  et $v\in  T$, voisin de  $y''$ manquant  $x''$.  Si
  $x''\neq x'$ alors  $\{ x, u, x'', y'', v\}$  induit un trou impair.
  Donc  $x''=x'$.   De  même,   $y''=y'$.   Mais  alors  $u$  et  $v$
  permettent  de  satisfaire  la conclusion~(\ref{rr.o.rrleap})~:  une
  contradiction.

  Pour tout $t\in  T$ on note $f(t)$ l'ensemble  des $\{t\}$-arêtes de
  $P$. On  note $t_1, t_2, \dots,  t_\alpha$ les éléments  de $T$.  On
  sait d'après ce qui précède que $|f(t_1) \cup f(t_2) \cup \dots \cup
  f(t_\alpha)|$ est le  nombre d'intervalles de longueur 1  de $P$. Ce
  nombre est impair et par la formule du crible on a~:

  \begin{eqnarray}
    |f(t_1) \cup f(t_2)  \cup \dots \cup f(t_\alpha)| & = & \sum_{i} |f(t_i)| \nonumber \\
    & & -  \sum_{i\neq j} |f(t_i) \cap f(t_j)| \nonumber \\
    & & \vdots \nonumber \\
    & & + (-1)^{(k+1)} \sum_{I\subset T, |I|=k} |\cap_{i\in I} f(t_i)| \nonumber \\
    & & \vdots \nonumber \\
    & & + (-1)^{(\alpha+1)} |f(t_1) \cap f(t_2) \cap \dots \cap f(t_\alpha)| \nonumber
  \end{eqnarray}

  Par hypothèse de  récurrence, on sait que pour  tout $I\subsetneq T$
  le nombre $|\cap_{i\in I} f(t_i)|$ est impair. De plus, le nombre de
  termes additionnés  dans le membre droit de  l'égalité ci-dessus est
  égal à  $2^\alpha -  1$ (nombre de  parties non vides  d'un ensemble
  avec  $\alpha$  éléments).   Donc  l'égalité  ci-dessus  peut  être
  réécrite modulo 2~:

  \[
  |f(t_1)  \cup  f(t_2) \cup  \dots  \cup  f(t_\alpha)|=  2^\alpha-2 +  (-1)^{(\alpha+1)}
  |f(t_1) \cap f(t_2) \cap \dots \cap f(t_\alpha)|
  \]

  Donc $|f(t_1) \cap f(t_2)  \cap \dots \cap f(t_\alpha)|$ est impair
  ce qui veut dire que $P$  possède un nombre impair de $T$-arêtes, et
  donc  que  $P^*$  possède   au  moins  un  sommet  $T$-complet~:  une
  contradiction.

  {\vspace{2ex}}

  {\noindent \it  Cas 1.1:  $T$ induit un  stable (notre  preuve).} Si
  $|T|=1$, alors $T\cup  V(P)$ induit un trou impair.   Donc $|T| \ge
  2$. Soit $t$ un sommet de  $T$.  Nous prétendons que $N(t) \cap V(P)
  = \{x,  x', y\}$ ou  $N(t) \cap P  = \{x,y',y\}$.  Notons  que, pour
  tout  $t\in T$,  les  sommets $T\setminus  \{t\}$-complets de  $P^*$
  manquent $t$.  On appelle  \emph{intervalle} de $P$ tout sous-chemin
  de $P$ dont les extrémités voient $t$ et dont les sommets intérieurs
  manquent $t$.  Puisque $x$ et $y$ voient $t$, les intervalles de $P$
  forment une partition de $P$  au sens des arêtes.  D'autre part nous
  savons qu'il y a un nombre impair de $T \setminus \{t\}$-arêtes dans
  $P$.  Donc il y a un intervalle de $P$ noté $\bp r \tp P \tp s \ep$
  qui  possède  un  nombre   impair  de  $T  \setminus  \{t\}$-arêtes.
  Supposons  que $x, r,  s, y$  apparaissent dans  cet ordre  sur $P$.
  Notons que $\bp r \tp P \tp  s \ep$ est de longueur au moins $2$ car
  sinon, l'un au  moins de $r$ ou de  $s$ serait  un sommet $T$-complet de
  $P^*$.  Donc  $V(\bp r \tp P  \tp s \ep) \cup\{t\}$  induit un trou
  pair et $\bp r \tp P \tp s \ep$ est de longueur paire.  Soit $r'$ le
  sommet $T  \setminus \{t\}$-complet de  $\bp r \tp  P \tp s  \ep$ le
  plus proche de $r$, et $s'$  le plus proche de $s$.  Notons que $\bp
  r' \tp P \tp s'  \ep$ contient toutes les $T \setminus \{t\}$-arêtes
  de  $\bp r  \tp  P \tp  s  \ep$ et  que ces  arêtes  sont en  nombre
  impair. Donc, par hypothèse de récurrence, $\bp r' \tp P \tp s' \ep$
  est  de  longueur impaire.   On  en  déduit  qu'exactement l'un  des
  chemins $\bp r \tp P \tp r' \ep$  et $\bp s' \tp P \tp s \ep$ est de
  longueur impaire.  On  suppose d'abord que $\bp r \tp  P \tp r' \ep$
  est de longueur impaire.  En particulier, $r\neq r'$ et $x \neq r$.

  Si $r'$ manque $y$, alors $P_1 = \bp r' \tp P \tp r \tp t \tp y \ep$
  est un chemin de $G$ de longueur impaire dont les extrémités sont $T
  \setminus  \{t\}$-complètes.   Or  $P_1$  n'a pas  de  $T  \setminus
  \{t\}$-arête.   Donc,   par  hypothèse  de   récurrence,  l'une  des
  conclusions~(\ref{rr.o.rrleap}) ou~(\ref{rr.o.rrhop}) est satisfaite
  pour $P_1$ et $T \setminus  \{t\}$. Mais cela est impossible car $t$
  n'a pas  de voisin dans $T  \setminus \{t\}$.  Donc  $r'$ doit voir
  $y$, ce qui implique $s' = s = y$ et $r' = y'$.

  Si $r$ manque $x$, alors $P_2 = \bp  r' \tp P \tp r \tp t \tp x \ep$
  est un chemin de $G$ de longueur impaire dont les extrémités sont $T
  \setminus  \{t\}$-complètes.   Or  $P_2$  n'a pas  de  $T  \setminus
  \{t\}$-arête.   Donc,   par  hypothèse  de   récurrence,  l'une  des
  conclusions~(\ref{rr.o.rrleap}) ou~(\ref{rr.o.rrhop}) est satisfaite
  pour $P_2$ et  $T \setminus \{t\}$. Mais cela  est encore impossible
  car $t$ n'a  pas de voisin dans $T \setminus  \{t\}$.  Donc $r$ doit
  voir $x$, ce qui  implique $r = x'$.  On a donc  bien $N(t) \cap P =
  \{x, x', y\}$. De même, si $\bp  s' \tp P \tp s \ep$ est de longueur
  impaire, alors  on montre $N(t) \cap P  = \{x, y', y\}$.   On a bien
  montré  que $N(t)  \cap  V(P) =  \{x, x',  y\}$  ou $N(t)  \cap P  =
  \{x,y',y\}$.

  Maintenant, puisque  $y'$ n'est  pas $T$-complet, il  y a  un sommet
  $u\in T$ tel  que $N(u) \cap P = \{x,x',y\}$,  et puisque $x'$ n'est
  pas $T$-complet  il y a un  sommet $v\in T$  tel que $N(v) \cap  P =
  \{x,y',y\}$.   Donc   $u$  et  $v$  permettent   de  satisfaire  la
  conclusion~(\ref{rr.o.rrleap})~: une contradiction.

  {\vspace{2ex}}
  
  {\noindent  \it Cas  1.2:  $T$  n'induit pas  un  stable (preuve  de
  Roussel et Rubio un peu  simplifiée).} Soit $Q = \overline{\bp u \tp
  \cdots \tp v \ep}$ un  plus long antichemin de $G[T]$.  Comme $G[T]$
  est  anticonnexe avec au  moins une  arête, $Q$  est de  longueur au
  moins $2$. De plus, par  le choix de $Q$,  $G[T\setminus \{u\}]$
  et $G[T \setminus  \{v\}]$ sont anticonnexes.  Donc on  sait que $P$
  possède un  nombre impair de $T\setminus \{u\}$-arêtes  et un nombre
  impair de  $T \setminus  \{v\}$-arêtes.  Notons qu'une  $T \setminus
  \{u\}$-arête et  une $T \setminus  \{v\}$-arête n'ont pas  de sommet
  commun car il  y aurait alors un sommet  $T$-complet dans $P^*$.  En
  particulier,    les    $T    \setminus    \{u\}$-arêtes    et    les
  $T\setminus\{v\}$-arêtes sont deux à deux distinctes.

  Supposons  que $Q$  soit de  longueur  paire.  Soit  $x_u x'_u$  une
  $T\setminus  \{u\}$-arête  de  $P$  et $y'_v  y_v$  une  $T\setminus
  \{v\}$-arête de $P$  telles que, sans perte de  généralité, $x, x_u,
  x'_u, y'_v, y_v, y$ apparaissent  dans cet ordre sur $P$.  Si $x'_u$
  manque $y'_v$  alors $\{x'_u, y'_v\}  \cup V(Q)$ induit  un antitrou
  impair.  Si  $x\neq x_u$ alors  $\{x_u, y'_v\} \cup V(Q)$  induit un
  antitrou impair.   Si $y_v \neq  y$ alors $\{x'_u, y_v\}  \cup V(Q)$
  induit un antitrou impair.  Donc $P  = \bp x_u \tp x'_u \tp y'_v \tp
  y_v \ep$, et  $V(P) \cup V(Q)$ induit un  antitrou impair.  Donc $Q$
  est de longueur impaire (au moins $3$).

  Supposons  que  $G[T\setminus\{u,v\}]$   ne  soit  pas  anticonnexe.
  Puisque $G[T]$  est anticonnexe,  il existe un  sommet $w$  dans une
  composante  connexe de  $\overline{G}[T \setminus  \{u,v\}]$  qui ne
  contient pas $Q^*$ et tel  que $w$ voit dans $\overline{G}$ au moins
  l'un  de  $u$  et  $v$.  Mais  alors  $V(Q)\cup\{w\}$  induit  dans
  $\overline{G}[T]$ un chemin plus long que $Q$ ou un trou impair~: une
  contradiction.  Donc $G[T \setminus \{u, v\}]$ est anticonnexe.

  Nous  savons  qu'il  y  a  un  nombre  impair  de  $T\setminus  \{u,
  v\}$-arêtes dans $P$.  Rappelons que $P$ possède un nombre impair de
  $T\setminus\{u\}$-arêtes,     et     un     nombre     impair     de
  $T\setminus\{v\}$-arêtes, qui sont différentes, et qui totalisent un
  nombre pair  de $T\setminus \{u, v\}$-arêtes.  Donc  $P$ possède au
  moins une  $T\setminus \{u,  v\}$-arête $x'' y''$  qui n'est  ni une
  $T\setminus  \{u\}$-arête  ni  une  $T \setminus  \{v\}$-arête.   On
  suppose sans perte de généralité que $x,x'',y'',y$ apparaissent dans
  cet ordre  sur $P$ et que  $y''\in P^*$.  Donc $y''$  manque $u$ ou
  $v$, par exemple  $v$.  Alors $y''$ voit $u$,  car sinon $V(Q) \cup
  \{y''\}$  induirait un antitrou  impair~; et  $x''$ manque  $u$, car
  sinon, $x''y''$ serait une  $T\setminus \{v\}$-arête~; et $x''$ voit
  $v$,  car sinon  $V(Q)\cup \{x''\}$  induirait antitrou  impair~; et
  $x''=x'$ car  sinon $V(Q)\cup\{x'', y'', x\}$  induirait un antitrou
  impair, et de  même $y''=y'$.  Donc $P=\bp x \tp x''  \tp y'' \tp y
  \ep$  et $V(Q)  \cup \{x'',  y''\}$  est un  antichemin de  longueur
  impaire  de $G$~:  la conclusion~(\ref{rr.o.rrhop})  est satisfaite,
  une contradiction.

  {\vspace{2ex}}

  \noindent  {\it Cas 2:  Il y  a un  sommet $T$-complet  dans $P^*$.}
  Soit $z$ un tel sommet.  Par hypothèse de récurrence, on applique le
  lemme au chemin $\bp x \tp P \tp  z \ep$ et à $T$.  Si on obtient la
  conclusion~(\ref{rr.o.rrhop}), alors $\bp x \tp P \tp z \ep$ est de
  longueur impaire, au moins 3, on le  note $\bp x \tp x' \tp P \tp z'
  \tp z \ep$ et  il existe deux sommets $u\in T$ et  $v\in T$ tels que
  $N(u) \cap V(P) = \{x, x', z\}$  et $N(v) \cap V(P) = \{x, z', z\}$.
  Mais alors $V(\bp x'  \tp P \tp z' \ep) \cup \{u,  v, y\}$ induit un
  trou impair.  Si  on obtient la conclusion~(\ref{rr.o.rrleap}) alors
  $\bp x  \tp P \tp  z \ep$  est de longueur  $3$ et ses  deux sommets
  intérieurs  sont les  extrémités  d'un antichemin  $Q$  de $G$  dont
  l'intérieur est dans  $T$.  Mais alors $V(Q) \cup  \{y\}$ induit un
  antitrou     impair.      Donc      on     obtient     l'une     des
  conclusions~(\ref{rr.o.rreven})  ou~(\ref{rr.o.rrodd}), et  le nombre
  de $T$-arêtes dans $\bp x \tp P  \tp z \ep$ et la longueur de $\bp x
  \tp P  \tp z \ep$  sont de même  parité.  On obtient  une conclusion
  similaire pour  $\bp z  \tp P \tp  y \ep$.   Au total, le  nombre de
  $T$-arêtes dans  $P$ et la  longueur de $P$  sont de même  parité et
  l'une des  conclusions~(\ref{rr.o.rreven}) ou~(\ref{rr.o.rrodd}) est
  satisfaite.
\end{preuve}

Nous  donnons  maintenant  deux  variantes  du  lemme  de  Roussel  et
Rubio. La première, qui n'a, à notre connaissance, aucune application,
montre  que si  $T$ est  stable  alors on  peut étendre  le lemme  aux
graphes      sans      trou      impair     et      abandonner      la
conclusion~(\ref{rr.o.rrhop})~:

\begin{lemme}
  Soient $G$ un graphe sans trou impair, $T$ un ensemble de sommets de
  $G$ et  $P$ un  chemin de $G$  tels que  $G[T]$ est stable,  $T$ est
  disjoint  de $V(P)$  et les  extrémités de  $P$  sont $T$-complètes.
  Alors ou bien~:
  \begin{outcomes}
  \item 
    $P$ est de longueur paire et possède un nombre pair de $T$-arêtes.
  \item 
    $P$  est  de longueur  impaire  et  possède  un nombre  impair  de
    $T$-arêtes.
  \item 
    $P$ est de longueur impaire, supérieure  ou égale à 3, et si on le
    note $P=  \bp x  \tp x'  \tp \cdots \tp  y' \tp  y \ep$,  alors il
    existe dans  $T$ deux  sommets non adjacents  $u$ et $v$  tels que
    $N(u) \cap  V(P) =  \{x, x', y\}$  et $N(v)  \cap V(P) =  \{x, y',
    y\}$.
  \end{outcomes}
\end{lemme}

\begin{preuve}
  Il  suffit  de  reprendre  la   preuve  du  lemme  principal  en  se
  retreignant au cas ``$T$ stable''.
\end{preuve}

La  deuxième variante  que  nous présentons  ci-dessous  est celle  de
l'article  original.   Elle  donne  une conclusion  moins  forte  mais
présente  l'avantage d'être  plus  symétrique et  de  bien montrer  la
similitude    qui    ne    saute    pas    aux    yeux    entre    les
conclusions~(\ref{rr.o.rrleap})  et~(\ref{rr.o.rrhop})  de la  version
principale. Nous avons recherché  des versions encore plus symétriques
du lemme,  en essayant par exemple  de remplacer le chemin  $P$ par un
ensemble  connexe  quelconque.   C'est  sans doute  possible  au  prix
d'artifices  divers  qui  n'ont   ---  nous  semble-t-il  ---  pas  de
meilleures motivations que le culte de la symétrie.

\begin{lemme}[Roussel et Rubio, \cite{roussel.rubio:01}]
  \label{rr.l.winutile}
  Soient $G$ un graphe de Berge,  $T$ un ensemble de sommets de $G$ et
  $P$  un chemin  de  $G$ tels  que  $G[T]$ est  anticonnexe, $T$  est
  disjoint  de $V(P)$  et les  extrémités de  $P$  sont $T$-complètes.
  Alors ou bien~:
  \begin{outcomes}
  \item 
    $P$ est de longueur paire et possède un nombre pair de $T$-arêtes.
  \item 
     $P$  est  de longueur  impaire  et  possède  un nombre  impair  de
    $T$-arêtes.
  \item 
    Il  existe un chemin  de longueur  impaire (au  moins 3)  dont les
    extrémités sont dans $T$ et dont l'intérieur est dans $V(P)$.
  \item  
    Il existe un antichemin de  longueur impaire (au moins 3) dont les
    extrémités sont dans $V(P)$ et dont l'intérieur est dans $T$.
  \end{outcomes}
\end{lemme}

Cette version  est une conséquence  directe de la  version principale.
Réciproquement, il  est assez facile de prouver  la version principale
en supposant connue la version ci-dessus.

\section{Restriction à des classes de graphes particuliers}

Le  lemme  de  Roussel  et  Rubio reste  évidemment  vrai  dans  toute
sous-classe des  graphes de Berge, mais dans  certaines d'entre elles,
on  obtient une  version  plus forte.   Nous  examinerons les  graphes
d'Artémis  pairs, les  graphes de  Meyniel et  les  graphes faiblement
triangulés. Au  chapitre suivant,  les trois lemmes  ci-dessous seront
utilisés pour démontrer l'existence de paires d'amis.

\index{Artémis~(graphe~d'---)!lemme~de~Roussel~et~Rubio}
\index{lemme~de~Roussel~et~Rubio!graphes~d'Artémis}
\index{lemme~de~Roussel~et~Rubio!graphes~d'Artémis~pairs}

\begin{lemme}
  \label{rr.l.wa}
  Soit $G$  un graphe d'Artémis  pair.  Soit $P$  un chemin de  $G$ et
  $T\subset V(G)$ disjoint de  $V(P)$ tels que $G[T]$ est anticonnexe
  et que les extrémités de  $P$ sont $T$-complètes.  Alors le nombre
  de  $T$-arêtes de  $P$ a  même parité  que la  longueur de  $P$.  En
  particulier, si  $P$ est de  longueur impaire, alors $P$  contient au
  moins un sommet $T$-complet.
\end{lemme}
\begin{preuve}
  On      applique      le      lemme~\ref{rr.l.w}.       Dans      la
  conclusion~(\ref{rr.o.rrleap}),   l'ensemble  $V(P)   \cup  \{u,v\}$
  induit  un  prisme  impair,   un  sous-graphe  interdit.   Dans  la
  conclusion~(\ref{rr.o.rrhop}),  si on note  $Q$ l'antichemin  de $G$
  dont l'intérieur est  dans $T$ et dont les  extrémités sont les deux
  sommets intérieurs de $P$, alors $V(P) \cup V(Q)$ induit un antitrou
  de $G$ de  taille au moins 6, encore  un sous-graphe interdit. Donc,
  comme    annoncé,    seules    les    conclusions~(\ref{rr.o.rreven})
  et~(\ref{rr.o.rrodd}) peuvent être satisfaites.
\end{preuve}

Le lemme suivant était déjà  connu de Meyniel sous une présentation un
peu  différente. Il montre  que dans  les graphes  de Meyniel  on peut
relaxer  l'hypothèse  d'anticonnexité  de  $T$  tout  en  donnant  une
conclusion beaucoup plus forte~:

\index{Meyniel~(graphe~de~---)!lemme~de~Roussel~et~Rubio}
\index{lemme~de~Roussel~et~Rubio!graphes de Meyniel}

\begin{lemme}[Meyniel, \cite{meyniel:76}]
  \label{rr.l.rrmeyniel}
  Soit  $G$ un  graphe  de Meyniel.   Soit  $P$ un  chemin  de $G$  de
  longueur impaire et $T \subset V(G)$ disjoint de $V(P)$ tels que les
  extrémités $P$  soient $T$-complètes.   Alors tout sommet  de $P$
  est $T$-complet.
\end{lemme}

\begin{preuve}
  Soit $t$ un sommet de  $T$. On appelle \emph{intervalle} de $P$ tout
  sous-chemin  de $P$  dont  les  extrémités voient  $t$  et dont  les
  sommets intérieurs manquent $t$.  Comme les extrémités de $P$ voient
  $t$, on sait que les intervalles de $P$ forment une partition de $P$
  au sens des  arêtes. De plus tout intervalle $P'$  est de longueur 1
  ou de longueur paire, sans quoi $V(P') \cup \{t\}$ induirait un trou
  impair. Donc, $P$ étant de longueur impaire, il contient au moins un
  intervalle de longueur~1.  Si $t$  manque un sommet de $P$, alors il
  existe dans $P$ des intervalles de longueur paire et des intervalles
  de longueur~1. Donc  on peut trouver sur $P$ des  sommets $x, y, z$
  vérifiant $x-y$ est un intervalle de  longueur~1 et $\bp y \tp P \tp
  z \ep$ est un intervalle de longueur paire. Mais alors $V(\bp x \tp
  P  \tp z  \ep) \cup  \{t\}$ induit  un cycle  impair avec  une seule
  corde, un sous-graphe interdit. Donc $t$ voit  tous les sommets de
  $P$.
\end{preuve}

Le lemme  de Roussel  et Rubio peut  également être renforcé  pour les
graphes   faiblement  triangulés.    Cette  fois,   on   peut  relaxer
l'hypothèses de parité de $P$.

\index{lemme~de~Roussel~et~Rubio!graphes~faiblement~triangulés}
\index{faiblement~triangulé~(graphe~---)!lemme~de~Roussel~et~Rubio}

\begin{lemme}
  \label{rr.l.rrwt}
  Soit $G$ un graphe faiblement triangulé.  Soit $P = \bp x \tp \cdots
  \tp y  \ep$ un chemin de  $G$ de longueur  au moins 3 et  $T \subset
  V(G)$ disjoint de $V(P)$ tels  $G[T]$ est anticonnexe et tel que les
  extrémités  de $P$ sont  $T$-complètes.  Alors  il existe  un sommet
  intérieur de $P$ qui est $T$-complet.
\end{lemme}

\begin{preuve}
  Remarquons d'abord  qu'aucun sommet $t\in T$ ne  manque deux sommets
  consécutifs de $P$, car sinon, $V(P) \cup \{t\}$ contient un trou de
  taille  au moins  5, un  sous-graphe  interdit.  Soit  alors $z$  un
  sommet de $V(P^*)$  qui voit un maximum de sommets  de $T$.  On peut
  supposer qu'il existe un sommet  $u \in T \setminus N(z)$ car sinon,
  $z$ permet de satisfaire la  conclusion du lemme.  Soient alors $x'$
  et $y'$ les  voisins de $z$ sur $P$, choisis de  sorte que $x, x',
  z, y',  y$ apparaissent dans cet  ordre sur $P$.   Alors, d'après la
  remarque  initiale, $u$  voit $x'$  et $y'$.   Quitte à  échanger les
  rôles de $x$ et $y$, on peut supposer $x' \neq x$.  D'après le choix
  de $z$, puisque $u$ voit $x'$  et manque $z$, il existe un sommet $v
  \in  T$  qui  voit $z$  et  qui  manque  $x'$.  Puisque  $G[T]$  est
  anticonnexe, il existe  un antichemin $Q$ de $G[T]$  reliant $u$ et
  $v$, et on  chosit $u$ et $v$ de manière à  minimiser la longueur de
  cet   antichemin.   D'après  la   remarque  initiale,   les  sommets
  intérieurs  $Q$ voient chacun  l'un de  $x'$ ou de  $z$, et  d'après le
  choix de $Q$, les sommets intérieurs de $Q$ voient tous $x'$ et $z$.
  Si  $x'$ manque  $x$  alors $V(Q)  \cup  \{ x,  x',  z\}$ induit  un
  antitrou de taille au moins $5$. Donc $x'$ voit $x$.  Si $z$ manque
  $y$, alors $V(Q) \cup \{z, x', y\}$ induit un antitrou de taille au
  moins 5.   Donc $z$ voit  $y$ et $y  = y'$. Mais alors  $V(Q) \cup
  \{x,  x', z,  y\}$ induit  un  antitrou de  taille au  moins 5~:  une
  contradiction.
\end{preuve}

Nous laissons  ouverte la question suivante~:  existe-t-il une version
du  lemme   de  Roussel  et   Rubio  pour  les   graphes  parfaitement
ordonnables~?  Cela serait  assez surprenant car on ne  connait pas de
caractérisation des graphes  parfaitement ordonnables par sous-graphes
induits interdits.

\section{Corollaires et variantes utiles}
\index{lemme~de~Roussel~et~Rubio!corollaires}

Ici, nous  présentons des corollaires  et des variantes de  la version
``Artémis'' du lemme  de Roussel et Rubio qui  révèleront leur utilité
au chapitre  suivant. Nous proposons  tout d'abord un  lemme technique
suggéré par  Bruce Reed.  Dans la  première version de  ce travail, ce
lemme était démontré et utilisé implicitement à de nombreuse reprises.
Le formuler à part permet donc de simplifier certaines preuves (toutes
les preuves où on l'invoque).

\begin{lemme}
  \label{rr.l.bruce}
  Soit $G$ un  graphe, et $\{a, b, c\}$ un  triangle de $G$. Supposons
  que  $G  \setminus \{a,b,c\}$  est  connexe,  que  chaque sommet  de
  $G\setminus \{a,b,c\}$ voit au plus un sommet de $\{a,b,c\}$, et que
  chaque  sommet  de $\{a,b,c\}$  possède  exactement  un voisin  dans
   $G\setminus \{a,b,c\}$. Alors, $G$ possède un prisme ou une pyramide
  de triangle $\{a,b,c\}$.
\end{lemme}

\begin{preuve}
  Soit $P$ un  chemin de $G\setminus \{a,b,c\}$ reliant  un voisin $u$
  de  $a$ et  un voisin  $v$  de $b$,  de longueur  minimale avec  ces
  propriétés. Soit $Q$ un chemin de $G\setminus \{a,b\}$ reliant $c$ à
  un sommet $w$ ayant des voisins dans $P$, de longueur minimale
  avec  ces   propriétés.   Notons  que   $P$  et  $Q$   existent  par
  connexité.  Notons que  $u$, $v$  et $w$  sont distincts  car chaque
  sommet  de  $G\setminus  \{a,b,c\}$   voit  au  plus  un  sommet  de
  $\{a,b,c\}$. On examine trois cas~:

  \begin{itemize}
  \item
    Si $w$  possède un  unique voisin $w'$  sur $P$, alors,  les trois
    chemins $\bp w' \tp w \tp Q \tp  c \ep$, $\bp w' \tp P \tp u \tp a
    \ep$ et $\bp w'  \tp P \tp v \tp b \ep$  induisent une pyramide de
    coin $w'$ et de triangle $\{a,b,c\}$.
  \item
    Si  $w$ possède exactement  deux voisins  adjacents $w',  w''$ sur
    $P$, alors  on suppose que $u,  w', w'', v$  apparaissent dans cet
    ordre sur $P$.  Les trois chemins $\bp w \tp Q \tp c \ep$, $\bp w'
    \tp  P \tp  u \tp  a \ep$  et  $\bp w''  \tp P  \tp v  \tp b  \ep$
    induisent  un prisme  de triangles  $\{a,b,c\}$  et $\{w,w',w''\}$
    (notons  que  $w\neq c$  car  $c$  voit  exactement un  sommet  de
    $G\setminus \{a,b,c\}$).
  \item
    Si $w$ possède au moins  deux voisins non adjacents sur $P$, alors
    on appelle $w'$ le voisin de $w$ sur $P$ le plus proche de $u$, et
    $w''$ le plus proche de $v$.  Les trois chemins $\bp w \tp Q \tp c
    \ep$, $\bp w \tp w' \tp P \tp u \tp a \ep$ et $\bp w \tp w'' \tp P
    \tp v \tp b \ep$ induisent une pyramide de coin $w$ et de triangle
    $\{a,b,c\}$  (notons que  $w\neq  c$ car  $c$  voit exactement  un
    sommet de $G\setminus \{a,b,c\}$).
  \end{itemize}
\end{preuve}

\begin{lemme}
  \label{rr.l.wh}
  Soit  $G$ un  graphe d'Artémis  pair,  soit $H$  un trou  de $G$  et
  $T\subset  V(G)$ tels  que $V(H)$  est disjoint  de $T$  et  que $H$
  comporte au moins deux sommets $T$-complets non adjacents. Alors $H$
  possède un nombre pair de $T$-arêtes.
\end{lemme}

\begin{preuve}
  Soient $x$ et $y$ deux sommets non-adjacents et $T$-complets de $H$.
  Comme  $H$ est  un trou  pair, les  deux chemins  de $H$  ayant pour
  extrémités  $x$  et $y$  sont  de  même  parité.  Donc,  d'après  le
  lemme~\ref{rr.l.wa}, leurs nombres de $T$-arêtes sont égaux
  modulo~2.
  Donc, au total, $H$ possède bien un nombre pair de $T$-arêtes.
\end{preuve}

On rappelle que les \emph{quasi-prismes} sont des graphes particuliers
définis  formellement  à   la  section~\ref{base.ss.pps}  (voir  aussi
figure~\ref{rr.fig.quasi-prisme}).

\begin{lemme}
  \label{rr.l.sgp}
  Soient  $G$  un  graphe   d'Artémis,  $S$  un  quasi-prisme  de  $G$
  d'extrémités $a$ et  $b$, et $T\subset V(G)$ un  ensemble de sommets
  de $G$ disjoint de $V(S)$,  tels que $G[T]$ est anticonnexe et $a,b$
  sont $T$-complets.   Alors dans  chaque triangle de  $S$, il y  a au
  moins deux sommets $T$-complets.
\end{lemme}

\begin{preuve}
  On note les sommets et les chemins du quasi-prisme $S$ comme indiqué sur
  la figure~\ref{rr.fig.quasi-prisme}.

  \begin{figure}[ht]
    \center
        \includegraphics{fig.rr.4}  \\  \rule{0cm}{3ex} \\ 
    \parbox{8.5cm}{
      \begin{itemize}
      \item
	$S_1$ est le chemin reliant $a$ et $a'$.
      \item
	$S_2$ est le chemin reliant $b$ et $b'$.
      \item
	$S_3$ est le chemin reliant $c$ et $c'$ sans passer par $d$.
      \item
	$S_4$ est le chemin reliant $d$ et $d'$ sans passer par $c$.
    \end{itemize}}
    \caption{Notations  pour  les  quasi-prismes d'extrémités  $a$  et
    $b$\label{rr.fig.quasi-prisme}}
  \end{figure}
  
  \begin{claim}
    \label{rr.c.sgp}
    Chaque  sommet $t\in  T$ voit  au moins  deux sommets  du triangle
    $\{a', c, d\}$
  \end{claim}
  
  \begin{preuveclaim}
    Supposons qu'il existe un sommet  $t\in T$ ayant au plus un voisin
    dans  $\{a',  c, d\}$.   Si  $t$  manque  entièrement le  triangle
    $\{a',c,d\}$ ou si  $a = a'$ et $N(t) \cap \{a',  c, d\} = \{a\}$,
    alors on  applique le lemme~\ref{rr.l.bruce} au  graphe induit par
    $V(S)  \cup \{t\}$.  Sinon,  $t$ a  un unique  voisin $v$  dans le
    triangle $\{a',c,d\}$, on note $v'$ l'unique voisin de $v$ dans $S
    \setminus \{a', c, d\}$,  et on applique le lemme~\ref{rr.l.bruce}
    au graphe  induit par $(V(S) \cup \{t\})  \setminus \{v'\}$.  Dans
    les deux cas  les hypothèses du lemme sont  bien satisfaites grâce
    au triangle $\{a',  c, d\}$, d'où une contradiction  car un graphe
    d'Artémis ne contient ni prisme ni pyramide.
  \end{preuveclaim}

  Supposons maintenant que le lemme  soit faux pour le triangle $\{a',
  c,  d\}$, c'est-à-dire  qu'il existe  deux  sommets $\alpha,\beta\in
  \{a',  c,  d\}$ et  deux  sommets $u,v\in  T$  tels  que $u$  manque
  $\alpha$ et $v$  manque $\beta$.  D'après~(\ref{rr.c.sgp}), $u$ voit
  $\beta$ et $v$ voit $\alpha$.   Puisque $G[T]$ est anticonnexe, il y
  a un antichemin $Q$ d'extrémités $u$  et $v$ dans $G[T]$, et on peut
  choisir  $u$  et  $v$  de  manière  à rendre  $Q$  aussi  court  que
  possible.  Cela  entraîne  que  les  sommets de  $Q^*$  voient  tous
  $\alpha$ et  $\beta$.  Mais alors  $V(Q) \cup \{\alpha,  b, \beta\}$
  induit un antitrou  de $G$~: contradiction.  Donc au  moins deux des
  sommets  $a', c$ et  $d$ sont  $T$-complets.  On  traite de  la même
  façon l'autre triangle du quasi-prisme.
\end{preuve}

\begin{figure}[p]
    \center
    \includegraphics{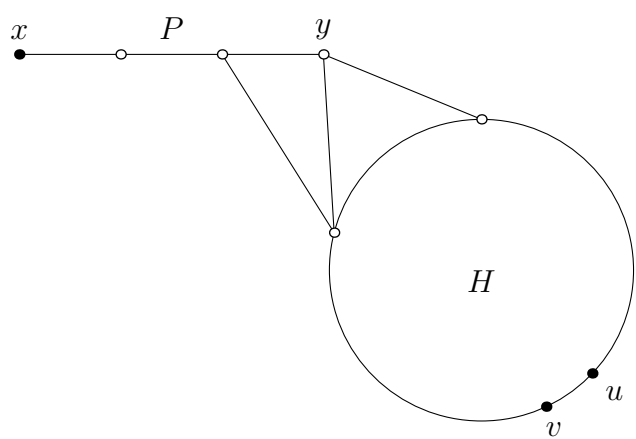} 
    \caption{Lemme~\ref{rr.l.sgt}. Les sommets ``pleins'' sont dans $C(T)$\label{rr.fig.sgt}}
  \end{figure}

\begin{lemme}
  \label{rr.l.sgt}
  Soient $G$ un graphe d'Artémis, $H$ un trou de $G$, $P$ un chemin de
  $G$  d'extrémités $x$  et  $y$,  et $T\subset  V(G)$  tels que  $T$,
  $V(H)$,  et  $V(P)$  sont  deux  à  deux  disjoints  et  $G[T]$  est
  anticonnexe. Soient $ab$ et $cd$  deux arêtes disjointes de $H$ tels
  que l'ensemble  des arêtes entre $P$  et $H$ est  $\{ay, by\}$.  Si
  $c$, $d$ et $x$ sont $T$-complets, alors au moins l'un de $a$ et $b$
  est $T$-complet.
\end{lemme}
\begin{preuve}
  On       pourra      consulter       la      figure~\ref{rr.fig.sgt}
  page~\pageref{rr.fig.sgt}.   On suppose,  quitte à  échanger  $c$ et
  $d$, que $a,c,d,b$ apparaissent dans  cet ordre sur $H$.  On appelle
  $P_1$  le  chemin  d'extrémités  $a$  et  $c$,  contenu  dans  $V(H)
  \setminus  \{b,d\}$, et  $P_2$  le chemin  d'extrémités  $b$ et  $d$
  contenu dans  $V(H) \setminus\{a,c\}$.  Supposons que  le lemme soit
  faux~: il  existe un  sommet $u\in  T$ qui manque  $a$ et  un sommet
  $v\in T$ qui manque $b$.

  \begin{claim}
    Le sommet $u$ voit $b$ et $y$. Le sommet $v$ voit $a$ et $y$.
  \end{claim}
  
  \begin{preuveclaim}
    Supposons que  $u$ ait au plus  un voisin dans $\{a,  b, y\}$.  Si
    $u$ manque entièrement le triangle $\{a,  b, y\}$ ou si $x = y$ et
    $N(u)  \cap   \{a,  b,  y\}   =  \{y\}$,  alors  on   applique  le
    lemme~\ref{rr.l.bruce} au  graphe induit par $V(P)  \cup V(H) \cup
    \{u\}$.  Sinon, $u$ a un unique voisin $u'$ dans le triangle $\{a,
    b, y\}$,  on note $u''$ l'unique  voisin de $u'$  dans $(V(P) \cup
    V(H))   \setminus    \{a,   b,    y\}$,   et   on    applique   le
    lemme~\ref{rr.l.bruce} au graphe induit  par $(V(P) \cup V(H) \cup
    \{u\}) \setminus  \{u''\}$.  Dans les  deux cas les  hypothèses du
    lemme sont bien satisfaites grâce  au triangle $\{a, b, y\}$, d'où
    une contradiction car un graphe d'Artémis ne contient ni prisme ni
    pyramide.
  \end{preuveclaim}

  Puisque  $G[T]$  est  anticonnexe,   il  existe  un  antichemin  $Q$
  d'extrémités  $u$ et  $v$ dans  $G[T]$.  On  choisit $u$  et  $v$ de
  manière  à  minimiser  la  longueur  de cet  antichemin.   Donc  les
  éventuels sommets intérieurs  de $Q$ voient tous $a$  et $b$.  Si un
  sommet  $w\in\{y,  c, d\}$  manque  $a$  et  $b$, alors  $V(Q)  \cup
  \{a,b,w\}$  induit  un  antitrou  de  longueur  au  moins~$5$~:  une
  contradiction.  Donc on doit  avoir $x=y$, $ac\in E(G), bd\in E(G)$.
  Mais alors $V(Q)\cup\{a,b,c,d,x\}$ induit un antitrou de longueur au
  moins~$7$~: une contradiction.
\end{preuve}

\begin{lemme}
  \label{rr.l.dpg}
  Soit $G$ un  graphe d'Artémis, $H$ un trou de $G$,  $P$ un chemin de
  $G$  d'extrémités $x$ et  $y$ et  $T\subset V(G)$  tels que  $T$ est
  disjoint de  $V(P) \cup  V(H)$, $G[H \cup  P]$ est connexe,  $x$ est
  $T$-complet et  il existe  deux sommets $T$-complets  adjacents $u,v
  \in V(H)$.   Alors ou bien  il existe un  sommet de $P$ qui  voit au
  moins l'un  de $u$ et $v$,  ou bien il existe  un sommet $T$-complet
  dans $V(H) \setminus \{u,v,x\}$.
\end{lemme}
\begin{preuve}
  On  pourra consulter  la figure~\ref{rr.fig.dpg}.   Notons  qu'on ne
  suppose pas $H$  et $P$ disjoints.  En fait  $V(P) \subset V(H)$ est
  même possible.  Suposons  que le lemme soit faux,  et considérons un
  contre-exemple avec $|H  \cup P|$ minimal.  On note  $H_u$ le chemin
  induit  par $V(H)  \setminus \{v\}$  et $H_v$  le chemin  induit par
  $V(H) \setminus \{u\}$.

  Soit $x'$  le sommet de  $P$ le plus  proche de $x$ ayant  un voisin
  dans $H$.  Soit  $u'$ le voisin de $x'$ sur $H_u$  le plus proche de
  $u$ et $v'$ le voisin de $x'$ sur $H_v$ le plus proche de $v$.

  Comme $H$, $P$ et $T$ forment un contre-exemple, $u' \not= u$ et $v'
  \not= v$.   Si $x \in  H$ alors d'après le  lemme~\ref{rr.l.wh}, $H$
  possède un nombre pair  de $T$-arêtes et donc $H\setminus \{u,v,x\}$
  possède  un  sommet $T$-complet  et  $H, P,  T$  ne  forment pas  un
  contre-exemple.   Donc  on  peut  supposer  $x \notin  H$,  ce  qui
  entraîne $V(\bp x \tp P \tp x' \ep) \cap V(H)= \emptyset$.

  Supposons $u'= v'$.  L'un des chemins $\bp x \tp P \tp x' \tp u' \tp
  H_u \tp u \ep$ et $\bp x \tp P  \tp x' \tp v' \tp H_v \tp v \ep$ est
  de  longueur  impaire.  D'après  le  lemme~\ref{rr.l.wa}, ce  chemin
  possède un  nombre impair de $T$-arêtes.   Donc il y a  au moins un
  sommet $T$-complet $x''$ dans $ (V(\bp  x \tp P \tp x' \ep) \cup V(H))
  \setminus \{u,v,x\}$.  Si $x''$ est  dans $H$, alors $H$, $P$ et $T$
  ne  forment pas un  contre-exemple.  Sinon,  c'est-à-dire si  $x'' \in
  V(\bp x \tp P \tp x' \ep)  \setminus \{x\}$, alors le trou $H$ et le
  chemin $\bp  x'' \tp P  \tp x' \ep$  forment un contre-exemple avec
  $|V(H)  \cup  V(\bp  x'' \tp  P  \tp  x'  \ep)|  < |H\cup  P|$~:  une
  contradiction.

  Supposons  $u' \neq  v'$ et  $u'v'\in E(G)$.   On peut  appliquer le
  lemme~\ref{rr.l.sgt} au  trou $H$ et au  chemin $\bp x \tp  P \tp x'
  \ep$~: au moins l'un de $u', v'$ est $T$-complet~: une contradiction.
  
  Supposons enfin $u'  \neq v'$ et $u'v' \notin  E$.  On définit $H'$,
  trou induit par  $V(\bp u \tp H_u  \tp u' \ep) \cup V(\bp  v \tp H_v
  \tp v'  \ep) \cup\{x'\}$.  Alors $H'$, $\bp  x \tp P \tp  x' \ep$ et
  $T$ forment un contre-exemple au lemme avec $|V(H') \cup V(\bp x \tp
  P \tp x' \ep)| < |V(H) \cup V(P)|$~: une contradiction.
\end{preuve}

  \begin{figure}[p]
    \center
    \includegraphics{fig.rr.6} 
    \caption{Lemme~\ref{rr.l.dpg}. Les sommets ``pleins'' sont dans $C(T)$\label{rr.fig.dpg}}
  \end{figure}

\chapter{Graphes d'Artemis}
\label{artemis.chap}

Ce     chapitre    est     consacré    à     la    preuve     de    la
conjecture~\ref{pair.conj.Artemis}~:   tout   graphe   d'Artémis   est
parfaitement contractile.   À l'exception notable de  l'usage du lemme
de Roussel  et Rubio,  l'essentiel des idées  était déjà  présent dans
l'étude  par Linhares  et Maffray  du cas  restreint aux  graphes sans
carré~\cite{linhares.maffray:evenpairsansc4}.  Notre  preuve donne un
algorithme de complexité $O(n^2m)$ pour colorier les graphes d'Artémis.
Nous spécialiserons notre preuve  à des classes restreintes de graphes
d'Artémis.  Nous verrons que pour les graphes de Meyniel, notre preuve
peut  être vue  comme une  généralisation  de la  preuve classique  de
l'existence   d'une  paire  d'amis.    Pour  les   graphes  faiblement
triangulés,   notre  preuve  donne   une  nouvelle   démonstration  de
l'existence  d'une  2-paire  (théorème~\ref{pair.t.deuxp}).   Dans  la
dernière        section,       nous        verrons        que       le
théorème~\ref{graphespar.t.sartemis} peut être vu comme un théorème de
décomposition  des graphes  d'Artémis. Nous  ferons le  point  sur les
liens  pouvant  exister entre  théorèmes  de  décomposition et  paires
d'amis.

Les  résultats  de ce  chapitre  ont été  soumis  au  {\it Journal  of
Combinatorial Theory, Series B}~\cite{nicolas:artemis}.

\begin{figure}[hb]
  \center
    \includegraphics{fig.artemis.3}\\
   \caption{Un graphe d'Artémis sans carré.\label{artemis.fig.sansc4}}
\end{figure}

\index{carré!graphe~d'Artémis~sans~---}  Avant  de  rentrer  dans  les
détails, nous donnons un aperçu des idées qui ont conduit à la preuve.
Ces idées, antérieures  au début de mes travaux,  sont dues à Frédéric
Maffray. Pour trouver une paire d'amis dans un graphe, pourquoi ne pas
considérer un plus court trou $H$  et deux sommets $u, v$ à distance 2
sur $H$,  c'est-à-dire voisins d'un  sommet $t$ de  $H$~?  Évidemment,
même dans une  classe aussi restreinte que les  graphes d'Artémis sans
carré, cette idée ne fonctionne  pas immédiatement, comme le montre le
graphe  représenté  figure~\ref{artemis.fig.sansc4}~: $\{a,b\}$  n'est
pas une paire d'amis à cause du sommet $a'$ qui est le deuxième sommet
d'un $P_4$  reliant $a$  et $b$.  Cependant,  on peut  légitimement se
plaindre  d'avoir manqué de  chance~: si  au lieu  de $H$  nous avions
choisi comme plus  court trou le trou $H'$  obtenu en remplaçant, dans
$H$, $a$ par  $a'$, alors, l'idée aurait fonctionné~:  $\{a', b\}$ est
bien une paire  d'amis de $G$.  En ce sens,  $a'$ est ``meilleur'' que
$a$.   Dans le cas  des graphes  sans carré,  Linhares et  Maffray ont
montré qu'en  fait le  voisinage de $t$  comporte deux cliques  $A$ et
$B$, $A$  étant un ensemble  de sommets ``pouvant remplacer  $a$'', et
$B$ étant un  ensemble de sommets ``pouvant remplacer  $B$''.  Sur $A$
(et aussi sur $B$), on peut alors définir une relation d'ordre partiel
(correspondant  à l'idée  déjà  évoquée de  sommet ``meilleur''  qu'un
autre).  Un  élément maximal de $A$  et un élément maximal  de $B$ (la
``meilleure'' paire) forment alors une paire d'amis de $G$.

Dans  notre  preuve,  le  sommet  $t$ est  remplacé  par  un  ensemble
anticonnexe  $T$. Dans  l'esprit du  lemme  de Roussel  et Rubio,  les
ensembles $A$ et  $B$ sont des ensembles de  sommets $T$-complets. Les
sommets de $V(H)  \setminus \{a, t, b\}$ sont  les sommets d'un chemin
appelé   $Z$  dont   une  extrémité   est  $A$-complète,   et  l'autre
$B$-complète.  Avant  la preuve  proprement dite, la  section suivante
introduit quelques notions préliminaires.

\section{Paires d'amis spéciales}

Soit  $G$ un  graphe.  On  note $C(T)$\index{C@$C(T)$}  l'ensemble des
sommets $T$-complets de $V(G)$.  Remarquons que par définition, $T$ et
$C(T)$  sont disjoints.  On  dit qu'un  ensemble non  vide $T\subseteq
V(G)$   est   \emph{intéressant}\index{intéressant~(ensemble~---)}  si
${G}[T]$ est anticonnexe  et si $G[C(T)]$ n'induit pas  une clique (ce
qui implique  d'ailleurs $|C(T)|\ge  2$ puisque le  graphe vide  et le
graphe réduit  à un sommet  sont des cliques).   Il se peut  très bien
qu'un graphe $G$ n'ait aucun  ensemble intéressant.  Dans ce cas, pour
tout  sommet  $v$, l'ensemble  $N(v)$  induit  une  clique (sans  quoi
$\{v\}$ serait  un ensemble intéressant).  Cela implique donc  que $G$
est une réunion disjointe de cliques sans aucune arête entre elles.

\index{spéciale~(paire~d'amis~---)}
\index{paire~d'amis!---~spéciale}

Rappelons   que   les   \emph{quasi-prismes}   ont   été   définis   à
section~\ref{base.ss.pps}.  Une paire d'amis $\{a,b\}$ de $G$ est dite
\emph{spéciale}   si  $G$  ne   contient  ancun   quasi-prisme  strict
d'extrémités  $a$  et  $b$  (voir  figure~\ref{artemis.fig.speciale}).
L'intérêt des  paires d'amis spéciales provient d'un  lemme facile que
voilà~:

\begin{figure}
  \center
  \includegraphics{fig.artemis.1}\\
  \rule{0cm}{1cm}\\
  \parbox{9cm}{  On a représenté  un graphe  d'artémis $G$.   La paire
  $\{a,b\}$ est une  paire d'amis de $G$ qui n'est  pas spéciale et la
  paire $\{a',b'\}$ est une paire d'amis spéciale de $G$.}
  \caption{Paires        d'amis         spéciales        et        non
  spéciales.\label{artemis.fig.speciale}}
\end{figure}

\begin{lemme}[Cf. \cite{everett.f.l.m.p.r:ep}]
  \label{artemis.l.sep}
  Soit $G$ un  graphe et $\{x,y\}$ une paire d'amis  de $G$.  Alors on
  a~:
  \begin{outcomes}
  \item
    Si $G$ ne contient pas de trou impair, alors $G/xy$ ne contient
    pas de trou impair.
  \item
    Si  $G$ ne  contient pas  d'antitrous, alors,  $G/xy$  ne contient
    aucun antitrou différent de $\overline{C_6}$.
  \item
    Si $G$ ne contient ni prismes ni quasi-prismes d'extrémités $x$ et $y$,
    alors, $G/xy$ ne contient aucun prisme.
  \item 
    Si $G$ est d'Artémis et si $\{x,y\}$ est une paire d'amis spéciale
    de $G$, alors $G/xy$ est d'Artémis.
  \end{outcomes}
\end{lemme}

\noindent Notre théorème principal est le suivant~:

\begin{theoreme}
  \label{artemis.t.main}
  Soit $G$ un graphe d'Artémis.  Si $T$ est un ensemble intéressant de
  $G$ alors $C(T)$ contient une paire d'amis spéciale de $G$.
\end{theoreme}

\noindent En admettant  ce théorème, la conjecture de  Reed et Everett
est très simple à démontrer~:

\index{Artémis~(graphe~d'---)!paire~d'amis}
\index{paire~d'amis!graphes~d'Artémis}
\begin{corollaire}
  \label{artemis.cor.main}
  Tout  graphe  d'Artémis différent  d'une  clique  possède une  paire
  d'amis spéciale. Tout graphe d'Artémis est parfaitement contractile.
\end{corollaire}

\begin{preuve}
  Si  $G$ est  une réunion  disjointe d'au  moins deux  cliques alors,
  toute paire $\{x, y\}$ de sommets non adjacents est une paire d'amis
  spéciale de $G$.  Si $G$  n'est pas une réunion disjointe de cliques
  alors il existe  au moins un ensemble intéressant,  et donc, d'après
  le  théorème~\ref{artemis.t.main}, une  paire d'amis  spéciale $\{x,
  y\}$.  Dans tous les cas,  on obtient une paire d'ami spéciale $\{x,
  y\}$.  D'après le lemme~\ref{artemis.l.sep}, $G/xy$ est d'Artémis et
  par  une  récurrence évidente,  on  voit  que  $G$ est  parfaitement
  contractile.
\end{preuve}

\section{Preuve du théorème 5.2}
\label{artemis.s.proof}

Soit  $G$  un graphe  d'Artémis  et  $T$  un ensemble  intéressant  de
$G$. Notons alors que $G$  n'est pas une réunion disjointe de cliques.
On démontre  le théorème par  récurrence sur $|V(G)|$. La  plus petite
valeur possible pour $|V(G)|$ est 3, correpondant au cas où $G$ est un
chemin  à 3  sommets,  graphe pour  lequel  le théorème  est vrai.  On
suppose donc  $|V(G)|\ge 4$,  et on suppose  que le théorème  est vrai
pour tous les graphes de taille strictement inférieure à $|V(G)|$.

Remarquons tout d'abord qu'il suffit de démontrer le théorème pour les
ensembles intéressants  maximaux pour l'inclusion.   Car si $T$  est un
ensemble intéressant de $G$, alors tout ensemble intéressant $T'$ avec
$T\subset T'$ vérifiera $C(T')\subseteq  C(T)$. Donc une paire d'amis
spéciale de $G$ dans $C(T')$ sera aussi dans $C(T)$.

Soit donc  $T$ un ensemble  intéressant maximal.  Remarquons  que pour
tout sommet $z\in V(G) \setminus (T \cup C(T))$, l'ensemble $N(z) \cap
C(T)$  induit une  clique de  $G$, car  sinon, $T\cup\{z\}$  serait un
ensemble  intéressant  plus  grand.   En  effet,  $T\cup\{z\}$  serait
anticonnexe (puisque $T$ l'est  et $z\not\in C(T)$) et $C(T\cup \{z\})
= C(T)\cap N(z)$.

Appelons                                                   \emph{chemin
sortant}\index{chemin!---~sortant~pour~$C(T)$}\index{sortant!chemin~---~pour~$C(T)$}
tout chemin $P$  reliant deux sommets non adjacents  de $C(T)$ et dont
les  sommets intérieurs  sont  dans $V(G)  \setminus  (T \cup  C(T))$.
D'après  le lemme~\ref{rr.l.wa},  il n'y  a pas  de chemin  sortant de
longueur impaire.  De  plus, si $P$ est un  chemin sortant de longueur
paire, alors sa  longueur est au moins $4$.   En effet l'unique sommet
intérieur $z$  d'un éventuel chemin  sortant de longueur 2  serait tel
que  $N(z)\cap   C(T)$  n'induit  pas  une   clique,  contredisant  la
maximalité de $T$. On distingue maintenant deux cas~:

\noindent{\bf\rule{0cm}{4ex}Cas 1~:} {\it Il n'y a aucun chemin sortant.}

\noindent\rule{0cm}{2.5ex}Soit $\{a,b\}$ une  paire d'amis spéciale du
graphe $G[C(T)]$.  Rappelons que $G[C(T)]$ n'est pas une clique ce qui
implique l'existence d'une telle paire, par hypothèse de récurrence si
$G[C(T)]$ n'est pas une  réunion disjointe de cliques, trivialement si
$G[C(T)]$ est une réunion disjointe de cliques.

Soit $P$ un chemin de $G$ d'extrémités $a$ et $b$.  Si $P$ contient un
sommet $t\in  T$ alors  $P = \bp  a \tp t  \tp b  \ep$, et $P$  est de
longueur $2$.  Si $V(P)  \cap T=\emptyset$, alors $V(P) \subset C(T)$,
car  sinon, $P$  contiendrait  un  chemin sortant.   Donc  $P$ est  de
longueur paire  et $\{a, b\}$ est  une paire d'amis de  $G$.  De plus,
s'il existe  un quasi-prisme strict  $S$ d'extrémités $a$ et  $b$ dans
$G$, alors $V(S)  \subset C(T)$ car sinon, tous les  sommets de $S$ se
trouvant sur un  chemin reliant $a$ et $b$, il  en résulterait que $S$
contient un chemin sortant.  Mais cela contredit que $\{a,b\}$ est une
paire d'amis  spéciale de $G[C(T)]$.  Finalement, $\{a,  b\}$ est bien
une paire d'amis spéciale de $G$.

\noindent{\bf\rule{0cm}{4ex}Cas 2~:} {\it Il existe un chemin sortant.}

\noindent\rule{0cm}{2.5ex}Soit $\bp \alpha \tp  z_1 \tp \cdots \tp z_n
\tp \beta \ep$  un plus court chemin sortant.   Sa longueur est $n+1$.
Notons que  d'après ce qui  précède, $n$ est  impair et $n\ge  3$.  On
note $Z= \bp z_1 \tp \cdots \tp z_n$.  On définit~:
\begin{eqnarray*}
A&=&\{v\in C(T)\mid vz_1\in E(G), vz_i\not\in E(G) \ (i=2,\ldots, n)\},\\
B&=&\{v\in C(T)\mid vz_n\in E(G), vz_i\not\in E(G) \ (i=1,\ldots, n-1)\}.
\end{eqnarray*}

Notons que $A$ est non vide,  car $\alpha\in A$, et que $A$ induit une
clique, car $A  \subseteq N(z_1) \cap C(T)$.  De  même, $B$ induit une
clique non vide. Par définition, $A  \cap B = \emptyset$.  De plus, il
n'y a  pas d'arête $uv$  avec $u\in A,  v\in B$, car alors  $\{u, z_1,
\ldots, z_n, v\}$ induirait un trou impair.

Nous  allons  maintenant  montrer par  les  lemmes~\ref{artemis.l.abi}
à~\ref{artemis.l.absep}  que  des  sommets  $a\in A$,  $b\in  B$  bien
choisis forment une paire d'amis spéciale de $G$.  À partir d'ici, les
raisonnements   sont  assez   proches  de   la  preuve   du   cas  sans
carré~\cite{linhares.maffray:evenpairsansc4}.

\begin{lemme}
  \label{artemis.l.abi}
  $C(T) \cap N(Z) \subseteq A \cup B \cup C(T\cup A\cup B)$.
\end{lemme}

\begin{preuve}
  Soit $w\in C(T)\cap N(Z)$. Il existe donc une arête $z_jw$ avec $z_j
  \in Z$ ($1\le j\le n$).   Supposons $w \notin C(T\cup A\cup B)$~: il
  existe donc  un sommet  $u\in A\cup B$  avec $uw\not\in  E(G)$. Sans
  perte de  généralité, on peut supposer  $u \in A$.  Soit  $i$ le plus
  petit entier $i$ tel que $z_iw  \in E(G)$.  Alors $\bp u \tp z_1 \tp
  \cdots \tp z_i  \tp w \ep$ est un chemin  sortant de longueur $i+1$,
  donc on doit avoir $i=n$, et $w\in B$.
\end{preuve}

\begin{lemme}
  \label{artemis.l.pq}
  Soit $P=  \bp u \tp  u' \tp \cdots  \tp v' \tp  v \ep$ un  chemin de
  longueur  impaire avec $u\in  A$ et  $v\in B$.   Alors $u'\in  A$ ou
  $v'\in B$.
\end{lemme}
\begin{preuve}
  Notons que  $P$ est  de longueur au  moins $3$  car il n'y  a aucune
  arête entre $A$ et $B$.  De plus $P$ ne contient aucun sommet de $T$
  ni aucun sommet de $C(T\cup A\cup B)$ (sinon, $P$ serait de longueur
  2).  Supposons que le lemme~\ref{artemis.l.pq} soit faux~: $u'\notin
  A$  et  $v'\notin  B$.   Nous  allons voir  que  cela  entraîne  une
  contradiction.

  \begin{claim}
    \label{artemis.c.eqct}
    Les  seules arêtes  entre  $Z$ et  $P  \cap C(T)$  sont $z_1u$  et
    $z_nv$.
  \end{claim}

  \begin{preuveclaim}
    En effet, si  $zw$ est une arête avec  $z \in Z$ et $w  \in P \cap
    C(T)$, alors,  puisque $w \notin C(T  \cup A\cup B)$  comme on l'a
    déjà  remarqué, et  d'après le  lemme~\ref{artemis.l.abi},  on a~:
    $w\in  A\cup B$.   Puisque $A$  est une  clique, le  cas  $w\in A$
    entraîne $w=u$  (et donc $zw=z_1u$  comme souhaité) ou  $w=u'$ (et
    donc $u'\in  A$, cas  exclu par hypothèse).   Le cas $w\in  B$ est
    similaire.
  \end{preuveclaim}

  On  marque les sommets  de $P$  qui ont  un voisin  dans $Z$  --- en
  particulier les sommets de  $V(P) \cap V(Z)$ sont marqués.  Appelons
  \emph{intervalle}  de $P$ tout  sous-chemin de  $P$, de  longueur au
  moins $1$,  dont les  extrémités sont marquées  et dont  les sommets
  intérieurs  ne  le  sont pas.   D'après~(\ref{artemis.c.eqct}),  les
  sommets  marqués de $V(P)  \setminus \{u,v\}$  sont tous  dans $V(P)
  \setminus C(T)$.   Puisque $u$ et $v$ sont  marqués, les intervalles
  de $P$ forment une partition  de $P$ au sens des arêtes.  Remarquons
  qu'au moins un  sommet intérieur de $P$ doit  être marqué, car sinon
  $V(Z) \cap V(P)  = \emptyset$ et $V(Z) \cup  V(P)$ induirait un trou
  impair. Donc \emph{$P$ possède  au moins deux intervalles}.  De plus,
  par le lemme de Roussel  et Rubio (lemme~\ref{rr.l.wa}), on sait que
  \emph{$P$ possède un nombre impair de $T$-arêtes.}  Donc~:

  \begin{claim}
    \label{artemis.c.intervalle}
    Il existe dans $P$ un intervalle $Q$ qui contient un nombre impair
    de $T$-arêtes.
  \end{claim}

  Notons  que comme  $P$ possède  au  moins deux  intervalles, $Q$  ne
  contient pas  à la  fois $u$  et $v$. Sans  perte de  généralité, on
  suppose $v \notin  V(Q)$.  Soient $w$ et $x$  les extrémités de $Q$.
  Appelons $w'$ (resp. $x'$), le sommet de $Q$ appartenant à $C(T)$ et
  choisi  aussi proche  que  possible  de $w$  (resp.   $x$). On  peut
  supposer que $u$, $w$, $w'$, $x'$, $x$, et $v$ apparaissent dans cet
  ordre sur  $P$.  On a  $w' \not= x'$  car $Q$ contient au  moins une
  $T$-arête.  Pour  que tout  soit bien clair,  notons que  $V(Q) \cap
  V(Z) = \emptyset$ car si un sommet $Q$ était dans $Z$ il serait, lui
  ainsi que  ses voisins dans $Q$, marqué.   Donc on aurait $Q  = wx =
  w'x'$.  Mais  alors, l'un de $w', x'$  serait à la fois  dans $Z$ et
  dans $C(T)$, deux ensembles  par définition disjoints.  Notons aussi
  que si  $w'= w$, alors ce  sommet est dans $P\cap  C(T)\cap N(Z)$ et
  donc,   d'après~(\ref{artemis.c.eqct}),  $w'=w=u$.    D'autre  part,
  $x'=x$ est impossible car $x$ serait alors dans $N(Z) \cap C(T)$, ce
  qui  d'après~(\ref{artemis.c.eqct})  entraînerait   $v  =  x$~:  une
  contradiction.  Il s'ensuit que $Q$ est de longueur au moins $2$.

  D'après la définition des  sommets marqués, il existe un sous-chemin
  $Z'$ de $Z$ d'extrémités $z_i$ et  $z_j$, tel que $z_iw$ et $z_jx$
  sont des arêtes  ($i<j$, $i=j$ ou $i>j$ sont  possible).  On choisit
  $Z'$  minimal avec  ces propriétés,  ce qui  fait que  ses éventuels
  sommets intérieurs manquent $w$ et  $x$. Par conséquent, $H = G[V(Q)
  \cup V(Z')]$ est  un trou, évidemment pair.  Notons  qu'il n'y a pas
  de sommets  de $C(T)$  dans $V(H)  \setminus V(\bp w'  \tp Q  \tp x'
  \ep)$, par définition de $w'$ et $x'$ et parce que $Z \subseteq V(G)
  \setminus(  T\cup  C(T))$.  En  fait,  les  $T$-arêtes  de $H$  sont
  exactement les $T$-arêtes de $Q$.

  Si $w'$  et $x'$ ne sont  pas des sommets consécutifs  de $P$, alors
  $H$ possède deux sommets non  adjacents de $C(T)$ et possède pourtant
  un   nombre    impair   de   $T$-arêtes~:    une   contradiction   au
  lemme~\ref{rr.l.wh}.  Donc $w'$ et $x'$ sont consécutifs sur $P$.

  Notons $k=  \max \{i,j\}$ ($k\ge 1$).   On définit un  chemin $Y$ en
  posant $Y= \bp z_{k+1}  \tp Z \tp z_n \tp v \ep$  si $k<n$ et $Y= v$
  si $k=n$.  Notons que $H\cup Y$  est connexe car $z_k$ est un sommet
  de $H$ adjacent à $Y$.  Nous  prétendons que chaque sommet $z \in Y$
  manque $w'$ et $x'$.  En effet, $v$ lui même manque ces deux sommets
  car $w',$ $x'$, $x$,  $v$ sont quatre sommets distincts apparaissant
  dans cet ordre  sur $P$~; et si  $z \in V(\bp z_{k+1} \tp  Z \tp z_n
  \ep)$, alors $z$ manque $x'$ car $x'$ n'est pas marqué.  De plus, si
  $z$ voit $w'$, alors $w'$ est marqué et $w'=w$, et donc $w'\in P\cap
  C(T)$.  Mais alors, comme $z\neq  z_1$ et $w' \neq v$, l'arête $zw'$
  contredit~(\ref{artemis.c.eqct}).   Donc il est  bien vrai  que $z$
  manque $w'$ et  $x'$.  Finalement, le triplet formé  du trou $H$, du
  chemin $Y$  et de  l'ensemble $T$ contredit  le lemme~\ref{rr.l.dpg}
  puisqu'il  n'y a pas  de sommet  de $C(T)$  dans $H  \setminus \{w',
  x'\}$.
\end{preuve}

Poursuivons  la preuve  du théorème~\ref{artemis.t.main}.   On définit
une relation  $<_A$ sur $A$ en  posant $u<_A u'$ si  et seulement s'il
existe un chemin de longueur impaire  de $u$ vers un sommet de $B$ tel
que $u'$ est le deuxième sommet  de ce chemin. Nous allons montrer que
$<_A$ est une relation d'ordre (partiel).

\begin{lemme}\label{artemis.l.antisym}
  La relation $<_A$ est antisymétrique.
\end{lemme}
\begin{preuve}
  Supposons  en  vue  d'une  contradiction qu'il  existe  des  sommets
  $u,v\in A$ tels que $u<_A v$ et $v<_A u$.  Donc il existe un chemin
  de longueur  impaire $P_u  = \bp  u_0 \tp \cdots  \tp u_p  \ep$ avec
  $u=u_0$, $v=u_1$, $u_p\in B$, $p\ge  3$, $p$ impair, et il existe un
  chemin de  longueur impaire $P_v =  \bp v_0 \tp \cdots  \tp v_q \ep$
  avec $v=v_0$,  $u=v_1$, $v_q\in B$,  $q\ge 3$, $q$ impair.   On peut
  avoir  $u_p=v_q$, et  sinon $u_p$  voit  $v_q$ puisque  $B$ est  une
  clique.

  \begin{claim}
    \label{artemis.c.nopupva}
    Aucun sommet  de $V(P_u) \cup V(P_v)  \setminus\{u,v\}$ n'est dans
    $A$.
  \end{claim}

  \begin{preuveclaim}
    En  effet,  un sommet  de  $V(P_u)  \cup V(P_v)  \setminus\{u,v\}$
    manque au moins l'un de $u$ et  $v$, et ne peut donc être dans $A$
    qui est une clique.
  \end{preuveclaim}

  \begin{claim}
    \label{artemis.c.nopupvb}
    Aucun  sommet de $V(P_u)  \cup V(P_v)  \setminus\{u_p,v_q\}$ n'est
    dans $B$.
  \end{claim}

  \begin{preuveclaim}
    En effet, puisque $B$ est une  clique, un sommet de $(P_u \cup P_v
    \setminus  \{u_p,v_q\})  \cap  B$  ne peut  qu'être  $u_{p-1}$  ou
    $v_{q-1}$.  Mais,  si $u_{p-1} \in B$  alors $\bp u_1  \tp P_u \tp
    u_{p-1} \ep$ est un chemin de longueur impaire de $A$ vers $B$, et
    le lemme~\ref{artemis.l.pq}  implique ou  bien $u_2\in A$  (ce qui
    contredit~(\ref{artemis.c.nopupva}))  ou bien  $u_{p-2}\in  B$ (ce
    qui contredit que  $B$ est une clique). Le  cas $v_{q-1}\in B$ est
    similaire.
  \end{preuveclaim}
  
  Soit  $r$ le plus  petit entier  tel que  $u_r \in  V(P_u) \setminus
  \{u_0, u_1\}$  a un voisin dans  $P_v$, et soit $s$  le plus petit
  entier tel  que $u_rv_s$ est une arête,  avec $2 \le s  \le q$.  De
  tels entiers existent puisque $v_q$ lui-même a un voisin dans $P_u$.
  Notons  que  les  sommets  $u_1,  \ldots,  u_r,  v_1,  \ldots,  v_s$
  induisent un trou $H$, et que $r$ et $s$ ont de ce fait même parité.

  \begin{claim}
    \label{artemis.c.rs}
    On peut supposer que ou bien~:
    \begin{outcomes}
    \item
      \label{artemis.c.rsoa}
      $r=p$ et $s=q$,
    \item
      \label{artemis.c.rsob}
      $r<p$, $s<q$ et $ \bp u_{r+1}  \tp P_u \tp u_p \ep = \bp v_{s+1}
      \tp P_v \tp v_q \ep$.
    \end{outcomes}
  \end{claim}
  
  \begin{preuveclaim}
    Soit $t$  le plus grand  entier tel que  $u_r v_t \in  E(G)$, avec
    $2\le s\le t\le q$.

    Si $t-s$ est pair  alors $\bp u_1 \tp P_u \tp u_r  \tp v_t \tp P_v
    \tp v_q$ est  un chemin de longueur impaire de  $A$ vers $B$.  Son
    dernier  sommet est  $u_2$, et  son avant-dernier  sommet  $w$ est
    $v_{q-1}$  (si   $t<  q$)  ou   $u_r$  (si  $t=q$).    D'après  le
    lemme~\ref{artemis.l.pq}, appliqué à ce chemin, nous avons ou bien
    $u_2\in A$ (mais cela contredit~(\ref{artemis.c.nopupva})) ou bien
    $w\in B$.  D'après~(\ref{artemis.c.nopupvb})  ce dernier cas n'est
    possible que si $w=u_r=u_p$ (et donc $t=q$). Dans ce cas, $\bp v_1
    \tp P_v \tp v_s \tp u_p  \ep$ est un chemin de longueur impaire de
    $A$  vers $B$,  et  d'après le  lemme~\ref{artemis.l.pq}, on  doit
    avoir  $v_2\in A$ (mais  cela contredit~(\ref{artemis.c.nopupva}))
    ou $v_s\in B$.  Ce dernier  cas n'est possible que si $v_s=v_q$ et
    la  conclusion~(\ref{artemis.c.rsoa})  de~(\ref{artemis.c.rs}) est
    bien satisfaite.

    Si $t-s$  est impair avec $t \ge  s+3$ alors $\bp v_1  \tp P_v \tp
    v_s \tp u_r \tp v_t \tp P_v \tp v_q \ep$ est un chemin de longueur
    impaire de  $A$ vers $B$.  Son  deuxième sommet est  $v_2$, et son
    avant-dernier sommet $w$ est $v_{q-1}$ (si $t< q$) ou $u_r$ (si $t
    =q$).  D'après  le lemme~\ref{artemis.l.pq} appliqué  à ce chemin,
    on      doit      avoir       $v_2\in      A$      (mais      cela
    contredit~(\ref{artemis.c.nopupva}))     ou     $w     \in     B$.
    D'après~(\ref{artemis.c.nopupvb})  ce dernier  cas  n'est possible
    que si  $w=u_r=u_p$ (et  donc $r=p$ et  $t=q$).  Donc $r$  et $t$
    sont impairs, ce qui est  impossible car $t-s$ est impair et $r-s$
    est pair.

    Le seul cas restant est $t=s+1$,  ce qui entraîne que $\bp u_0 \tp
    P_u \tp  u_r \tp  v_{s+1} \tp P_v  \tp v_q  \ep$ est un  chemin de
    longueur impaire, qui peut jouer le rôle de $P_u$.  Cette fois, la
    conclusion~(\ref{artemis.c.rsob})    de~(\ref{artemis.c.rs})   est
    satisfaite.
  \end{preuveclaim}

  On  appelle \emph{intervalle}  tout chemin  de $H$,  de  longueur au
  moins  $1$, dont  les extrémités  voient $z_1$  et dont  les sommets
  intérieurs manquent $z_1$.
  
  \begin{claim}
    \label{artemis.c.H3}
    Les intervalles de $H$ sont au nombre de 3 au moins et ils forment
    une partition de $H$ au sens des arêtes.
  \end{claim}

  \begin{preuveclaim}
    Soit   $P$    un   chemin   défini    de   la   sorte~:    Si   la
    conclusion~(\ref{artemis.c.rsoa})    de~(\ref{artemis.c.rs})   est
    satisfaite,  on pose  $P=  \bp z_2  \tp  Z \tp  z_n  \ep$.  Si  la
    conclusion~(\ref{artemis.c.rsob})    de~(\ref{artemis.c.rs})   est
    satisfaite,  on définit  $P$  comme étant  un chemin  d'extrémités
    $z_2$  et  $u_{r+1}$,  contenu   dans  $\{z_2,  \dots,  z_n\}  \cup
    \{u_{r+1}, u_{r+2},  \dots, u_p\}$  ($P$ existe car  l'ensemble où
    ses sommets  sont tenus de  se trouver est connexe).   On remarque
    que  le lemme~\ref{rr.l.dpg}  peut s'appliquer  au triplet  $H, P,
    \{z_1\}$~: $H$ est un trou, $P$ est un chemin, $\{z_1\}$ induit un
    graphe  anticonnexe  disjoint  de  $H\cup  P$,  $G[H\cup  P]$  est
    connexe,  tout sommet  $z\in V(P)$  manque  $u$ et  $v$ (si  $z\in
    Z[z_2, z_n]$  parce que $u, v  \in A$, si $z\in  V(\bp u_{r+1} \tp
    P_u \tp u_p \ep)$ parce que $P_u$ et $P_v$ sont des chemins et $r,
    s \ge 2$), et enfin $u,  v, z_2$ sont dans $C(\{z_1\})$.  Donc le
    lemme~\ref{rr.l.dpg}  montre qu'un  sommet de  $V(H) \setminus\{u,
    v\}$ est  dans $C(\{z_1\})=N(z_1)$.  Donc $z_1$  possède au moins
    $3$ voisins dans $H$.
  \end{preuveclaim}

  \begin{claim}
    \label{artemis.c.HT}
    $H$ possède un nombre pair de $T$-arêtes.
  \end{claim}

  \begin{preuveclaim}
    Rappelons   que   $u$   et   $v$   sont  dans   $C(T)$.    Si   la
    conclusion~(\ref{artemis.c.rsoa})    de~(\ref{artemis.c.rs})   est
    satisfaite,  alors  $u_r$  et  $v_s$  sont  dans  $C(T)$.   Si  la
    conclusion~(\ref{artemis.c.rsob})    de~(\ref{artemis.c.rs})   est
    satisfaite, alors on applique le lemme~\ref{rr.l.sgt} au trou $H$,
    au chemin $\bp  u_{r+1} \tp P_u \tp u_p \ep$  et à l'ensemble $T$.
    À cause  des arêtes  $u v$ et  $u_r v_s$,  le lemme~\ref{rr.l.sgt}
    implique  qu'au moins  l'un de  $u_r$  et $v_s$  est dans  $C(T)$.
    Donc,      que      la     conclusion      de~(\ref{artemis.c.rs})
    soit~(\ref{artemis.c.rsoa})   ou~(\ref{artemis.c.rsob}),  on  peut
    affirmer que $H$  contient au moins trois sommets  de $C(T)$.  Par
    le lemme~\ref{rr.l.wh}, on obtient la conclusion souhaitée.
  \end{preuveclaim}

  \begin{figure}
    \begin{center}
      \includegraphics{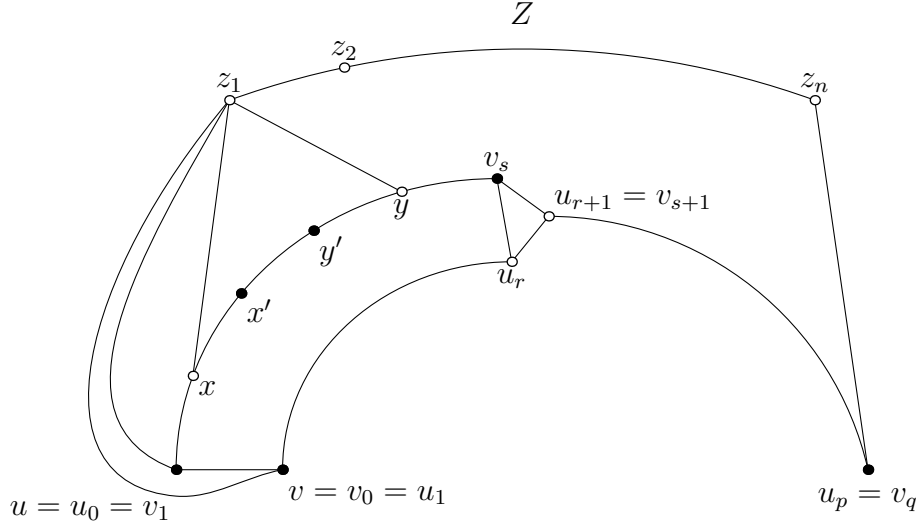}
    \caption{Lemme~\ref{artemis.l.antisym}, sommets $x$, $x'$, $y'$ et
    $y$.       Les       sommets       ``pleins''      sont       dans
    $C(T)$.\label{artemis.fig.antisym}}
    \end{center}
  \end{figure}

  À      partir       d'ici,      on      peut       consulter      la
  figure~\ref{artemis.fig.antisym}.  Remarquons que  $\bp u \tp v \ep$
  est   un    intervalle   qui   contient    une   $T$-arête.    Donc,
  d'après~(\ref{artemis.c.H3})  et~(\ref{artemis.c.HT}),  il existe  un
  intervalle $Q$ de $H$, différent de  $\bp u \tp v \ep$, qui contient
  un  nombre  impair  de  $T$-arêtes.   On  appelle  $x$  et  $y$  les
  extrémités de $Q$,  et on appelle $x'$ le sommet  de $Q\cap C(T)$ le
  plus proche de $x$ et $y'$  le sommet de $Q\cap C(T)$ le plus proche
  de $y$,  de manière que $u, x,  x', y', y, v$  apparaissent dans cet
  ordre sur $H$.  Puisque $H$  contient au moins trois intervalles, on
  peut supposer  qu'au moins l'un de  $x$ et $y$, par  exemple et sans
  perte  de généralité  $y$,  est différent  de  $u$ et  de $v$.   Par
  conséquent, $y$ manque l'un de $u$ et $v$.

  Si $Q$ est de longueur $1$, on a  $Q= \bp x \tp y \ep =\bp x' \tp y'
  \ep$,  et $y  \in N(z_1)  \cap C(T)$.   Le lemme~\ref{artemis.l.abi}
  implique $y \in  A \cup B \cup C(T  \cup A \cup B)$.  Mais  $y \in A
  \cup C(T\cup A\cup B)$ est impossible  car $y$ manque l'un de $u$ et
  $v$, et $y\in B$ est impossible car $y \in N(z_1$).  Donc $Q$ est de
  longueur  au moins  $2$, $V(Q)  \cup \{z_1\}$  induit un  trou $H_1$
  contenant un nombre impair  de $T$-arêtes, et le lemme~\ref{rr.l.wh}
  appliqué à la paire $H_1, T$ implique  que $\bp x' \tp Q \tp y' \ep$
  est de  longueur $1$. Donc $x'$  et $y'$ sont les  seuls sommets de
  $C(T)$ dans $H_1$.

  \begin{figure}
    \begin{center}
      \includegraphics{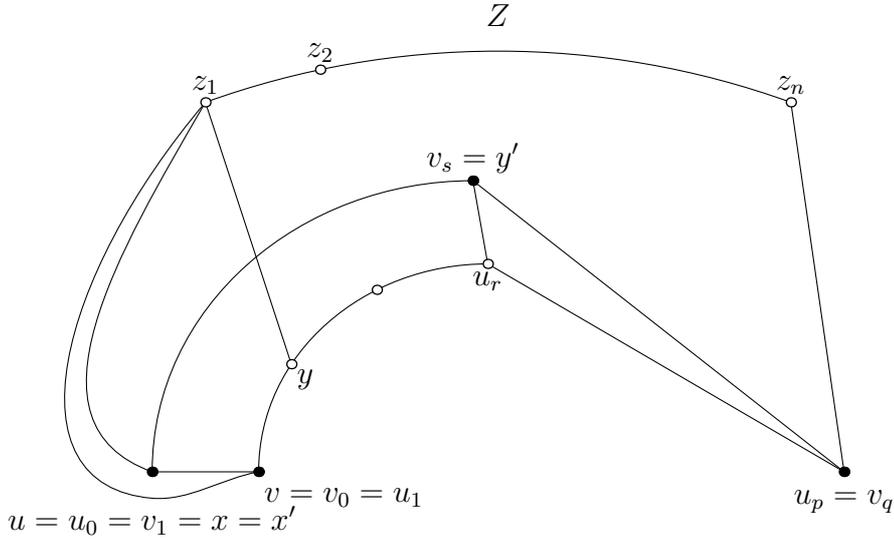}
      \caption{Lemme~\ref{artemis.l.antisym},  fin de  la  preuve. Les
	sommets             ``pleins''            sont            dans
	$C(T)$.\label{artemis.fig.antisym2}}
    \end{center}
  \end{figure}

  Supposons $x\neq u$. Le lemme~\ref{rr.l.dpg} appliqué au trou $H_1$,
  au chemin réduit à $u$ et  à l'ensemble $T$ montre que $u$ doit voir
  $x'$, ce qui implique $x=x'=v_2$.  On a alors $x\in N(Z) \cap C(T)$,
  ce qui contredit le  lemme~\ref{artemis.l.abi}. Donc $x=u$ et comme
  $x$  est   dans  $C(T)$,   $x'=x=u$.   On  applique   maintenant  le
  lemme~\ref{rr.l.dpg} au trou $H_1$, au chemin $\bp u_p \tp z_n \tp Z
  \tp  z_2 \ep$, et  à l'ensemble  $T$.  Ceci  montre qu'il  existe un
  sommet de $\bp u_p \tp z_n \tp  Z \tp z_2 \ep$ qui voit l'un de $x'$
  et $y'$. Ce  sommet ne peut pas être un $z_i,  (i=2, \dots, n$), car
  ces $z_i$ manquent $x'=x=u$ par définition de $Z$ et manquent $y'$ à
  cause du lemme~\ref{artemis.l.abi}. Ce  sommet ne peut donc être que
  $u_p$,  ce qui  entraîne  $s=2$  et $u_p$  voit  $v_s$.  Cela  n'est
  possible que si  $s=q-1$ ou $s=q$. Mais $s=q$  est impossible car on
  aurait alors $x'\in  A$, $y'\in B$ et $x'y'$  serait une arête entre
  un  sommet de  $A$ et  un sommet  de $B$~:  une  contradiction. Donc
  $s=q-1$.  À  partir de maintenant,  la situation ressemble  plutôt à
  celle représentée  figure~\ref{artemis.fig.antisym2}. Le chemin $\bp
  u \tp z_1 \tp y \tp P_u \tp u_r \ep$ est de longueur paire car $H_1$
  est un trou pair.  Donc le chemin  $\bp u \tp z_1 \tp y \tp P_u \tp
  u_r \tp  u_p \ep$  est un chemin  de longueur impaire,  reliant deux
  sommets   de  $C(T)$   et  ne   comportant  aucune   $T$-arête~:  une
  contradiction au lemme~\ref{rr.l.wa}.
\end{preuve}

\begin{lemme}
  La relation $<_A$ est transitive.
\end{lemme}

\begin{preuve}  
  Soient $u,  v, w$ trois  sommets de  $A$ tels que  $u <_A v  <_A w$.
  Puisque $v  <_A w$, il  existe un  chemin $P =  \bp v_0 \tp  v_1 \tp
  \cdots  \tp v_q  \ep$ avec  $v_0=v$,  $v_1=w$, $v_q\in  B$, $q$  est
  impair, $q\ge 3$.

  Si $u$ n'a pas de voisin sur  $\bp v_2 \tp P \tp v_q \ep$ alors $\bp
  u \tp  v_1 \tp  P \tp  v_1 \tp v_q  \ep$ est  un chemin  de longueur
  impaire  de  $u$  vers  $v_q$   qui  montre  que  $u  <_A  w$  comme
  souhaité. On peut donc supposer  que $u$ possède un voisin $v_i$ sur
  $\bp  v_2 \tp  P  \tp v_q  \ep$, et  on  choisit $i$  le plus  grand
  possible avec  cette propriété ($2\le i\le  q$).  On a  $i<q$ car il
  n'y a pas d'arêtes entre $A$ et $B$.

  Si $i$ est  impair ($3 \le i \le  q-2$), alors $\bp u \tp  v_i \tp P
  \tp  v_q \ep$  est un  chemin impair  de $A$  vers $B$.   D'après le
  lemme~\ref{artemis.l.pq}, nous avons alors $v_i\in A$ ou $v_{q-1}\in
  B$.  Le premier de ces cas  est impossible car $A$ est une clique et
  donc $v_{q-1}\in B$.   Mais alors, $\bp v_0 \tp u \tp  v_i \tp P \tp
  v_{q-1} \ep$ est un chemin de  longueur impaire de $v$ vers $B$, qui
  montre   $v   <_A   u$,   contredisant   l'antisymétrie   de   $<_A$
  (lemme~\ref{artemis.l.antisym}).

  Si $i$ est pair ($2\le i \le q-1$)  alors $\bp v \tp u \tp v_i \tp P
  \tp v_q$ est un chemin de  longueur impaire de $v$ vers un sommet de
  $B$ qui  montre $v <_A  u$ contredisant à nouveau  l'antisymétrie de
  $<_A$ (lemme~\ref{artemis.l.antisym}).
\end{preuve}

Les deux  lemmes précédents  montrent bien que  $<_A$ est  une relation
d'ordre (partiel). On définit de  même une relation sur $B$~: si $u,
u' \in B$, on pose $u <_B u'$ si et seulement s'il existe un chemin de
longueur impaire  de $u$  vers un sommet  de $A$  tel que $u'$  est le
deuxième sommet  de ce  chemin.  Les ensembles  $A$ et $B$  jouant des
rôles  symétriques,  on voit  que  la  relation  $<_B$ est  aussi  une
relation d'ordre.  Notons que $A$ et $B$ étant des ensembles finis, on
sait qu'il existe dans $A$ et dans $B$ des éléments maximaux.

\begin{lemme}
  \label{artemis.l.absep}
  Soit $a$ un sommet maximal de $(A, <_A)$ et $b$ un sommet maximal de
  $(B, <_B)$.  Alors $\{a,b\}$ est une paire d'amis spéciale de $G$.
\end{lemme}

\begin{preuve}
  Supposons qu'il existe  un chemin de longueur impaire  $Q= \bp a \tp
  a' \tp \cdots \tp b' \tp b  \ep$. Ce chemin est de longueur au moins
  3   car   il  n'y   a   aucune  arête   entre   $A$   et  $B$.    Le
  lemme~\ref{artemis.l.pq} montre que $a'\in A$ (et donc $a<_A a'$) ou
  $b'\in  B$  (et  donc  $b<_B  b'$),  dans tous  les  cas  on  a  une
  contradiction à  la maximalité  de $a$ et  $b$.  Donc  $\{a,b\}$ est
  bien une paire d'amis de $G$.

  Supposons  que $\{a,b\}$  ne soit  pas spéciale,  c'est-à-dire qu'il
  existe un quasi-prisme  strict $S$ d'extrémités $a$ et  $b$ (on note
  les   sommets  et   les  chemins   du  quasi-prisme   comme  indiqué
  figure~\ref{rr.fig.quasi-prisme}).   On  pourra  se  reporter  à  la
  figure~\ref{artemis.fig.final}  page~\pageref{artemis.fig.final} qui
  représente certains  des sommets que  nous définirons au long  de la
  preuve qui  suit.  Par  symétrie on peut  supposer que $S_1$  est de
  longueur au  moins $1$, et  on appelle $a_1$  le voisin de  $a'$ sur
  $S_1$.  Puisqu'aucun  sommet de $S$  ne voit à  la fois $a$  et $b$,
  aucun sommet de  $S$ n'est dans $T$.  Notons que $\bp  a \tp S_1 \tp
  a' \tp c \tp S_3 \tp c' \tp b'  \tp S_2 \tp b \ep$ et $\bp a \tp S_1
  \tp a'  \tp d  \tp S_4 \tp  d' \tp b'  \tp S_2  \tp b \ep$  sont des
  chemins  de longueur  paire car  $\{a,  b\}$ est  une paire  d'amis.
  Appelons $H$ le trou induit par $V(S_3) \cup V(S_4)$.

  \begin{claim}
    \label{artemis.c.acdcd}
    Aucun sommet de $S \setminus\{a,b\}$ n'est dans $N(Z) \cap C(T)$.
  \end{claim}

  \begin{preuveclaim}
    Supposons qu'il existe  un sommet $u\in V(S)\setminus\{a,b\}$ dans
    $N(Z)  \cap C(T)$.   D'après  le lemme~\ref{artemis.l.abi},  $u\in
    A\cup B\cup C(T \cup A\cup B)$.  Mais $u \in C(T\cup A\cup B)$ est
    impossible car aucun  sommet de $S$ ne voit à la  fois $a$ et $b$.
    Donc $u \in A  \cup B$.  Si $u\in A$, $u$ doit  être le voisin de
    $a$  sur  $S_1$  (car  $A$  induit une  clique).   Donc  $(V(S_1)
    \setminus\{a\})  \cup  V(S_3) \cup  V(S_2)$  induit  un chemin  de
    longueur  impaire  $P_u$  de  $u\in  A$ vers  $b$  (rappelons  que
    $\{a,b\}$ est  une paire  d'amis).  Puisque le  voisin de  $u$ sur
    $P_u$  n'est  pas  dans  $A$  (car  $A$  induit  une  clique),  le
    lemme~\ref{artemis.l.pq} implique  que le voisin de  $b$ sur $P_u$
    est  dans $B$.   Mais cela  contredit  la maximalité  de $b$  pour
    l'ordre $<_B$.   Si $u\in B$, $u$  doit être le voisin  de $b$ sur
    $S_2$ (si $S_2$  est de longueur au moins $1$)  ou l'un de $c',d'$
    (si $S_2$ est de longueur $0$), mais dans l'un ou l'autre des cas,
    un  argument similaire  à  la  situation $u\in  A$  conduit à  une
    contradiction.
  \end{preuveclaim}

  D'après le lemme~\ref{rr.l.sgp} appliqué à $S$ et $T$, on sait que~:

  \begin{claim}
     \label{artemis.c.ald}
    Au moins deux des sommets $a',  c, d$ et deux des sommets $b', c',
    d'$ sont dans $C(T)$.
  \end{claim}

  Puisque  $a  \neq a'$  (souvenons  nous  que  $b=b'$ est  possible),
  de~(\ref{artemis.c.acdcd}) et~(\ref{artemis.c.ald})  on peut déduire
  les deux fait suivants~:

  \begin{claim}
    \label{artemis.c.ctnz1}
    Si  l'un  de   $a',c,d$  est  dans  $N(Z)$,  alors   il  est  dans
    $N(Z)\setminus C(T)$  et les deux autres  sont dans $C(T)\setminus
    N(Z)$.
  \end{claim}

  \begin{claim}
    \label{artemis.c.ctnz2}
    Si  l'un  de  $c',  d'$   est  dans  $N(Z)$,  alors  il  est  dans
    $N(Z)\setminus  C(T)$, l'autre est  dans $C(T)\setminus  N(Z)$, et
    $b'\in C(T)$.
  \end{claim}

  Notons  que  $a' \notin  Z$  car sinon  l'un  de  $c,d$ serait  dans
  $N(Z)\cap C(T)$.

  On définit un chemin $P$ comme suit~: soit $a''$ le sommet de $N(Z)$
  le plus  proche de $a_1$ sur  $\bp a \tp  S_1 \tp a_1 \ep$,  et soit
  $b''$ le  sommet de  $N(Z)$ le  plus proche de  $b'$ sur  $S_2$ (les
  sommets $a''$ et $b''$ existent à cause de $a$ et $b$).  Soient $z_i
  \in V(Z)  \cap N(a'')$ et $z_j\  \in V(Z) \cap N(b'')$  tels que $\bp
  z_i \tp  Z \tp z_j  \ep$ est aussi  court que possible  (notons que
  $i<j$, $i=j$ ou  $i>j$ sont possibles).  Posons $P=  \bp a_1 \tp S_1
  \tp a''  \tp z_i  \tp Z \tp  z_j \tp  b'' \tp S_2  \tp b'  \ep$.  On
  vérifie que  $P$ est un chemin  de $G$ d'extrémités $a_1$  et $ b'$.
  D'après le lemme~\ref{rr.l.dpg} appliqué à $H$ , $P$ et $\{a'\}$, il
  existe un  sommet $z$  de $P$  qui voit l'un  de $c$  ou $d$,  et on
  choisit  un  tel  $z$  aussi   proche  que  possible  de  $a_1$  sur
  $P$. Quitte à échanger $S_3$ et  $S_4$, on suppose que $z$ voit $c$.
  La définition  des quasi-prismes entraîne  $z \not\in V(\bp  a_1 \tp
  S_1 \tp  a'') \cup V(S_2)$,  et donc $z\in  V(\bp z_i \tp Z  \tp z_j
  \ep)$.  D'après~(\ref{artemis.c.ctnz1}),  on a $c\in  N(z) \setminus
  C(T)$  et $a', d\in  C(T) \setminus  N(Z)$.  Donc  $a'$ n'a  pas de
  voisins  dans $V(P) \setminus  \{a_1\}$, car  un tel  voisin devrait
  obligatoirement se  trouver dans  $\bp z_i \tp  Z \tp z_j  \ep$ sous
  peine  de  contredire  la   définition  des  quasi-prismes,  ce  qui
  entraînerait $a'\in N(Z)$~: une contradiction. Autrement dit~:

  \begin{claim}
    \label{artemis.c.pa1}
    $\bp a' \tp a_1 \tp P \tp b' \ep$ est un chemin de $G$.
  \end{claim}

  On note $b_1$ le voisin de $b'$ sur $P$ ($b_1$ peut être $z_j$ ou un
  sommet de  $S_2$).  D'après le lemme~\ref{rr.l.dpg}  appliqué à $H$,
  $\bp b_1 \tp  P \tp z \ep$  et $\{b'\}$, il existe un  sommet $y$ du
  chemin $\bp b_1  \tp P \tp z  \ep$ qui voit l'un de  $c', d'$.  Sous
  peine de contredire la définition des quasi-prismes, on a $y \in \bp
  z  \tp  P   \tp  z_j  \ep  =   \bp  z  \tp  Z  \tp   z_j  \ep$  Donc,
  d'après~(\ref{artemis.c.ctnz2}), on sait que~:

  \begin{claim}
    \label{artemis.c.cpdpt}
    Exactement l'un de $c',d'$ est dans $C(T) \setminus N(Z)$, l'autre
    est dans $N(Z)\setminus C(T)$, et $b'\in C(T)$.
  \end{claim}

  \noindent De plus~:
  
  \begin{claim}
    \label{artemis.c.ccnz}
    $cc'\notin E(G)$.
  \end{claim}

  \begin{preuveclaim}
    Supposons $cc'\in E(G)$.  Si $c'\in  C(T)$ alors $\bp a' \tp c \tp
    c' \ep$  est un chemin  sortant de longueur 2~:  une contradiction.
    Si  $c'\not\in C(T)$, alors  d'après~(\ref{artemis.c.cpdpt}), $\bp
    a' \tp  c \tp  c' \tp b'  \ep$ est  un chemin sortant  de longueur
    impaire~: une contradiction.
  \end{preuveclaim}

  Appelons $H_1$ le cycle induit par  $V(\bp a_1 \tp P \tp z \ep) \cup
  \{a',    c\}$,     qui    possède    au     moins    $4$    sommets.
  D'après~(\ref{artemis.c.pa1}) et d'après la définition de $z$, $H_1$
  est un trou.  On pose $S'_4 = \bp d \tp S_4 \tp d' \tp c'$.

  \begin{claim}
    \label{artemis.c.nowx}
    Il n'y aucune arête entre $\bp a_1 \tp P \tp z \ep$ et $S'_4$.
  \end{claim}
  
  \begin{preuveclaim}
    Supposons qu'il existe  une arête $xw$ avec $x\in  V(\bp a_1 \tp P
    \tp z  \ep)$ et $w  \in V(S'_4)$.   On a $w  \neq d$ car  $d \notin
    N(Z)$ et $d$ manque tous  les sommet de $V(S_1) \setminus \{a'\}$.
    Aucun sommet de $V(S'_4) \setminus \{d\}$ n'est adjacent à $a'$ ou
    $c$   à   cause   de    la   définition   des   quasi-prismes   et
    de~(\ref{artemis.c.ccnz}).   Mais  alors  le  lemme~\ref{rr.l.dpg}
    appliqué  à $H_1$,  $S'_4 \setminus\{d\}$  et $\{d\}$  entraîne une
    contradiction.
  \end{preuveclaim}

  On sait qu'il existe un sommet $z'$ de $\bp z \tp P \tp z_j \ep$ qui
  a un voisin  $\delta \in V(S'_4)$ (à cause de  $y$). On choisit $z'$
  aussi proche que possible de $z$ sur $\bp z \tp Z \tp z_j \ep$ et on
  choisit $\delta\in V(S'_4) \cap  N(z')$ aussi proche que possible de
  $d$ sur $S'_4$.   Par ces choix,  $\bp z \tp Z \tp  z' \tp \delta
  \tp S'_4  \tp d$ est  un chemin  de $G$.  On  a $\delta \neq  d$ car
  $z'\in Z$ et $d\notin N(Z)$.   Considérons le cycle $H_2$ induit par
  $V(\bp a_1 \tp P \tp z' \tp \delta \tp S'_4 \tp d \ep ) \cup\{a'\}$,
  qui possède au  moins $4$ sommets.  Supposons que  $H_2$ possède une
  corde.  Les définitions de $S, P, z', \delta$ et le fait que $a' \tp
  a_1 \tp P \tp  b'$ et $\bp z \tp P \tp z' \tp  \delta \tp S'_4 \tp d
  \ep$  soient des  chemins de  $G$  impliquent que  les seules  cordes
  possibles dans $H_2$ sont des arêtes de type $wx$ avec $w\in V(S_4)$
  et  $x \in V(z_i  \tp Z  \tp z  \ep)$.  Mais  de telles  arêtes sont
  interdites d'après~(\ref{artemis.c.nowx}).  Donc  $H_2$ est un trou
  pair.  Le sommet  $c$ possède trois voisins sur  $H_2$ ($a'$, $d$ et
  $z$).  Donc, d'après le  lemme~\ref{rr.l.wh}, $H_2$ possède un nombre
  pair de $c$-arêtes, l'une d'elle étant $a'd$.  Manifestement, il n'y
  a pas  de $c$-arête sur $\bp d  \tp S'_4 \tp \delta  \ep$.  De plus,
  $\delta z'$  ne peut pas  être une $c$-arête, car  cela impliquerait
  $\delta    =    c'$   et    $cc'\in    E(G)$~:   une    contradiction
  avec~(\ref{artemis.c.ccnz}).    Donc   les   $c$-arêtes  de   $H_2$
  différentes de $a'd$ (et il y  en a un nombre impair) sont dans $\bp
  z \tp P \tp z' \ep$.  Appelons $z''$ le voisin de $c$ le plus proche
  de $z'$ sur $\bp z \tp P  \tp z' \ep$, de sorte que les $c$-arêtes
  de $H_2$  différentes de $a'd$ se trouvent  dans $\bp z \tp  P \tp z''
  \ep$.  D'après le  lemme~\ref{rr.l.wa} appliqué à $z \tp  P \tp z''$
  et $\{c\}$, on obtient~:

  \begin{claim}
    $\bp z \tp P \tp z'' \ep$ est de longueur impaire.
  \end{claim}

  D'après le lemme~\ref{rr.l.wh}, le nombre de $T$-arêtes de $H_1$ est
  pair ou  égal à $1$. Dans  ce dernier cas, les  sommets $V(H_1) \cap
  C(T)$  sont $a'$  et  $a_1$  puisque $c  \notin  C(T)$.  D'après  le
  lemme~\ref{rr.l.dpg} appliqué  à $H_1$,  $(\bp b' \tp  P \tp  z \ep)
  \setminus \{z\}$  et $T$, certains sommets  de $(\bp b' \tp  P \tp z
  \ep)  \setminus \{z\}$  voient l'un  de  $a'$ et  $a_1$.  Mais  cela
  contredit~(\ref{artemis.c.pa1}).  Donc~:

  \begin{claim}
    $H_1$ possède un nombre pair de $T$-arêtes.
  \end{claim}

  \begin{figure}[htb]
    \begin{center}
       \includegraphics{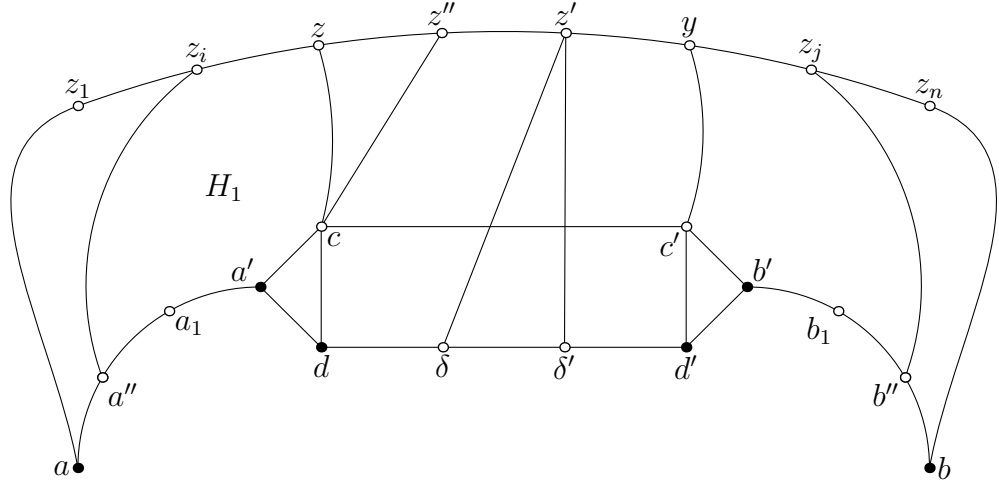}\\
    \end{center}
    \caption{Lemme~\protect{\ref{artemis.l.absep}} dans  le cas $d'\in
      C(T)$.       Les      sommets      ``pleins''     sont      dans
      $C(T)$.\label{artemis.fig.final}}
  \end{figure}

  D'après~(\ref{artemis.c.cpdpt}), l'un de  $c', d'$ est dans $C(T)$~:
  on le note  $\gamma'$. On définit un chemin  $S''_4$ comme suit.  Si
  $\delta = c'$ (ce qui entraîne $\gamma'=d'$), on pose $S''_4= \bp c'
  \tp  d' \ep$  et $\delta'  = c'$.   Si $\delta  \in S_4$  alors soit
  $\delta'$ le  voisin de  $z'$ sur  $S_4$ qui est  le plus  proche de
  $d'$. On pose alors $S''_4 =  \bp \delta' \tp S'_4 \tp \gamma' \ep$.
  On pose  $Q_1 = \bp a'  \tp c \tp z''  \tp P \tp z'  \tp \delta' \tp
  S''_4 \tp  \gamma' \ep$ et $Q_2  = \bp a' \tp  a_1 \tp P  \tp z' \tp
  \delta'  \tp S''_4  \tp \gamma'  \ep$.   Ces chemins  sont bien  des
  chemins         de         $G$        d'après~(\ref{artemis.c.pa1}),
  (\ref{artemis.c.ccnz}),   (\ref{artemis.c.nowx})   et  d'après   les
  définitions de $z'$, $z''$ et $\delta'$.  Leurs extrémités sont $a'$
  et $\gamma'$, qui sont tous  les deux dans $C(T)$. Puisque $H_1$ est
  un trou pair  et puisque $\bp z  \tp P \tp z'' \ep$  est de longueur
  impaire, les  chemins $Q_1$ et $  Q_2$ ont des  longueurs de parités
  différentes.  Le  nombre de  $T$-arêtes de $Q_1$  et celui  de $Q_2$
  sont de même parité, car  $H_1$ possède un nombre pair de $T$-arêtes
  dont $cz$ et $a'c$ ne font  pas partie (car $c \notin C(T)$), et car
  $\bp z \tp P \tp z'' \ep$  ne contient pas de $T$-arêtes (car $Z$ ne
  contient pas de  sommets de $C(T)$). Donc l'un  des chemins $Q_1$ et
  $Q_2$ contredit le lemme~\ref{rr.l.wa}.
\end{preuve}

\noindent La conclusion du lemme~\ref{artemis.l.absep} montre bien que
le théorème~\ref{artemis.t.main} est vrai.

\section{Conséquences algorithmiques}

Nous allons donner ici une série d'algorithmes fondés sur la preuve du
théorème~\ref{artemis.t.main}.   Le dernier  d'entre eux  permettra de
colorier les graphes d'Artémis  en temps polynomial~: $O(n^2m)$. Notre
version originale de l'algorithme avait pour complexité $O(n^6)$, mais
Bruce Reed  a remarqué de  nombreuses améliorations.  Nous  lui devons
notamment l'analyse  ``amortie'' de l'algorithme,  c'est-à-dire l'idée
de  comptabiliser   globalement  tous  les   appels  récursifs.   Nous
rappelons que nous donnerons également un algorithme de reconnaissance
des  graphes d'Artémis  au  chapitre~\ref{reco.chap}.  Cet  algorithme
n'utilisera  pas   l'existence  d'une  paire   d'amis  mais  détectera
directement les prismes dans les graphes de Berge.

\begin{algorithme}
  \label{artemis.a.t}
  \begin{itemize}
  \item[\sc Entrée~:] 
    Un graphe $G$.
  \item[\sc Sortie~:]  Un ensemble intéressant maximal de  $G$, ou une
    paire d'amis spéciale de $G$, ou ``$G$ est une clique''.
  \item[\sc Calcul~:] \mbox{}
    \begin{itemize}
      \item[Phase  1~:  recherche  d'un  sommet  $t$  non  simplicial.]
	\mbox{}

	Si  $G$ n'est  pas  connexe, stopper  et  retourner une  paire
	$\{a,b\}$  avec  $a$  et   $b$  dans  2  composantes  connexes
	distinctes de  $G$. Si $G$  est connexe, calculer le  degré de
	chaque  sommet.  Si  tous  les sommets  sont  de degré  $n-1$,
	stopper et retourner ``$G$ est une clique''. Sinon, considérer
	un sommet $u$ de degré strictement inférieur à $n-1$.  Marquer
	tous  les  voisins  de   $u$.   Parmi  les  sommets  marqués,
	rechercher un sommet ayant un voisin non marqué. Nommer $t$ un
	tel sommet.

      \item[Phase 2~: Initialisation des marques.] \mbox{}

	Donner à $t$ la marque ``$T$'', donner aux voisins de $t$ la
	marque ``$C(T)$''. Considérer tous les autres sommets comme
	non-marqués.

      \item[Phase 3~: recherche de l'ensemble intéressant maximal $T$.]
      \mbox{}

      Tant qu'il existe un sommet $u$ non marqué faire~:

      Si  $N(u) \cap  C(T)$ est  une clique,  donner à  $u$  la marque
      ``clique''. Si $N(u)  \cap C(T)$ n'est pas une  clique, donner à
      $u$  la marque  ``$T$''. Pour  chaque non-voisin  de  $u$ marqué
      ``C(T)'',  effacer la  marque ``C(T)''  et considérer  $u$ comme
      non-marqué. Fin du ``faire''.

      Stopper et retourner l'ensemble des sommets marqués ``$T$''.
    \end{itemize}

  \item[\sc Complexité~:] $O(\max (n+m,m(n-k)))$ où $k$ est le nombre de
    sommets de $C(T)$ pour  l'ensemble $T$ retourné (si aucun ensemble
    $T$ n'est retourné, on adopte la convention $k=n$ et la complexité
    est alors $O(n+m)$).
  \end{itemize}
\end{algorithme}

\begin{figure}[htb]
  \begin{center}
    \includegraphics{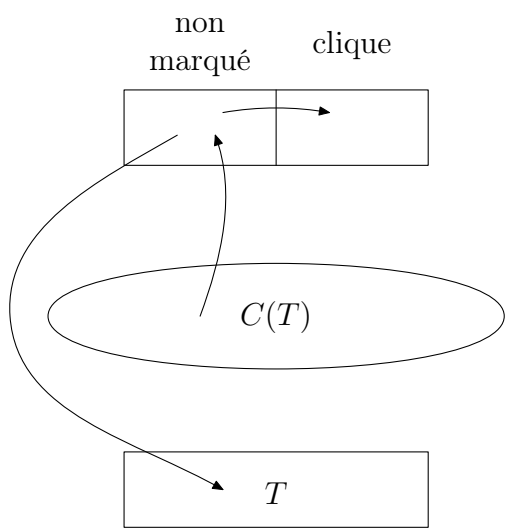}\\
  \end{center}
  \caption{Algorithme~\protect{\ref{artemis.a.t}}~:   changements   de
  marques  possibles  pour  les  sommets  lors  de  l'augmentation  de
  $T$.\label{artemis.fig.trans}}
\end{figure}

\begin{preuve}
  Il est  clair que notre algorithme est  conforme aux spécifications.
  Reste à analyser sa complexité.

  La phase~1 nécessite au  pire $n+m$ étapes  de calcul.  La  phase~2
  également. Pour la phase~3,  chaque itération de la boucle nécessite
  au plus $m$  pas de calculs.  En effet, décider  si $N(u) \cap C(T)$
  est une  clique nécessite au plus  $m$ pas (c'est  un processsus qui
  s'arête dès qu'on rencontre  une non-arête).  L'effaçage des marques
  ``$C(T)$''  pose problème  dans la  mesure où  on travaille  sur des
  non-arêtes. Mais on  peut remarquer que le nombre  de pas est majoré
  par le degré du sommet $t$ initial car l'ensemble $C(T)$ ne fait que
  décroître à partir de sa valeur initiale $N(t)$.

  Donc, chaque itération ``tant que'' nécessite au plus $m$ étapes. La
  figure~\ref{artemis.fig.trans} montre  à l'aide de  flèches la façon
  dont les sommets peuvent changer de marque. Il apparaît qu'un sommet
  peut être  au plus  une fois  non marqué. De  plus, les  $k$ sommets
  restant à la fin dans $C(T)$  ne sont à aucun moment non-marqués, et
  n'ont  de ce fait  jamais compté  pour une  itération de  la boucle.
  Donc, la  boucle nécessite bien  au pire $n-k$ itérations.   D'où la
  complexité annoncée.
\end{preuve}

\begin{algorithme}
  \label{artemis.a.zon}
  \begin{itemize}
    \item[\sc  Entrée~:]  Un  graphe  connexe  $G$  et  $T$,  ensemble
    intéressant maximal de $G$.
    \item[\sc  Sortie~:] Si  $C(T)$  possède un  chemin sortant,  ``il
    existe un chemin sortant'', sinon ``pas de chemin sortant''.
    \item[\sc Calcul~:] 
      Pour chaque sommet $u$ de $G\setminus (T \cup C(T))$ faire~:

      À partir  de $u$ et dans  $G\setminus T$, marquer  les sommets à
      l'aide  d'un parcours  en  largeur d'abord.   Si l'ensemble  des
      sommets  marqués de  $C(T)$  n'est pas  une  clique, stopper  et
      retourner ``il existe un chemin sortant''. Fin du ``faire''.
 
      Si aucun sommet $u$ n'a  permis de stopper, stopper et retourner
      ``pas de chemin sortant''
    \item[\sc  Complexité~:]  $O(m(n-k))$  où  $k$ est  le  nombre  de
    sommets de $C(T)$.
  \end{itemize}
\end{algorithme}

\begin{preuve}
  Il   est    clair   que   notre   algorithme    est   conforme   aux
  spécifications. Considérer  tous les  sommets $u$ de  $G\setminus (T
  \cup  C(T))$ nécessite  au  pire $O(n-k)$  étapes.   Le parcours  en
  largeur  d'abord  nécessite  au   pire  $O(m)$  étapes.   Tester  si
  l'ensemble  des  sommets marqués  de  $C(T)$  n'est  pas une  clique
  nécessite au pire $O(m)$ étapes. D'où la complexité annoncée.
\end{preuve}

\begin{algorithme}
  \label{artemis.a.z}
  \begin{itemize}
    \item[\sc  Entrée~:] Un  graphe $G$  et $T$,  ensemble intéressant
    maximal de $G$ tel qu'il existe un chemin sortant pour $C(T)$.
    \item[\sc Sortie~:]  Un plus court chemin sortant.
    \item[\sc Calcul~:] Pour tout sommet $u \in C(T)$, on calcule dans
      $G \setminus (T \cup (N(u) \cap  C(T)) )$ un plus court chemin à
      partir de $u$ et dont  l'autre extrémité $v$ est dans $C(T)$ (si
      un  tel  chemin existe).   On  retourne  un  chemin de  longueur
      minimale parmi ces chemins.
    \item[\sc Complexité~:] $O(nm)$.
  \end{itemize}
\end{algorithme}

\begin{preuve}
  Il est  clair que notre  algorithme est conforme  aux spécifications
  (notons que par hypothèse,  il existe un chemin sortant). Considérer
  tous les sommets de $C(T)$  nécessite au pire $n$ étapes.  Le calcul
  du  plus court  chemin  nécessite  au pire  $O(m)$  étapes, d'où  la
  complexité annoncée.
\end{preuve}

\begin{algorithme}
  \label{artemis.a.ep}
  \begin{itemize}
    \item[\sc  Entrée~:]  Un   graphe  $G$  d'Artémis,  $T$,  ensemble
      intéressant  maximal de $G$  et $Z$  chemin sortant  de longueur
      minimale pour $C(T)$.
    \item[\sc Sortie~:] Une paire d'amis spéciale de $G$.
    \item[\sc Calcul~:]  On note $z_1$ l'une des  extrémités du chemin
    sortant, et $z_k$ l'autre  extrémité. On calcule les ensembles $A$
    et $B$~:  $A$ est l'ensemble des  sommets de $C(T)$ qui  de $Z$ ne
    voit que $z_1$, $B$ que $z_k$.

    Pour tout sommet  $a$ de $A$ calculer un plus  court chemin de $a$
    vers $B$  dans $G\setminus(T\cup  (N(a)\setminus A))$.  Si  un tel
    chemin n'existe  pas, retenir ``$a$  est maximal''. De  même, pour
    tout  sommet  $b$  de  $B$   calculer  un  plus  court  chemin  de
    $B\setminus  b$  vers  $A$ dans  $G\setminus(T\cup  (N(b)\setminus
    B))$.  Si un tel chemin n'existe pas, retenir ``$b$ est maximal''.
    On retourne une  paire $\{a,b\}$ où $a$ est  un élément maximal de
    $A$, et $b$ de $B$.
    \item[\sc Complexité~:] $O(nm)$.
  \end{itemize}
\end{algorithme}

\begin{preuve}
  D'après la preuve du théorème~\ref{artemis.t.main}, un élément $a\in
  A$ est maximal pour la  relation $<_A$ si et seulement s'il n'existe
  pas  de chemin  $P$ de  $a$  vers $B$  dans $G\setminus(T\cup  (N(a)
  \setminus A))$.   En effet, d'après  le lemme~\ref{artemis.l.pq}, un
  tel chemin $P$ choisi de  longueur minimale est de longueur impaire,
  et son unique  sommet $a'$ voyant $a$ vérifie $a'  >_A a$.  Donc, si
  $a$ est maximal,  un tel chemin n'existe pas,  et réciproquement, si
  un tel chemin n'existe pas,  $a$ est maximal.  Notons que d'après la
  preuve du  théorème~\ref{artemis.t.main}, il existe  bien un élément
  maximal pour $<_A$, que  l'algorithme détectera forcément.  De même,
  l'algorithme détectera  bien un sommet  $b$ maximal pour  $<_B$.  Et
  d'après le lemme~\ref{artemis.l.absep}, $\{a,b\}$ est bien une paire
  d'amis spéciale  de $G$.  Donc  notre algorithme est conforme  à ses
  spécifications.

  Calculer  les  ensembles  $A$   et  $B$  nécessite  au  plus  $O(m)$
  étapes. Le calcul des plus courts  chemins prend lui aussi au plus $m$
  étapes, et  au pire ce calcul  est à faire $n$  fois.  La complexité
  est donc bien $O(nm)$.
\end{preuve}

\begin{algorithme}
  \label{artemis.a.ept}
  \begin{itemize}
  \item[\sc Entrée~:] Un graphe d'Artémis $G$.
  \item[\sc Sortie~:]  Une paire d'amis  spéciale de $G$ ou  ``$G$ est
    une clique''.
  \item[\sc Calcul~:] On invoque l'algorithme~\ref{artemis.a.t}.  S'il
    retourne ``$G$ est une clique'',  on retourne la même chose.  S'il
    retourne  une  paire  d'amis  spéciale  de  $G$,  on  la  retourne
    également.   S'il retourne  un ensemble  intéressant  maximal $T$,
    alors  on invoque  l'algorithme~\ref{artemis.a.zon}.   S'il répond
    ``Pas  de   chemin  sortant'',  alors   on  invoque  récursivement
    l'algorithme~\ref{artemis.a.ept} sur le graphe $G[C(T)]$, qui nous
    retournera  une  paire  d'amis  spéciale  de  $G[C(T)]$  que  l'on
    retourne.  S'il retourne ``il existe un chemin sortant'', alors on
    invoque l'algorithme~\ref{artemis.a.z} qui nous fournit un chemin
    sortant de longueur minimale, puis l'algorithme~\ref{artemis.a.ep}
    qui  nous retourne  une paire  d'amis  spéciale de  $G$, que  l'on
    retourne.
  \item[\sc Complexité~:] $O(nm)$.
  \end{itemize}
\end{algorithme}

\begin{preuve}
  Dans  la preuve  du théorème~\ref{artemis.t.main},  on a  montré que
  s'il  y a  un  ensemble intéressant  maximal  $T$ et  pas de  chemin
  sortant pour  $C(T)$, alors une  paire d'amis spéciale  de $G[C(T)]$
  est  également  une  paire  d'amis  spéciale  de  $G$.   Donc  notre
  algorithme est bien conforme à ses spécifications.

  On   va   montrer   par    induction   sur   $n$   qu'un   appel   à
  l'algorithme~\ref{artemis.a.ept} nécessite au pire $O(nm)$ étapes de
  calcul. Supposons que  cela soit vrai pour tous  les graphes avec au
  plus $n-1$ sommets.  On appelle  l'algorithme sur un graphe avec $n$
  sommets. Si le  graphe est une clique, ou s'il  a un chemin sortant,
  alors   il   n'y   a   pas   d'appel  récursif,   et   l'appel   aux
  algorithmes~\ref{artemis.a.z}   et~\ref{artemis.a.ep}  nécessite  au
  pire $O(nm)$  étapes.  S'il  n'y a pas  de chemin sortant,  alors on
  note     $k$     la    taille     de     $C(T)$.     L'appel     aux
  algorithmes~\ref{artemis.a.t}  et~\ref{artemis.a.zon}  nécessite  au
  pire  $O((n-k)m)$  étapes,  et  par hypothèse  d'induction,  l'appel
  récursif   à  l'algorithme~\ref{artemis.a.ept}  nécessite   au  pire
  $O(km)$ étapes, ce qui totalise bien $O(nm)$ étapes.
\end{preuve}

\index{coloration!graphes d'Artémis}
\index{Artémis~(graphe~d'---)!coloration}
\begin{algorithme} \label{artemis.a.color}
  \begin{itemize}
    \item[\sc Entrée~:] Un graphe d'Artémis $G$.
    \item[\sc Sortie~:] Une coloration optimale de $G$ et une clique
    de taille maximale de $G$.
  \item[\sc Calcul~:] \mbox{}

    Coloration~:   on  invoque   l'algorithme~\ref{artemis.a.ept}.  Si
    celui-ci   répond  ``$G$   est   une  clique'',   on  colore   $G$
    trivialement.  Sinon,  tant qu'on  n'obtient  pas  une clique,  on
    contracte       la      paire      d'amis       retournée      par
    l'algorithme~\ref{artemis.a.ept},   et   on   appelle  à   nouveau
    l'algorithme~\ref{artemis.a.ept} sur le graphe contracté. Quand on
    arrive à une  clique de taille $\omega$, on  donne à chaque sommet
    de $G$ la couleur de son représentant dans la clique.

    Clique maximale~: on marque tous les sommets de la clique, puis on
    redécontracte chaque paire d'amis $\{a,b\}$, en ne laissant marqué
    que celui  des deux  sommets qui avec  les autres  sommets marqués
    induit  une clique taille  $\omega$.  Une  fois toutes  les paires
    décontractés, on retourne l'ensemble des sommets marqués.
    \item[\sc Complexité~:] $O(n^2m)$.
  \end{itemize}
\end{algorithme}

\begin{preuve}
  D'après    les    propriétés    des    paires    d'amis    spéciales
  (lemme~\ref{artemis.l.sep}),  les   contractions  de  paires  d'amis
  spéciales conduisent bien à  une clique. D'après les explications de
  la  section  d'introduction  aux graphes  parfaitement  contractiles
  (section~\ref{pair.s.pc}), le théorème  de Fonlupt et Uhry notamment
  (théorème~\ref{pair.t.fonlupt}),  l'algorithme  conduit  bien à  une
  coloration optimale de $G$ et à une clique de taille maximale de
  $G$. L'ensemble des sommets marqués retournés à la fin de
  l'exécution  induit bien une clique de taille maximale de~$G$.

  Au pire, il faut $O(n)$  contractions de paires d'amis pour parvenir
  à une clique, d'où la complexité annoncée.
\end{preuve}

Nous avons donc bien un algorithme de coloration des graphes d'Artémis
de compléxité $O(n^2m)$. On peut facilement modifier cet algorithme de
manière à lui faire accepter  n'importe quel graphe en entrée. Le seul
problème est que  sur un graphe qui n'est pas  d'Artémis, il risque de
ne pas exister d'élément maximal  pour les ordres $<_A$ et $<_B$. Mais
alors on peut toujours retourner une paire quelconque de sommets (tant
mieux  si c'est  une paire  d'amis). La  phase finale  de l'algorithme
risque  alors de  contracter des  paires qui  ne sont  pas  des paires
d'amis, ce qui aura pour  conséquence que la coloration finale ne sera
pas forcément  optimale, et que l'ensemble  retourné à la  fin ne sera
pas forcément une  clique.  Mais il se peut très  bien que malgré tout
l'algorithme  retourne une  coloration avec  $\omega$ couleurs  et une
clique de taille  $\omega$, auquel cas on aura  la certitude que c'est
bien une coloration optimale et  une plus grande clique (sans être sûr
que le  graphe de  départ est d'Artémis).   Notre algorithme  est donc
robuste\index{robuste~(algorithme~---)}~:   ou  bien   il   donne  une
coloration optimale et une clique de taille maximale, ou bien il donne
il donne  autre chose,  et on a  la certitude  que le graphe  donné en
entrée n'était pas d'Artémis.  Dans ce dernier cas, une trace complète
de l'exécution est un objet  de taille polynomiale qui certifie que le
graphe  en entrée  n'était pas  d'Artémis.  Une  question intéressante
serait de savoir si à partir  de cet objet complexe, on peut retrouver
facilement un trou impair, un antitrou  ou un prisme (en tout cas cela
peut    être    fait   directement    grâce    aux   algorithmes    du
chapitre~\ref{reco.chap}).

\section{Conséquences pour des sous-classes d'Artémis}

Nous  analysons ici  notre preuve  dans le  cadre restreint  des trois
classes     de    graphes    parfaitement     contractiles    définies
pages~\pageref{pair.ss.meyniel} à~\pageref{pair.ss.po}.

\subsection*{Graphes de Meyniel}

\label{artemis.ss.meyniel}

Notre preuve du  théorème~\ref{artemis.t.main} ressemble beaucoup à la
preuve classique de l'existence d'une paire d'amis dans les graphes de
Meyniel.  Nous redonnons cette preuve avec nos notations.

\index{Meyniel~(graphe~de~---)!paire~d'amis}
\index{paire~d'amis!graphes~de~Meyniel}

\begin{theoreme}[Meyniel, \cite{meyniel:87}]
  Soit $G$ un  graphe de Meyniel.  Si $G$ n'est  pas une clique, alors
  $G$ possède une paire d'amis.
\end{theoreme}

\begin{preuve}
  On procède par récurrence.  Si $G$  a moins de trois sommets, si $G$
  est une  clique ou  une réunion disjointe  de cliques, le  lemme est
  trivial.  Sinon, soit $T$  un ensemble intéressant quelconque de $G$.
  Par  hypothèse  de  récurrence,   il  existe  une  paire  d'amis  de
  $G[C(T)]$. Mais  d'après le lemme~\ref{rr.l.rrmeyniel},  cette paire
  d'amis est aussi une paire d'amis de $G$.
\end{preuve}

Notons  que  si  on  exécute l'algorithme~\ref{artemis.a.ept}  sur  un
graphe  de Meyniel,  la contraction  de la  paire d'amis  retournée ne
redonne    pas     toujours    un    graphe     de    Meyniel    (voir
figure~\ref{pair.fig.meyniel}).  Pour  les  graphes de  quasi-Meyniel\index{quasi-Meyniel~(graphe~de~---)},
nous ne savons pas ce qu'il en est.

\subsection*{Graphes faiblement triangulés}
\label{artemis.ss.ft}

Si on spécialise notre preuve aux graphes faiblement triangulés, alors
on  obtient une  nouvelle  preuve d'un  ancien résultat~:  l'existence
d'une 2-paire.

\begin{theoreme}[Hayward, Ho\`ang, Maffray, \cite{hayward.hoang.m:90}]
  Soit  $G$   un  graphe  faiblement  triangulé  qui   n'est  pas  une
  clique. Alors $G$ possède une $2$-paire.
\end{theoreme}

\begin{preuve}
  Si  $G$ est  une réunion  disjointe d'au  moins deux  cliques, alors
  toute paire de sommets non  adjacents forme une $2$-paire. Sinon, on
  considère   un  ensemble  intéressant   maximal  $T$.    D'après  le
  lemme~\ref{rr.l.rrwt},  il   n'y  a  pas  de   chemin  sortant  pour
  $C(T)$. Par  hypothèse de récurrence, on trouve  une $2$-paire $\{u,
  v\}$ de  $G[C(T)]$. Mais  $\{u, v\}$ est  aussi une 2-paire  de $G$,
  puisque  tout  chemin  entre $u$  et  $v$  passant  par $T$  est  de
  longueur~2 et puisqu'il n'y a pas de chemin sortant pour $C(T)$.
\end{preuve}

\index{coloration!graphes~faiblement~triangulés}
\index{faiblement~triangulé~(graphe~---)!coloration}
\index{faiblement~triangulé~(graphe~---)!existence~d'une~2-paire}
\index{2-paire!graphes~faiblement~triangulés}
\index{paire@2-paire!graphes~faiblement~triangulés}

La   preuve  du  théorème   ci-dessus  montre   que  si   on  applique
l'algorithme~\ref{artemis.a.ept}  à  un  graphe faiblement  triangulé,
alors  l'algorithme  de  recherche  des  éléments  maximaux  pour  les
ensembles  $A$  et  $B$ (algorithme~\ref{artemis.a.ep})  n'est  jamais
appelé puisqu'il n'y  a pas de chemin sortant. Il  est en fait inutile
d'appeler   les   algorithmes   de   recherche   de   chemin   sortant
(algorithmes~\ref{artemis.a.zon}  et~\ref{artemis.a.z}).  On  est donc
obligé      quoi     qu'il     arrive      d'appeler     récursivement
l'algorithme~\ref{artemis.a.ept}.  À  chaque fois, on  doit rechercher
un ensemble  intéressant maximal (algorithme~\ref{artemis.a.t}).  Donc
la complexité totale de l'algorithme~\ref{artemis.a.ept} est inchangée
même si on  le restreint aux graphes faiblement  triangulés. On notera
cependant  que la  paire d'amis  retournée est  toujours  une 2-paire,
puisque toute  paire de sommets non-adjacents  d'une réunion disjointe
de  cliques  est  une  deux-paire.   Rappelons également  à  titre  de
comparaison  que  l'algorithme  le  plus  performant à  ce  jour  pour
colorier   les  graphes  faiblement   triangulés  a   pour  complexité
$O(n+m^2)$              (Hayward,              Spinrad              et
Sritharan~\cite{hayward.spinrad.sritharan:optwt}).

\subsection*{Graphes parfaitement ordonnables}

\index{parfaitement~ordonnable~(graphe~---)!coloration}
\index{coloration!graphes~parfaitement~ordonnables}

Notre preuve de la  parfaite contractilité des graphes d'Artémis donne
une  preuve totalement nouvelle  pour le  cas particulier  des graphes
parfaitement  ordonnables.   À  notre  connaissance, cette  preuve  ne
généralise  aucun  ancien résultat.   Notons  que  si nous  appliquons
l'algorithme~\ref{artemis.a.ept} à  un graphe parfaitement ordonnable,
alors  nous ne pouvons  pas garantir  que la  contraction de  la paire
d'amis retournée donne un nouveau graphe parfaitement ordonnable.

\section{Théorèmes de décomposition et paires d'amis}

Nous donnons  ici des théorèmes  de décomposition pour les  classes de
graphes intéressantes du point de vue des paires d'amis. Ces théorèmes
sont tous  énoncés presque explicitement (et bien  sûr démontrés) dans
la      preuve       du      théorème      fort       des      graphes
parfaits~\cite{chudvovsky.r.s.t:spgt}, et nous ne pouvons donc pas les
revendiquer  pour nôtres. Nous  avons tout  d'abord besoin  d'un petit
lemme technique~:

\index{line-graphe~de~subdivision~bipartie~de~$K_4$!contient~un~prisme~long}
\begin{lemme}
  \label{artemis.l.ltech}
  Soit $G$  le line-graphe d'une subdivision bipartie  de $K_4$. Alors
  ou  bien  $G$  est   $L(K_{3,3}\setminus  e)$  et  $G$  contient  un
  prisme~($\overline{C_6}$), ou bien $G$ contient un prisme long.
\end{lemme}

\begin{preuve}
  On suppose que $G$ est le line-graphe d'un graphe $R$ et que $R$ est
  une  subdivision  bipartie  d'un  $K_4$  dont on  note  les  sommets
  $a,b,c,d$. On choisit une arête  du $K_4$ subdivisée en un chemin de
  longueur maximale --- sans perte  de généralité, on suppose que $ab$
  est choisie.  L'arête $ab$ du $K_4$ est donc subdivisée en un chemin
  de longueur au moins 2. Dans  $R$, il existe un unique chemin de $a$
  vers $b$  qui ne  passe pas  par $c,d$~: on  note l'ensemble  de ses
  arête $E_0$.  Il  existe un unique chemin de $a$  vers $b$ qui passe
  par  $c$ sans  passer par  $d$~: on  note l'ensemble  de  ses arêtes
  $E_c$. Enfin on note $E_d$  l'ensemble des arêtes de l'unique chemin
  de $R$ passant par $d$ sans passer par $c$. On vérifie que $E_0 \cup
  E_c \cup E_d$, qui peut être vu comme un ensemble de sommets de $G$,
  induit un prisme  de $G$.  Si ce prisme  n'est pas $\overline{C_6}$,
  on   satisfait  une  conclusion   du  lemme.    Si  ce   prisme  est
  $\overline{C_6}$, alors  les seules arête subdivisées  du $K_4$ sont
  $ab$ et $cd$,  l'arête $ab$ n'est subdivisée qu'une  seule fois, et,
  d'après  le choix  initial de  $ab$, l'arête  $cd$ n'est  elle aussi
  subdivisée qu'une seule fois.  Donc, $R = K_{3,3}\setminus e$, et on
  satisfait l'autre conclusion du lemme.
\end{preuve}

Voici un théorème qui résume les propriétés des graphes d'Artémis~:

\index{Artémis~(graphe~d'---)!paire~d'amis}
\index{paire~d'amis!graphes~d'Artémis}
\index{décomposition!graphes d'Artémis}
\index{décomposition!par~2-joint}
\index{décomposition!par~partition~antisymétrique}
\index{Artémis~(graphe~d'---)!décomposition}
\index{Artémis~(graphe~d'---)!caractérisation}
\begin{theoreme}
  \label{artemis.t.srtuctart}
  Soit $G$ un graphe. Les quatre propriétés suivantes sont équivalentes~:
  \begin{outcomes}
  \item
    \label{artemis.o.srtuctart1}
    $G$ est d'Artémis.
  \item
    \label{artemis.o.srtuctart2}
    $G$ ne contient ni trou impair, ni antitrou long, ni prisme.
  \item
    \label{artemis.o.srtuctart3}
    Pour tout sous-graphe induit $G'$ de $G$,  pour tout ensemble $T$
    anticonnexe de  $G'$, ou  bien $C(T)$ est  une clique, ou  bien il
    existe une paire d'amis spéciale de $G'$ dans $C(T)$.
  \item
    \label{artemis.o.srtuctart4}
    Pour tout sous-graphe induit $G'$ de $G$, ou bien $G'$ possède une
    partition antisymétrique paire, ou  bien $G'$ est biparti, ou bien
    $G'$ est une clique.
  \end{outcomes}
\end{theoreme}

\begin{preuve}
  Les                           propriétés~(\ref{artemis.o.srtuctart1})
  et~(\ref{artemis.o.srtuctart2}) sont équivalentes par définition. La
  propriété       ~(\ref{artemis.o.srtuctart2})       implique      la
  propriété~(\ref{artemis.o.srtuctart3})           d'après          le
  théorème~\ref{artemis.t.main}. La  réciproque est plus  facile~: les
  antitrous longs  et les  trous impairs n'ont  pas de  paire d'amis~;
  quant aux prismes, ils peuvent en avoir, mais le voisinage d'un coin
  n'en contient  jamais. 

  Notons qu'un  graphe d'Artémis ne contient jamais  de line-graphe de
  subdivision       bipartie      de      $K_4$       (d'après      le
  lemme~\ref{artemis.l.ltech},  ils contiennent  des  prismes), ni  de
  double  diamant   (qui  contient  un   $\overline{C_6}$).   Donc  le
  théorème~\ref{graphespar.t.sartemis}
  page~\pageref{graphespar.t.sartemis}       montre       que       la
  propriété~(\ref{artemis.o.srtuctart2})          implique          la
  propriété~(\ref{artemis.o.srtuctart4}).   La   réciproque  est  plus
  directe~: on  vérifie facilement que les antitrous  longs, les trous
  impairs  et  les prismes  ne  sont pas  bipartis,  ne  sont pas  des
  cliques, et n'ont pas de partition antisymétrique paire.
\end{preuve}

\index{décomposition~et~paires~d'amis}  Le théorème  ci-dessus suggère
que, de fait, les paires  d'amis s'héritent d'une certaine manière par
partition antisymétrique paire, bien que nous n'arrivions pas à savoir
pourquoi et comment. Il serait  donc intéressant d'avoir une preuve de
l'existence de paire d'amis dans les graphes d'Artémis par le théorème
de   décomposition.    S'il  est   bien   vrai   que  les   partitions
antisymétriques paires préservent les paires d'amis, alors le théorème
de décomposition suivant est  un encouragement à rechercher des paires
d'amis dans les graphes bipartisans~:


\index{bipartisan~(graphe~---)!décomposition}
\index{décomposition!graphes bipartisans}
\begin{theoreme}
  \label{artemis.t.srtuctbip}
  Soit $G$ un graphe. Les trois propriétés suivantes sont équivalentes~:
  \begin{outcomes}
  \item
    \label{artemis.o.srtuctbip1}
    $G$ est bipartisan.
  \item
    \label{artemis.o.srtuctbip2}
    $G$ et $\overline{G}$ ne contiennent ni trou impair,  ni $L(K_{3,3}
    \setminus e)$, ni double-diamant,  ni prisme long.
  \item
    \label{artemis.o.srtuctbip3}
    Pour tout sous-graphe induit $G'$ de $G$, ou bien $G'$ possède une
    partition antisymétrique paire, ou  bien $G$ est biparti, ou bien
    $G$ est complémentaire de biparti.
  \end{outcomes}
\end{theoreme}

\begin{preuve}
  Les                           propriétés~(\ref{artemis.o.srtuctbip1})
  et~(\ref{artemis.o.srtuctbip2})  sont  équivalentes par  définition.
  Le lemme~\ref{artemis.l.ltech} montre que les graphes bipartisans
  sont sans line-graphe de subdivision bipartie de $K_4$. 
  Les                              théorème~\ref{graphespar.t.sartemis}
  et~\ref{graphespar.t.sbipar}      montrent     alors      que     la
  propriété~(\ref{artemis.o.srtuctbip2})          implique          la
  propriété~(\ref{artemis.o.srtuctbip3}).   La   réciproque  est  plus
  directe~: on  vérifie facilement que les trous  impairs, les prismes
  longs, $L(K_{3,3}  \setminus e)$  et le double  diamant ne  sont pas
  bipartis, ne sont  pas complémentaires de bipartis, et  n'ont pas de
  partition antisymétrique paire.
\end{preuve}

Pour  les  graphes d'Artémis  pairs,  la  situation  est un  peu  plus
compliquée.  Il est facile de  vérifier que ces graphes ne contiennent
pas  de  line-graphe  de  subdivision  bipartite de  $K_4$~:  de  tels
sous-graphes contienne  toujours un prisme impair (ce  fait bien connu
est   démontré   un    peu   plus   loin~:   lemme~\ref{reco.l.techK4}
page~\pageref{reco.l.techK4}).   Les  théorèmes  de décomposition  des
graphes    de    Berge    (théorèmes~\ref{graphespar.t.spgtprismepair}
et~\ref{graphespar.t.sartemis})   indiquent  alors  que   les  graphes
d'Artémis pairs sont  soit bipartis, soit des cliques,  soit un prisme
pair avec 9  sommets, soit ont une partition  antisymétrique paire, ou
bien ont un deux-joint.  Toutefois, il doit y avoir un théorème un peu
plus fort qui nécessite une nouvelle définition.
 
\index{Artémis~pair~(graphe~d'---)!décomposition}
\index{décomposition!graphes~d'Artémis~pairs}  Soit  $(X_1,  X_2)$  un
2-joint d'un graphe $G$ (on reprend les notations habituelles pour les
2-joints,   voir   page~\pageref{graphespar.p.2joint}).   On   appelle
\emph{chemin sortant}  pour le 2-joint  tout chemin de  $G[X_i]$ ayant
une  extrémité dans  $A_i$,  l'autre extrémité  dans  $B_i$, tous  les
autres sommets hors de $A_i, B_i$  ($i = 1, 2$).  On dit qu'un 2-joint
est   \emph{pair}\index{2-joint!pair}\index{pair!2-joint~---}   (resp.
\emph{impair}\index{2-joint!impair}\index{impair!2-joint~---}) si tous
ses chemins sortants  sont pairs (resp. impairs).  Il  est clair qu'un
2-joint $(X_1,  X_2)$ d'un  graphe de Berge  $G$ est toujours  pair ou
impair, sans quoi $G$ contient un trou impair.

\index{décomposition!par~2-joint~pair}
Il doit être possible de  montrer qu'un graphe d'Artémis pair qu'on ne
peut  décomposer  que par  2-joint  peut  toujours  se décomposer  par
2-joint pair.  En effet, les 2-joints pairs apparaissent naturellement
dans  les  prismes pairs,  structure  à  partir  de laquelle  on  peut
décomposer le graphe  tout entier dans la preuve  du théorème fort des
graphes  parfaits  (voir  théorème~\ref{graphespar.t.spgtprismepair}).
Pour donner un théorème  de décomposition des graphes d'Artémis pairs,
il  devrait donc  suffire  de recopier  le  chapitre de  la preuve  du
théorème  fort  des  graphes  parfaits  consacré au  prisme  pair,  en
remplaçant partout ``2-joint'' par ``2-joint pair''.  Ce travail n'est
peut-être pas insurmontable,  mais n'a d'intérêt que si  on parvient à
utiliser le théorème de décomposition obtenu.

\chapter{Problèmes de reconnaissance}
\label{reco.chap}

Reconnaître  une  classe  de   graphes,  c'est  donner  un  algorithme
permettant de décider  si un graphe appartient ou non  à la classe. Il
est  toujours intéressant  de disposer  de ce  genre  d'algorithme. En
effet, est-il vraiment utile  de disposer d'algorithmes performants de
coloration pour une classe si l'on  est même pas capable de décider si
l'algorithme fonctionne pour un graphe donné~?

Nous débuterons  ce chapitre  par un état  de l'art en  rappelant tout
d'abord à  titre d'illustration  un algorithme naïf  de reconnaissance
des graphes  faiblement triangulés et en donnant  des indications très
succintes sur  la résolution récente du problème  de la reconnaissance
des  graphes de  Berge.  Nous  donnerons ensuites  des  algorithmes en
temps polynomial pour détecter  certains types de sous-graphes induits
(prismes, prismes pairs,  prismes impairs, line-graphes de subdivision
de  $K_4$).  À  chaque fois,  nous  serons obligés  de restreindre  le
champs  d'application de  nos  algorithmes à  des  classes de  graphes
particuliers  (graphes sans  pyramide ou  sans trou  impair  selon les
cas).   Nous  justifierons  cette  restriction  en  prouvant  que  nos
problèmes de détection deviennent  NP-complets lorsqu'on essaye de les
résoudre  pour  tout  graphe.   En  conclusion,  nous  parviendrons  à
reconnaître  en temps  polynomial les  classes de  graphes suivantes~:
graphes d'Artémis,  d'Artémis pairs et bipartisans.  Nous donnerons un
algorithme  de  coloration  en   temps  polynomial  pour  les  graphes
d'Artémis  pairs,   dont  nous  ne   sommes  pas  parvenu   à  prouver
l'optimalité. Toutefois, nous montrons qu'il donne bien une coloration
optimale, à condition que la  conjecture de Everett et Reed soit vraie
(conjecture~\ref{pair.conj.EverettReed}).

Les résultats  de ce chapitre ont  été soumis au {\it  SIAM Journal of
Discrete Mathematics}~\cite{nicolas:reco}.

\section{État de l'art}

Calculer des plus courts chemins pour trouver des sous-graphes induits
paraît assez naïf  en première analyse. Par exemple,  si on cherche un
prisme dans un  graphe, si on connait les  deux triangles $\{a_1, a_2,
a_3\}$ et $\{  b_1, b_2, b_3\}$ d'un prisme du graphe,  alors il n'y a
aucune raison  {\it a priori} pour  que des plus  courts chemins $P_i$
reliant $a_i$ à $b_i$ forment un prisme.  Pourtant, un algorithme bien
connu  de reconnaissance des  graphes faiblement  triangulés fontionne
sur ce  principe. Maria Chudnovsky et  Paul Seymour ont  montré que ce
type  de techniques  peut aider  à reconnaître  les graphes  de Berge,
moyennant de multiples précautions.

\subsection{Graphes faiblement triangulés}
\label{reco.s.wt}

Le lemme et l'algorithme ci-dessous  sont faciles et bien connus.  Ils
ne sont publiés nulle part pour  autant que nous ayons pu le vérifier,
et  on peut  considérer qu'ils  appartiennent au  folklore.   Ils sont
néanmoins utiles et donnent une bonne  idée de ce qui va suivre.  Pour
chaque  type  de  sous-graphe  à  détecter, nous  donnerons  un  lemme
affirmant qu'en  calculant des plus  courts chemins entre  des sommets
d'un sous-graphe donné de ce  type (un prisme par exemple), on obtient
à nouveau un  sous-graphe induit de même type.   Ce lemme permettra la
mise  en \oe uvre  d'un algorithme  qui ``devine''  un nombre  fixé de
sommets du sous-graphe  puis qui en retrouve le  reste par des calculs
de plus courts chemins.

On rappelle qu'un \emph{long trou} dans un graphe est un cycle induit
de taille au moins 5.
 
\begin{lemme}
  Soit $G$  un graphe et $H$ un  long trou de $G$  de taille minimale.
  Soit $\bp  a \tp u \tp  b \ep$ un chemin  de $H$.  Soit  $P$ un plus
  court chemin de $a$ vers $b$ ne contenant ni voisin de $u$ ni voisin
  commun de $a$  et $b$.  Alors $\{u\} \cup V(P)$  induit un long trou
  de $G$ de taille minimale.
\end{lemme}

\begin{preuve}
Clair.
\end{preuve}

\index{reconnaissance!graphes~faiblement~triangulés}
\index{faiblement~triangulé~(graphe~---)!reconnaissance}
\begin{algorithme} \label{reco.a.ft}
  \begin{itemize}
  \item[\sc Entrée :] Un graphe $G$.
  \item[\sc Sortie :]  Si $G$ possède un long  trou, alors l'algorithme
    retourne  un long  trou  de  $G$ de  taille  minimale.  Sinon,  il
    retourne ``Pas de long trou''.
  \item[\sc Calcul :]\mbox{}
    \begin{itemize}
    \item
      Considérer tous  les triplets $(a,  u, b)$ de $G$.   Pour chacun
      d'eux calculer  par l'algorithme \ref{base.a.pluscourt}  un plus
      court  chemin  $P$  de  $a$   vers  $b$  parmi  les  sommets  de
      $V\backslash (N(u) \cup (N(a)\cap N(b)))$.  Si $\{u\} \cup V(P)$
      induit un trou  de $G$, de taille inférieure  à un éventuel trou
      déjà rencontré, alors stocker ce trou en mémoire.
    \item
      Une  fois tous  les triplets  énumérés,  retourner l'eventuel
      trou de taille minimale.  S'il n'y en a pas, retourner ''pas
      de trou''.
    \end{itemize}
  \item[\sc Complexité :] $O(n^5)$.
  \end{itemize}
\end{algorithme}

\begin{preuve} Clair. \end{preuve}

L'algorithme~\ref{reco.a.ft} peut  être utilisé à  des fins multiples.
Par exemple  pour reconnaître  les graphes faiblement  triangulés~: il
suffit de l'exécuter  sur $G$ puis sur $\overline{G}$.   En fait, on a
déjà vu  page~\pageref{pair.t.deuxp} qu'il  existe un algorithme  dû à
Hayward,  Spinrad  et Sritharan~\cite{hayward.spinrad.sritharan:optwt}
pour    reconnaître    les    graphes   faiblement    triangulés    en
$O(n+m^2)$. Notons  que cet algorithme ne détecte  pas directement les
trous.             Récemment,            S.~Nikolopoulos            et
L.~Pallios~\cite{nikolopoulos.palios:hole} ont  donné un algorithme de
complexité     $O(n+m^2)$      qui     détecte     directement     les
trous.\index{détection!trous}

\subsection{Reconnaissance des graphes de Berge}

Comme  on l'a  dit plus  haut, le  problème de  la  reconnaissance des
graphes  parfaits  est  resté   ouvert  jusqu'en  2002.   Puis,  Maria
Chudnovsky, Paul  Seymour, Gérard Cornuéjols, Xinming  Liu et Kristina
Vu{\v  s}kovi{\'c}  ont  simultanément  proposé  deux  algorithmes  de
reconnaissance     des      graphes     de     Berge      en     temps
polynomial~\cite{chudnovsky.c.l.s.v:cleaning,chudnovsky.seymour:reco,cornuejols.liu.vuskovic:reco}.
Ce  problème semble après  coup moins  difficile que  la preuve  de la
conjecture forte  des graphes  parfaits, et il  a pourtant  été résolu
postérieurement et  sans que le  théorème de structure des  graphes de
Berge n'apporte aucune  aide (voir section~\ref{graphespar.s.poq} pour
des détails  sur ce  point).  Il faut  noter que les  deux algorithmes
reconnaissent les  graphes \emph{de  Berge}, et ne  reconnaissent donc
les  graphes parfaits  que  dans la  mesure  où le  théorème fort  des
graphes parfaits est vrai.  Il faut également noter que le problème de
la reconnaissance des graphes sans trou impair est encore ouvert.

Pour  fournir un  aperçu des  méthodes  employées, il  nous faut  deux
définitions.  Soient $G$ un graphe et $C$ un trou de $G$. On dit qu'un
sommet  $v$ de  $G$ est  \emph{majeur} sur  $C$ si  l'ensemble  de ses
voisins  sur $C$  n'est  inclu dans  aucun  $P_3$ de  $C$. On  appelle
\emph{raquette}\index{raquette} pour $C$ un  sommet majeur sur $C$ qui
voit exactement trois sommets de $C$. On remarquera que si $v$ est une
raquette  pour un trou  impair $C$  de $G$  de taille  minimale, alors
l'ensemble $V(C) \cup \{v\}$ induit une pyramide.
 
Les deux algorithmes  reposent sur une première phase  commune dite de
\emph{nettoyage}\index{nettoyage}               du              graphe
\cite{chudnovsky.c.l.s.v:cleaning}.  Cette  première phase permet soit
de conclure en temps polynomial que le graphe n'est pas de Berge, soit
de fournir une liste d'au  plus $O(n^8)$ sous-ensembles de $V(G)$ avec
la propriété suivante~:  si $C$ est un trou  impair de taille minimale
de $G$,  et si  $C$ est  sans raquette, alors  l'un des  ensembles est
disjoint de $C$  et contient tous les sommets majeurs  de $C$.  On dit
alors que le trou $C$  est \emph{nettoyé} (de ses voisins majeurs), et
on dit que le graphe $G$ est \emph{nettoyé}.

La phase  de nettoyage  permet donc de  ramener la  reconnaissance des
graphes  de Berge  à deux  problèmes plus  simples~: la  détection des
trous  impairs possèdant  une raquette,  puis la  détection  des trous
impairs nettoyés de taille minimale.   Ces problèmes ont été résolus de
deux  manières  différentes.  D'une  part,  Chudnosky  et Seymour  ont
proposé un algorithme direct de  recherche de pyramides dans un graphe
quelconque,   puis   de   trou    impair   dans   un   graphe   nettoyé
\cite{chudnovsky.seymour:reco}.   Cet  article   a  été  notre  source
d'inspiration principale pour  ce chapitre.  D'autre part, Cornuéjols,
Liu et  Vu\v skovi\'c \cite{cornuejols.liu.vuskovic:reco}  ont proposé
un   algorithme  de   reconnaissance  fondé   sur  leur   théorème  de
\index{décomposition~et~reconnaissance}
décomposition       des      graphes       sans       trou      impair
(théorème~\ref{graphespar.t.sansti}
page~\pageref{graphespar.t.sansti}).    Nous    avons   aussi   essayé
d'adapter  cette deuxième  méthode à  des problèmes  non  résolus mais
cette fois sans succès.

Nous  mentionnons pour  mémoire l'algorithme  suivant qui  utilise des
techniques de plus courts chemins.

\index{pyramide!détection~des~---}
\index{détection!pyramides}
\begin{algorithme}[Chudnovsky et Seymour~\cite{chudnovsky.seymour:reco}.]\label{reco.a.py}
  \begin{itemize}
  \item
    [\sc Entrée :] Un graphe $G$.
  \item
    [\sc  Sortie :] Si  $G$ est  sans pyramide,  l'algorithme retourne
    ``Pas de pyramide''. Sinon, il retourne une pyramide de $G$ de
    taille minimale.
  \item
    [\sc Complexité :] $O(n^9)$.
  \end{itemize}
\end{algorithme}

Enfin  nous citons  le principal  résultat sur  la  reconnaissance des
graphes de Berge~:

\index{reconnaissance!graphes~de~Berge}
\index{Berge~(graphe de ---)!reconnaissance}
\begin{algorithme}[\cite{chudnovsky.c.l.s.v:cleaning,chudnovsky.seymour:reco}] \label{reco.a.berge}
  \begin{itemize}
  \item[]{\bf (Chudnovsky, Seymour \& Cornuéjols, Liu, Vu\v skovi\'c)}
  \item
    [\sc Entrée :] Un graphe $G$.
  \item
    [\sc  Sortie  :]  Si  $G$  est  de  Berge,  l'algorithme  retourne
    {\sc Berge}. Sinon, il retourne {\sc Non Berge}.
  \item
    [\sc Complexité :] $O(n^9)$.
  \end{itemize}
\end{algorithme}

\section{Détection de sous-graphes}
\index{détection!sous-graphes|(}

Nous allons  maintenant voir une  série d'algorithmes de  détection de
sous-graphes.   Nous utiliserons intensivement  la technique  des plus
courts  chemins mise en  \oe uvre  par Chudnovsky  et Seymour  pour la
détection  des  pyramides.   Signalons   aussi  que  les  articles  de
Conforti, Cornuéjols  et Vu\v skovi\'c  sur les graphes de  Berge sans
carré    et     la    décomposition    des     graphes    sans    trou
impair~\cite{conforti.c.v:square,conforti.c.v:dstrarcut}    contiennent
de  nombreux résultats  sur la  manière dont  un sommet  peut  voir un
prisme. Les lemmes  qui suivent utilisent des résultats  du même type,
mais beaucoup plus simples.

\subsection{Prisme ou Pyramide}
\index{prisme!détection~des~---}
\index{pyramide!détection~des~---}

Nous  donnons ici  deux algorithmes  de détection  des prismes  ou des
pyramides.   Appliqués à  des graphes  sans pyramide,  ces algorithmes
deviennent  simplement des  algorithmes de  détection de  prismes.  Le
premier d'entre eux s'inspire directement des techniques utilisées par
Chudnovsky et  Seymour pour la  détection des pyramides,  mais possède
une complexité  bien meilleure car lorsque  l'on cherche simultanément
un prisme ou une pyramide dans un graphe, la technique des plus courts
chemins  se  simplifie.   Le  deuxième  algorithme  est  un  peu  plus
performant et  ne repose  plus aussi directement  sur les  plus courts
chemins.

On       dit       qu'un       sous-graphe      de       $G$       est
\emph{optimal}\index{optimal!sous-graphe ---} si son nombre de sommets
est inférieur ou  égal au nombre de sommets de tout  prisme de $G$, et
inférieur ou égal au nombre de sommets de toute pyramide de $G$.

\begin{lemme} \label{reco.l.pyraopt}
  Soit $G$  un graphe quelconque.   Soient $P_1$, $P_2$,  $P_3$, trois
  chemins de $G$ formant une pyramide optimale $F$ de triangle $\{b_1,
  b_2, b_3\}$. On note $a_1$ le coin  de $F$ et on suppose que pour $i
  = 1, 2, 3$, $P_i$ a pour extrémités $a_1$ et $b_i$.

  Soit  $P'_1$ un  plus  court chemin  de  $a_1$ vers  $b_1$ dont  les
  sommets intérieurs manquent $b_{2}$ et $b_{3}$.  Alors ou bien~:

  \begin{outcomes}
  \item
    L'ensemble  $V(P'_{1})  \cup  V(P_{2})  \cup V(P_{3})$  induit  un
    prisme optimal de coin $a_1$ et de triangle $\{b_1, b_2, b_3\}$
  \item
    Les  chemins  $P'_{1}$,  $P_{2}$,  $P_{3}$  forment  une  pyramide
    optimale de coin $a_1$ et de triangle $\{b_1, b_2, b_3\}$.
  \end{outcomes}

  On a un résultat similaire pour $P_2$ et $P_3$.
\end{lemme}

\begin{preuve}
  Si $a_1$ voit $b_1$, le lemme  est trivial. On suppose donc $a_1 b_1
  \notin E$.  Notons qu'alors $(P'_1)$ est de longueur au moins~2.  On
  note $H$  le trou induit par  $V(P_2) \cup V(P_3)$.  On  note $c$ le
  sommet de ${P'_1}^*$, le plus  proche  de $b_1$, et ayant un
  voisin dans  $H$.  Notons que $c$  existe car $a_1 \in  V(H)$ et car
  $P'_1$ est  de longueur au  moins~2.  Restent à examiner  quatre cas
  conduisant à une contradiction ou à ce qu'avance le lemme~:

  \begin{itemize}
  \item 
    Si  $c$ ne  voit  de $H$  que  $a_1$, alors  $P'_1$, $P_2$,  $P_3$
    forment une pyramide optimale de coin $a_1$ et de triangle $\{b_1,
    b_2, b_3\}$.
  \item 
    Si  $c$ ne voit  de $H$  que $a_1$  et un  voisin de  $a_1$, alors
    $V(P'_1) \cup V(P_2) \cup V(P_3)$ induit un prisme optimal de coin
    $a_1$ et de triangle $\{b_1, b_2, b_3\}$.
  \item
    Si $c$  manque $a_1$  alors le graphe  $F'$ induit par  $V(H) \cup
    V(\bp b_1 \tp  P'_1 \tp c\ep)$ et le  triangle $\{b_1, b_2, b_3\}$
    satisfont toutes les  hypothèses du lemme~\ref{rr.l.bruce}, ce qui
    montre que  $F'$ contient  un prisme ou  une pyramide.   Mais cela
    contredit  l'optimalité  de  $F$  car $|F'|<|F|$  (il  manque  les
    sommets de $P'_1$ entre $c$ et $a_1$).

  \item
    Si $c$  voit $a_1$  et au  moins un non-voisin  de $a_1$  sur $H$,
    alors il  existe un chemin $Q$  de $H$ de longueur  au moins~2, et
    dont les extrémités sont $a_1$ et un voisin de $c$.  Mais alors le
    graphe $F'$  induit par $(V(H) \cup V(P'_1))  \setminus V(Q^*)$ et
    le triangle $\{b_1, b_2, b_3\}$ satisfont toutes les hypothèses du
    lemme~\ref{rr.l.bruce}, ce qui montre  que $F'$ contient un prisme
    ou  une pyramide.   Mais cela  contredit l'optimalité  de  $F$ car
    $|F'|<|F|$ (il manque les sommets de $Q^*$).
  \end{itemize}

  Si $c$ voit  de $H$  un autre
  sommet que $a_1$ alors le graphe $F'$ induit par $V(H) \cup V( P'_1)
  \setminus  \{a_1\}$ et  le  triangle $\{b_1,  b_2, b_3\}$  satisfont
  toutes les  hypothèses du lemme~\ref{rr.l.bruce}, ce  qui montre que
  $F'$  contient  un  prisme  ou  une pyramide.  Mais  cela  contredit
  l'optimalité  de $F$  car $|F'|<|F|$~:  il manque  à $F'$  le sommet
  $a_1$. Donc, $c$  ne voit de $H$ que le  sommet $a_1$, exactement ce
  qu'on voulait montrer.
\end{preuve}

\begin{lemme} \label{reco.l.prismeopt}
  Soit $G$  un graphe quelconque.   Soient $P_1$, $P_2$,  $P_3$, trois
  chemins de  $G$ formant un  prisme optimal $F$ de  triangles $\{a_1,
  a_2, a_3\}$ et  $\{b_1, b_2, b_3\}$.  On suppose pour $i  = 1, 2, 3$
  que $P_i$ a pour extrémités $a_i$ et $b_i$.

  Soient~:
  \begin{itemize}
  \item
    $P'_1$ un plus  court chemin de $a_1$ vers  $b_1$ dont les sommets
    intérieurs manquent $b_2$ et $b_3$.
  \item
    $P'_2$ un plus  court chemin de $a_1$ vers  $b_2$ dont les sommets
    intérieurs manquent $b_1$ et $b_3$.
  \item
    $P'_3$ un plus  court chemin de $a_1$ vers  $b_3$ dont les sommets
    intérieurs manquent $b_1$ et $b_2$.
  \end{itemize}
  
   Alors~:

  \begin{outcomes}
  \item 
    \label{reco.o.prismeopt.1} 
    Les chemins $P'_1$, $P_2$, $P_3$  forment un prisme optimal de $G$
    de coin $a_1$ et de triangle $\{b_1, b_2, b_3\}$.
  \item 
    \label{reco.o.prismeopt.2}  
    Ou bien $V(P_1) \cup V(P'_2) \cup V(P_3)$ induit un prisme optimal
    de coin $a_1$  et de triangle $\{b_1, b_2,  b_3\}$, ou bien $P_1$,
    $P'_2$, $P_3$  forment une pyramide  optimale de coin $a_1$  et de
    triangle $\{b_1, b_2, b_3\}$.
  \item 
    \label{reco.o.prismeopt.3}  
    Ou bien $V(P_1) \cup V(P_2) \cup V(P'_3)$ induit un prisme optimal
    de coin $a_1$  et de triangle $\{b_1, b_2,  b_3\}$, ou bien $P_1$,
    $P_2$, $P'_3$ forment  une pyramide optimale de $G$  de coin $a_1$
    et de triangle $\{b_1, b_2, b_3\}$.
  \end{outcomes}
\end{lemme}

\begin{preuve}
  On utilisera sans mention explicite le fait que $(P_i)$ est au moins
  aussi long que $(P'_i)$.
  
  On commence par prouver la conclusion~(\ref{reco.o.prismeopt.1}). On
  note $H$  le trou de  $G$ induit par  $V(P_2) \cup V(P_3)$.   Si les
  sommets intérieurs de $P'_1$ manquent entièrement $H$, alors $P'_1$,
  $P_2$,  $P_3$  forment un  prisme  optimal  de  $G$.  Donc  on  peut
  supposer l'existence  d'un sommet intérieur  $c$ de $P'_1$  qui voit
  $H$  et choisi  aussi proche  que possible  de $b_1$.   Montrons que
  l'existence de ce sommet entraîne une contradiction~: le graphe $F'$
  induit  par $V(H)  \cup V(P'_1)  \setminus \{a_1\}$  et  le triangle
  $\{b_1,   b_2,   b_3\}$   satisfont   toutes   les   hypothèses   du
  lemme~\ref{rr.l.bruce}, ce qui montre que $F'$ contient un prisme ou
  une  pyramide.    Mais  cela  contredit  l'optimalité   de  $F$  car
  $|F'|<|F|$.  Ceci prouve la conclusion~(\ref{reco.o.prismeopt.1}).

  Prouvons la  conclusion~(\ref{reco.o.prismeopt.2}).  On note  $H$ le
  trou  de $G$  induit par  $V(P_1) \cup  V(P_3)$. Comme  $a_1$ manque
  $b_2$, on sait que $P'_2$  possède des sommets intérieurs.  Soit $c$
  le sommet  intérieur de  ${P'_2}$ qui voit  $V(H)$, et  choisi aussi
  proche que  possible de $b_2$.  Notons  que $c$ existe  car $a_1 \in
  H$.  Restent à examiner quatre cas conduisant à une contradiction ou
  à ce qu'avance le lemme~:

  \begin{itemize}
  \item 
    Si  $c$ ne  voit  de $H$  que  $a_1$, alors  $P_1$, $P'_2$,  $P_3$
    forment une pyramide optimale de coin $a_1$ et de triangle $\{b_1,
    b_2, b_3\}$.
  \item 
    Si  $c$ ne voit  de $H$  que $a_1$  et un  voisin de  $a_1$, alors
    $V(P_1) \cup V(P'_2) \cup V(P_3)$ induit un prisme optimal de coin
    $a_1$ et de triangle $\{b_1, b_2, b_3\}$.
  \item
    Si $c$  manque $a_1$  alors le graphe  $F'$ induit par  $V(H) \cup
    V(\bp b_2 \tp  P'_2 \tp c\ep)$ et le  triangle $\{b_1, b_2, b_3\}$
    satisfont toutes les  hypothèses du lemme~\ref{rr.l.bruce}, ce qui
    montre que  $F'$ contient  un prisme ou  une pyramide.   Mais cela
    contredit  l'optimalité  de  $F$  car $|F'|<|F|$  (il  manque  les
    sommets de $P'_2$ entre $c$ et $a_1$).

  \item
    Si $c$  voit $a_1$  et au  moins un non-voisin  de $a_1$  sur $H$,
    alors il  existe un chemin $Q$  de $H$ de longueur  au moins~2, et
    dont les extrémités sont $a_1$ et un voisin de $c$.  Mais alors le
    graphe $F'$  induit par $(V(H) \cup V(P'_2))  \setminus V(Q^*)$ et
    le triangle $\{b_1, b_2, b_3\}$ satisfont toutes les hypothèses du
    lemme~\ref{rr.l.bruce}, ce qui montre  que $F'$ contient un prisme
    ou  une pyramide.   Mais cela  contredit l'optimalité  de  $F$ car
    $|F'|<|F|$ (il manque les sommets de $Q^*$).
  \end{itemize}

  Ceci prouve la  conclusion~(\ref{reco.o.prismeopt.2}).  La preuve de
  la conclusion~(\ref{reco.o.prismeopt.3}) est semblable à celle de la
  conclusion~(\ref{reco.o.prismeopt.2}).
\end{preuve}

Nous  donnons un lemme  qui est  une variante  du précédent,  mais qui
concerne  cette  fois  les  prismes longs  (c'est-à-dire  les  prismes
différents de $\overline{C_6}$). Un prisme long \emph{minimal} dans un
graphe $G$ est un prisme long  dont le nombre de sommets est inférieur
ou égal à celui de tout autre prisme long du graphe.

\begin{lemme} \label{reco.l.prismelongopt}
  Soit $G$ un graphe sans pyramide.  Soient $P_1$, $P_2$, $P_3$, trois
  chemins  de $G$  formant un  prisme  long minimal  $F$ de  triangles
  $\{a_1, a_2, a_3\}$ et $\{b_1, b_2, b_3\}$.  On suppose pour $i = 1,
  2, 3$ que $P_i$ a pour extrémités $a_i$ et $b_i$. Soit $i\in \{1, 2,
  3\}$.  Soit $P'_i$ un plus court chemin de $a_1$ vers $b_i$ dont les
  sommets intérieurs  manquent $b_{i+1}$ et  $b_{i+2}$ (l'addition des
  indices s'entend modulo 3).

   \noindent   Alors,  l'ensemble   $V(P'_i)   \cup  V(P_{i+1})   \cup
   V(P_{i+2})$ induit un  prisme long optimal de $G$  de coin $a_1$ et
   de triangle $\{b_1, b_2, b_3\}$.
\end{lemme}

\begin{preuve} 
  La   preuve  est   en  tout   point   similaire  à   la  preuve   du
  lemme~\ref{reco.l.prismeopt}. Notons  qu'ici, le graphe  $G$ est par
  hypothèse sans pyramide ce qui fait qu'on ne peut jamais trouver de
  pyramide.
\end{preuve}

Nous proposons l'algorithme suivant pour détecter les pyramides ou les
prismes~:

\index{détection!prismes ou pyramides}
\index{prisme!détection~des~---}

\begin{algorithme}\label{reco.a.pypri}
  \begin{itemize}
  \item[\sc Entrée :] Un graphe $G$.

  \item[\sc Sortie :] Si le  graphe est sans pyramide et sans prisme~:
    ``Ni  prisme, ni  pyramide''.  Sinon~:  un prisme  optimal  ou une
    pyramide optimale de $G$.

  \item[\sc Calcul :] \mbox{}
    \begin{itemize}

    \item
      Considérer successivement tous  les quadruplets $(a_1, b_1, b_2,
      b_3)\in V^4$.
      
      Calculer pour  $i=1, 2, 3$ un  plus court chemin  $P_i$ de $a_1$
      vers  $b_i$ dont  les sommets  intérieurs manquent  le  reste du
      quadruplet.   Si  l'ensemble $V(P_1)  \cup  V(P_2) \cup  V(P_3)$
      induit un prisme ou une pyramide, alors le stocker en mémoire.

    \item
      Si  aucun  quadruplet  n'a  donné  de  prisme  ou  de  pyramide,
      retourner ``Ni prisme, ni pyramide''. Sinon, retourner un prisme
      ou  une  pyramide  de  taille  minimale parmi  ceux  stockés  en
      mémoire.

    \end{itemize}
  \item[\sc Complexité :] $O(n^4m)$.
  \end{itemize}
\end{algorithme}

\begin{preuve}
  Si $G$  est sans  prisme et  sans pyramide, alors  il est  clair que
  l'algorithme répond ``Ni prisme, ni pyramide''.
  
  Si  $G$ possède  une pyramide  ou un  prisme, alors  $G$  possède un
  prisme optimal ou  une pyramide optimale, que l'on  note $F$.  À une
  certaine  étape, l'algorithme  considère un  quadruplet  $(a_1, b_1,
  b_2, b_3)$,  tel que $\{b_1, b_2,  b_3\}$ est un triangle  de $F$ et
  $a_1  \notin \{b_1,  b_2, b_3\}$  est un  coin de  $F$. L'algorithme
  calcule alors  des plus  courts chemins $P_i$  de $a_1$  vers $b_i$.
  Par trois  applications consécutives des lemmes~\ref{reco.l.pyraopt}
  et~\ref{reco.l.prismeopt},  on  voit que  $V(P_1)  \cup V(P_2)  \cup
  V(P_3)$ induit  un prisme ou une  pyramide optimale de  $G$, de coin
  $a_1$  et  de  triangle  $\{b_1,  b_2,  b_3\}$.   Donc  l'algorithme
  retourne  ce sous-graphe, ou  un sous-graphe  \emph{ad hoc}  de même
  taille.
  
  Dans le pire des cas, tester tous les quadruplets nécessite $O(n^4)$
  étapes.  Calculer  les  plus  courts  chemins  et  vérifier  que  le
  sous-graphe obtenu est un  prisme ou une pyramide nécessite $O(m)$
  étapes.  Donc, au pire, l'algorithme se termine en $O(n^4m)$ étapes.
\end{preuve}

En remplaçant partout dans  l'énoncé ci dessus ``prisme'' par ``prisme
long'' et en invoquant  le lemme~\ref{reco.l.prismelongopt} à la place
du  lemme~\ref{reco.l.prismeopt},  on arrive  à  détecter les  prismes
longs dans les graphes sans pyramide~:

\index{long!détection~des~prismes~---}
\index{détection!prismes longs}
\index{prisme!détection~des~---~longs}

\begin{algorithme}\label{reco.a.prilong}
  \begin{itemize}
  \item[\sc Entrée :] Un graphe $G$ sans pyramide.

  \item[\sc Sortie :]  Si le graphe est sans  prisme long~: ``Pas de
    prisme long''.  Sinon~: un prisme long minimal de $G$.

  \item[\sc Calcul :] Similaire à l'algorithme précédent.

  \item[\sc Complexité :] $O(n^6)$.
  \end{itemize}
\end{algorithme}

Voici  maintenant l'algorithme le  plus rapide  que nous  ayons trouvé
pour  détecter les  prismes ou  les pyramides.  Notons que  l'usage du
lemme  suggéré par  Bruce  Reed (lemme~\ref{rr.l.bruce})  a permis  de
raccourcir très sensiblement sa preuve.

\index{détection!prismes ou pyramides}
\index{prisme!détection~des~---}
\index{pyramide!détection~des~---}

\begin{algorithme}\label{reco.a.rapidepypri}
  \begin{itemize}
  \item[\sc Entrée :] Un graphe $G$.

  \item[\sc  Sortie  :]  Si  le  graphe est  sans  pyramide  et  sans
    prisme~: ``Ni prisme ni pyramide''.  Sinon, ``Il y a un prisme ou
    une pyramide''.

  \item[\sc Calcul :] \mbox{}

    Pour tous les triangles $T = \{b_1, b_2, b_3\}$ de $G$ faire :
    
    \begin{enumerate}
      \setcounter{enumi}{-1}
    \item \label{reco.mic.rapidepypri0} 
      Calculer $X_1$, ensemble des sommets  de $G$ qui voient $b_1$ et
      manquent $b_2,  b_3$. Calculer de même $X_2,  X_3$. Calculer $X$
      l'ensemble des sommets de $G$ qui ne voient aucun sommet de $T$.
      Calculer les composantes connexes de $X$.

    \item \label{reco.mic.rapidepypri1} 
      Si pour  $i=1$, $2$ ou $3$, il
      existe  un  sommet  de $X_i$  voyant  à  la  fois un  sommet  de
      $X_{i+1}$ et un sommet  de $X_{i+2}$, alors stopper et retourner
      ``Il y  a un prisme  ou une pyramide''. (L'addition  des indices
      s'entend modulo 3.)

    \item
      \label{reco.mic.rapidepypri2} 
      On  marque ainsi composantes  connexes de  $X$. Donner  à chaque
      composante $H$ la marque $i$ s'il existe un sommet de $H$ voyant
      un  sommet de  $X_i$.  Si  une composante  $H$ reçoit  les trois
      marques $1,  2, 3$, stopper et  retourner ``Il y a  un prisme ou
      une pyramide''.

    \item
      \label{reco.mic.rapidepypri3}
      Pour tout $(i,j) \in \{1,  2, 3\}^2$, pour chaque composante $H$
      de $X$  possèdant les deux  marques $i,j$ et pour  chaque sommet
      $x$  de $X_i$  possèdant un  voisin dans  $H$, donner  à  $x$ la
      marque $j$. Pour tout $i\in  \{1, 2, 3\}$, et pour chaque sommet
      $x$ de  $X_i$ ayant un voisin  dans $X_j$ $(i\neq  j)$, donner à
      $x$  la marque  $j$. Si  un sommet  de $X_1  \cup X_2  \cup X_3$
      possède deux  marques différentes, stopper et répondre  ``Il y a
      un prisme ou une pyramide''.
    \end{enumerate}

    Si aucun triangle n'a permis de stopper pour retourner ``Il y a un
    prisme  ou une  pyramide'', stopper  et retourner  ``Ni  prisme ni
    pyramide''.
      
  \item
    [\sc Complexité :] $O(n^3(n+m))$.
  \end{itemize}
\end{algorithme}

\begin{preuve} Montrons d'abord~:

  \begin{claim}
    \label{reco.c.ppirep}
    Si $G$ possède  un prisme ou une pyramide  $F$, alors l'algorithme
    retourne ``Il y a un prisme ou une pyramide''.
  \end{claim}
  
  \begin{preuveclaim}
    Soit $\{b_1,  b_2, b_3\}$  un triangle de  $F$. On note  $c_1$ le
    voisin de  $b_1$ dans $F\setminus  \{b_2, b_3\}$ et on  définit de
    même $c_2$, $c_3$.  Pour $i  = 1, 2, 3$, notons que l'algorithme
    met $c_i$ dans $X_i$.

    Si $F$  n'a que 6  sommets, alors on  constate que l'un  de $c_1,
    c_2,  c_3$ voit  les deux  autres.  Donc  l'algorithme  stoppe à
    l'étape~\ref{reco.mic.rapidepypri1} et retourne ``Il y a un prisme
    ou une pyramide''.

    Si  $F$ a plus  de 6  sommets, et  si $F$  n'est pas  une pyramide
    possèdant  un  chemin de  longueur~1,  alors  on  constate que  $F
    \setminus \{b_1, b_2, b_3, c_1,  c_2, c_3\}$ est un graphe connexe
    qui voit $X_1,  X_2, X_3$.  Donc, l'algorithme met  les sommets de
    $F \setminus \{b_1, b_2, b_3, c_1, c_2, c_3\}$ dans une composante
    de   $X$    qui   reçoit   les   3   marques,    puis   stoppe   à
    l'étape~\ref{reco.mic.rapidepypri2} et retourne ``Il y a un prisme
    ou une pyramide''.
  
    Si $F$  a plus de  6 sommets  et si $F$  est une pyramide  avec un
    chemin de  longueur~1, alors  on suppose à  une symétrie  près que
    $c_1$ est le  coin de la pyramide. Si $c_2 c_1  \in E$ et $c_3
    c_1 \in  E$ alors  $F$ a six  sommets~: l'algorithme  aurait déjà
    stopper.  On  peut donc  supposer à une  symétrie près  que $c_2$
    manque $c_1$. Donc, dans le chemin de la pyramide reliant $b_2$ à
    $c_1$,  il  existe  des  sommets  entre  $c_2$  et  $c_1$,  qui
    appartiennent  à une composante  connexe $H$  de $X$  recevant les
    marques $1$ et $2$. Le sommet  $c_1$ est alors un sommet de $X_1$
    recevant la  marque~$2$ (car il a  un voisin dans  $H$). Si $c_1$
    voit $c_3$, alors $c_1$ reçoit la marque~3. Sinon, il existe des
    sommets  sur le chemin  de $F$  reliant $b_3$  à $c_1$,  qui sont
    entre  $c_3$ et  $c_1$, et  qui appartiennent  à  une composante
    connexe $H'$ de $X$ recevant  les marques $1$ et $3$. Donc, $c_1$
    reçoit   encore  la  marque~3.   Finalement,  $c_1$   reçoit  les
    marques~$2$     et~$3$,     donc     l'algorithme     stoppe     à
    l'étape~\ref{reco.mic.rapidepypri3} et retourne ``Il y a un prisme
    ou une pyramide''.
  \end{preuveclaim}

  \begin{claim}
    \label{reco.c.repipp}
    Si l'algorithme  retourne ``Il  y a un  prisme ou  une pyramide'',
    alors $G$ possède un prisme ou une pyramide.
  \end{claim}

  \begin{preuveclaim}
    Soit  $T  =  \{b_1,  b_2,  b_3\}$ le  triangle  à  l'étude  duquel
    l'algorithme stoppe.

    Si  l'algorithme   stoppe  à  l'étape~\ref{reco.mic.rapidepypri1},
    alors il existe à une symétrie  près un sommet $c_1$ de $X_1$, qui
    voit deux sommets  $c_{2}\in X_2$ et $c_{3} \in  X_3$.  Si $c_{2}$
    voit $c_{3}$, alors l'ensemble  $\{c_1, c_2, c_3, b_1, b_2, b_3\}$
    induit  un prisme.  Si  $c_{2}$ manque  $c_{3}$, alors  l'ensemble
    $\{c_1, c_2, c_3, b_1, b_2, b_3\}$ induit une pyramide.

    Si  l'algorithme   stoppe  à  l'étape~\ref{reco.mic.rapidepypri2},
    alors l'ensemble $X$ possède une composante connexe $H$ comprenant
    des  sommets marqués  $1$,  des sommets  marqués  $2$ et  d'autres
    marqués $3$.  Il existe dans $G$ trois sommets $c_1 \in X_1$, $c_2
    \in X_2$ et $c_3  \in X_3$, tel que pour $i= 1,  2, 3$, $c_i$ voit
    un sommet de $H$.   On applique alors le lemme~\ref{rr.l.bruce} au
    graphe induit par $\{b_1, b_2, b_3,  c_1, c_2, c_3\} \cup H$ et au
    triangle $\{b_1, b_2, b_3\}$.  On constate que $G$ possède bien un
    prisme ou une pyramide.

    Si  l'algorithme   stoppe  à  l'étape~\ref{reco.mic.rapidepypri3},
    alors, à une symétrie près, il existe un sommet $c_1$ de $X_1$ qui
    a des voisins marqués~$2$  et~$3$. Notons qu'il est impossible que
    $c_1$ ait des voisins dans  $X_2$ et $X_3$, car alors l'algorithme
    aurait  stoppé à l'étape~\ref{reco.mic.rapidepypri1}.   On suppose
    donc sans  perte de  généralité que $c_1$  n'a pas de  voisin dans
    $X_2$. Donc, $c_1$  a reçu la marque~2 en  raison d'une composante
    connexe  $H$ de  $X$  comprenant des  sommets  marqués~$1$ et  des
    sommets marqués~$2$. Notons qu'il n'y  a aucune arête entre $H$ et
    $X_3$, car  sinon $H$  aurait reçu les  3 marques  et l'algorithme
    aurait stoppé à l'étape~\ref{reco.mic.rapidepypri2}. Soit alors un
    chemin $P$  de $H$ reliant  un voisin $d_1$  de $c_1$ à  un sommet
    $d_2$ ayant un  voisin $c_2$ dans $X_2$, de  longueur minimal avec
    ces propriétés.   Si $c_1$  a reçu la  marque~$3$ en  raison d'une
    arête entre  $c_1$ et  un sommet $c_3$  de $X_3$, alors  les trois
    chemins $\bp b_1  \tp c_1\ep$, $\bp b_2 \tp c_2 \tp  d_2 \tp P \tp
    d_1 \tp  c_1 \ep$  et $\bp b_3  \tp c_3  \tp c_1 \ep$  forment une
    pyramide de triangle $\{b_1, b_2, b_3\}$ et de coin $c_1$.  Sinon,
    $c_1$ a reçu  la marque~3 en raison d'une  composante connexe $H'$
    de  $X$   comprenant  des  sommets  marqués~$1$   et  des  sommets
    marqués~$3$. Notons  qu'il n'y a  aucune arête entre $H$  et $H'$,
    car  ce sont des  composantes connexes  de $X$,  ni entre  $H'$ et
    $X_2$  car sinon  $H'$  aurait les  3  marques et  l'algorithme
    aurait stoppé à l'étape~\ref{reco.mic.rapidepypri2}. Soit alors un
    chemin $Q$  de $H'$ reliant un  voisin $e_1$ de $c_1$  à un sommet
    $e_3$ ayant un  voisin $c_3$ dans $X_3$, de  longueur minimal avec
    ces  propriétés.  Alors les  trois chemins  $\bp b_1  \tp c_1\ep$,
    $\bp b_2 \tp  c_2 \tp d_2 \tp P  \tp d_1 \tp c_1 \ep$  et $\bp b_3
    \tp c_3 \tp e_3 \tp Q \tp e_1 \tp c_1 \ep$ forment une pyramide de
    triangle $\{b_1, b_2, b_3\}$ et de coin $c_1$.
  \end{preuveclaim}

  D'après~(\ref{reco.c.ppirep}) et~(\ref{reco.c.repipp}), l'algorithme
  est conforme  à ses spécifications.  Énumérer tous  les triangles de
  $G$ nécessite au pire $O(n^3)$  étapes. Reste à vérifier que chacune
  des autres étapes  peut être réalisée par un  parcours de l'ensemble
  des  arêtes de  $G$.   Pour l'étape~\ref{reco.mic.rapidepypri0},  il
  suffit  pour chaque  arête de  vérifier si  l'une des  extrémité est
  $b_i$, et de marquer  l'autre extrémité ``$X_i$''. Les sommets ayant
  une et une  seule marque ``$X_i$'' seront les  élément de $X_i$. Les
  sommets n'ayant  reçu aucune marque  seront les éléments de  $X$ (on
  risque de devoir parcourir $V(G)$ pour les prendre en compte).  Pour
  les                                étapes~\ref{reco.mic.rapidepypri1}
  et~\ref{reco.mic.rapidepypri2},  un   parcours  des  arêtes  suffit.
  Notons qu'on peut  vérifier au fur et à mesure  des marquages si une
  composante  $H$ a  trois marques,  sans qu'il  y ait  besoin  de les
  reparcourir  après coup.   Pour l'étape~\ref{reco.mic.rapidepypri3},
  il suffit  pour chaque  arête $xy$  avec $x\in X_i$  et $y\in  H$ de
  marquer $x_i$. On  vérifie au fur et à mesure  si sommet $x$ possède
  deux marques.

  Donc, dans le pire des cas, l'algorithme nécessite $O(n^3(n+m))$ étapes.
\end{preuve}

On  constate  que les  deux  algorithmes  ci-dessus sont  relativement
performants  comparés   à  l'algorithme~\ref{reco.a.py}.   Ce  dernier
semble d'ailleurs difficile à améliorer.

\subsection{Prismes pairs}

Nous donnons ici un algorithme de détection des prismes pairs dans les
graphes sans trou impair. Nous  n'avons pas trouvé d'utilisation à cet
algorithme, mais il est potentiellement intéressant pour au moins deux
raisons~: les prismes pairs jouent un rôle important dans la preuve de
la   conjecture    forte   des   graphes   parfaits    ---   voir   le
théorème~\ref{graphespar.t.spgtprismepair}.   De  plus,  il peut  être
intéressant, en vue d'algorithmes futurs, de savoir détecter les seuls
sous-graphes à la  fois autorisés dans les graphes  d'Artémis pairs et
interdits dans les graphes d'Artémis.

On peut ici encore appliquer les méthodes de plus courts chemins, mais
cette fois,  le cheminement est  moins direct. On doit  considérer des
sommets ``au  milieu des chemins'', pour des  raisons qui apparaîtront
dans  la  preuve  du lemme~\ref{reco.l.prismepair}.   L'algorithme  de
Chudnovsky et Seymour~\cite{chudnovsky.seymour:reco} pour la détection
des pyramides utilise  cette même ruse des sommets  ``au milieu'', qui
ont pour conséquence fâcheuse mais sans doute inévitable de multiplier
par $n^3$ la complexité de notre algorithme.

\index{milieu~(---~d'un~chemin)}
Soit $P$  un chemin de longueur  paire et d'extrémités $a$  et $b$. On
appelle  \emph{milieu} de  $P$ l'unique  sommet $m$  de  $P$ vérifiant
$\lg(\bp a \tp P  \tp m\ep ) = \lg(\bp m \tp P  \tp b\ep )$.  Soit $F$
un  prisme pair  formé  par les  chemins  $P_1$, $P_2$,  $P_3$, et  de
triangles $\{a_1, a_2, a_3\}$, $\{b_1,  b_2, b_3\}$, de sorte que pour
$i=1,2,3$,  le chemin  $P_i$ soit  d'extrémités $a_i$  et  $b_i$. Soit
$m_i$  le sommet  au milieu  du  chemin $P_i$.   On dit  alors que  le
$9$-uplet  $(a_1, a_2,  a_3, b_1,  b_2, b_3,  m_1, m_2,  m_3)$  est la
\emph{trame}\index{trame!d'un~prisme~pair} de $F$.

\begin{lemme}
\label{reco.l.prismepair}
Soit $G$  un graphe  sans trou impair  et soit  $F$ un prisme  pair de
taille minimale de  $G$.  Supposons que $F$ soit  de trame $(a_1, a_2,
a_3, b_1,  b_2, b_3, m_1, m_2,  m_3)$ et formé par  les chemins $P_1$,
$P_2$, $P_3$ avec $a_i, m_i, b_i \in V(P_i)$ ($i= 1, 2, 3$).  Soit $R$
un chemin de $G$ d'extrémités  $a_1, m_1$, dont les sommets intérieurs
manquent $a_2,  a_3, b_2$ et $b_3$,  et de longueur  minimale avec ces
propriétés.

\noindent Alors, $\bp a_1 \tp R \tp  m_1 \tp P_1 \tp b_1 \ep$ est un
chemin de $G$  que l'on note $R_1$, et $R_1$,  $P_2$, $P_3$ forment un
prisme pair de $G$ de taille minimale.
\end{lemme}

\begin{preuve}
  Soit $k$ la longueur de  $P_1$, qui est paire par hypothèse.  Notons
  que $|E(R)| \le k/2$.  On note $H$ le trou de $G$ induit par $V(P_2)
  \cup  V(P_3)$.    Si  les  sommets  intérieurs   de  ${R}$  manquent
  entièrement $H$, alors soit $R_1$ un plus court chemin de $a_1$ vers
  $b_1$ contenu dans $V(R) \cup V(\bp m_1 \tp P_1 \tp b_1 \ep)$.  On a
  $|E(R_1)|\le k$ et les chemins $R_1$, $P_2$, $P_3$ forment un prisme
  $F'$ avec $|V(F')|  \le |V(F)|$.  Puisque $G$ est  sans trou impair,
  $R_1$ est de longueur paire (sinon $V(R_1) \cup V(P_2)$ induirait un
  trou  impair),  et $F'$  est  un  prisme  pair minimal.  L'inégalité
  ci-dessus est donc  une égalité, ce qui montre  qu'en fait $R_1$ est
  égal à $\bp a_1  \tp R \tp m_1 \tp P_1 \tp  b_1 \ep$~: la conclusion
  du lemme est satisfaite.

  On peut  donc supposer qu'il existe  un sommet intérieur  $c$ de $R$
  ayant des voisins dans $H$.   Nous allons montrer que l'existence de
  ce sommet  entraîne une contradiction.  On choisit  $c$ aussi proche
  que possible  de $m_1$. On va  voir que l'existence  de $c$ entraine
  une contradiction.  Soit $S$ un  plus court chemin de $c$ vers $b_1$
  contenu dans $V(\bp c \tp R \tp  m_1 \ep) \cup V(\bp m_1 \tp P_1 \tp
  b_1 \ep)$.   On a  $|E(S)|< k$ car  $|E(R)|\le k/2$ et  $c\neq a_1$.
  D'après la définition  de $c$, aucun sommet de  $S\setminus b_1$ n'a
  de voisin  dans $H$.  On définit  deux sous-chemins de  $H$~: $H_2 =
  H\backslash a_3$  et $H_3  = H \backslash  a_2$.  Pour $i=2,  3$, on
  note $c_i$ le sommet de $H_i$ voisin de $c$ le plus proche de $a_i$.

  \begin{claim}
    On peut supposer que $c_2c_3$ est une arête de $P_2$.
  \end{claim}

  \begin{preuveclaim}
    Si  $c_2 =  c_3$ alors  $V(S)\cup V(P_2)  \cup V(P_3)$  induit une
    pyramide de  triangle $\{b_1, b_2,  b_3\}$ et de coin  $c_2$~: une
    contradiction.  Donc $c_2 \neq c_3$.  Si $c_2$ manque $c_3$, alors
    $V(S) \cup  V(P_2)\cup V(P_3)$  contient une pyramide  de triangle
    $\{b_1, b_2, b_3\}$ et de coin $c$~: une contradiction.  Donc $c_2
    c_3$ est  bien une arête  de $G$, et,  quitte à échanger  $P_2$ et
    $P_3$, on peut supposer que $c_2 c_3$ est une arête de $P_2$.
  \end{preuveclaim}

  Donc $V(S)\cup V(P_2)\cup V(P_3)$ induit un prisme $F'$ de triangles
  $\{b_1,  b_2, b_3\}$  et $\{c,  c_2, c_3\}$  avec  $|V(F')|< |V(F)|$
  puisque  $|E(S)|<  k$.  Donc  $F'$  est  un  prisme impair,  ce  qui
  implique que $\bp c_3 \tp P_2 \tp b_2 \ep$ est un chemin de longueur
  impaire, et donc $\bp a_2 \tp P_2 \tp c_2 \ep$ un chemin de longueur
  paire.  Soit  $R'$ un  plus court chemin  de $c$ vers  $a_1$ contenu
  dans $V(\bp  c \tp R  \tp m_1  \ep) \cup V(\bp  m_1 \tp P_1  \tp a_1
  \ep)$.   On  a  $|E(R')|<   k$  car  $|E(R)|\le  k/2$.   D'après  la
  définition de $c$,  aucun sommet de $R'\setminus a_1$  n'a de voisin
  dans $H$.   Donc $R'$ est de  longueur paire car  sinon, $V(R') \cup
  V(\bp  a_2 \tp P_2  \tp c_2  \ep)$ induirait  un trou  impair.  Donc
  $V(R')\cup V(P_2)\cup  V(P_3)$ induit  un prisme $F''$  de triangles
  $\{a_1, a_2,  a_3\}$ et  $\{c, c_2, c_3\}$,  et $F''$ est  un prisme
  pair vérifiant  $|V(F'')|< |V(F)|$ car  $|E(R')|< k$.  Ceci  est une
  contradiction.
\end{preuve}

On  peut maintenant  donner  un algorithme  de  détection des  prismes
pairs~:

\index{détection!prismes~pairs}
\index{prisme!détection~des~---~pairs}
\begin{algorithme}\label{reco.a.prismepair}
  \begin{itemize}
  \item[\sc Entrée :] Un graphe $G$ sans trou impair.

  \item[\sc Sortie :]  Si le graphe est sans  prisme pair~: ``Pas de
    prisme pair''.  Sinon~: un prisme pair minimal de $G$.

  \item[\sc Calcul :] \mbox{}
    \begin{itemize}

    \item
      Considérer  successivement  tous  les  9-uplets  $(a_1$,  $a_2$,
      $a_3$, $b_1$, $b_2$, $b_3$, $m_1$, $m_2$, $m_3) \in V(G)^9$.
      
      Calculer pour  $i=1, 2, 3$ un  plus court chemin  $R_i$ de $a_i$
      vers  $m_i$ dont  les sommets  intérieurs manquent  le  reste du
      9-uplet, et   un plus court  chemin $S_i$ de  $m_i$ vers $b_i$
      dont les  sommets intérieurs manquent  le reste du  9-uplet.  Si
      l'ensemble $\cup_{i=1,2,3}(V(R_i) \cup V(S_i))$ induit un prisme
      pair, le stocker en mémoire.

    \item
      Si  aucun 9-uplet  n'a donné  de  prisme pair,  ``Pas de  prisme
      pair''. Sinon, retourner un prisme pair de taille minimale parmi
      ceux stockés en mémoire.

    \end{itemize}
  \item[\sc Complexité :] $O(n^{11})$.
  \end{itemize}
\end{algorithme}

\begin{preuve}
  Si $G$  est sans  prisme pair, alors  il est clair  que l'algorithme
  répond ``Pas de prisme pair''.
  
  Si  $G$ possède un  prisme pair,  alors $G$  possède un  prisme pair
  minimal,  que l'on  note $F$.   À une  certaine  étape, l'algorithme
  considère un 9-uplet $(a_1, a_2,  a_3, b_1, b_2, b_3, m_1, m_2, m_3)
  \in V(G)^9$ qui est la  trame de $F$. L'algorithme calcule alors des
  plus courts  chemins $R_i, S_i$.  Par  six applications consécutives
  du lemme~\ref{reco.l.prismepair}, on voit que $\cup_{i=1,2,3}(V(R_i)
  \cup V(S_i))$ induit un prisme pair.  Donc l'algorithme retourne ce
  sous-graphe, ou un sous-graphe \emph{ad hoc} de même taille.
  
  Dans le  pire des cas,  tester tous les 9-uplets  nécessite $O(n^9)$
  étapes.  Calculer  les  plus  courts  chemins  et  vérifier  que  le
  sous-graphe obtenu est un  prisme ou une pyramide nécessite $O(n^2)$
  étapes.   Donc,  au pire,  l'algorithme  se  termine en  $O(n^{11})$
  étapes.
\end{preuve}


\subsection{Line-graphes de  subdivisions de $K_4$}
\index{line-graphe~de~subdivision~bipartie~de~$K_4$!détection~des~---}
\index{détection!line-graphes~de~subdivision~bipartie~de~$K_4$}

Les  line-graphes de subdivisions  biparties de  $K_4$ jouent  un rôle
important dans la preuve du théorème des graphes parfaits --- voir par
exemple le théorème~\ref{graphespar.t.spgtlsbk4}.  Détecter ce type de
sous-graphe peut donc  être intéressant en soi.  Ensuite,  et cela fut
notre motivation initiale, l'algorithme donné ci-dessous est essentiel
pour la détection des prismes  impairs, et donc pour la reconnaissance
des graphes d'Artémis pairs.

\index{line-graphe~de~subdivision~bipartie~de~$K_4$!figure}
\begin{figure}
  \center  \includegraphics{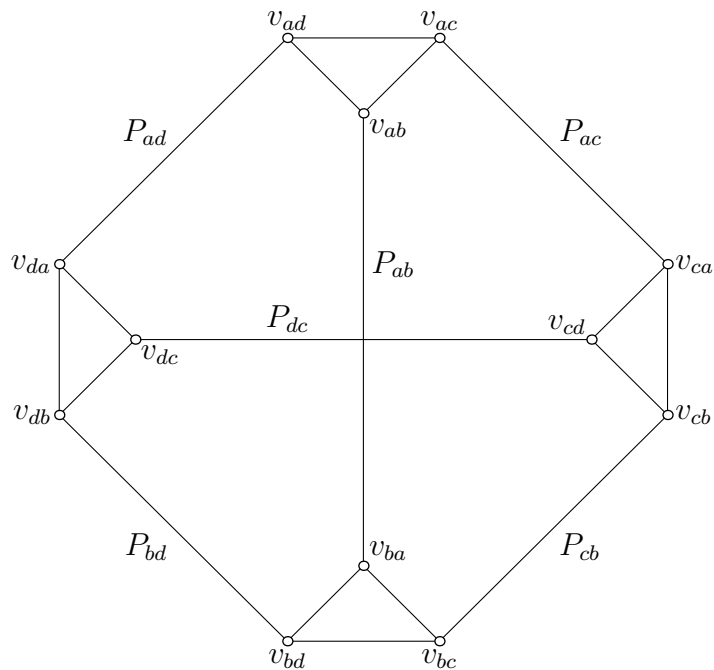}
  \caption{Line-graphe d'une  subdivision bipartie de  $K_4$\label{reco.f.lK4}}
\end{figure}

Pour ce qui suit,  on pourra consulter la figure~\ref{reco.f.lK4}.  Si
$F$ est le line-graphe d'une  subdivision bipartie $R$ de $K_4$, alors
on convient de noter $a$, $b$,  $c$ et $d$ les quatre sommets du $K_4$
correspondant. On ordonne arbitrairement $a$, $b$, $c$ et $d$~: $a < b
< c  < d$.  À chaque arête  $ij$, $i<j$, de $K_4$,  correspond donc un
chemin de $F$  qu'on note $P_{ij}$.  On note  $v_{ij}$ et $v_{ji}$ les
extrémités de  $P_{ij}$, de  sorte que pour  $i =  a,b,c,d$ l'ensemble
$T_i = \{ v_{ij}; j\in \{a, b, c, d\} \text{ , } j\neq i\}$ induise un
triangle.  Finalement,  on  obtient  que  les  six  chemins  $P_{ab}$,
$P_{ac}$, $P_{ad}$, $P_{bc}$, $P_{bd}$, $P_{cd}$, et les douze sommets
$v_{ab}$, $v_{ac}$, $v_{ad}$,  $v_{ba}$, $v_{bc}$, $v_{bd}$, $v_{ca}$,
$v_{cb}$, $v_{cd}$, $v_{da}$, $v_{db}$ et $v_{dc}$ vérifient~:

\vspace{1ex}

\begin{itemize}
\item
  Les chemins $P_{ab}$,  $P_{ac}$, $P_{ad}$, $P_{bc}$,
  $P_{bd}$, $P_{cd}$ sont disjoints.
\item 
  $F =  P_{ab} \cup  P_{ac} \cup P_{ad}  \cup P_{bc} \cup  P_{bd} \cup
  P_{cd}$
\item 
  Pour $i  \in\{ a,b,c,d \}$,  $j \in \{a,  b, c, d  \}$, $i <  j$, le
  chemin $P_{ij}$ a pour extrémités $v_{ij}$ et $v_{ji}$.
\item
  Les  ensembles $\{  v_{ab}, v_{ac},  v_{ad}\}$, $\{  v_{ba}, v_{bc},
  v_{bd}\}$,  $\{ v_{ca},  v_{cb},  v_{cd}\}$ et  $\{ v_{da},  v_{db},
  v_{dc}\}$ sont des triangles de $F$.
\item
  Les  seules arêtes  de  $F$  sont celles  des  chemins et  triangles
  précédement décrits.
\end{itemize}

\vspace{2ex}

On remarquera que certains des chemins peuvent être de longueur nulle,
ce qui correspond  à une arête de $K_4$  non subdivisée.  Par exemple,
$P_{ab}$  pourrait être  de  longueur nulle,  auquel  cas les  sommets
$v_{ab}$ et  $v_{ba}$ seraient confondus.  Modulo le  choix initial de
l'ordre  sur  $\{a,b,c,d\}$, les  sommets  et  chemins ci-dessus  sont
définis  sans  ambiguïté.   Les  sommets  $v_{ij}$  sont  appelés  les
\emph{coins} de $F$.  Les  chemins $P_{ij}$ sont appelés \emph{chemins
de base} de  $F$.  Chaque sommet $i$ du $K_4$  donne un triangle $T_i$
de $F$.  On  appelle ces triangle \emph{triangles de  base} de $F$. Il
faut  noter que  $F$  peut posséder  d'autres  triangles, lorsque  des
arêtes du $K_4$ de départ ne sont pas subdivisées.

Les trois lemmes à suivre sont  bien connus mais ne sont publiés nulle
part sous la forme exacte dont nous avons besoin.  C'est pourquoi nous
donnons  leur preuve  pour ne  rien  laisser dans  l'ombre.  Le  lemme
suivant montre  en quoi la détection des  line-graphes de
subdivision  bipartie  de $K_4$  peut  aider  à  détecter des  prismes
impairs.

\index{line-graphe~de~subdivision~bipartie~de~$K_4$!contient~un~prisme~impair}
\begin{lemme}
  \label{reco.l.LSBK4->prismeimpair}
  Soit $G$ le line-graphe  d'une subdivision bipartie de $K_4$. Alors
  $G$ contient un prisme impair.
\end{lemme}

\begin{preuve}
  On suppose sans perte de généralité que $a$ et $b$ sont du même côté
  de  la  bipartition  du  $K_4$.   Donc l'arête  $ab$  du  $K_4$  est
  subdivisée en un  chemin de longueur paire.  Les  arêtes incidente à
  $a$ et $b$ donnent, dans le line-graphe $G$, les deux triangles d'un
  prisme, et  le chemin de $a$  vers $b$ dans la  subdivision de $K_4$
  donne, dans  le line-graphe, un  chemin de longueur  impaire.  Comme
  $G$ est de Berge, on sait que $F$ est un prisme impair.
\end{preuve}

Le  lemme suivant  nous  permettra par  la  suite de  montrer que  les
subdivisions de $K_4$ que  nous détecterons sont bien des subdivisions
\emph{biparties}.   On  dit  qu'une   subdivision  $R$  de  $K_4$  est
\index{triviale}  \index{subdivision!---~triviale}  \emph{triviale} si
$R=K_4$.      Le    line-graphe     de     $K_4$    est     représenté
figure~\ref{base.fig.k4} page~\pageref{base.fig.k4}.

\begin{lemme} 
  \label{reco.l.techK4}
  Soit $R$  une subdivision quelconque de  $K_4$.  Alors on  a l'un et
  l'un seulement des cas suivants~:
  \begin{itemize}
  \item
    $R$ est triviale, c'est-à-dire $R = K_4$.
  \item 
    $L(R)$ contient un trou impair.
  \item
    $R$ est une subdivision bipartie de $K_4$.
  \end{itemize}
\end{lemme}
  
\begin{preuve}
  Si $R$  est triviale, alors $R= K_4$  et il est clair  que $L(R)$ ne
  contient pas de  trou impair et que $R$ n'est  pas biparti.  On peut
  donc supposer  $R \neq  K_4$.  Il est  clair également que  les deux
  dernières  suppositions ne  peuvent être  satisfaites simultanément,
  car les line-graphes de graphes  bipartis sont des graphes de Berge.
  Il reste  juste à prouver  qu'on ne peut avoir  simultanément $L(R)$
  sans trou impair  et $R$ non biparti.  Supposons  donc, en vue d'une
  contradiction, qu'il en soit ainsi.

  Appelons $a$, $b$,  $c$ et $d$ les quatre sommets  du $K_4$, et pour
  $i<j \in \{a,  b, c, d\}$ notons $C_{ij}$  la subdivision de l'arête
  $ij$.  Puisque $L(R)$ est sans trou impair et $R$ non biparti, alors
  $R$  contient un cycle  impair (non  nécessairement induit).   Si ce
  cycle impair  n'est pas un triangle,  alors il donne  dans $L(R)$ un
  trou  impair~: une contradiction.   Si c'est  un triangle,  alors on
  peut supposer sans perte de  généralité que ce triangle est $\{a, b,
  c\}$. Comme $R\neq K_4$, on sait  qu'au moins une arête du $K_4$ est
  subdivisée, par  exemple $ad$ sans perte de  généralité.  Mais alors
  l'un  des  ensembles  $E(C_{ad})  \cup  \{cd\}  \cup  E(C_{cd})$  et
  $E(C_{ad}) \cup  \{bd\} \cup \{bc\} \cup C_{cd}$  est l'ensemble des
  arêtes d'un  cycle impair de $R$  de longueur au moins  5, qui donne
  dans $L(R)$ un trou impair, encore une contradiction.
\end{preuve}

En combinant les deux lemmes précédents, on obtient~:

\begin{lemme} 
  \label{reco.l.techK4syn}
  Soit $R$  une subdivision quelconque de  $K_4$.  Alors ou bien~:
  \begin{itemize}
  \item
    $R$ est triviale, c'est-à-dire $R = K_4$.
  \item 
    $L(R)$ contient un trou impair.
  \item
    $L(R)$ contient un prisme impair.
  \end{itemize}
\end{lemme}

Pour des raisons qui apparaîtront dans la preuve du prochain lemme, on
a besoin  de prendre en compte  les sommets au milieu  des chemins des
line-graphes.   Or ces chemins  peuvent être  de longueur  impaire, et
n'avoir  pas de ce  fait de  sommet central.  On est  donc amené  à la
définition suivante~:

Soit $F$ le line-graphe d'une  subdivision de $K_4$ avec les notations
ci-dessus. Soient $m_{ab}$,  $m_{ac}$, $m_{ad}$, $m_{bc}$, $m_{bd}$ et
$m_{cd}$ six  sommets de $F$. Si  pour tout $i<j$,  le sommet $m_{ij}$
est   dans   $P_{ij}$  et   vérifie   $d_{P_{ij}}(v_{ij},  m_{ij})   -
d_{P_{ij}}(v_{ji},  m_{ij}) \in  \{-1, 0,  1\}$, alors  on dit  que le
18-uplet $(v_{ab},  v_{ac}, \dots, v_{cd}, m_{ab},  \dots, m_{cd}) \in
V^{18}$                             est                            une
\emph{trame}\index{trame!d'un~line~graphe~de~subdivision~de~$K_4$}   de
$F$.  Les six sommets $m_{ab}$, $m_{ac}$, $m_{ad}$, $m_{bc}$, $m_{bd}$
et $m_{cd}$ s'interprètent comme des sommets ``proches des milieux des
chemins''   $P_{ab}$,  $P_{ac}$,   $P_{ad}$,  $P_{bc}$,   $P_{bd}$  et
$P_{cd}$.  Voici maintenant le lemme principal~:

\begin{lemme} \label{reco.l.subK4}
  Soit $G$ un graphe sans  pyramide. Soit $F$ un sous-graphe induit de
  $G$ qui est une subdivision non triviale de $K_4$. On suppose $F$ de
  taille minimale et de trame $(v_{ab}, v_{ac}, \dots, v_{cd}, m_{ab},
  \dots,  m_{cd})$. Soit  $P$ un  plus court  chemin de  $v_{ab}$ vers
  $m_{ab}$  dont les  sommets  intérieurs manquent  tous les  $v_{ij}$
  autres que $v_{ab}$.  Alors~:

  \noindent L'ensemble  $(F \backslash  v_{ab} P_{ab} m_{ab})  \cup P$
  induit  un line-graphe d'une  subdivision non  triviale de  $K_4$ de
  taille minimale.
  
  \noindent On  a un résultat  similaire pour les cinq  autres couples
  $(v_{ij}, m_{ij})$ et les six couples $(v_{ji}, m_{ij})$.
\end{lemme}

\begin{figure}
  \center
  \includegraphics{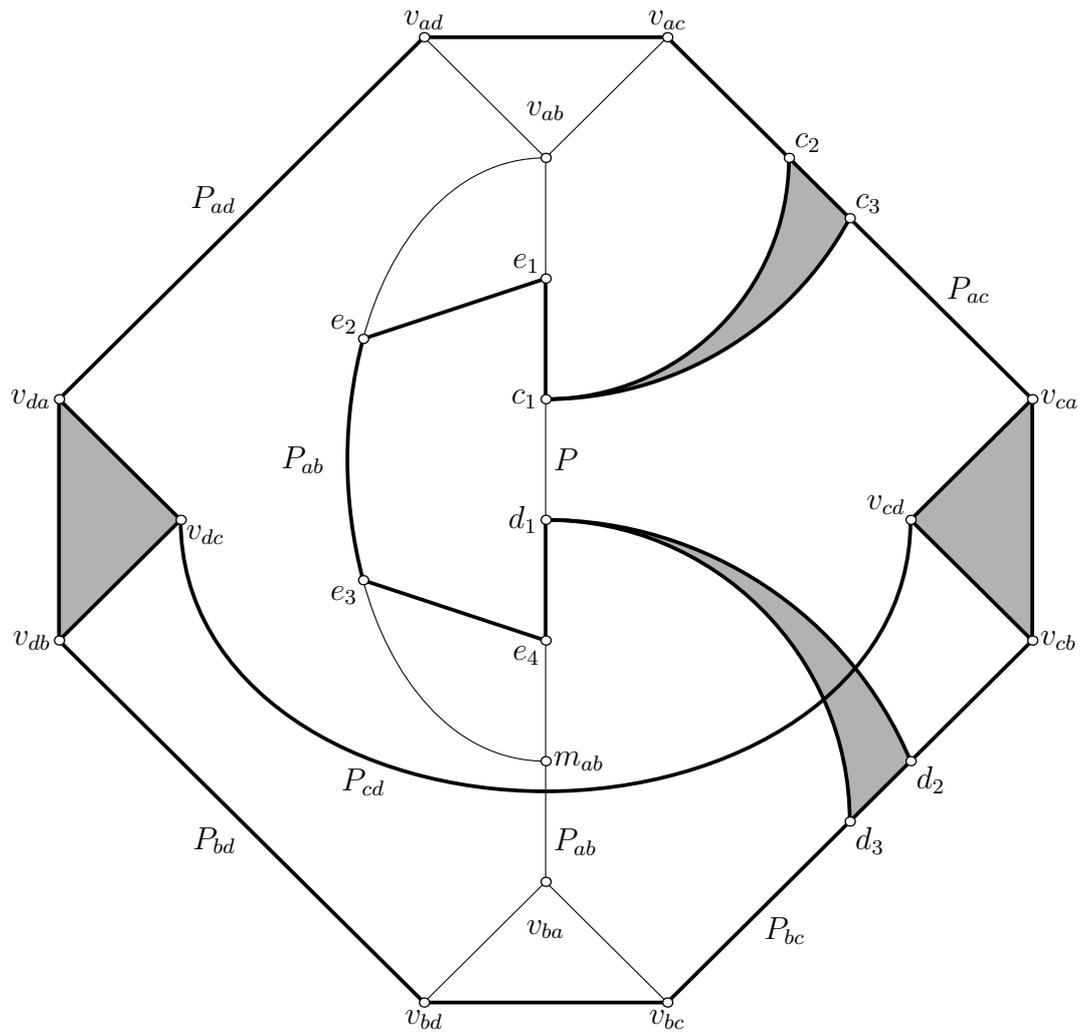}
  \caption{$F$ et $P$ (preuve  de \ref{reco.l.subK4}).  Les chemins en
    gras et les triangles grisés constituent le sous-graphe induit $E$
    défini avant (\ref{reco.c.subK4.EF}).} 
  \label{reco.f.preuveLK4}
\end{figure}

\begin{preuve}
  On     pourra     consulter     la     figure~\ref{reco.f.preuveLK4}
  page~\pageref{reco.f.preuveLK4}.   On   note  $F'  =   F  \backslash
  P_{ab}$.  Si $v_{ab} = m_{ab}$,  ou si $v_{ab}$ voit $m_{ab}$, alors
  la conclusion  du lemme est  trivialement satisfaite.  On  peut donc
  supposer que $v_{ab}$ et $m_{ab}$  sont des sommets distincts et non
  adjacents, ce qui implique en outre $m_{ab} \neq v_{ba}$.

  \begin{claim} \label{reco.c.subK4.sP}
    Si  les sommets  intérieurs de  $P$ manquent  entièrement  $F'$, la
    conclusion du lemme est satisfaite.
  \end{claim}

  \begin{preuveclaim}
    Soit $u$  le sommet de $\bp v_{ab}  \tp P \tp m_{ab}  \ep$ le plus
    proche  de $v_{ab}$  et ayant  des  voisins dans  $\bp m_{ab}  \tp
    P_{ab}  \tp v_{ba}  \ep$. Soit  $u'$ le  voisin de  $u$  dans $\bp
    m_{ab}  \tp P_{ab}  \tp v_{ba}  \ep$  le plus  proche   de
    $v_{ba}$.  Alors,  $\bp v_{ab} \tp P  \tp u \tp u'  \tp P_{ab} \tp
    v_{ba}  \ep$ est  un chemin  sans corde  que l'on  note  $P'$.  On
    constate  que  $V(P')  \cup  V(F')$ induit  le  line-graphe  d'une
    subdivision non triviale de $K_4$.   Donc $V(P') \cup V(F')$ a une
    taille supérieure  ou égale  à celle de  $F$, ce qui  est possible
    seulement si  $u=m_{ab}$ car $P$  est un plus court  chemin.  Mais
    dans ce cas,  $(V(F) \setminus V(\bp v_{ab} \tp  P_{ab} \tp m_{ab}
    \ep)  \cup  V(P)$  induit  le line-graphe  d'une  subdivision  non
    triviale de  $K_4$ de taille  minimale, et la conclusion  du lemme
    est bien satisfaite.
  \end{preuveclaim}

  On peut donc supposer l'existence d'un sommet $c_1$ de $P$ ayant des
  voisins  dans   $F'$,  et  choisi  aussi  proche   que  possible  de
  $v_{ab}$. On définit de même  $d_1$, sommet de $P$ ayant des voisins
  dans $F'$,  et choisi aussi  proche que possible de  $m_{ab}$.  Nous
  allons   voir  que   l'existence   de  ces   sommets  entraîne   une
  contradiction.


  \begin{claim} \label{reco.c.subK4.Nc_1}
    $N(c_1) \cap V(F')$  est une arête de $F'$  et $N(d_1) \cap V(F')$
    est une arête de $F'$
  \end{claim}

  \begin{preuveclaim}
    On  pose $H  = P_{ac}  \cup P_{bc}  \cup P_{bd}  \cup  P_{ad}$. On
    suppose d'abord  que le  sommet $c_1$ a  des voisins dans  le trou
    $H$.

    On définit deux chemins de  $H$~: $H_{ac} = H \setminus v_{ad}$ et
    $H_{ad}  = H  \setminus v_{ac}$.  On  définit $c_2$  le sommet  de
    $H_{ac}$  le plus  proche  de  $v_{ac}$ et  voisin  de $c_1$.   On
    définit $c_3$ le sommet de  $H_{ad}$ le plus proche de $v_{ad}$ et
    voisin de $c_1$.  Si  $c_2=c_3$, alors l'ensemble $V(H) \cup V(\bp
    c_2 \tp  c_1 \tp P  \tp v_{ab} \ep)$  induit une pyramide  de coin
    $c_2$  et   de  triangle   $\{  v_{ab},  v_{ac},   v_{ad}\}$~:  une
    contradiction.   Donc $c_2  \neq  c_3$.  Si  $c_2$ manque  $c_3$,
    alors on  suppose sans perte  de généralité que  $v_{ac}$, $c_2$,
    $c_3$ et $v_{ad}$ apparaissent dans cet ordre sur $H$.  Et on voit
    que les trois chemins $\bp c_1 \tp P \tp v_{ab} \ep$, $\bp c_1 \tp
    c_2 \tp H_{ac} \tp v_{ac} \ep$  et $\bp c_1 \tp c_3 \tp H_{ad} \tp
    v_{ad} \ep$  induisent une pyramide  de coin $c_1$ et  de triangle
    $\{ v_{ab}, v_{ac}, v_{ad}\}$~: une contradiction.  Donc le sommet
    $c_1$  a pour seuls  voisins dans  $H$ les  deux sommets  $c_2$ et
    $c_3$  qui sont  en  outre  adjacents. On  suppose  sans perte  de
    généralité que $c_2 c_3$ est une arête de $P_{ac} \cup P_{bc}$.

    Si $c$  possède, en plus de  ses voisins sur $H$,  des voisins sur
    $P_{cd}$, alors on  note $c_4$ le voisin de  $c_1$ sur $P_{cd}$ le
    plus proche de  $v_{dc}$.  On voit que les chemins  $\bp c_1 \tp P
    \tp v_{ab} \ep$,  $\bp c_1 \tp c_4 \tp P_{cd}  \ep \cup P_{ad}$ et
    $\bp c_1 \tp c_2 \tp H_{ac} \tp v_{ac} \ep$ induisent une pyramide
    de coin  $c_1$ et de  triangle $\{v_{ab}, v_{ac},  v_{ad}\}$~: une
    contradiction.  Donc $c_1$ manque $P_{cd}$.

    Finalement, si $c_1$  voit des sommets de $H$ alors  il ne voit de
    $F'$        qu'une        arête,        ainsi        que        le
    prévoit~(\ref{reco.c.subK4.Nc_1}). On peut donc supposer que $c_1$
    manque $H$.  Dans ce cas, d'après~(\ref{reco.c.subK4.sP}), on peut
    supposer que $c_1$ a des voisins dans $P_{cd}^*$.

    Si $c_1$  n'a qu'un seul  voisin $c\in P_{cd}$, alors  les chemins
    $\bp c \tp c_1 \tp P \tp v_{ab} \ep$, $\bp c \tp P_{cd} \tp v_{cd}
    \ep \cup P_{ac}$ et $\bp c  \tp P_{cd} \tp v_{dc} \ep \cup P_{ad}$
    induisent  une pyramide  de coin  $c$ et  de triangle  $\{ v_{ab},
    v_{ac},  v_{ad}\}$~:  une contradiction.   Donc  $c_1$ a  plusieurs
    voisins dans  $P_{cd}^*$. On  note $c_2$ le  voisin de  $c_1$ dans
    $P_{cd}^*$ le plus proche de  $v_{cd}$ et $c_3$ le voisin de $c_1$
    dans $P_{cd}^*$ le  plus proche de $v_{dc}$. Si  $c_2$ et $c_3$ ne
    sont pas  adjacents, alors les chemins  $\bp c_1 \tp  P \tp v_{ab}
    \ep$, $\bp c_1  \tp c_2 \tp P_{cd} \tp v_{cd}  \ep \cup P_{ac}$ et
    $\bp c_1 \tp c_3 \tp  P_{cd} \tp v_{dc} \ep \cup P_{ad}$ induisent
    une  pyramide de  coin $c_1$  et de  triangle $\{  v_{ab}, v_{ac},
    v_{ad}\}$~: une contradiction.  Donc  les sommets $c_2$ et $c_3$ se
    comportent  encore ainsi  que l'annonce~(\ref{reco.c.subK4.Nc_1}).
    La preuve pour $d_1$ est similaire.
  \end{preuveclaim}

  On note  donc $c_2 ,  c_3$ les deux  voisins de $c_1$ dans  $F'$, et
  $d_2, d_3$ les deux voisins de $d_1$ dans $F'$. Avant de poursuivre,
  il faut  vérifier que  les line-graphes que  nous allons  trouver ne
  dégénèrent pas en $L(K_4)$~:

  \begin{claim} \label{reco.c.subK4.degenere}
    Soit  $J$ un sous-graphe  de $F'  \cup P$  qui est  le line-graphe
    d'une  subdivision  \emph{quelconque}   de  $K_4$.   Si  l'un  des
    triangles de $J$ contient $c_1$, alors $J$ est non trivial.
  \end{claim}

  \begin{preuveclaim}
    En général, pour vérifier  qu'un line-graphe $J$ d'une subdivision
    de $K_4$ est non trivial, il suffit de vérifier que $J$ possède un
    coin qui n'appartient qu'à un seul triangle de base de $J$.  C'est
    bien le cas ici~: $J$ ne peut pas être le line-graphe de $K_4$ car
    $c_1$ ne peut appartenir qu'à un seul triangle de $J$.
  \end{preuveclaim}

  \begin{claim} \label{reco.c.subK4.c_2c_3}
    On peut supposer  que $c_2 c_3$ est une arête  de $P_{ac}$ et $d_2
    d_3$ une arête de $P_{cb} \cup P_{bd}$.
  \end{claim}

  \begin{preuveclaim}
    Si $c_2  c_3$ est une  arête de $P_{cd}$, alors  l'ensemble $V(\bp
    c_1  \tp P  \tp v_{ab}  \ep)  \cup V(P_{ac})  \cup V(P_{ad})  \cup
    V(P_{cd})$  induit le line-graphe  d'une subdivision  bipartie de
    $K_4$ de  taille strictement  inférieure à celle  de $F$,  qui est
    non-triviale     d'après~(\ref{reco.c.subK4.degenere})~:    une
    contradiction.   Si $c_2  c_3$ est  une arête  de  $P_{bc}$, alors
    l'ensemble $V(\bp c_1 \tp P  \tp v_{ab} \ep) \cup V(F')$ induit le
    line-graphe  d'une  subdivision   bipartie  de  $K_4$  de  taille
    strictement  inférieure  à  celle  de $F$,  qui  est  non-triviale
    d'après~(\ref{reco.c.subK4.degenere})~:une contradiction.

    Donc $c_2 c_3$ est une  arête de $P_{ac} \cup P_{ad}$.  On montre
    de même que $d_2 d_3$ est  une arête de $P_{bc} \cup P_{bd}$.  Par
    symétrie et  sans perte de  généralité, on peut supposer  que $c_2
    c_3$ est une arête de $P_{ac}$.
  \end{preuveclaim}

  On suppose  sans perte de  généralité que $v_{ac}$, $c_2$,  $c_3$ et
  $v_{ad}$ apparaissent dans cet  ordre sur $H$. On suppose également
  que $d_2$, $d_3$ et $v_{ad}$ apparaissent dans cet ordre sur $H$.

  \begin{claim} \label{reco.c.subK4.c_1d_1}
    $c_1$ et $d_1$ sont des sommets distincts et non adjacents.
  \end{claim}

  \begin{preuveclaim}
    Par  (\ref{reco.c.subK4.Nc_1}) et  (\ref{reco.c.subK4.c_2c_3}), on
    sait  que $c_1$  et $d_1$  sont distincts.   Si $c_1$  voit $d_1$,
    alors l'ensemble  $\{c_1, d_1\} \cup V(H')$  induit un line-graphe
    d'une  subdivision   bipartie  de  $K_4$   de  taille  strictement
    inférieure       à      $F$,       qui       est      non-triviale
    d'après~(\ref{reco.c.subK4.degenere})~:   une  contradiction.   On
    sait donc que $c_1$ manque $d_1$.
  \end{preuveclaim}

  On note  alors $e_1$  le sommet de  $\bp c_1  \tp P \tp  v_{ab} \ep$
  ayant des voisins  dans $(\bp m_{ab} \tp P_{ab}  \tp v_{ab} \ep)^*$,
  et  choisi aussi  proche que  possible de  $c_1$. On  note  $e_4$ le
  sommet de  $\bp d_1  \tp P  \tp m_{ab} \ep$  ayant des  voisins dans
  $(\bp m_{ab} \tp  P_{ab} \tp v_{ab} \ep)^*$, et  choisi aussi proche
  que possible de $d_1$.  On choisit alors $e_2$ dans $(\bp m_{ab} \tp
  P_{ab}  \tp v_{ab}  \ep)^*$, voisin  de $e_1$,  et $e_3$  dans $(\bp
  m_{ab} \tp P_{ab} \tp v_{ab}  \ep)^*$, voisin de $d_1$.  Les sommets
  $e_2$  et $e_3$  sont choisis  aussi  proches que  possible l'un  de
  l'autre.

  \begin{claim} \label{reco.c.subK4.e_1}
    $e_1 \neq v_{ab}$.
  \end{claim}

  \begin{preuveclaim}
    Sinon, les chemins $\bp v_{ab} \tp P \tp c_1 \ep$, $\bp v_{ab} \tp
    v_{ac} \tp P_{ac}  \tp c_2 \ep$ et $\bp v_{ab}  \tp P_{ab} \tp e_3
    \tp e_4 \tp P \tp d_1 \tp  d_2 \tp H_{ac} \tp c_3 \ep$ forment une
    pyramide de coin $v_{ab}$ et de triangle $\{c_1, c_2, c_3\}$.
  \end{preuveclaim}

  Les  critères pour  le  choix des  sommets  ci-après impliqués,  les
  propriétés  (\ref{reco.c.subK4.c_1d_1})  et (\ref{reco.c.subK4.e_1})
  et le  fait, signalé  en début de  preuve, que $m_{ab}$  et $v_{ba}$
  sont distincts,  montrent que  $\bp c_1  \tp P \tp  e_1 \tp  e_2 \tp
  P_{ab} \tp e_3 \tp e_4 \tp P  \tp d_1 \ep$ est un chemin de $G$ dont
  les sommets  intérieurs manquent  $F'$. On note  ce chemin  $P'$. On
  constate que l'ensemble $E = V(P') \cup V(F')$ induit le line-graphe
  d'une  subdivision   non  triviale   de  $K_4$  (non   triviale  par
  (\ref{reco.c.subK4.degenere})).

  \begin{claim}  \label{reco.c.subK4.EF}
    La  taille  de $E$  est  strictement  inférieure  à la  taille  de
    $F$.
  \end{claim}

  \begin{preuveclaim}
    Il  suffit de  montrer que  les  six chemins  de base  de $E$  ont
    strictement moins d'arêtes que les six chemins de base de $F$.  On
    note $\alpha$ le nombre d'arêtes  du chemin $\bp v_{ab} \tp P_{ab}
    \tp  m_{ab} \ep$.   On note  $\beta$  le nombre  d'arêtes de  $\bp
    v_{ba} \tp P_{ab} \tp m_{ab} \ep$.  On note $\delta$ le nombre des
    arêtes $e$ de $F'$ appartenant à des chemins de base $F$.

    Le  nombre d'arêtes des  chemins de  $E$ est  inférieur ou  égal à
    $\delta +2 \alpha  -3 $.  Cette borne correspond  au pire des cas,
    celui où $e_4$ et $m_{ab}$  sont confondus, où il existe un unique
    sommet de $P_{ab}$  entre $c_1$ et $d_1$, où  $e_1$ voit $v_{ab}$,
    où $e_2$  voit $v_{ab}$ et  où $P$ et  $\bp v_{ab} \tp  P_{ab} \tp
    m_{ab}  \ep$  sont de  même  longueur.   Dans  ce cas,  le  nombre
    d'arêtes du chemin de $E$ allant  de $c_1$ à $d_1$ est bien égal à
    $2 \alpha -3 $.

    Le nombre d'arêtes des chemins de $F$ est égal à $\alpha + \beta +
    \delta =  2 \alpha -  \varepsilon + \delta$ avec  $\varepsilon \in
    \{0,1\}$  car,  d'après  la  définition  de  $m_{ab}$,  il  existe
    $\varepsilon$ tel que $\alpha = \beta + \varepsilon$.  Et donc le
    nombre  d'arêtes  des  chemins  de  base de  $E$  est  strictement
    inférieur à celui des chemins de  base de $F$.  Donc la taille de
    $E$ est strictement inférieure à celle de $F$.
  \end{preuveclaim}

  Nous arrivons donc à une contradiction.
\end{preuve}

On peut maintenant donner  un algorithme de détection des line-graphes
de subdivision non triviale de  $K_4$ dans les graphes sans pyramide.
Il  faut  noter  que  cette  restriction  aux  subdivisions  \emph{non
triviales}  n'est  aucunement   indispensable  pour  que  l'algorithme
fonctionne correctement.  Nous avons choisi ce point de vue restrictif
pour pouvoir facilement traiter ensuite les line-graphes de subdivision
\emph{bipartie} de  $K_4$ (rappelons que $K_4$ est  un graphe quelque
peu pathologique, en ce qu'il est la seule sudivision non bipartie de
$K_4$ dont le line-graphe est de Berge, voir figure~\ref{base.fig.k4}
page~\pageref{base.fig.k4}).

\begin{algorithme}\label{reco.a.subK4}
  \begin{itemize}
  \item[\sc Entrée :] Un graphe $G$ sans pyramide.
    
  \item[\sc Sortie :] Si le graphe est sans line-graphe de subdivision
    non triviale  de $K_4$~: ``Pas de  LGSNTK4''.  Sinon, l'algorithme
    retourne un line-graphe de subdivision non triviale de $K_4$.
  \item[\sc Calcul :] \mbox{}
    \begin{itemize}
    \item
      Considérer   successivement   tous   les  18-uplets   $(v_{ab}$,
      $v_{ac}$,  $\dots$, $v_{cd}$,  $m_{ab}$,  $\dots$, $m_{cd})  \in
      V^{18}$.
      
      Calculer pour tout $i\in\{a, b,  c, d\}$ et tout $j\in\{a, b, c,
      d\}$ tel  que $i<j$  un plus court  chemin $P_{ij}$  de $v_{ij}$
      vers $m_{ij}$  dont les sommets intérieurs manquent  le reste du
      18-uplet,  et un  plus court  chemin $Q_{ij}$  de  $v_{ji}$ vers
      $m_{ij}$  dont  les  sommets  intérieurs manquent  le  reste  du
      18-uplet.

      Si $\cup_{i,j} (V(P_{ij}) \cup V(Q_{ij}))$ induit le line-graphe
      d'une  subdivision  non  triviale  de $K_4$,  alors  stocker  ce
      sous-graphe en mémoire et considérer les 18-uplets suivants.
     
    \item
      Retourner l'un des sous-graphes  stockés de taille minimale.  Si
      aucun  18-uplet  n'a donné  de  line-graphe  de subdivision  non
      triviale de $K_4$, retourner ``Pas de LGSNTK4''.
    \end{itemize}

  \item[\sc Complexité :] $O(n^{20})$.
  \end{itemize}
\end{algorithme}

\begin{preuve} 
  Si $G$  est sans line-graphe  de subdivision non triviale  de $K_4$,
  alors  il  est clair  que  l'algorithme  n'en  trouve pas  et  qu'il
  retourne ``Pas de LGSNTK4''.
  
  Si  $G$ contient  un line-graphe  d'une subdivision  non-triviale de
  $K_4$, alors $G$  en contient un de taille  minimale, que l'on note
  $F$.

  À une certaine étape,  l'algorithme considère un 18-uplet $(v_{ab}$,
  $v_{ac}$,  $\dots$,   $v_{cd}$,  $m_{ab}$,  $\dots$,   $m_{cd})  \in
  V^{18}$, qui est une trame de $F$.  Notons que les douze plus courts
  chemins   recherchés   par   l'algorithme   existent.    Par   douze
  applications consécutives du  lemme~\ref{reco.l.subK4}, on voit que
  $\cup_{i,j}  (V(P_{ij}) \cup  V(Q_{ij}))$ est  un  line-graphe d'une
  subdivision  non  triviale  de  $K_4$  de  taille  minimale.   Donc
  l'algorithme  retourne   ce  sous-graphe,  ou   peut-être  un  autre
  sous-graphe de même taille et ayant les mêmes coins.
  
  Dans  le   pire  des  cas,  tester  tous   les  18-uplets  nécessite
  $O(n^{18})$ étapes. Calculer les plus courts chemins et vérifier que
  le sous-graphe obtenu est  le line-graphe de subdivision bipartie de
  $K_4$ prend $O(n^2)$ étapes.  Donc, au pire, l'algorithme se termine
  en $O(n^{20})$ étapes.
\end{preuve}

Nous   pouvons   maintenant  donner   un   algorithme  détectant   les
line-graphes de subdivision \emph{bipartie} de $K_4$ dans les graphes
sans trous impairs.

\begin{algorithme}\label{reco.a.bergeSansLK4}
  \begin{itemize}
  \item[\sc Entrée :] Un graphe $G$ sans trou impair.
  \item[\sc Sortie :] Si le graphe est sans line-graphe de subdivision
    bipartie  de  $K_4$~:  ``Pas  de LGSBK4''.   Sinon,  l'algorithme
    retourne  un  line-graphe de  subdivision  bipartie  de $K_4$  de
    taille minimale.
  \item[\sc Calcul :] Exécuter l'algorithme~\ref{reco.a.subK4}.
  \item[\sc Complexité :] $O(n^{20})$.
  \end{itemize}
\end{algorithme}

\begin{preuve} 
  D'après  le  lemme~\ref{reco.l.techK4} on  sait  que,  dans $G$,  les
  line-graphes  de  subdivision  \emph{non-triviales}  de  $K_4$  sont
  exactement  les  line-graphes  de  subdivision  \emph{bipartie}  de
  $K_4$.    Ces  dernières   configurations  sont   donc  correctement
  détectées par l'algorithme~\ref{reco.a.subK4}.
\end{preuve}

\subsection{Prismes impairs}

Nous pouvons désormais fournir  un algorithme en temps polynomial pour
détecter les prismes impairs.  Il  faut se montrer prudent car dans ce
cas  la technique  des plus  courts  chemins fonctionne mal~:

\begin{figure}[ht]
  \center
  \includegraphics{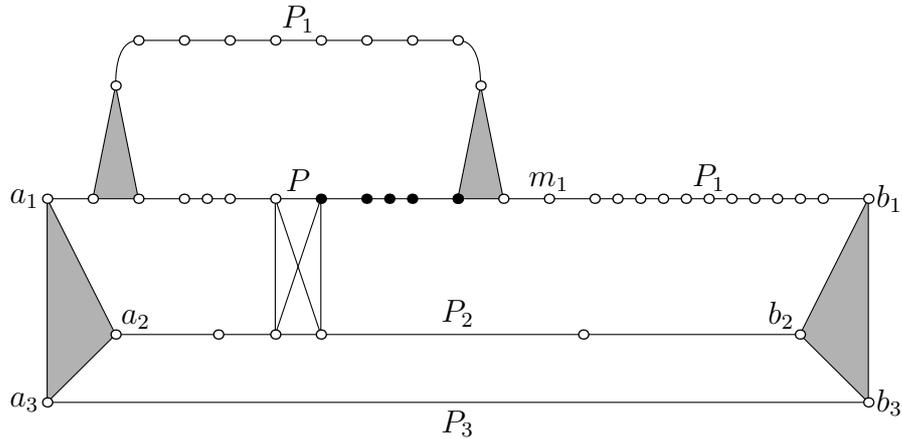}
  \caption{Graphe      contenant      six      prismes
    impairs\label{reco.f.prismeimp}}
  \end{figure}


Le  graphe  $G$  représenté figure~\ref{reco.f.prismeimp}  comporte  5
cliques maximales de taille au  moins trois. On vérifie que les quatre
cliques  grisées   sont  deux  à  deux   ennemies  (voir  sous-section
\ref{pair.ss.cliquesennemies}    page~\pageref{pair.ss.cliquesennemies}
pour de plus amples explications  sur les cliques ennemies).  Tous les
chemins sortants les  reliant à la clique non  grisée sont de longueur
paire.  On en  déduit que $G$ est le  line-graphe d'un graphe biparti,
et  qu'il est  de Berge  (pour s'en  convaincre, on  peut  utiliser le
lemme~\ref{pair.l.caracennemie}, qui  aura finalement servi  à quelque
chose~!).  Pour  toute paire de  cliques grisées, il existe  un unique
prisme impair de $G$ qui  ait ces cliques pour triangles.  On constate
que les  chemins $P_1$, $P_2$ et  $P_3$ induisent un  prisme impair de
$G$  de taille  minimale.  Pourtant,  si  on remplace  $P_1$, ou  même
simplement le  sous-chemin $\bp a_1 \tp  P_1 \tp m_1 \ep$  par un plus
court chemin $P$, on n'a pas la garantie d'obtenir à nouveau un prisme
impair.   Un  algorithme naïf  par  plus  courts  chemins risque  donc
d'échouer. En fait un tel algorithme (calculer les plus courts chemins
pour chaque sextuplet $(a_1, a_2, a_3, b_1, b_2, b_3)\in V^6$) executé
sur $G$ finirait bien par  détecter un prisme impair~: celui contenant
les  deux triangles  les  plus en  haut  sur le  dessin.   Mais il  ne
trouverait pas un prisme de taille  minimale de $G$, ce qui fait qu'on
ne voit comment prouver qu'il fonctionne (s'il fonctionne~!).

On  peut cependant  remarquer que  $G$ contient  le  line-graphe d'une
subdivision bipartie de $K_4$~:  le sous-graphe induit en oubliant les
sommets ``pleins''.  Si ce type de configuration est interdite, il y a
donc  encore un  espoir  que  les techniques  de  plus courts  chemins
fonctionnent.  Le  lemme suivant montre  que cet espoir est  fondé. On
remarquera qu'il n'est  même plus nécessaire de partir  d'un prisme de
taille  minimale,  ni  de  plus  courts chemins,  ce  qui  montre  que
l'interdiction  des  line-graphes  de  subdivision de  $K_4$  est  une
contrainte  très forte  dont on  a déjà  eu une  illustration  dans la
preuve   du   théorème   fort   des   graphes   parfaits   (voir   les
théorèmes~\ref{graphespar.t.spgtprismepair}
à~\ref{graphespar.t.sartemis}).

\begin{lemme} \label{reco.l.prismeimp}
  Soit  $G$ un  graphe  sans  trou impair  et  sans line-graphe  de
  subdivision  bipartie de $K_4$.   Soit $F$  un prisme  de triangles
  $\{a_1, a_2,  a_3\}$ et  $\{ b_1,  b_2, b_3 \}$  et formé  par trois
  chemins $P_i$  ($i=1, 2, 3$)  ayant pour extrémités $a_i$  et $b_i$.
  Soit $P$ un  chemin de $a_1$ vers $b_1$  dont les sommets intérieurs
  manquent $a_2$, $a_3$, $b_2$ et $b_3$. Alors~:

  \begin{center}
    $P$, $P_2$, $P_3$ forment un prisme de $G$ de même parité que $F$.
  \end{center}

\end{lemme}

\begin{preuve}
  On  pose $H  = P_2  \cup  P_3$.  Si  les sommets  de $P^*$  manquent
  entièrement le trou  $H$, alors $P$, $P_2$, $P_3$  forment un prisme
  impair de $G$  (impair car $G$ est sans trou  impair).  On peut donc
  supposer qu'il existe un sommet de $P^*$ qui a des voisins dans $H$.

  On  note $c_1$  le sommet  de $P^*$  ayant des  voisins dans  $H$ et
  choisi aussi proche  que possible de $a_1$.  On  définit les chemins
  $H_2 =  H \setminus \{  b_3 \}$  et $H_3 =  H \setminus \{  b_2 \}$.
  Pour $i=2,  3$, on note $c_i$  le sommet de $H_i$,  voisin de $c_1$,
  choisi aussi proche que possible de $b_i$.

  \begin{claim}
    On peut supposer que $c_2 c_3$ est une arête de $P_2^*$.
  \end{claim}

  \begin{preuveclaim}
    Si $c_2 = c_3$, alors les trois chemins $\bp c_2 \tp c_1 \tp P \tp
    a_1 \ep$, $\bp  c_2 \tp P_2 \tp  a_2 \ep$ et $\bp c_2  \tp P_2 \tp
    b_2 \tp  b_3 \tp  P_3 \tp  a_3 \ep$ forment  une pyramide  de coin
    $c_2$ et  de triangle $\{a_1,  a_2, a_3\}$~: une  contradiction. On
    peut donc supposer  $c_2 \neq c_3$. Si $c_2$ et  $c_3$ ne sont pas
    adjacents, alors les  trois chemins $\bp c_1 \tp  P \tp a_1 \ep$,
    $\bp c_1 \tp c_2 \tp P_2 \tp  a_2 \ep$ et $\bp c_1 \tp c_3 \tp P_2
    \tp b_2 \tp b_3 \tp P_3  \tp a_3 \ep$ forment une pyramide de coin
    $c_1$ et de triangle $\{a_1, a_2, a_3\}$~: une contradiction.
  \end{preuveclaim}

  \begin{claim}
    $c_1$ manque $b_1$.
  \end{claim}

  \begin{preuveclaim}
    Sinon, les trois  chemins $\bp c_1 \tp b_1 \ep$,  $\bp c_1 \tp c_3
    \tp P_2 \tp b_2 \ep$ et $\bp c_1 \tp P \tp a_1 \tp a_3 \tp P_3 \tp
    b_3 \ep$ forment une pyramide de coin $c_1$ et de triangle $\{b_1,
    b_2, b_3\}$.
  \end{preuveclaim}

  On  note $a'_1$ le  voisin de  $a_1$ dans  $P_1$.  On  définit alors
  $d_1$  le sommet de  $\bp a'_1  \tp P  \tp c_1  \ep$ le  plus proche
  possible de $c_1$ et ayant  des voisins dans $P_1$.  On définit $d_2
  \in P_1$, voisin de $d_1$ aussi proche que possible de $a_1$ et $d_3
  \in P_1$, voisin de $d_1$ aussi proche que possible de $b_1$.

  \begin{claim}
    $d_2 d_3$ est une arête de $P_1$.
  \end{claim}

  \begin{preuveclaim}
    Si $d_2 = d_3$, alors les trois chemins $\bp d_2 \tp d_1 \tp P \tp
    c_1 \ep$, $\bp d_2 \tp P_1 \tp a_1 \tp a_2 \tp P_2 \tp c_2 \ep$ et
    $\bp d_2 \tp P_1 \tp b_1 \tp  b_2 \tp P_2 \tp c_3 \ep$ forment une
    pyramide de  coin $d_2$  et de triangle  $\{c_1, c_2,  c_3\}$~: une
    contradiction. On peut donc supposer  $d_2 \neq d_3$.  Si $d_2$ et
    $d_3$ ne sont pas adjacents, alors les trois chemins $\bp d_1 \tp
    d_2 \tp P_1 \tp a_1 \ep$, $\bp d_1 \tp d_3 \tp P_1 \tp b_1 \tp b_3
    \tp P_3 \tp a_3 \ep$ et $\bp d_1 \tp P \tp c_1 \tp c_2 \tp P_2 \tp
    a_2 \ep$ forment une pyramide de coin $d_1$ et de triangle $\{a_1,
    a_2, a_3\}$~: une contradiction.
  \end{preuveclaim}

  On remarque maintenant que les quatre triangles $\{a_1, a_2, a_3\}$,
  $\{b_1, b_2, b_3\}$, $\{c_1, c_2, c_3\}$, $\{d_1, d_2, d_3\}$ et les
  six chemins $P_3$,  $\bp a_2 \tp P_2 \tp c_2 \ep$,  $\bp a_1 \tp P_1
  \tp d_2 \ep$,  $\bp b_2 \tp P_2  \tp c_3 \ep$, $\bp b_1  \tp P_1 \tp
  d_3 \ep$, $\bp c_1 \tp P \tp d_1 \ep$ induisent le line-graphe d'une
  subdivision \emph{quelconque} de $K_4$  qui n'est pas le line-graphe
  de    $K_4$,   car    $a_3    \neq   b_3$.     Donc,   d'après    le
  lemme~\ref{reco.l.techK4},   $G$  contient   le   line-graphe  d'une
  subdivision \emph{bipartie} de $K_4$, ce qui est contradictoire.

\end{preuve}

Voici maintenant  un algorithme de détection des  prismes impairs dans
les graphes sans trou impair.
 
\index{prisme!détection~des~---~impairs}
\index{détection!prismes~impairs}

\begin{algorithme}\label{reco.a.prismeimp}
  \begin{itemize}
  \item[\sc  Entrée  :]  Un  graphe   $G$  sans  trou  impair.

  \item[\sc Sortie :] Si $G$  est sans prisme impair~: ``Pas de prisme
  impair''. Sinon~: ``Il y a un prisme impair''

  \item[\sc  Calcul :] \mbox{}
    \begin{itemize}
    \item 
      À  l'aide  de  l'algorithme  \ref{reco.a.subK4}, tester  si  $G$
      possède un line-graphe de subdivision non triviale de $K_4$.  Si
      tel  est  le  cas,  stopper  et  répondre ``Il  y  a  un  prisme
      impair''. Sinon poursuivre.
    \item 
      Considérer successivement  tous les sextuplets  $(a_1, a_2, a_3,
      b_1, b_2, b_3)\in V^6$.  Calculer pour $i=1, 2, 3$ un plus court
      chemin  $P_i$ de $a_i$  vers $b_i$  dont les  sommets intérieurs
      manquent  $a_{i+1}$,  $a_{i+2}$,  $b_{i+1}$  et  $b_{i+2}$  (les
      indices sont calculés modulo 3).  Si $P_1$, $P_2$, $P_3$ forment
      un prisme impair, alors  stopper l'algorithme et répondre ``Il y
      a un prisme impair''.
   
    \item 
      Si  aucun  sextuplet n'a  donné  de  prisme,  répondre ``Pas  de
      prisme impair''.
    \end{itemize}
  \item[\sc Complexité :] $O(n^{20})$.
  \end{itemize}
\end{algorithme}

\begin{preuve}
  Si $G$  est sans prisme  impair, alors $G$  ne peut pas  contenir de
  line-graphe d'une  subdivision non triviale de $K_4$  car d'après le
  lemme~\ref{reco.l.techK4syn}, un tel  line-graphe serait trivial, ou
  contiendrait un trou impair, ou contiendrait un prisme impair~: dans
  tous les  cas une contradiction.   Donc l'algorithme donne  la bonne
  réponse~: ``Pas de prisme impair''.
  
  Si $G$ possède un prisme  impair et un line-graphe d'une subdivision
  non triviale de  $K_4$, alors cela est détecté  à la deuxième étape,
  et l'algorithme  retourne ``Il  y a un  prisme impair'' comme  il se
  doit.

  Si $G$ possède  un prisme impair et ne  possède aucun line-graphe de
  subdivision \emph{non triviale} de $K_4$, alors, en particulier, $G$
  ne  possède  aucun line-graphe  de  subdivision \emph{bipartie}  de
  $K_4$, et  on peut appliquer  le lemme~\ref{reco.l.prismeimp}.  Soit
  alors $F$ un  prisme impair quelconque de $G$,  de triangles $\{a_1,
  a_2,  a_3\}$  et  $\{b_1,   b_2,  b_3\}$.   À  une  certaine  étape,
  l'algorithme  considère le sextuplet  formé par  ces six  sommets et
  calcule  pour  $i=1,  2,  3$  les plus  courts  chemins  $P_i$  (qui
  existent).      Par    trois     applications     consécutives    du
  lemme~\ref{reco.l.prismeimp}, on  voit que $P_1 \cup  P_2 \cup P_3$
  induit un  prisme impair, qui sera donc  bien détecté~: l'algorithme
  retourne ``Il y a un prisme impair''.
\end{preuve}
\index{détection!sous-graphes|)}

\section{Des problèmes NP-complets}

Nous avons résolu  dans ce chapitre un certain  nombre de problèmes de
détection de  sous-graphes induits dans  les graphes sans  pyramide ou
dans  les graphes  sans trou  impair.  Nous allons  voir dans  cette
section que si l'on cherche à les étendre aux graphes quelconques, ces
problèmes deviennent  NP-complets.  Nous montrerons  que les problèmes
de détection  des prismes, des  prismes pairs, des prismes  impairs et
des line-graphes de subdivision de $K_4$ sont tous NP-complets.

Notre  preuve s'inspire  directement de  la construction  de Bienstock
dans  son article  sur la  NP-complétude  de la  recherche des  paires
d'amis.  Dans~\cite{bienstock:evenpair},  Bienstock parvient à réduire
le problème 3-{\sc sat} à la recherche de paire d'amis dans un graphe.
Il parvient également à réduire 3-{\sc sat} au problème suivant :

\begin{probleme}[$\Pi'$]
  \begin{itemize}
  \item[\sc Instance  :] Un graphe $G$  et deux sommets $a$  et $b$ de
    $G$.
  \item[\sc Question  :] Y a-t-il  un trou de  $G$ passant par  $a$ et
  $b$~?
  \item[\sc      Complexité       :]      NP-complet,      (Bienstock,
  \cite{bienstock:evenpair}).
  \end{itemize}
\end{probleme}

En  adaptant  la preuve  de  Bienstock,  nous  allons montrer  que  le
problème $\Pi'$  demeure NP-complet même  si l'on se restreint  au cas
des graphes sans  triangle.  Pour alléger notre propos,  il aurait été
préférable  d'utiliser les  \emph{résultats} de  Bienstock  plutôt que
d'adapter sa  preuve, mais nous  n'y sommes pas parvenus.   Notons que
notre  construction est  plus simple  que  celle de  Bienstock car  ce
dernier prouvait  des propriétés plus difficiles  (NP-complétude de la
recherche de paires  d'amis). On notera qu'{\it a priori}  il n'y a pas
de rapport  directe entre la recherche  de prismes et  la recherche de
paires d'amis. Il est donc  remarquable que ces deux structures soient
liés  à la  fois  par des  théorèmes  (paire d'amis  dans les  graphes
d'Artémis,  où  les  prismes  sont  interdits) et  par  la  preuve  de
NP-complétude à suivre.

\begin{probleme}[$\Pi$]
  \begin{itemize}
  \item[\sc Instance  :]  
    Un graphe  $G$  sans triangle. Deux sommets  $a$  et  $b$ de  $G$,
    non adjacents et de degré~2.
  \item[\sc Question :] 
    Y a-t-il un trou de $G$ passant par $a$ et $b$ ?
  \item[\sc Complexité :] NP-complet. 
  \end{itemize}
\end{probleme}

\begin{preuve} 
  Pour  un   rappel  de  la   définition  du  problème   3-{\sc  sat},
  voir~\ref{base.ss.npcomplet}  page~\pageref{base.ss.npcomplet}.   Il
  est clair  que $\Pi$ est un  problème NP.  Soit $f$  une instance de
  {\sc $3$-sat}, consistant en $m$  clauses $C_1, \ldots, C_m$ sur $n$
  variables $x_1, \ldots, x_n$.  Nous allons construire un graphe sans
  triangle $G_f$  dont la  taille sera majorée  par un polynôme  en la
  taille de  $f$ et  contenant deux sommets  distingués $a$ et  $b$ de
  degré deux, tels  que $G_f$ contient un trou passant  par $a$ et $b$
  si  et  seulement  si  $f$  est  satisfaisable.   Ceci  prouvera  le
  théorème.  Les étapes de  la construction de $G_f$ sont représentées
  figure~\ref{reco.f.bienstock}.

  \begin{figure}[!htb]
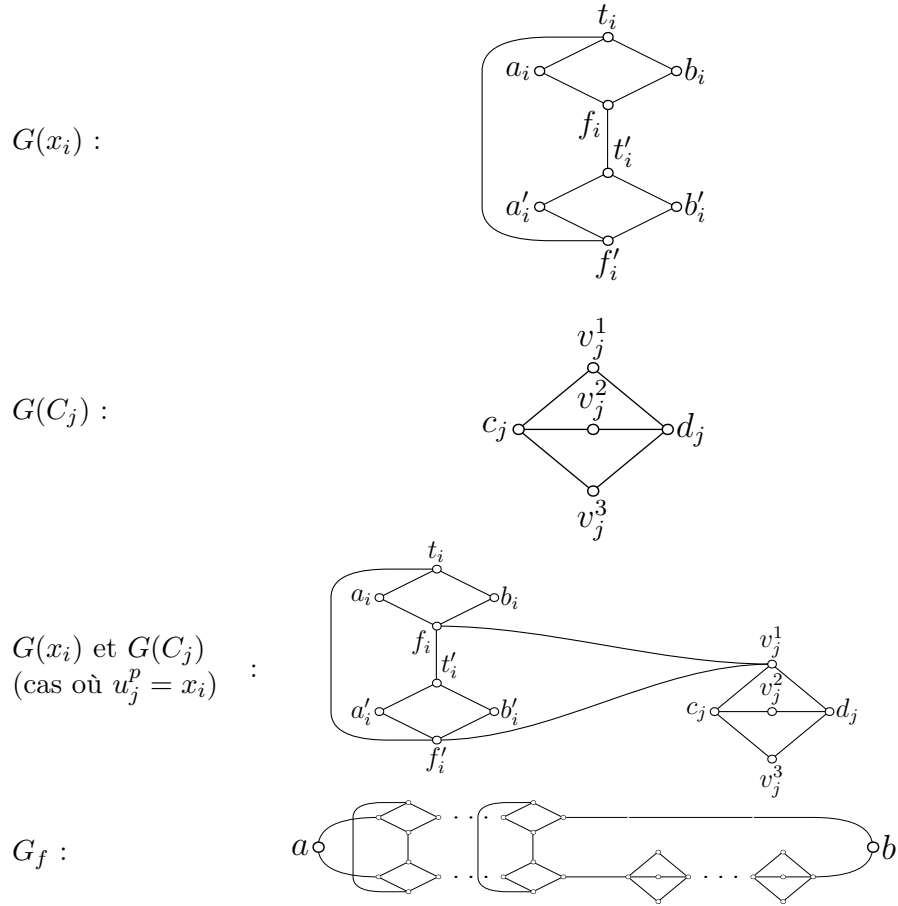

    \center 
    \begin{tabular}{lc}
      $G(x_i)$ : &\parbox[c]{3cm}{\includegraphics[width= 3cm]{fig.reco.4}}\\
      $G(C_j)$ : &\parbox[c]{3cm}{\includegraphics[width= 3cm]{fig.reco.5}}\\
      \parbox[c]{3cm}{$G(x_i)$ et $G(C_j)$  \\ 
	(cas où $u^p_j=x_i$)} :
      &\parbox[c]{7cm}{\includegraphics[width=    7cm]{fig.reco.6}}\\
      $G_f$        :        &\parbox[c]{8cm}{\includegraphics[width=
	  8cm]{fig.reco.7}}\\
    \end{tabular}
    \caption{Étapes de la construction de $G_f$}
    \label{reco.f.bienstock}
  \end{figure}

  \vspace{1ex}  
  \begin{itemize}
  \item
    Pour  chaque variable  $x_i$  ($i=1, \ldots,  n$),  on définit  le
    graphe  $G(x_i)$ avec  huit  sommets $a_i$,  $b_i$, $t_i$,  $f_i$,
    $a'_i$, $b'_i$, $t'_i$, $f'_i$ et dix arêtes $a_i t_i$, $a_i f_i$,
    $b_i t_i$, $b_i f_i$, $a'_i t'_i$, $a'_i f'_i$, $b'_i t'_i$, $b'_i
    f'_i$, $t_i f'_i$, $t'_i f_i$.

  \item 
    Pour chaque clause $C_j$ ($j=1, \ldots, m$), avec $C_j = u_j^1\vee
    u_j^2\vee u_j^3$, où chaque $u_j^p$  ($p=1, 2, 3$) est un littéral
    de $\{x_1, \ldots, x_n, \overline{x}_1, \ldots, \overline{x}_n\}$,
    on définit un graphe $G(C_j)$  avec cinq sommets $c_j, d_j, v_j^1,
    v_j^2, v_j^3$  et six arêtes de  sorte que chaque  $c_j, d_j$ voit
    chaque $v_j^1, v_j^2, v_j^3$.

  \item 
    Pour  $p=1, 2,  3$, si  $u_j^p=x_i$  alors on  ajoute deux  arêtes
    $v_j^p    f_i$    et    $v_j^p    f'_i$.     Si    au    contraire
    $u_j^p=\overline{x}_i$ alors on ajoute  deux arêtes $v_j^p t_i$ et
    $v_j^p t'_i$.
  \end{itemize}

  \vspace{2ex}

  Le  graphe  $G_f$ est  obtenu  à  partir  de l'union  disjointe  des
  $G(x_i)$ et des $G(C_j)$ de la manière suivante.  Pour $i=1, \ldots,
  n-1$, on  ajoute les arêtes  $b_i a_{i+1}$ et $b'_i  a'_{i+1}$.  On
  rajoute  l'arête $b'_n  c_1$.  Pour  $j=1, \ldots,  m-1$,  on ajoute
  l'arête $d_j c_{j+1}$.  On introduit les deux sommets distingués $a$
  et $b$,  puis on  ajoute les arêtes  $a a_1$,  $a a'_1$, $b  d_m$, $b
  b_n$.  La taille de $G_f$ est  un polynôme en la taille $f$ (en fait
  la taille de $G_f$ est  $O(n+m)$). On constate que $G_f$ ne contient
  pas de triangle, et que $a$ et $b$ sont de degré~2.

  On  va  montrer  que  si  l'on  dispose  d'un  algorithme  en  temps
  polynomial pour résoudre $\Pi$,  alors on dispose d'un algorithme en
  temps polynomial pour résoudre 3-{\sc sat}. Pour ce faire, il suffit
  de montrer qu'une  instance $f$ de 3-{\sc sat}  est satisfaisable si
  et seulement si le graphe $G_f$  contient un trou passant par $a$ et
  $b$.
  
  \begin{claim}
    \label{reco.c.bienstock1}
    Si $f$ est satisfaisable, alors $G_f$ contient un trou passant par
    $a$ et $b$.
  \end{claim}
   
  \begin{preuveclaim}
    Soit $\xi\in\{0, 1\}^n$ un  vecteur booléen qui satisfait $f$.  On
    trouve un trou  dans $G$ en sélectionnant les  sommets comme suit.
    On sélectionne $a$ et $b$.   Pour $i=1, \ldots, n$, on sélectionne
    $a_i$,  $b_i$,   $a'_i$,  $b'_i$;   de  plus,  si   $\xi_i=1$,  on
    sélectionne $t_i,  t'_i$, tandis que si  $\xi_i=0$, on sélectionne
    $f_i, f'_i$.  Pour $j=1,  \ldots, m$, puisque $\xi$ satisfait $f$,
    au moins l'un  des trois littéraux de $C_j$ est  égal à $1$.  Donc
    il  existe $p\in\{1,  2, 3\}$  tel que  $u_j^p=1$.  On sélectionne
    alors  $c_j$, $d_j$  et $v_j^p$.   Il  est clair  que les  sommets
    sélectionnés  appartiennent à  un cycle  $Z$ qui  contient  $a$ et
    $b$. On vérifie que $Z$ est sans corde.  Le seul point délicat est
    que $Z$ ne contient aucune corde entre un sommet de $G(C_j)$ et un
    sommet de $G(x_i)$,  car une telle corde serait  ou bien une arête
    $t_i v_j^p$ (ou $t'_iv_j^p$) avec $u_j^p = x_i$ et $\xi_i = 1$, ou
    bien, symétriquement, une arête $f_i v_j^p$ (ou $f'_i v_j^p$) avec
    $u_j^p  = \overline{x}_i$  et $\xi_i=0$.   Dans tous  les  cas, on
    aurait une contradiction avec la méthode de sélection.
  \end{preuveclaim}
   
  Réciproquement, supposons  que $G_f$ contienne un trou  $Z$ avec $a,b
  \in V(Z)$.  Alors  $Z$ contient $a_1$ et $a'_1$  puisque ce sont les
  seuls voisins de $a$ dans $G_f$.

  \begin{claim}\label{clm:zgxi}
    \label{reco.c.bienstock2}
    Pour  $i=1, \ldots,  n$, $Z$  contient exactement  six  sommets de
    $G_i$: quatre d'entre eux sont $a_i, a'_i, b_i, b'_i$, et les deux
    autres sont ou bien $t_i$ et $t'_i$, ou bien $f_i$ et $f'_i$.
  \end{claim}

  \begin{preuveclaim} 
    On étudie  d'abord le cas  $i=1$.  Puisque $a,  a_1 \in Z$  et que
    $a_1$  n'a que  trois voisins  $a, t_1,  f_1$, exactement  l'un de
    $t_1, f_1$ est dans $Z$.   De même exactement l'un de $t'_1, f'_1$
    est dans $Z$.  Si $t_1$ et $f'_1$ sont dans $Z$, alors les sommets
    $a, a_1, a'_1,  t_1, f'_1$ sont tous dans $Z$  et ils induisent un
    trou qui ne contient pas  $b$~: une contradiction.  De même on n'a
    pas simultanément  $t'_1$ et $f_1$ dans $Z$.   Donc, par symétrie,
    on peut supposer que $t_1, t'_1$ sont dans $Z$ tandis que $f_1$ et
    $f'_1$ n'y  sont pas.   S'il existe un  sommet $u_j^p  \in G(C_j)$
    (avec $1\le j\le m$, $1\le p\le  3$) dans $Z$ et que ce sommet est
    adjacent  à $t_1$,  alors,  puisque ce  sommet  $u_j^p$ est  aussi
    adjacent à  $t'_1$, on  voit que les  sommets $a, a_1,  a'_1, t_1,
    t'_1,  u_j^p$ sont  tous  dans $Z$  et  induisent un  trou qui  ne
    contient pas $b$~: une contradiction.  Donc le voisin de $t_1$ sur
    $Z$ qui  est différent de $a_1$  n'est pas dans  l'un des $G(C_j)$
    ($1\le j\le m$), donc ce voisin est $b_1$.  De même, $b'_1 \in Z$.
    On  a donc  prouvé~(\ref{reco.c.bienstock2}) pour  $i=1$.  Puisque
    $b_1\in Z$ et  qu'exactement l'un de $t_1$ et  $f_1$ est dans $Z$,
    puisque  $b_1$ est de  degré $3$,  on conclut  que $a_2$  est dans
    $Z$. De même,  on conclut que $b_2$ est dans  $Z$.  Pour $i=2$, la
    preuve  est  essentiellement  identique   au  cas  $i=1$,  et  par
    récurrence on parvient à $i=n$.
  \end{preuveclaim}

  \begin{claim}
    \label{reco.c.bienstock3}
    Pour  $j=1, \ldots,  m$, $Z$  contient $c_j$,  $d_j$ et  un sommet
    exactement parmi $v_j^1, v_j^2, v_j^3$.
  \end{claim}

  \begin{preuveclaim} 
    On        étudie        d'abord        le        cas        $j=1$.
    D'après~(\ref{reco.c.bienstock2}),   $b'_n   \in   Z$  et un   sommet
    exactement parmi $t'_n, f'_n$  est dans $Z$.  Donc, puisque $b'_n$
    est de degré $3$, $c_1$  est dans $Z$.  Par conséquent, exactement
    un sommet parmi  $v_1^1, v_1^2, v_1^3$ est dans  $Z$, par exemple,
    sans perte de  généralité, $v_1^1$.  Le voisin de  $v_1^1$ dans $Z$
    qui est différent de $c_1$ ne  peut pas être un sommet de $G(x_i)$
    ($1\le i\le n$), car ce serait alors $t_i$ (ou $f_i$). Or, on sait
    d'après~(\ref{reco.c.bienstock2}) que  les deux voisins  de $t_i$
    (ou $f_i$) dans $Z$ sont  des sommets de $G(x_i)$. Le sommet $t_i$
    (ou $f_i$) aurait alors  trois voisins dans $Z$~: une contradiction
    avec l'état  de trou.  Finalement, l'autre voisin  de $u_1^1$ dans
    $Z$  est  $d_1$,  et  on a  démontré~(\ref{reco.c.bienstock3})  pour
    $j=1$.  Puisque  $d_1$ est de  degré $4$ et qu'exactement  l'un de
    $v_1^1, v_1^2, v_1^3$ est dans $Z$, on constate que $c_2$ est dans
    $Z$.  La  preuve est essentiellement identique   pour $j=2$, et
    par récurrence on parvient à $j=m$.
  \end{preuveclaim}

  À partir  de $Z$, on  construit maintenant un vecteur  booléen $\xi$
  comme suit. Pour $i=1, \ldots,  n$, si $Z$ contient $t_i$ et $t'_i$,
  on pose $\xi_i = 1$; si  $Z$ contient $f_i$ et $f'_i$ on pose $\xi_i
  =  0$.   D'après~(\ref{reco.c.bienstock2}),  cela  a bien  un  sens.

  Soit    $C_j$    ($1\le    j\le    m$)   une    clause    de    $f$.
  D'après~(\ref{reco.c.bienstock3})  et sans  perte de  généralité, on
  peut  supposer que  $v_j^1$ est  dans $Z$.   Si $u_j^1  =  x_i$ avec
  $i\in\{1, .., n\}$,  alors la construction de $G_f$  entraîne que ni
  $f_i$ ni  $f'_i$ ne sont dans  $Z$. Donc, $t_i$ et  $t'_i$ sont dans
  $Z$, et donc  $\xi_i=1$ et la clause $C_j$  est satisfaite par $x_i$.
  De même, si $u_j^1 =  \overline{x}_i$ avec $i\in\{1, .., n\}$, alors
  la construction  de $G_f$  entraîne que ni  $t_i$ ni $t'_i$  ne sont
  dans $Z$.  Donc, $f_i$ et $f'_i$ sont dans $Z$, donc $\xi_i=0$ et la
  clause $C_j$ est satisfaite par $\overline{x}_i$.  Finalement, $\xi$
  satisfait $f$.
\end{preuve}

La  NP-complétude de  $\Pi$  entraîne directement  celle de  plusieurs
problèmes.  Nous  pourrions peut-être écrire un  problème générique de
détection de sous-graphes induits formés  de cliques et de chemins les
reliant,  mais  cela  serait  d'un  formalisme  un  peu  lourd.   Nous
préférons  nous restreindre  ici aux  problèmes pertinents  pour notre
étude.   Nous espérons  convaincre  le lecteur  que,  pour des  usages
futurs,  le  problème  $\Pi$  est  assez facile  à  réduire  à  divers
problèmes.

\index{prisme!détection~des~---} 
\index{prisme!détection~des~---~pairs} 
\index{prisme!détection~des~---~impairs} 
\index{détection!prismes~pairs} 
\index{détection!prismes~impairs} 
\index{détection!prismes} 
\index{line-graphe~de~subdivision~bipartie~de~$K_4$!détection~des~---}
\index{détection!line-graphes~de~subdivision~bipartie~de~$K_4$}

\begin{theoreme}
  \label{thm:prismsnpc}
  Les cinq problèmes suivants sont NP-complets~:
  \begin{enumerate}
  \item
    \begin{itemize}
    \item[\sc  Instance :]  Un  graphe $G$  contenant exactement  deux
      triangles.
    \item[\sc Question :] $G$ contient-il un prisme ?
    \end{itemize}

  \item
    \begin{itemize}
    \item[\sc Instance  :]  
      Un graphe  $G$  contenant exactement deux triangles.
    \item[\sc Question :] 
      $G$ contient-il un prisme impair ?
    \end{itemize}

  \item
    \begin{itemize}
    \item[\sc Instance  :]  
      Un graphe  $G$  contenant exactement deux triangles.
    \item[\sc Question :] 
      $G$ contient-il un prisme pair ?
    \end{itemize}

  \item
    \begin{itemize}
    \item[\sc Instance  :]  
      Un graphe  $G$  contenant exactement quatre triangles.
    \item[\sc Question :] 
      $G$ contient-il le line-graphe d'une subdivision  de $K_4$ ?
    \end{itemize}

  \item
    \begin{itemize}
    \item[\sc Instance  :]  
      Un graphe  $G$  contenant exactement quatre triangles.
    \item[\sc Question :] 
      $G$ contient-il le line-graphe d'une subdivision bipartie de $K_4$ ?
    \end{itemize}

  \end{enumerate}
\end{theoreme}

\begin{preuve}
  Pour chacun des  cinq problèmes, on va supposer  que l'on connaît un
  algorithme de décision  en temps polynomial. Puis, on  va en réduire
  polynomialement  le problème  considéré  au problème  $\Pi$, ce  qui
  prouvera le théorème. Soit donc $(G, a, b)$ une instance du problème
  $\Pi$.  On  rappelle que par  définition $G$ est sans  triangle, que
  $a$ et $b$ sont non adjacents et de degré~2.

    \begin{figure}[!htb]
      \center
      \includegraphics{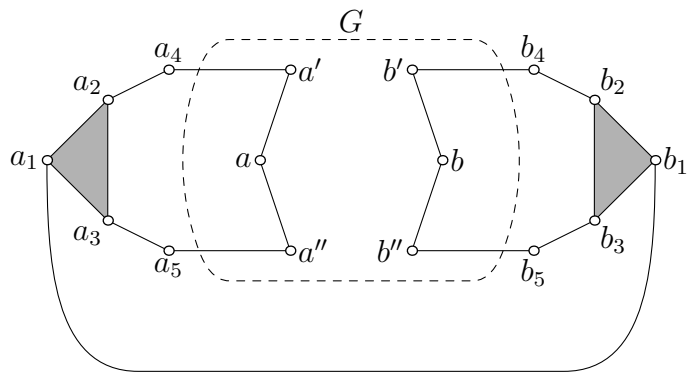}
      \caption{$G$  et $G'$,  réduction  de $\Pi$  à  la recherche  de
      prismes\label{reco.f.Piprism}}
    \end{figure}

  \begin{enumerate}

  \item 
    À  partir  de  $G$,  on  construit  le  graphe  $G'$  ainsi  (voir
    figure~\ref{reco.f.Piprism})~:  on remplace le  sommet $a$  par le
    graphe avec cinq sommets $a_1$, $a_2$, $a_3$, $a_4$, $a_5$ et cinq
    arêtes $a_1 a_2$, $a_1 a_3$,  $a_2 a_3$, $a_2 a_4$, $a_3 a_5$.  On
    relie $a_4$ à l'un des voisins  de $a$ et $a_5$ à l'autre voisin de
    $a$.  On  remplace de même le  sommet $b$ par le  graphe avec cinq
    sommets  $b_1$, $b_2$,  $b_3$, $b_4$,  $b_5$ et  cinq  arêtes $b_1
    b_2$, $b_1 b_3$, $b_2 b_3$,  $b_2 b_4$, $b_3 b_5$.  On relie $b_4$
    à l'un  des voisins de  $b$ et $b_5$  à l'autre voisin de  $b$.  On
    relie $a_1$  et $b_1$.  Comme $G$  est sans triangle,  on voit que
    $G'$  possède un prisme  si et  seulement si  $G$ possède  un trou
    passant  par $a$  et $b$.   Il est  en outre  évident que  $G'$ ne
    possède que deux triangles. Donc,  si l'on sait détecter un prisme
    en  temps  polynomial,  alors  on  sait  résoudre  $Pi$  en  temps
    polynomial.
    
  \item  
    À  partir de  $G$, on  construit le  graphe $G'$  comme ci-dessus.
    Puis on construit huit graphes auxiliaires $G(i,j,k)$ avec $(i, j,
    k)  \in \{0,1\}^3$  de la  manière suivante~:  selon que  $i=1$ ou
    $i=0$, on  subdivise ou non l'arête  $a_2 a_4$ en  lui ajoutant un
    sommet.   De même,  selon que  $j=1$  ou $j=0$  pour l'arête  $a_3
    a_5$. De  même, selon que $k=1$  ou $k=0$ pour  l'arête $a_1 b_1$.
    On constate alors  que $G$ possède un trou passant  par $a$ et $b$
    si  le  graphe  obtenu  par  réunion  disjointe  des  $G_{i,j,k}$
    posssède un prisme impair.

  \item 
    À  partir  de  $G$,  on  construit les  huit  graphes  auxiliaires
    $G(i,j,k)$ comme  ci-dessus. On vérifie  alors que $G$  possède un
    trou passant  par $a$ et $b$  si et seulement si  le graphe obtenu
    par réunion disjointe des $G_{i,j,k}$ posssède un prisme pair.

    \begin{figure}[!htb]
      \center
      \includegraphics{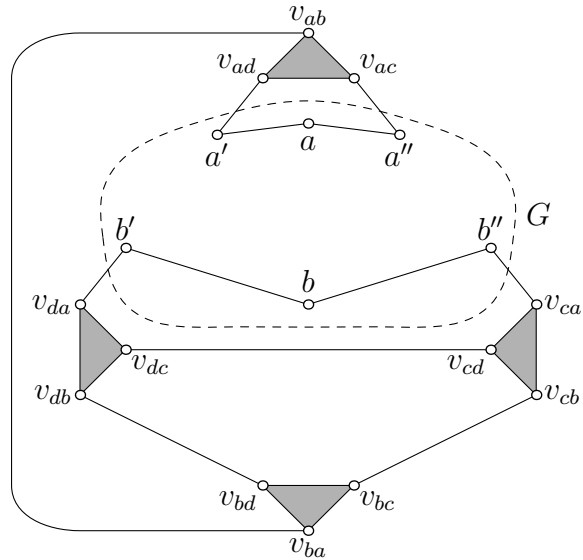}
      \caption{Réduction  de $\Pi$  à  la recherche  de
      line-graphes de subdivision de $K_4$\label{reco.f.PiLGBK4}}
    \end{figure}

  \item
    À  partir  de $G$,  on  construit un  graphe  $G'$  de la  manière
    suivante   (voir  figure~\ref{reco.f.PiLGBK4})~:  on   enlève  les
    sommets $a$ et $b$ et  on ajoute douze sommets $v_{ab}$, $v_{ac}$,
    $v_{ad}$,   $v_{ba}$,  $v_{bc}$,  $v_{bd}$,   $v_{ca}$,  $v_{cb}$,
    $v_{cd}$, $v_{da}$, $v_{db}$ et  $v_{dc}$. On ajoute des arêtes de
    sorte  que  chaque ensemble  $\{  v_{ab},  v_{ac}, v_{ad}\}$,  $\{
    v_{ba}, v_{bc},  v_{bd}\}$, $\{  v_{ca}, v_{cb}, v_{cd}\}$  et $\{
    v_{da}, v_{db}, v_{dc}\}$ soit  un triangle.  On ajoute les arêtes
    $v_{ab}  v_{ba}$,  $v_{dc}  v_{cd}$,  $v_{bd} v_{db}$  et  $v_{bc}
    v_{cb}$.  On relie $v_{ad}$ à l'un des voisins de $a$ dans $G$, et
    $v_{ac}$ à  l'autre voisin de $a$  dans $G$.  On  relie $v_{da}$ à
    l'un des voisins  de $b$ dans $G$ et $v_{ca}$  à l'autre voisin de
    $b$ dans $G$.  On vérifie que  $G$ possède un trou passant par $a$
    et  $b$ si  et seulement  si  $G'$ contient  le line-graphe  d'une
    subdivision bipartie de $K_4$.

  \item
    À partir de $G$, on construit un graphe $G'$ comme ci-dessus. Puis
    on construit quatre graphes  auxiliaires $G(i,j)$ avec $(i, j) \in
    \{0,1\}^2$ de la  manière suivante~: selon que $i=1$  ou $i=0$, on
    subdivise ou  non l'arête $v_{ad}  a'$ en lui ajoutant  un sommet.
    De même  selon que $j=1$ ou  $j=0$ pour l'arête  $v_{ac} a''$.  On
    constate alors que  $G$ possède un trou passant par  $a$ et $b$ si
    et  seulement  si  le  graphe  obtenu par  réunion  disjointe  des
    $G_{i,j}$  posssède  un  line-graphe  de subdivision  bipartie  de
    $K_4$.

  \end{enumerate}\mbox{}
\end{preuve}

\section{Reconnaissance et coloration}
\index{reconnaissance!lien avec la coloration}

Nous   pouvons  résumer  les   résultats  de   ce  chapitre   dans  le
tableau~\ref{reco.tab.res}.  Chaque ligne  correspond à un problème de
détection que l'on analyse dans  une classe de graphes correspondant à
une colonne.  Pour trois problèmes  de détection dans les graphes sans
pyramide  (signalés  par  ``?''),   nous  ne  sommes  pas  parvenus  à
déterminer de complexité. Si nous devions risquer une conjecture, nous
dirions que ces problèmes ont la même complexité que le problème de la
détection des trous impairs (conjecture peu aventureuse).

\index{détection!résumé~des~résultats}
\begin{table}[ht]
  \center
  \begin{tabular}{l|c|c|c|}
    \  & Graphes  &  Graphes sans   &
    Graphes sans  \\
    \ & généraux& pyramide & trou impair \\ \hline
    \rule{0ex}{3ex}Prisme ou Pyramide & $n^5$ & $n^5$ & $n^5$  \\
    Pyramide   & $n^9$ \cite{chudnovsky.seymour:reco} & $1$ & $1$  \\
    Prisme & NPC & $n^5$ & $n^5$  \\
    LGSNT$K_4$ & NPC & $n^{20}$ & $n^{20}$  \\
    LGSB$K_4$ & NPC & ? & $n^{20}$  \\
    Prisme impair & NPC & ? & $n^{20}$  \\
    Prisme pair & NPC & ? & $n^{11}$  \\
  \hline
\end{tabular}
  \caption{Complexité de divers problèmes de détection\label{reco.tab.res}}
\end{table}

Grâce  à  nos  algorithmes  de  détection,  nous  pouvons  donner  des
algorithmes  de reconnaissance  pour les  classes de  graphes étudiées
dans les chapitres précédents~:

\index{Artémis~(graphe~d'---)!reconnaissance}
\index{reconnaissance!graphes~d'Artémis}

\begin{algorithme}\label{reco.a.artemis}
  \begin{itemize}
  \item[\sc Entrée :] Un graphe $G$.

  \item[\sc Sortie  :] Si  le graphe  est d'Artémis~:  {\sc Oui}.  Sinon,
  {\sc Non}.

  \item[\sc Calcul :] \mbox{}
    \begin{itemize}
    \item
      Vérifier que  $G$ ne  possède pas d'antitrou  long à  l'aide de
      l'algorithme~\ref{reco.a.ft} exécuté sur $\overline{G}$.

    \item
      Vérifier    que    $G$    est    de   Berge    à    l'aide    de
      l'algorithme~\ref{reco.a.berge}.
    \item
      Vérifier que  $G$ est sans prisme  et sans pyramide  à l'aide de
      l'algorithme~\ref{reco.a.pypri}.
    \item 
      Si  toutes  les vérifications  ont  donné {\sc Oui},  retourner
      {\sc Oui}. Sinon, retourner {\sc Non}.
    \end{itemize}
  \item
    [\sc Complexité :] $O(n^9)$.
  \end{itemize}
\end{algorithme}

\begin{preuve}
  Quasiment évidente.   Il faut juste  noter que les graphes  de Berge
  n'ont   pas   de   pyramide.    Donc,   si   $G$   est   de   Berge,
  l'algorithme~\ref{reco.a.pypri} vérifie simplement  que $G$ est sans
  prisme.
\end{preuve}

\index{Artémis~(graphe~d'---)!coloration}
\index{coloration!graphes~d'Artémis}  L'algorithme  de  reconnaissance
des  graphes  d'Artémis  permet  également  de  colorier  les  graphes
d'Artémis en temps polynomial.  On énumère toutes les paires d'amis de
$G$, que l'on contracte (rappelons qu'on sait vérifier qu'une paire de
sommets  d'un graphe  de Berge  est  une paire  d'amis).  Pour  chaque
paire, on  vérifie que  le graphe obtenu  est bien d'Artémis  (grâce à
l'algorithme ci-dessus).  Si c'est le cas, on recommence sur le graphe
contracté     juqu'à      obtenir     une     clique.       Par     le
théorème~\ref{artemis.t.main},  on  sait  que cela  fonctionnera.   Au
pire, il faudra appeler $n^3$ fois l'algorithme~\ref{reco.a.artemis}~:
car dans  un graphe à $n$  sommets, il y  a au pire $O(n^2)$  paires à
tester, et au  pire, il faudra faire $n$  contractions avant d'obtenir
une  clique.    La  complexité   de  cet  algorithme   est  finalement
$O(n^{12})$,   ce    qui   est   assez   médiocre    par   rapport   à
l'algorithme~\ref{artemis.a.color} qui est de complexité $O(n^2m)$.

\index{reconnaissance!graphes~d'Artémis~pairs}
\index{Artémis~pair~(graphe~d'---)!reconnaissance}
\begin{algorithme}\label{reco.a.grenoble}
  \begin{itemize}
  \item[\sc Entrée :] Un graphe $G$.

  \item[\sc  Sortie  :]  Si   le  graphe  est  d'Artemis  pair~:  {\sc
  Oui}. Sinon, {\sc Non}.

  \item[\sc Calcul :] \mbox{}
    \begin{itemize}
          
    \item 
      À l'aide de  l'algorithme~\ref{reco.a.artemis}, verifier que $G$
      est d'Artémis.  Si tel est  le cas poursuivre, sinon  stopper et
      répondre {\sc Non}.

    \item 
      À l'aide de l'algorithme~\ref{reco.a.prismeimp} vérifier que $G$
      est sans prisme impair.  Si tel est le cas, retourner {\sc Oui}
      sinon, stopper et répondre {\sc Non}.
    \end{itemize}
  \item[\sc Complexité :] $O(n^{20})$. 
\end{itemize}
\end{algorithme}

\begin{preuve}
  Il faut simplement remarquer  que les graphes d'Artémis ne possèdent
  pas  de  trou  impair,   condition  requise  pour  pouvoir  exécuter
  l'algorithme~\ref{reco.a.prismeimp}.
\end{preuve}

Si    les    conjectures   de    Everett    et   Reed    (conjecture~%
\ref{pair.conj.EverettReed2}            et~\ref{pair.conj.EverettReed},
page~\pageref{pair.conj.EverettReed})  devaient  être prouvées,  alors
cela  impliquerait  l'existence  d'un  algorithme par  contraction  de
paires d'amis pour colorier  les graphes parfaitement contractiles. Il
faut bien noter que ce fait, certes plausible, n'a rien d'évident même
une fois connu  l'algorithme de recherche d'une paire  d'amis dans les
graphes  de Berge.   Si la  conjecture~\ref{pair.conj.EverettReed} est
vraie, alors  les graphes d'Artémis pairs sont  exactement les graphes
parfaitement contractiles.   Dans ce  cas, tout graphe  d'Artémis pair
différent d'une  clique possède une  paire d'amis dont  la contraction
redonne     un      graphe     d'Artémis     pair      (d'après     la
conjecture~\ref{pair.conj.EverettReed2}).  Le principe de l'algorithme
est donc simple~:  pour toute paire d'amis du  graphe, vérifier que sa
contraction donne  à nouveau  un graphe parfaitement  contractile.  Si
les conjectures de Reed et Everett sont vraies, cette vérification est
possible  grâce  à  l'algorithme~\ref{reco.a.grenoble}.  Une  fois  la
paire trouvée, recommencer sur le graphe contracté jusqu'à obtenir une
clique.  Au  pire, il  faut $O(n^3) \times  O(n^{20})$ étapes,  ce qui
donne une complexité de $O(n^{23})$. En résumé~:

\index{Artémis~pair~(graphe~d'---)!coloration}
\index{coloration!graphes d'Artémis~pairs}
\begin{algorithme}  \label{reco.a.colorpc}
  \begin{itemize}
  \item[\sc Entrée :] Un graphe d'Artémis pair $G$.

  \item[\sc Sortie :]  Une coloration de $G$, qui  est optimale si les
  conjectures~\ref{pair.conj.EverettReed2}
  et~\ref{pair.conj.EverettReed} sont vraies.

  \item[\sc Calcul :] Voir ci-dessus.
  
  \item[\sc Complexité :] $O(n^{23})$. 
\end{itemize}
\end{algorithme}

Voilà  une motivation  supplémentaire pour  chercher une  preuve  à la
conjecture  de Everett  et Reed,  car  quelle que  soit cette  preuve,
algorithmique ou non, on sait  déjà qu'elle permettra de colorier tous
les graphes  parfaitement contractiles  en temps polynomial.   On peut
toutefois  espérer,  comme  pour  les graphes  d'Artémis,  une  preuve
donnant un algorithme de coloration plus efficace.

On  peut aussi  donner  un algorithme  de  reconnaissance des  graphes
bipartisans~:

\index{reconnaissance!graphes bipartisans}
\index{bipartisan~(graphe~---)!reconnaissance}
\begin{algorithme}\label{reco.a.bipartisan}
  \begin{itemize}
  \item[\sc Entrée :] Un graphe $G$.

  \item[\sc Sortie  :] Si  le graphe  est bipartisan~:  {\sc Oui}.  Sinon,
  {\sc Non}.

  \item[\sc Calcul :] \mbox{}
    \begin{itemize}
     \item
      Vérifier    que    $G$    est    de   Berge    à    l'aide    de
      l'algorithme~\ref{reco.a.berge}.
    \item
      Vérifier que  $G$ est sans prisme  long à l'aide de 
      l'algorithme~\ref{reco.a.prilong}.
    \item 
      Vérifier  que $G$  ne possède  ni double  diamant  ni $L(K_{3,3}
      \setminus e)$ en étudiant  systématiquement tous les 8-uplets de
      sommets de $G$.
    \item 
      Si  toutes  les vérifications  ont  donné {\sc Oui},  retourner
      {\sc Oui}. Sinon, retourner {\sc Non}.
    \end{itemize}
  \item
    [\sc Complexité :] $O(n^9)$.
  \end{itemize}
\end{algorithme}

\begin{preuve}
  Clair.
\end{preuve}

\section{Questions ouvertes}

La question  ouverte la plus évidente  (non pas à résoudre,  mais à se
poser) est  celle de la  reconnaissance des graphes sans  trou impair.
Notons que Conforti, Cornuéjols, Kapoor et Vu\v skovi\'c ont résolu un
problème  similaire~:   la  reconnaissance  des   graphes  sans  trou
pair~\cite{conforti.c.k.v:eh1,conforti.c.k.v:eh2}.  Mais  de nombreux
autres problèmes demeurent.

\subsection*{Reconnaissance rapide des graphes d'Artémis}

\index{Artémis~(graphe~d'---)!reconnaissance~rapide}
\index{reconnaissance!graphes~d'Artémis}

\index{nettoyage}  \index{raquette} Les sections  précédentes semblent
indiquer que  la reconnaissance des graphes d'Artémis  est plus facile
que  celle  des  graphes   de  Berge.   Tout  d'abord,  l'interdiction
simultanée des  prismes et des pyramides facilite  leur détection.  La
section        de       l'article        de        Chudnovsky       et
Seymour~\cite{chudnovsky.seymour:reco}  consacrée à  la  détection des
pyramides est  donc simplifiée, et leur algorithme  accéléré.  De plus
pour  la  partie   ``nettoyage'',  on  n'a  plus  à   se  soucier  des
``raquettes''  des  plus  petits  trous  impairs:  elles  donnent  des
pyramides.  L'interdiction  des antitrous facilite la  preuve du lemme
sur les  composantes anticonnexes de sommets majeurs  d'un trou impair
de taille  minimale (lemme 4.1 de~\cite{chudnovsky.c.l.s.v:cleaning}).
Malgré  toutes ces simplifications,  certaines étapes  essentielles de
l'algorithme  de  reconnaissance des  graphes  de  Berge  ne sont  pas
accélérées   dans  le   cas  des   graphes  d'Artémis,   notamment  le
spectaculaire  lemme  5.1  de~\cite{chudnovsky.c.l.s.v:cleaning},  qui
permet de ``deviner''  en $O(n^9)$ les sommets majeurs  d'un trou dans
l'étape du nettoyage.

\subsection*{Plus courts chemins impairs}

\index{chaîne!algorithme~de~plus~courte~---~de~longueur~paire}
\index{chemin!plus~court~---~de~longueur~impaire}    Il    existe   un
algorithme  de  plus   courte  \emph{chaîne}  impaire,  découvert  par
J.~Edmonds         (nous        l'avons         déjà        mentionné,
algorithme~\ref{pair.a.edmonds}  page~\pageref{pair.a.edmonds}).  Mise
à part son application directe  pour la détection de paire d'amis dans
les  line-graphes  (algorithme~\ref{pair.a.eplg}),  nous  n'avons  pas
réussi  à  utiliser cet  algorithme  pour  résoudre  des problèmes  de
reconnaissance.  La difficulté, c'est  qu'à cause des cordes, une plus
courte  chaîne impaire  ne donne  pas  toujours un  plus court  chemin
impair.

On sait grâce à Bienstock~\cite{bienstock:evenpair} que le problème de
décider si une  paire de sommets d'un graphe  quelconque est une paire
d'amis est  NP-complet.  Donc le  problème de la recherche  d'un plus
court chemin  impair entre deux sommets donnés  d'un graphe quelconque
est lui aussi clairement NP-difficile.  Mais comme on l'a déjà dit, le
problème de la paire d'amis devient polynomial dans le cas particulier
des  graphes de Berge.   András Seb{\H  o} m'a  donc posé  la question
suivante~:
 
\begin{question}
  Existe-t-il un  algorithme de  plus court \emph{chemin}  de longueur
  impaire dans les graphes sans trou impair~?
\end{question}

Notons que  s'il existait un jour  un algorithme de  recherche de plus
court  trou impair,  alors  cet algorithme  permettrait facilement  de
répondre positivement  à la question  ci-dessus.  Mais réciproquement,
une  réponse  positive  à  la  réponse  ci-dessus  n'impliquerait  pas
nécessairement qu'on  sache détecter les trous impairs  dans un graphe
quelconque.  Signalons la  thèse  d'Éric Tannier~\cite{tannier:these},
dont  le dernier  chapitre est  consacré à  une vue  d'ensemble  de la
littérature sur les problèmes de recherche de chemins avec contraintes
de parité.

\subsection*{Trigraphes de Berge}
\index{Berge~(trigraphe~de~---)}
\index{trigraphe}

Dans sa thèse, Maria Chudnovsky  a introduit la notion de trigraphe de
Berge         que        nous        avons         déja        évoquée
page~\pageref{graphespar.ss.trigraphes}. Rappelons qu'un trigraphe $G$
est un triplet $(V, E_1, E_2)$.  Les ensembles $E_1$ et $E_2$ sont des
sous-ensembles disjoints de l'ensemble  des paires de $V$.  L'ensemble
$E_1$ représente l'ensemble des arêtes ``obligatoires'' de $G$, tandis
que   $E_2$  est   l'ensemble  des   arêtes   ``optionnelles''.   Plus
formellement, on appelle \emph{réalisation} de $G$ tout graphe $G'=(V,
E)$ vérifiant $E_1 \subseteq E  \subseteq E_1 \cup E_2$.  Un trigraphe
est dit \emph{de Berge} si toutes ses réalisations sont des graphes de
Berge.  On  peut naturellement s'interroger sur  la reconnaissance des
trigraphes de Berge, problème  que nous soupçonnons NP-difficile, sans
avoir pu le prouver.

Nous  avons trouvé deux  articles (\cite{figuereido.k.v:1joinsandwich}
et   \cite{cerrioli.e.f.k:homogeneoussandwich})    portant   sur   des
questions  similaires,  appelées par  leurs  auteures (Maria  Cerioli,
Hazel Everett, Celina de  Figueiredo, Sulamita Klein et Kristina Vu{\v
s}kovi{\'c}),   ``problèmes  de   sandwichs''.    Reformulé  dans   le
vocabulaire  des trigraphes,  un problème  de sandwich  consiste  à se
demander, étant donné  un trigraphe $G$ et un graphe  $H$, si l'une au
moins des réalisations de $G$ contient $H$. Reconnaître les trigraphes
de Berge, c'est  donc se poser le problème du  sandwich pour les trous
impairs et les antitrous impairs.

\subsection*{Lecture algorithmique de la preuve de la conjecture forte
des graphes parfaits}

Les  algorithmes  de  détection  de  sous-graphes  présentés  dans  ce
chapitre étaient à l'origine motivés par la reconnaissance des graphes
d'Artémis. Par une curieuse  coïncidence, ils conduisent à détecter la
majorité  des sous-graphes utilisés  dans la  preuve de  la conjecture
forte des graphes parfaits. On peut donc se demander s'il est possible
de lire cette preuve comme un algorithme en temps polynomial qui étant
donné  un  graphe de  Berge  montre qu'il  est  basique,  ou alors  le
décompose.  Répondre  à cette question  nécessiterait probablement des
vérifications  fastidieuses pour un  résultat assez  maigre, puisqu'on
sait   directement   détecter   les   2-joints   et   les   partitions
antisymétriques  dans  n'importe  quel graphe.  Mentionnons  toutefois
qu'un tel travail permettrait  peut-être de résoudre (pour les graphes
de Berge) une question ouverte pour autant que nous le sachions~:

\index{partition~antisymétrique!détection}
\index{détection!partition~antisymétrique}
\begin{question}
  Existe-t-il un algorithme en temps polynomial qui décide si un
  graphe possède une partition antisymétrique paire~?
\end{question}


\appendix
\addcontentsline{toc}{chapter}{Annexes}

\addcontentsline{toc}{chapter}{\numberline {A}Bibliographie}

\addcontentsline{toc}{chapter}{\numberline {B}Index}
\small
\input{theseArxiv.ind}
\newpage
\thispagestyle{empty}
\mbox{}
\voffset-1.1cm
\textheight25cm
\leftmargin0cm
\newpage
\normalsize
\thispagestyle{empty}
\noindent\hspace{-.7cm}
\begin{tabular}{l}
  \parbox{12cm}{
    \noindent
	{\bf  Résumé~:} Ce  travail  a pour  motivation une  meilleure
	compréhension des  graphes parfaits.  La preuve en  2002 de la
	conjecture   des  graphes   parfaits  de   Claude   Berge  par
	Chudnovsky, Robertson,  Seymour et  Thomas a jeté  une lumière
	nouvelle  sur ce  domaine de  la combinatoire,  mais  a laissé
	plusieurs  questions en  suspens,  notamment l'existence  d'un
	algorithme   combinatoire  de  co\-lo\-ra\-tion   des  graphes
	parfaits.

	\vspace{.5ex}
	
	\noindent
	Une paire  d'amis d'un graphe  est une paire de  sommets telle
	que  tous les  chemins  les reliant  sont  de longueur  paire.
	Comme l'ont montré Fonlupt et Urhy, la contraction d'une paire
	amis préserve  le nombre  chromatique du graphe,  et appliquée
	récursivement,   permet   dans   certains  cas   de   colorier
	optimalement  le  graphe.   Nous  prouvons une  conjecture  de
	Everett et  Reed affirmant que cette  approche fonctionne pour
	une classe de graphes  parfaits~: les graphes d'Artémis.  Nous
	en déduisons  un algorithme de coloration  des graphes Artémis
	de complexité $O(n^2m)$.

	\vspace{.5ex}
	
	\noindent
	Nous  donnons   un  algorithme   de  complexité  $O(n^9)$   pour  la
  reconnaissance  des  graphes   d'Artémis.  D'autres  algorithmes  de
  reconnaissance  sont  donnés,  tous   fondés  sur  des  routines  de
  détection de sous-graphes dans  des graphes de Berge.  Nous montrons
  que ces problèmes de détection  sont NP-complets si on cherche à les
  étendre aux graphes quelconques.

  \vspace{.8ex}

  \noindent
  {\bf   Mots-clef~:}  Graphe   parfait,  algorithme,   paire  d'amis,
  coloration, reconnaissance, graphe parfaitement contractile.

  \vspace{.9cm} 
  }
  \\\hline
  \parbox{12cm}{

  \vspace{.3cm} 

  \noindent
  {\bf Abstract:}  This work is motivated  by the desire  for a better
  understanding of  perfect graphs.  The  proof of the  Claude Berge's
  perfect graph  conjecture in 2002 by  Chudnovsky, Robertson, Seymour
  and Thomas has shed a new light on this field of combinatorics.  But
  some questions  are still unsettled, particulary the  existence of a
  combinatorial algorithm for the coloring of perfect graphs.

  \vspace{.5ex}

  \noindent
  An even pair of  a graph is a pair of vertices  such that every path
  joining them  has even  length. As proved  by Fonlupt and  Uhry, the
  contraction of an even pair preserves the chromatic number, and when
  applied  recursively may  lead to  an optimal  coloring. We  prove a
  conjecture of Everett  and Reed saying that this  method works for a
  class  of perfect  graphs: Artemis  graphs. This  yields  a coloring
  algorithm for Artemis graphs with complexity $O(n^2m)$.

  \vspace{.5ex}

  \noindent
  We  give  an  $O(n^9)$  algorithm  for the  recognition  of  Artemis
  graphs. Other  recognition algorithms are  also given, each  of them
  based on subgraph detection routines for Berge graphs.  We show that
  these subgraph  detection problems are NP-complete  when extended to
  general graphs.

  \vspace{1ex}

  \noindent
  {\bf  Key-words:}  Perfect graph,  algorithm,  even pair,  coloring,
  recognition, perfectly contractile graph.
  }

  \vspace{.9cm}

  \\
  \hline
  \parbox{12cm}{
  \vspace{.3cm}

  \noindent
  {\bf  Discipline~:}  Mathématiques Informatique.

  \noindent
  {\bf  Formation doctorale~:}  Recherche Opérationnelle, Combinatoire
  et Optimisation.  
  }
\end{tabular}
\end{document}